# Iwasawa Theory for Elliptic Curves

Ralph Greenberg

University of Washington

## 1. Introduction.

The topics that we will discuss have their origin in Mazur's synthesis of the theory of elliptic curves and Iwasawa's theory of $\mathbb{Z}_p$-extensions in the early 1970s. We first recall some results from Iwasawa's theory. Suppose that $F$ is a finite extension of $\mathbb{Q}$ and that $F_\infty$ is a Galois extension of $F$ such that $\mathrm{Gal}(F_\infty/F) \cong \mathbb{Z}_p$, the additive group of $p$-adic integers, where $p$ is any prime. Equivalently, $F_\infty = \bigcup_{n \geq 0} F_n$, where, for $n \geq 0$, $F_n$ is a cyclic extension of $F$ of degree $p^n$ and $F = F_0 \subset F_1 \subset \cdots \subset F_n \subset F_{n+1} \subset \cdots$. Let $h_n$ denote the class number of $F_n$, $p^{e_n}$ the exact power of $p$ dividing $h_n$. Then Iwasawa proved the following result.

**Theorem 1.1.** *There exist integers $\lambda$, $\mu$, and $\nu$, which depend only on $F_\infty/F$, such that $e_n = \lambda n + \mu p^n + \nu$ for $n \gg 0$.*

The idea behind the proof of this result is to consider the Galois group $X = \mathrm{Gal}(L_\infty/F_\infty)$, where $L_\infty$ is the maximal abelian extension of $F_\infty$ which is unramified at all primes of $F_\infty$ and such that $\mathrm{Gal}(L_\infty/F_\infty)$ is a pro-$p$ group. In fact, $L_\infty = \bigcup_{n \geq 0} L_n$, where $L_n$ is the $p$-Hilbert class field of $F_n$ for $n \geq 0$. Now $L_\infty/F$ is Galois and $\Gamma = \mathrm{Gal}(F_\infty/F)$ acts by inner automorphisms on the normal subgroup $X$ of $\mathrm{Gal}(L_\infty/F)$. Thus, $X$ is a $\mathbb{Z}_p$-module and $\Gamma$ acts on $X$ continuously and $\mathbb{Z}_p$-linearly. It is natural to regard $X$ as a module over the group ring $\mathbb{Z}_p[\Gamma]$, but even better over the completed group ring

$$\Lambda = \mathbb{Z}_p[[\Gamma]] = \varprojlim \mathbb{Z}_p[\mathrm{Gal}(F_n/F)],$$

where the inverse limit is defined by the ring homomorphisms induced by the restriction maps $\mathrm{Gal}(F_m/F) \to \mathrm{Gal}(F_n/F)$ for $m \geq n \geq 0$. The ring $\Lambda$ is sometimes called the "Iwasawa algebra" and has the advantage of being a complete local ring. More precisely, $\Lambda \cong \mathbb{Z}_p[[T]]$, where $T$ is identified with $\gamma - 1 \in \Lambda$. Here $\gamma \in \Gamma$ is chosen so that $\gamma|_{F_1}$ is nontrivial, and 1 is the identity element of $\Gamma$ (and of the ring $\Lambda$). Then $\gamma$ generates a dense subgroup of $\Gamma$ and the action of $T = \gamma - 1$ on $X$ is "topologically nilpotent." This allows one to consider $X$ as a $\Lambda$-module.

Iwasawa proves that $X$ is a finitely generated, torsion $\Lambda$-module. There is a structure theorem for such $\Lambda$-modules which states that there exists a "pseudo-isomorphism"

$$X \sim \bigoplus_{i=1}^{t} \Lambda/(f_i(T)^{a_i}),$$

where each $f_i(T)$ is an irreducible element of $\Lambda$ and the $a_i$'s are positive integers. (We say that two finitely generated, torsion $\Lambda$-modules $X$ and $Y$ are pseudo-isomorphic when there exists a $\Lambda$-homomorphism from $X$ to $Y$ with finite kernel and cokernel. We then write $X \sim Y$.) It is natural to try to recover $\mathrm{Gal}(L_n/F_n)$ from $X = \mathrm{Gal}(L_\infty/F_\infty)$.

Suppose that $F$ has only one prime lying over $p$ and that this prime is totally ramified in $F_\infty/F$. (Totally ramified in $F_1/F$ suffices for this.) Then one can indeed recover $\mathrm{Gal}(L_n/F_n)$ from the $\Lambda$-module $X$. We have

$$\mathrm{Gal}(L_n/F_n) \cong X/(\gamma^{p^n} - 1)X.$$

The isomorphism is induced from the restriction map $X \to \mathrm{Gal}(L_n/F_n)$. Here is a brief sketch of the proof: $\mathrm{Gal}(F_\infty/F_n)$ is topologically generated by $\gamma^{p^n}$; one verifies that $(\gamma^{p^n} - 1)X$ is the commutator subgroup of $\mathrm{Gal}(L_\infty/F_n)$; and one proves that the maximal abelian extension of $F_n$ contained in $L_\infty$ is precisely $F_\infty L_n$. (This last step is where one uses the fact that there is only one prime of $F_n$ lying over $p$.) Then one notices that $\mathrm{Gal}(L_n/F_n) \cong \mathrm{Gal}(F_\infty L_n/F_\infty)$. If $F$ has more than one prime over $p$, one can still recover $\mathrm{Gal}(L_n/F_n)$ for $n \gg 0$, somehow taking into account the inertia subgroups of $\mathrm{Gal}(L_\infty/F_n)$ for primes over $p$. (Primes not lying over $p$ are unramified.) One can find more details about the proof in [Wa2].

The invariants $\lambda$ and $\mu$ can be obtained from $X$ in the following way. Let $f(T)$ be a nonzero element of $\Lambda$: $f(T) = \sum\limits_{i=0}^{\infty} c_i T^i$, where $c_i \in \mathbb{Z}_p$ for $i \geq 0$. Let $\mu(f) \geq 0$ be defined by: $p^{\mu(f)} | f(T)$, but $p^{\mu(f)+1} \nmid f(T)$ in $\Lambda$. Thus, $f(T)p^{-\mu(f)}$ is in $\Lambda$ and has at least one coefficient in $\mathbb{Z}_p^\times$. Define $\lambda(f) \geq 0$ to be the smallest $i$ such that $c_i p^{-\mu(f)} \in \mathbb{Z}_p^\times$. (Thus, $f(T) \in \Lambda^\times$ if and only if $\lambda(f) = \mu(f) = 0$.) Let $f(T) = \prod\limits_{i=1}^{t} f_i(T)^{a_i}$. The ideal $(f(T))$ of $\Lambda$ is called the "characteristic ideal" of $X$. Then it turns out that the $\lambda$ and $\mu$ occurring in Iwasawa's theorem are given by $\lambda = \lambda(f)$, $\mu = \mu(f)$. For each $i$, there are two possibilities: either $f_i(T)$ is an associate of $p$, in which case $\mu(f_i) = 1$, $\lambda(f_i) = 0$, and $\Lambda/(f_i(T)^{a_i})$ is an infinite group of exponent $p^{a_i}$, or $f_i(T)$ is an associate of a monic polynomial of degree $\lambda(f_i)$, irreducible over $\mathbb{Q}_p$, and "distinguished" (which means that the nonleading coefficients are in $p\mathbb{Z}_p$), in which case $\mu(f_i) = 0$ and $\Lambda/(f_i(T)^{a_i})$ is isomorphic to $\mathbb{Z}_p^{\lambda(f_i)a_i}$ as a group. Then, $\lambda = \Sigma a_i \lambda(f_i)$, $\mu = \Sigma a_i \mu(f_i)$. The invariant $\lambda$ can be described more simply as $\lambda = \mathrm{rank}_{\mathbb{Z}_p}(X/X_{\mathbb{Z}_p\text{-tors}})$, where $X_{\mathbb{Z}_p\text{-tors}}$ is the torsion subgroup of $X$. Equivalently, $\lambda = \dim_{\mathbb{Q}_p}(X \otimes_{\mathbb{Z}_p} \mathbb{Q}_p)$.

The invariants $\lambda = \lambda(F_\infty/F)$ and $\mu = \mu(F_\infty/F)$ are difficult to study. Iwasawa found examples of $\mathbb{Z}_p$-extensions $F_\infty/F$ where $\mu(F_\infty/F) > 0$. In his examples there are infinitely many primes of $F$ which decompose completely in $F_\infty/F$. In these lectures, we will concentrate on the "cyclotomic $\mathbb{Z}_p$-extension" of $F$ which is defined as the unique subfield $F_\infty$ of $F(\mu_{p^\infty})$ with $\Gamma = \mathrm{Gal}(F_\infty/F) \cong \mathbb{Z}_p$. Here $\mu_{p^\infty}$ denotes the $p$-power roots of unity. It is easy to show that all nonarchimedean primes of $F$ are finitely decomposed in $F_\infty/F$. More precisely, if $v$ is any such prime of $F$, then the corresponding decomposition subgroup $\Gamma(v)$ of $\Gamma$ is of finite index. If $v \nmid p$, then the inertia subgroup is trivial, i.e., $v$ is unramified. (This is true for any $\mathbb{Z}_p$-extension.) If $v|p$, then the corresponding inertia subgroup of $\Gamma$ is of finite index. Iwasawa has conjectured that $\mu(F_\infty/F) = 0$ if $F_\infty/F$ is the cyclotomic $\mathbb{Z}_p$-extension. In the case where $F$ is an abelian extension of $\mathbb{Q}$, this has been proved by Ferrero and Washington. (See [FeWa] or [Wa2].)

On the other hand, $\lambda(F_\infty/F)$ can be positive. The simplest example is perhaps the following. Let $F$ be an imaginary quadratic field. Then all $\mathbb{Z}_p$-extensions of $F$ are contained in a field $\widetilde{F}$ such that $\mathrm{Gal}(\widetilde{F}/F) \cong \mathbb{Z}_p^2$. (Thus, there are infinitely many $\mathbb{Z}_p$-extensions of $F$.) Letting $F_\infty/F$ still be the cyclotomic $\mathbb{Z}_p$-extension, one can verify that $\widetilde{F}/F_\infty$ is unramified if $p$ is a prime that splits completely in $F/\mathbb{Q}$. Thus in this case, $F_\infty \subseteq \widetilde{F} \subseteq L_\infty$ and hence $X = \mathrm{Gal}(L_\infty/F_\infty)$ has a quotient $\mathrm{Gal}(\widetilde{F}/F_\infty) \cong \mathbb{Z}_p$. Therefore, $\lambda(F_\infty/F) \geq 1$ if $p$ splits in $F/\mathbb{Q}$. Notice that, since $\widetilde{F}/F$ is abelian, the action of $T = \gamma - 1$ on $\mathrm{Gal}(\widetilde{F}/F_\infty)$ is trivial. Thus, $X/TX$ is infinite. Now if one considers the $\Lambda$-module $Y = \Lambda/(f_i(T)^{a_i})$, where $f_i(T)$ is irreducible in $\Lambda$, then $Y/TY$ is infinite if and only if $f_i(T)$ is an associate of $T$. Therefore, if $F$ is an imaginary quadratic field in which $p$ splits and if $F_\infty$ is the cyclotomic $\mathbb{Z}_p$-extension of $F$, then $T|f(T)$, where $f(T)$ is a generator of the characteristic ideal of $X$. One can prove that $T^2 \nmid f(T)$. (This is an interesting exercise. It is easy to show that $X/TX$ has $\mathbb{Z}_p$-rank 1. One must then show that $X/T^2X$ also has $\mathbb{Z}_p$-rank 1. See [Gr1] for a more general "semi-simplicity" result.)

In contrast, suppose that $F$ is again imaginary quadratic, but that $p$ is inert in $F/\mathbb{Q}$. Then $F$ has one prime over $p$, which is totally ramified in the cyclotomic $\mathbb{Z}_p$-extension $F_\infty/F$. As we sketched earlier, it then turns out that $X/TX$ is finite and isomorphic to the $p$-primary subgroup of the ideal class group of $F$. In particular, it follows that if $p$ does not divide the class number of $F$, then $X = TX$. Nakayama's Lemma for $\Lambda$-modules then implies that $X = 0$ and hence $\lambda(F_\infty/F) = 0$ for any such prime $p$. In general, for arbitrary $n \geq 0$, the restriction map $X \to \mathrm{Gal}(L_n/F_n)$ induces an isomorphism

$$X/\theta_n X \overset{\sim}{\to} \mathrm{Gal}(L_n/F_n),$$

where $\theta_n = \gamma^{p^n} - 1 = (1 + T)^{p^n} - 1$. We can think of $X/\theta_n X$ as $X_{\Gamma_n}$, the maximal quotient of $X$ on which $\Gamma_n$ acts trivially. Here $\Gamma_n = \mathrm{Gal}(F_\infty/F_n)$. It is interesting to consider the duals of these groups. Let

$$S_n = \mathrm{Hom}(\mathrm{Gal}(L_n/F_n), \mathbb{Q}_p/\mathbb{Z}_p), \qquad S_\infty = \mathrm{Hom}_{\mathrm{cont}}(X, \mathbb{Q}_p/\mathbb{Z}_p).$$

Then we can state that $S_n \cong S_\infty^{\Gamma_n}$, where the isomorphism is simply the dual of the map $X_{\Gamma_n} \overset{\sim}{\to} \mathrm{Gal}(L_n/F_n)$. Here $S_\infty^{\Gamma_n}$ denotes the subgroup of $S_\infty$ consisting of elements fixed by $\Gamma_n$. The map $S_n \to S_\infty^{\Gamma_n}$ will be an isomorphism if $F$ is any number field with just one prime lying over $p$, totally ramified in $F_\infty/F$. But returning to the case where $F$ is imaginary quadratic and $p$ splits in $F/\mathbb{Q}$, we have that $S_\infty^\Gamma$ is infinite. (It contains $\mathrm{Hom}(\mathrm{Gal}(\widetilde{F}/F_\infty), \mathbb{Q}_p/\mathbb{Z}_p)$ which is isomorphic to $\mathbb{Q}_p/\mathbb{Z}_p$.) Thus, $S_\infty^{\Gamma_n}$ is always infinite, but $S_n$ is finite, for all $n \geq 0$. The groups $S_n$ and $S_\infty$ are examples of "Selmer groups," by which we mean that they are subgroup of Galois cohomology groups defined by imposing local restrictions. In fact, $S_n$ is the group of cohomology classes in $H^1(G_{F_n}, \mathbb{Q}_p/\mathbb{Z}_p)$ which are unramified at all primes of $F_n$, and $S_\infty$ is the similarly defined subgroup of $H^1(G_{F_\infty}, \mathbb{Q}_p/\mathbb{Z}_p)$. Here, for any field $M$, we let $G_M$ denote the absolute Galois group of $M$. Also, the action of the Galois groups on $\mathbb{Q}_p/\mathbb{Z}_p$ is taken to be trivial. As is customary, we will denote the Galois cohomology group $H^i(G_M, *)$ by $H^i(M, *)$. We will denote $H^i(\mathrm{Gal}(K/M), *)$ by $H^i(K/M, *)$ for any Galois extension $K/M$. We always require cocycles to be continuous. Usually, the group indicated by $*$ will be a $p$-primary group which is given the discrete topology. We will also always understand $\mathrm{Hom}(\ ,\ )$ to refer to the set of continuous homomorphisms.

Now we come to Selmer groups for elliptic curves. Suppose that $E$ is an elliptic curve defined over $F$. We will later recall the definition of the classical Selmer group $\mathrm{Sel}_E(M)$ for $E$ over $M$, where $M$ is any algebraic extension of $F$. Right now, we will just mention the exact sequence

$$0 \to E(M) \otimes (\mathbb{Q}/\mathbb{Z}) \to \mathrm{Sel}_E(M) \to \text{Ш}_E(M) \to 0,$$

where $E(M)$ denotes the group of $M$-rational points on $E$ and $\text{Ш}_E(M)$ denotes the Shafarevich-Tate group for $E$ over $M$. We denote the $p$-primary subgroups of $\mathrm{Sel}_E(M)$, $\text{Ш}_E(M)$ by $\mathrm{Sel}_E(M)_p$, $\text{Ш}_E(M)_p$. The $p$-primary subgroup of the first term above is $E(M) \otimes (\mathbb{Q}_p/\mathbb{Z}_p)$. Also, $\mathrm{Sel}_E(M)_p$ is a subgroup of $H^1(M, E[p^\infty])$, where $E[p^\infty]$ is the $p$-primary subgroup of $E(\overline{\mathbb{Q}})$. As a group, $E[p^\infty] \cong (\mathbb{Q}_p/\mathbb{Z}_p)^2$, but the action of $G_F$ is quite nontrivial. Let $F_\infty/F$ denote the cyclotomic $\mathbb{Z}_p$-extension. We will now state a number of theorems and conjectures, which constitute part of what we call "Iwasawa Theory for $E$." Some of the theorems will be proved in these lectures. We always assume that $F_\infty$ is the cyclotomic $\mathbb{Z}_p$-extension of $F$.

**Theorem 1.2 (Mazur's Control Theorem).** *Assume that $E$ has good, ordinary reduction at all primes of $F$ lying over $p$. Then the natural maps*

$$\mathrm{Sel}_E(F_n)_p \to \mathrm{Sel}_E(F_\infty)_p^{\Gamma_n}$$

*have finite kernel and cokernel, of bounded order as $n$ varies.*

The natural maps referred to are those induced by the restriction maps $H^1(F_n, E[p^\infty]) \to H^1(F_\infty, E[p^\infty])$. One should compare this result with the remarks made above concerning $S_n$ and $S_\infty^{\Gamma_n}$. We will discuss below the cases where $E$ has either multiplicative or supersingular reduction at some primes of $F$ lying over $p$. But first we state an important conjecture of Mazur.

**Conjecture 1.3.** *Assume that $E$ has good, ordinary reduction at all primes of $F$ lying over $p$. Then $\mathrm{Sel}_E(F_\infty)_p$ is $\Lambda$-cotorsion.*

Here $\Gamma = \mathrm{Gal}(F_\infty/F)$ acts naturally on the group $H^1(F_\infty, E[p^\infty])$, which is a torsion $\mathbb{Z}_p$-module, every element of which is killed by $T^n$ for some $n$. Thus, $H^1(F_\infty, E[p^\infty])$ is a $\Lambda$-module. $\mathrm{Sel}_E(F_\infty)_p$ is invariant under the action of $\Gamma$ and is thus a $\Lambda$-submodule. We say that $\mathrm{Sel}_E(F_\infty)_p$ is $\Lambda$-cotorsion if

$$X_E(F_\infty) = \mathrm{Hom}(\mathrm{Sel}_E(F_\infty)_p, \mathbb{Q}_p/\mathbb{Z}_p)$$

is $\Lambda$-torsion. Here $\mathrm{Sel}_E(F_\infty)_p$ is a $p$-primary group with the discrete topology. Its Pontryagin dual $X_E(F_\infty)$ is an abelian pro-$p$ group, which we regard as a $\Lambda$-module. It is not hard to prove that $X_E(F_\infty)$ is finitely generated as a $\Lambda$-module (and so, $\mathrm{Sel}_E(F_\infty)_p$ is a "cofinitely generated" $\Lambda$-module). In the case where $E$ has good, ordinary reduction at all primes of $F$ over $p$, one can use theorem 1.2. For $X_E(F) = \mathrm{Hom}(\mathrm{Sel}_E(F)_p, \mathbb{Q}_p/\mathbb{Z}_p)$ is known to be finitely generated over $\mathbb{Z}_p$. (In fact, the weak Mordell-Weil theorem is proved by showing that $X_E(F)/pX_E(F)$ is finite.) Write $X = X_E(F_\infty)$ for brevity. Then, by theorem 1.2, $X/TX$ is finitely generated over $\mathbb{Z}_p$. Hence, $X/\mathfrak{m}X$ is finite, where $\mathfrak{m} = (p, T)$ is the maximal ideal of $\Lambda$. By a version of Nakayama's Lemma (valid for profinite $\Lambda$-modules $X$), it follows that $X_E(F_\infty)$ is indeed finitely generated as a $\Lambda$-module. (This can actually be proved for any prime $p$, with no restriction on the reduction type of $E$.) Here is one important case where the above conjecture can be verified.

**Theorem 1.4.** *Assume that $E$ has good, ordinary reduction at all primes of $F$ lying over $p$. Assume also that $\mathrm{Sel}_E(F)_p$ is finite. Then $\mathrm{Sel}_E(F_\infty)_p$ is $\Lambda$-cotorsion.*

This theorem is an immediate corollary of theorem 1.2, using the following exercise: if $X$ is a $\Lambda$-module such that $X/TX$ is finite, then $X$ is a torsion $\Lambda$-module. The hypothesis on $\mathrm{Sel}_E(F)_p$ is equivalent to assuming that both the Mordell-Weil group $E(F)$ and the $p$-Shafarevich-Tate group $Ш_E(F)$ are finite. A much deeper case where conjecture 1.3 is known is the following. The special case where $E$ has complex multiplication had previously been settled by Rubin [Ru1].

**Theorem 1.5 (Kato-Rohrlich).** *Assume that $E$ is defined over $\mathbb{Q}$ and is modular. Assume also that $E$ has good, ordinary reduction or multiplicative reduction at $p$ and that $F/\mathbb{Q}$ is abelian. Then $\mathrm{Sel}_E(F_\infty)_p$ is $\Lambda$-cotorsion.*

The case where $E$ has multiplicative reduction at a prime $v$ of $F$ lying over $p$ is somewhat analogous to the case where $E$ has good, ordinary reduction at $v$. In both cases, the $G_{F_v}$-representation space $V_p(E) = T_p(E) \otimes \mathbb{Q}_p$ has an unramified 1-dimensional quotient. (Here $T_p(E)$ is the Tate-module for $E$; $V_p(E)$ is a 2-dimensional $\mathbb{Q}_p$-vector space on which the local Galois group $G_{F_v}$ acts, where $F_v$ is the $v$-adic completion of $F$.) It seems reasonable to believe that the analogue of Theorem 1.2 should hold. This was first suggested by Manin [Man] for the case $F = \mathbb{Q}$.

**Conjecture 1.6.** *Assume that $E$ has good, ordinary reduction or multiplicative reduction at all primes of $F$ lying over $p$. Then the natural maps*

$$\operatorname{Sel}_E(F_n)_p \to \operatorname{Sel}_E(F_\infty)_p^{\Gamma_n}$$

*have finite kernel and cokernel, of bounded order as $n$ varies.*

For $F = \mathbb{Q}$, this is a theorem. In this case, Manin showed that it would suffice to prove that $\log_p(q_E) \neq 0$, where $q_E$ denotes the Tate period for $E$, assuming that $E$ has multiplicative reduction at $p$. But a recent theorem of Barré-Sirieix, Diaz, Gramain, and Philibert [B-D-G-P] shows that $q_E$ is transcendental when the $j$-invariant $j_E$ is algebraic. Since $j_E \in \mathbb{Q}$, it follows that $q_E p^{-\operatorname{ord}(q_E)}$ is not a root of unity and so $\log_p(q_E) \neq 0$. For arbitrary $F$, one would need to prove that $\log_p(N_{F_v/\mathbb{Q}_p}(q_E^{(v)})) \neq 0$ for all primes $v$ of $F$ lying over $p$ where $E$ has multiplicative reduction. Here $F_v$ is the $v$-adic completion of $F$, $q_E^{(v)}$ the corresponding Tate period. This nonvanishing statement seems intractable at present.

If $E$ has supersingular reduction at some prime $v$ of $F$, then the "control theorem" undoubtedly fails. In fact, $\operatorname{Sel}_E(F_\infty)_p$ will not be $\Lambda$-cotorsion. More precisely, let

$$r(E, F) = \sum_{pss} [F_v : \mathbb{Q}_p],$$

where the sum varies over the primes $v$ of $F$ where $E$ has potentially supersingular reduction. Then one can prove the following result.

**Theorem 1.7.** *With the above notation, we have*

$$\operatorname{corank}_\Lambda(\operatorname{Sel}_E(F_\infty)_p) \geq r(E, F).$$

This result is due to P. Schneider. He conjectures that equality should hold here. (See [Sch2].) This would include for example a more general version of conjecture 1.3, where one assumes just that $E$ has potentially ordinary or potentially multiplicative reduction at all primes of $F$ lying over $p$. As a consequence of theorem 1.7, one finds that

$$\operatorname{corank}_{\mathbb{Z}_p}(\operatorname{Sel}_E(F_\infty)_p^{\Gamma_n}) \geq r(E, F)p^n$$

for $n \geq 0$. This follows from the fact that $\Lambda/\theta_n \Lambda \cong \mathbb{Z}_p^{p^n}$. (The ring $\Lambda/\theta_n \Lambda$ is just $\mathbb{Z}_p[\operatorname{Gal}(F_n/F)]$.) One uses the fact that there is a pseudo-isomorphism from $X_E(F_\infty)$ to $\Lambda^r \oplus Y$, where $r = \operatorname{rank}_\Lambda(X_E(F_\infty))$, which is the $\Lambda$-corank of $\operatorname{Sel}_E(F_\infty)_p$, and $Y$ is the $\Lambda$-torsion submodule of $X_E(F_\infty)$. However, it is reasonable to make the following conjecture. We continue to assume that $F_\infty/F$ is the cyclotomic $\mathbb{Z}_p$-extension, but make no assumptions on the reduction type for $E$ at primes lying over $p$. The conjecture below follows from results of Kato and Rohrlich when $F$ is abelian over $\mathbb{Q}$ and $E$ is defined over $\mathbb{Q}$ and modular.

**Conjecture 1.8.** *The $\mathbb{Z}_p$-corank of $\mathrm{Sel}_E(F_n)_p$ is bounded as $n$ varies.*

If this is so, then the map $\mathrm{Sel}_E(F_n)_p \to \mathrm{Sel}_E(F_\infty)_p^{\Gamma_n}$ must have infinite cokernel when $n$ is sufficiently large, provided that we assume that $E$ has potentially supersingular reduction at $v$ for at least one prime $v$ of $F$ lying over $p$. Of course, assuming that the $p$-Shafarevich-Tate group is finite, the $\mathbb{Z}_p$-corank of $\mathrm{Sel}_E(F_n)_p$ is just the rank of the Mordell-Weil group $E(F_n)$. If one assumes that $E(F_n)$ does indeed have bounded rank as $n \to \infty$ then one can deduce the following nice consequence: $E(F_\infty)$ *is finitely generated. Hence, for some $n \geq 0$, $E(F_\infty) = E(F_n)$.* This is proved in Mazur's article [Maz1]. The crucial step is to show that $E(F_\infty)_{\mathrm{tors}}$ is finite. We refer the reader to Mazur (proposition 6.12) for a detailed proof of this helpful fact. (We will make use of it later. See also [Im] or [Ri].) Using this, one then argues as follows. Let $t = |E(F_\infty)_{\mathrm{tors}}|$. Choose $m$ so that $\mathrm{rank}(E(F_m))$ is maximal. Then, for any $P \in E(F_\infty)$, we have $kP \in E(F_m)$ for some $k \geq 1$. Then $g(kP) = kP$ for all $g \in \mathrm{Gal}(F_\infty/F_m)$. That is, $g(P) - P$ is in $E(F_\infty)_{\mathrm{tors}}$ and hence $t(g(P) - P) = O_E$. This means that $tP \in E(F_m)$. Therefore, $tE(F_\infty) \subseteq E(F_m)$, from which it follows that $E(F_\infty)$ is finitely generated.

On the other hand, let us assume that $E$ has good, ordinary reduction or multiplicative reduction at all primes $v$ of $F$ lying over $p$. Assume also that $\mathrm{Sel}_E(F_\infty)_p$ is $\Lambda$-cotorsion, as is conjectured. Then one can prove conjecture 1.8 very easily. Let $\lambda_E$ denote the $\lambda$-invariant of the torsion $\Lambda$-module $X_E(F_\infty)$. That is, $\lambda_E = \mathrm{rank}_{\mathbb{Z}_p}(X_E(F_\infty)) = \mathrm{corank}_{\mathbb{Z}_p}(\mathrm{Sel}_E(F_\infty)_p)$. We get the following result.

**Theorem 1.9.** *Under the above assumptions, one has*

$$\mathrm{corank}_{\mathbb{Z}_p}(\mathrm{Sel}_E(F_n)_p) \leq \lambda_E.$$

*In particular, the rank of the Mordell-Weil group $E(F_n)$ is bounded above by $\lambda_E$.*

This result follows from the fact that the maps $\mathrm{Sel}_E(F_n)_p \to \mathrm{Sel}_E(F_\infty)_p$ have finite kernel. This turns out to be quite easy to prove, as we will see in section 3. Also, the rank of $E(F_n)$ is the $\mathbb{Z}_p$-corank of $E(F_n) \otimes (\mathbb{Q}_p/\mathbb{Z}_p)$, which is of course bounded above by $\mathrm{corank}_{\mathbb{Z}_p}(\mathrm{Sel}_E(F_n)_p)$. (Equality holds if $\mathrm{III}_E(F_n)_p$ is finite.) Let $\lambda_E^{M\text{-}W}$ denote the maximum of $\mathrm{rank}(E(F_n))$ as $n$ varies, which is just $\mathrm{rank}(E(F_\infty))$. Let $\lambda_E^{\mathrm{III}} = \lambda_E - \lambda_E^{M\text{-}W}$. We let $\mu_E$ denote the $\mu$-invariant of the $\Lambda$-module $X_E(F_\infty)$. If necessary to avoid confusion, we might write $\lambda_E = \lambda_E(F_\infty/F)$, $\mu_E = \mu_E(F_\infty/F)$, etc. Then we have the following analogue of Iwasawa's theorem.

**Theorem 1.10.** *Assume that $E$ has good, ordinary reduction at all primes of $F$ lying over $p$. Assume that $\mathrm{Sel}_E(F_\infty)_p$ is $\Lambda$-cotorsion and that $\mathrm{III}_E(F_n)_p$ is finite for all $n \geq 0$. Then there exist $\lambda, \mu,$ and $\nu$ such that $|\mathrm{III}_E(F_n)_p| = p^{e_n}$, where $e_n = \lambda n + \mu p^n + \nu$ for all $n \gg 0$. Here $\lambda = \lambda_E^{\mathrm{III}}$ and $\mu = \mu_E$.*

As later examples will show, each of the invariants $\lambda_E^{M-W}$, $\lambda_E^{\mathrm{III}}$, and $\mu_E$ can be positive. Mazur first pointed out the possibility that $\mu_E$ could be positive, giving the following example. Let $E = X_0(11)$, $p = 5$, $F = \mathbb{Q}$, and $F_\infty = \mathbb{Q}_\infty =$ the cyclotomic $\mathbb{Z}_5$-extension of $\mathbb{Q}$. Then $\mu_E = 1$. (In fact, $(f_E(T)) = (p)$.) There are three elliptic curves/$\mathbb{Q}$ of conductor 11, all isogenous. In addition to $E$, one of these elliptic curves has $\mu = 2$, another has $\mu = 0$. In general, suppose that $\phi : E_1 \to E_2$ is an $F$-isogeny, where $E_1$, $E_2$ are defined over $F$. Let $\Phi : \mathrm{Sel}_{E_1}(F_\infty)_p \to \mathrm{Sel}_{E_2}(F_\infty)_p$ denote the induced $\Lambda$-module homomorphism. It is not hard to show that the kernel and cokernel of $\Phi$ have finite exponent, dividing the exponent of $\ker(\phi)$. Thus, $\mathrm{Sel}_{E_1}(F_\infty)_p$ and $\mathrm{Sel}_{E_2}(F_\infty)_p$ have the same $\Lambda$-corank. If they are $\Lambda$-cotorsion, then the $\lambda$-invariants are the same. The characteristic ideals of $X_{E_1}(F_\infty)$ and $X_{E_2}(F_\infty)$ differ only by multiplication by a power of $p$. If $F = \mathbb{Q}$, then it seems reasonable to make the following conjecture. For arbitrary $F$, the situation seems more complicated. We had believed that this conjecture should continue to be valid, but counterexamples have recently been found by Michael Drinen.

**Conjecture 1.11.** *Let $E$ be an elliptic curve defined over $\mathbb{Q}$. Assume that $\mathrm{Sel}_E(\mathbb{Q}_\infty)_p$ is $\Lambda$-cotorsion. Then there exists a $\mathbb{Q}$-isogenous elliptic curve $E'$ such that $\mu_{E'} = 0$. In particular, if $E[p]$ is irreducible as a $(\mathbb{Z}/p\mathbb{Z})$-representation of $G_{\mathbb{Q}}$, then $\mu_E = 0$.*

Here $E[p] = \ker(E(\overline{\mathbb{Q}}) \xrightarrow{p} E(\overline{\mathbb{Q}}))$. P. Schneider has given a simple formula for the effect of an isogeny on the $\mu$-invariant of $\mathrm{Sel}_E(F_\infty)_p$ for arbitrary $F$ and for odd $p$. (See [Sch3] or [Pe2].) Thus, the above conjecture effectively predicts the value of $\mu_E$ for $F = \mathbb{Q}$.

Suppose that $\mathrm{Sel}_E(F_\infty)_p$ is $\Lambda$-cotorsion. Let $f_E(T)$ be a generator of the characteristic ideal of $X_E(F_\infty)$. Then $\lambda_E = \lambda(f_E)$ and $\mu_E = \mu(f_E)$. We have

$$X_E(F_\infty) \sim \prod_{i=1}^{t} \Lambda/(f_i(T)^{a_i})$$

where the $f_i(T)$'s are irreducible elements of $\Lambda$, and the $a_i$'s are positive. If $(f_i(T)) = (p)$, then it is possible for $a_i > 1$. However, in contrast, it seems reasonable to make the following "semi-simplicity" conjecture.

**Conjecture 1.12.** *Let $E$ be an elliptic curve defined over $F$. Assume that $\mathrm{Sel}_E(F_\infty)_p$ is $\Lambda$-cotorsion. The action of $\Gamma = \mathrm{Gal}(F_\infty/F)$ on $X_E(F_\infty) \otimes_{\mathbb{Z}_p} \mathbb{Q}_p$ is completely reducible. That is, $a_i = 1$ for all $i$'s such that $f_i(T)$ is not an associate of $p$.*

Assume that $E$ has good, ordinary reduction at all primes of $F$ lying over $p$. Theorem 1.2 then holds. In particular, $\mathrm{corank}_{\mathbb{Z}_p}(\mathrm{Sel}_E(F)_p)$, which is equal to $\mathrm{rank}_{\mathbb{Z}_p}(X_E(F_\infty)/TX_E(F_\infty))$, would equal the power of $T$ dividing $f_E(T)$, assuming the above conjecture. Also, the value of $\lambda_E^{M\text{-}W}$ would be equal to the number of roots of $f_E(T)$ of the form $\zeta - 1$, where $\zeta$ is a $p$-power root of unity, if we assume in addition the finiteness of $\mathrm{III}_E(F_n)_p$ for all $n$. For conjecture 1.12 would imply that this number is equal to the $\mathbb{Z}_p$-rank of $X_E(F_\infty)/\theta_n X_E(F_\infty)$ for $n \gg 0$.

In section 4 we will introduce some theorems due to B. Perrin-Riou and to P. Schneider which give a precise relationship between $\mathrm{Sel}_E(F)_p$ and the behavior of $f_E(T)$ at $T = 0$. These theorems are important because they allow one to study the Birch and Swinnerton-Dyer conjecture by using the so-called "Main Conjecture" which states that one can choose the generator $f_E(T)$ so that it satisfies a certain interpolation property. We will give the statement of this conjecture for $F = \mathbb{Q}$, which was formulated by B. Mazur in the early 1970s (in the same paper [Maz1] where he proves theorem 1.2 and also in [M-SwD]).

**Conjecture 1.13.** *Assume that $E$ is an elliptic curve defined over $\mathbb{Q}$ which has good, ordinary reduction at $p$. Then the characteristic ideal of $X_E(\mathbb{Q}_\infty)$ has a generator $f_E(T)$ with the properties:*

*(i) $f_E(0) = (1 - \beta_p p^{-1})^2 L(E/\mathbb{Q}, 1)/\Omega_E$*

*(ii) $f_E(\phi(T)) = (\beta_p)^n L(E/\mathbb{Q}, \phi, 1)/\Omega_E \tau(\phi)$ if $\phi$ is a finite order character of $\Gamma = \mathrm{Gal}(\mathbb{Q}_\infty/\mathbb{Q})$ of conductor $p^n > 1$.*

We must explain the notation. First of all, fix embeddings of $\overline{\mathbb{Q}}$ into $\mathbb{C}$ and into $\overline{\mathbb{Q}}_p$. $L(E/\mathbb{Q}, s)$ denotes the Hasse-Weil $L$-series for $E$ over $\mathbb{Q}$. $\Omega_E$ denotes the real period for $E$, so that $L(E/\mathbb{Q}, 1)/\Omega_E$ is conjecturally in $\mathbb{Q}$. (If $E$ is modular, then $L(E/\mathbb{Q}, s)$ has an analytic continuation to the complex plane, and, in fact, $L(E/\mathbb{Q}, 1)/\Omega_E \in \mathbb{Q}$.) Let $\widetilde{E}$ denote the reduction of $E$ at $p$. The Euler factor for $p$ in $L(E/\mathbb{Q}, s)$ is $((1 - \alpha_p p^{-s})(1 - \beta_p p^{-s}))^{-1}$, where $\alpha_p$, $\beta_p \in \overline{\mathbb{Q}}$, $\alpha_p \beta_p = p$, $\alpha_p + \beta_p = 1 + p - |\widetilde{E}(\mathbb{F}_p)|$. Choose $\alpha_p$ to be the $p$-adic unit under the fixed embedding $\overline{\mathbb{Q}} \to \overline{\mathbb{Q}}_p$. Thus, $\beta_p p^{-1} = \alpha_p^{-1}$. For every complex-valued, finite order Dirichlet character $\phi$, $L(E/\mathbb{Q}, \phi, s)$ denotes the twisted Hasse-Weil $L$-series. In the above interpolation property, $\phi$ is a Dirichlet character whose associated Artin character factors through $\Gamma$. Using the fixed embeddings chosen above, we can consider $\phi$ as a continuous homomorphism $\phi : \Gamma \to \overline{\mathbb{Q}}_p^\times$ of finite order, i.e., $\phi(\gamma) = \zeta$, where $\zeta$ is a $p$-power root of unity in $\overline{\mathbb{Q}}_p$. Then $\phi(T) = \phi(\gamma - 1) = \zeta - 1$, which is in the maximal ideal of $\overline{\mathbb{Z}}_p$. Hence $f_E(\phi(T)) = f_E(\zeta - 1)$ converges in $\overline{\mathbb{Q}}_p$. The complex number $L(E/\mathbb{Q}, \phi, 1)/\Omega_E$ should be algebraic. In (ii), we regard it as an element of $\overline{\mathbb{Q}}_p$, as well as the Gaussian sum $\tau(\phi)$. For $p > 2$, conjecture 1.13 has been proven by Rubin when $E$ has complex multiplication. (See [Ru2].) If $E$ is a modular elliptic curve with good, ordinary reduction at $p$, then the existence of some power series satisfying the stated interpolation property $(i)$ and $(ii)$ was proven by Mazur and Swinnerton-Dyer in the early 1970s. We will denote it by $f_E^{\mathrm{anal}}(T)$. (See [M-SwD] or [M-T-T].) Conjecturally, this power series should be in $\Lambda$. This is proven in [St] if $E[p]$ is irreducible as a $G_\mathbb{Q}$-module. In general, it is only known to be $\Lambda \otimes_{\mathbb{Z}_p} \mathbb{Q}_p$. That is, $p^t f_E^{\mathrm{anal}}(T) \in \Lambda$ for some $t \geq 0$. Kato then proves that the characteristic ideal at least contains $p^m f_E^{\mathrm{anal}}(T)$ for some $m \geq 0$. Rohrlich proves that $L(E/\mathbb{Q}, \phi, 1) \neq 0$ for all but finitely many characters $\phi$ of $\Gamma$, which is equivalent to the statement $f_E^{\mathrm{anal}}(T) \neq 0$ as an element of $\Lambda \otimes_{\mathbb{Z}_p} \mathbb{Q}_p$. One can use Kato's theorem to prove conjecture 1.13 when $E$ admits a cyclic $\mathbb{Q}$-isogeny of degree $p$, where $p$ is odd and the kernel of the isogeny satisfies a certain condition (namely, the hypotheses in proposition 5.10 in these notes). This will be discussed in [GrVa].

Continuing to assume that $E/\mathbb{Q}$ is modular and that $p$ is a prime where $E$ has good, ordinary reduction, the so-called $p$-adic $L$-function $L_p(E/\mathbb{Q}, s)$ can be defined in terms of $f_E^{\mathrm{anal}}(T)$. We first define a canonical character

$$\kappa : \Gamma \to 1 + 2p\mathbb{Z}_p$$

induced by the cyclotomic character $\chi : \mathrm{Gal}(\mathbb{Q}(\mu_{p^\infty})/\mathbb{Q}) \xrightarrow{\sim} \mathbb{Z}_p^\times$ composed with the projection map to the second factor in the canonical decomposition $\mathbb{Z}_p^\times = \mu_{p-1} \times (1 + p\mathbb{Z}_p)$ for odd $p$, or $\mathbb{Z}_2^\times = \{\pm 1\} \times (1 + 4\mathbb{Z}_2)$ for $p = 2$. Thus, $\kappa$ is an isomorphism. For $s \in \mathbb{Z}_p$, define $L_p(E/\mathbb{Q}, s)$ by

$$L_p(E/\mathbb{Q}, s) = f_E^{\mathrm{anal}}(\kappa(\gamma)^{s-1} - 1).$$

The power series converges since $\kappa(\gamma)^{s-1} - 1 \in p\mathbb{Z}_p$. (Note: Let $t \in \mathbb{Z}_p$. The continuous group homomorphism $\kappa^t : \Gamma \to 1 + p\mathbb{Z}_p$ can be extended uniquely to a continuous $\mathbb{Z}_p$-linear ring homomorphism $\kappa^t : \Lambda \to \mathbb{Z}_p$. We have $\kappa^t(T) = \kappa(\gamma)^t - 1$ and $\kappa^t(f(T)) = f(\kappa(\gamma)^t - 1)$ for any $f(T) \in \Lambda$. Thus, $L_p(E/\mathbb{Q}, s)$ is $\kappa^{s-1}(f_E^{\mathrm{anal}}(T))$.) The functional equations for the Hasse-Weil $L$-series give a simple relation between the values $L(E/\mathbb{Q}, 1)$ and $L(E/\mathbb{Q}, \phi^{-1}, 1)$ occurring in the interpolation property for $f_E^{\mathrm{anal}}(T)$. Since $f_E^{\mathrm{anal}}(T)$ is determined by its interpolation property, one can deduce a simple relation between $f_E^{\mathrm{anal}}(T)$ and $f_E^{\mathrm{anal}}((1+T)^{-1} - 1)$. Omitting the details, one obtains a functional equation for $L_p(E/\mathbb{Q}, s)$:

$$L_p(E/\mathbb{Q}, 2 - s) = w_E \langle N_E \rangle^{s-1} L_p(E/\mathbb{Q}, s)$$

for all $s \in \mathbb{Z}_p$. Here $w_E$ is the sign which occurs in the functional equation for the Hasse-Weil $L$-series $L(E/\mathbb{Q}, s)$, $N_E$ is the conductor of $E$, and $\langle N_E \rangle$ is the projection of $N_E$ to $1 + 2p\mathbb{Z}_p$ as above.

The final theorem we will state is motivated by conjecture 1.13 and the above functional equation for the $p$-adic $L$-function $L_p(E/\mathbb{Q}, s)$. .The functional equation is in fact equivalent to the relation between $f_E^{\mathrm{anal}}(T)$ and $f_E^{\mathrm{anal}}((1+T)^{-1} - 1)$ mentioned above. In particular, $f_E^{\mathrm{anal}}(T^\iota)/f_E^{\mathrm{anal}}(T)$ should be in $\Lambda^\times$, where $T^\iota = (1 + T)^{-1} - 1$. The analogue of this statement is true for $f_E(T)$. More generally (for any $F$), we have:

**Theorem 1.14.** *Assume that $E$ is an elliptic curve defined over $F$ with good, ordinary reduction or multiplicative reduction at all primes of $F$ lying over $p$. Assume that $\mathrm{Sel}_E(F_\infty)_p$ is $\Lambda$-cotorsion. Then the characteristic ideal of $X_E(F_\infty)$ is fixed by the involution $\iota$ of $\Lambda$ induced by $\iota(\gamma) = \gamma^{-1}$ for all $\gamma \in \Gamma$.*

A proof of this result can be found in [Gr2] using the Duality Theorems of Poitou and Tate. There it is dealt with in a much more general context—that of Selmer groups attached to "ordinary" $p$-adic representations.

We will prove theorem 1.2 completely in the following two sections. Our approach is quite different than the approach in Mazur's article and in Manin's more elementary expository article. We first prove that, when $E$ has good, ordinary or multiplicative reduction at primes over $p$, the $p$-primary subgroups of $\mathrm{Sel}_E(F_n)$ and of $\mathrm{Sel}_E(F_\infty)$ have a very simple and elegant description. This is the main content of section 2. Once we have this, it is quite straightforward to prove theorem 1.2 and also a conditional result concerning conjecture 1.6 which we do in section 3. In this approach we avoid completely the need to study the norm map for formal groups over local fields, which is crucial in the approach in [Maz1] and [Man]. We also can use our description of the $p$-Selmer group to determine the $p$-adic valuation of $f_E(0)$, under the assumption that $E$ has good, ordinary reduction at primes over $p$ and that $\mathrm{Sel}_E(F)_p$ is finite. Section 4 is devoted to this comparatively easy special case of results of B. Perrin-Riou and P. Schneider found in [Pe1], [Sch1]. Their results give an expression involving a $p$-adic height determinant for the $p$-adic valuation of $(f_E(T)/T^r)|_{T=0}$, where $r = \mathrm{rank}(E(F))$, under suitable hypotheses. Finally, in section 5, (which is by far the longest section of this article) we will discuss a variety of examples to illustrate the results of sections 3 and 4 and also how our description of the $p$-Selmer group can be used for calculation. We also include in section 5 a number of remarks taken from [Maz1] (some of which are explained quite differently here) as well as various results which don't seem to be in the existing literature. Throughout this article, we have tried to include $p = 2$ in all of the main results. Perhaps surprisingly, this turns out to not be so complicated.

We will have very little to say about the case where $E$ has supersingular reduction at some primes over $p$. In recent years, this has become a very lively aspect of Iwasawa theory. We just refer the reader to [Pe4] as an introduction. In [Pe4], one finds the following concrete application of the theory described there: *Suppose that $E/\mathbb{Q}$ has supersingular reduction at $p$ and that $\mathrm{Sel}_E(\mathbb{Q})_p$ is finite. Then $\mathrm{Sel}_E(\mathbb{Q}_n)_p$ has bounded $\mathbb{Z}_p$-corank as $n$ varies.* This is, of course, a special case of conjecture 1.8. In the case where $E$ has good, ordinary reduction over $p$, theorem 1.4 gives the same conclusion. Another topic that we will not pursue is the behavior of the $p$-Selmer group in other $\mathbb{Z}_p$-extensions—for example, the anti-cyclotomic $\mathbb{Z}_p$-extension of an imaginary quadratic field. The analogues of conjectures 1.3 and 1.8 can in fact be false. We refer the reader to [Be], [BeDa1, 2], and [Maz4] for a discussion of this topic. We also will not pursue the analytic side of Iwasawa theory— questions involving the properties of $p$-adic $L$-functions and the $p$-adic version of a Birch and Swinnerton-Dyer conjecture. For this, one can learn something from the articles [M-SwD], [B-G-S], and [M-T-T]. Many of the ideas we discuss here can be extended to a far more general context. For an introduction to this, we refer the reader to [CoSc] and to [Gr2,3].

The author is grateful to the Fondazione Centro Internazionale Matematico Estivo and to Carlo Viola for the invitation to give lectures in Cetraro. This article is an extensively expanded version of those lectures, based considerably on research which was partially supported by the National Science Foundation. The author is also grateful for the support and hospitality of the American Institute of Mathematics during the Winter of 1998, when many of the results and examples described in section 5 were obtained. We want to thank Karl Rubin for many valuable discussions and for his help in the details of several examples, Ted McCabe for carrying out numerous calculations of $p$-adic $L$-functions which allowed us to verify the main conjecture in many cases, and Ken Kramer for explaining his results about elliptic curves with 2-power isogenies. We are also grateful to John Coates for many helpful remarks and to Y. Hachimori, K. Matsuno and T. Ochiai for finding a number of mistakes in the text.

## 2.    Kummer Theory for $E$.

Let $E$ be an elliptic curve defined over a number field $F$. If $M$ is any algebraic extension of $F$, Kummer theory for $E$ over $M$ leads quite naturally to the classical definition of the Selmer group $\mathrm{Sel}_E(M)$. The main objective of this section is to give a simplified description of its $p$-primary subgroup $\mathrm{Sel}_E(M)_p$ under the hypothesis that $E$ has either good, ordinary reduction or multiplicative reduction at all primes of $F$ lying over $p$. We will assume that $M$ is either a finite extension or a $\mathbb{Z}_p$-extension of $F$.

Kummer theory for the multiplicative group $M^\times$ is quite familiar. Regarding $M$ as a subfield of $\overline{F}$, a fixed algebraic closure of $F$ (or $\mathbb{Q}$), we can define the Kummer homomorphism

$$k : M^\times \otimes (\mathbb{Q}/\mathbb{Z}) \to H^1(M, \overline{F}^\times_{\mathrm{tors}})$$

as follows. Let $a \in M^\times$. Let $\alpha = a \otimes (m/n + \mathbb{Z}) \in M^\times \otimes (\mathbb{Q}/\mathbb{Z})$. Choose $b \in \overline{F}^\times$ such that $b^n = a^m$, using the fact that $\overline{F}^\times$ is a divisible group. Then one defines $k(\alpha)$ to be the class of the 1-cocycle $\phi_\alpha$ given by $\phi_\alpha(g) = g(b)/b$ for all $g \in G_M = \mathrm{Gal}(\overline{F}/M)$. The values of $\phi_\alpha$ are in $\overline{F}^\times_{\mathrm{tors}}$, the group of roots of unity in $\overline{F}$. The Kummer homomorphism is an isomorphism. Injectivity is easy to verify. Surjectivity is a consequence of Hilbert's Theorem 90, which asserts that $H^1(M, \overline{F}^\times) = 0$.

Since $E(\overline{F})$ is divisible, one can imitate the above definition, obtaining an exact sequence

$$0 \to E(M) \otimes (\mathbb{Q}/\mathbb{Z}) \xrightarrow{k} H^1(M, E(\overline{F})_{\mathrm{tors}}) \to H^1(M, E(\overline{F})) \to 0.$$

If $\alpha = a \otimes (m/n + \mathbb{Z}) \in E(M) \otimes (\mathbb{Q}/\mathbb{Z})$, then $k(\alpha)$ is the class of the 1-cocycle $\phi_\alpha$ given by $\phi_\alpha(g) = g(b) - b$ for all $g \in G_M$. Here $b \in E(\overline{F})$ satisfies $nb = ma$ on $E(\overline{F})$. However, in general, $H^1(M, E(\overline{F}))$ is nonzero. We will fix a prime $p$ and concentrate on the $p$-primary subgroups of the above groups. We let $\kappa = \kappa_M$ denote the corresponding Kummer homomorphism:

$$\kappa : E(M) \otimes (\mathbb{Q}_p/\mathbb{Z}_p) \to H^1(M, E[p^\infty]).$$

If $\eta$ is any prime of $M$, we define $M_\eta$ to be the union of the $\eta$-adic completions of all finite extensions of $F$ contained in $M$. Thus, if $\eta$ lies over the prime $v$ of $F$, then $M_\eta$ is an algebraic extension of $F_v$. By fixing an embedding $\overline{F} \to \overline{F}_v$ extending the embedding $M \to M_\eta$, one can identify $G_{M_\eta}$ with a subgroup of $G_M$, which of course is just the decomposition subgroup for some prime of $\overline{F}$ lying over $\eta$. We will let $\kappa_\eta$ denote the Kummer homomorphism for $E$ over $M_\eta$:

$$\kappa_\eta : E(M_\eta) \otimes (\mathbb{Q}_p/\mathbb{Z}_p) \to H^1(M_\eta, E[p^\infty]).$$

This is defined exactly as above. Now we can give the classical definition of the $p$-primary subgroup of the Selmer group for $E$ over $M$.

$$\mathrm{Sel}_E(M)_p = \ker\big(H^1(M, E[p^\infty]) \to \prod_\eta H^1(M_\eta, E[p^\infty])/\mathrm{Im}(\kappa_\eta)\big)$$

where $\eta$ runs over all primes of $M$ and the map is induced by $\phi \to (\phi\big|_{G_\eta})_\eta$ for any 1-cocycle $\phi$. We will denote the class of a 1-cocycle $\phi$ by $[\phi]$. Thus $[\phi]$ is in $\mathrm{Sel}_E(M)_p$ if and only if $[\phi\big|_{G_{M_\eta}}] \in \mathrm{Im}(\kappa_\eta)$ for all $\eta$. Obviously, $\mathrm{Im}(\kappa) \subseteq \mathrm{Sel}_E(M)_p$. The corresponding quotient $\mathrm{Sel}_E(M)_p/\mathrm{Im}(\kappa)$ is, by definition, $\mathrm{III}_E(M)_p$.

Faltings has proved that $E$ is determined up to $F$-isogeny by the $G_F$-representation space $V_p(E) = T_p(E) \otimes \mathbb{Q}_p$, where $T_p(E)$ denotes the $p$-adic Tate module for $E$. More precisely, the $G_F$-module $E[p^\infty] \cong V_p(E)/T_p(E)$ determines $E$ up to an $F$-isogeny of degree prime to $p$. Now $\mathrm{Sel}_E(M)_p$ is not changed by such $F$-isogenies, and hence one might hope to define it in a way which involves only the $G_F$-module $E[p^\infty]$. To do this, it suffices to give such a description of the subgroup $\mathrm{Im}(\kappa_\eta)$ of $H^1(M_\eta, E[p^\infty])$ for all primes $\eta$ of $M$. We will now proceed to do this under the assumption that $E$ has good, ordinary or multiplicative reduction at all primes of $F$ over $p$.

Assume at first that $M$ is a finite extension of $F$. Then $\eta | v$ for some prime $v$ of $F$, and $\eta | l$ for some prime $l$ of $\mathbb{Q}$ (possible $l = \infty$). If $l$ is a finite prime, then we have a theorem of Lutz: $E(M_\eta) \cong \mathbb{Z}_l^{[M_\eta:\mathbb{Q}_l]} \times U$ as a group, where $U = E(M_\eta)_{\mathrm{tors}}$ is finite. Now $\mathbb{Z}_l \otimes (\mathbb{Q}_p/\mathbb{Z}_p) = 0$ if $l \neq p$, whereas $\mathbb{Z}_p \otimes (\mathbb{Q}_p/\mathbb{Z}_p) \cong \mathbb{Q}_p/\mathbb{Z}_p$. Also, $U \otimes (\mathbb{Q}_p/\mathbb{Z}_p) = 0$. If $l = \infty$, then $M_\eta \cong \mathbb{R}$ or $\mathbb{C}$. In this case, $E(M_\eta) \cong T^{[M_\eta:\mathbb{R}]} \times U$, where $T = \mathbb{R}/\mathbb{Z}$ and $|U| \leq 2$. Since $T$ is divisible, we have $T \otimes (\mathbb{Q}_p/\mathbb{Z}_p) = 0$. We then obtain the following result.

**Proposition 2.1.** *If $\eta \nmid p$, then* $\mathrm{Im}(\kappa_\eta) = 0$. *If $\eta | p$, then*

$$\mathrm{Im}(\kappa_\eta) \cong (\mathbb{Q}_p/\mathbb{Z}_p)^{[M_\eta : \mathbb{Q}_p]}.$$

The first assertion can also be explained by using the fact that, for $\eta \nmid p$, $H^1(M_\eta, E[p^\infty])$ is a finite group. But $E(M_\eta) \otimes (\mathbb{Q}_p/\mathbb{Z}_p)$, and hence $\mathrm{Im}(\kappa_\eta)$ are divisible groups. Even if $M_\eta$ is an infinite extension of $F_v$, it is clear from the above that $\mathrm{Im}(\kappa_\eta) = 0$ if $\eta \nmid p$.

Assume that $E$ has good, ordinary reduction at $v$, where $v$ is a prime of $F$ lying over $p$. Then, considering $E[p^\infty]$ as a subgroup of $E(\overline{F}_v)$, we have the reduction map $E[p^\infty] \rightarrow \widetilde{E}[p^\infty]$, where $\widetilde{E}$ is the reduction of $E$ modulo $v$. Define $C_v$ by

$$C_v = \ker\left(E[p^\infty] \rightarrow \widetilde{E}[p^\infty]\right).$$

Now $E[p^\infty] \cong (\mathbb{Q}_p/\mathbb{Z}_p)^2$, $\widetilde{E}[p^\infty] \cong \mathbb{Q}_p/\mathbb{Z}_p$ as groups. It is easy to see that $C_v \cong \mathbb{Q}_p/\mathbb{Z}_p$. (In fact, $C_v = \mathcal{F}(\overline{\mathfrak{m}})[p^\infty]$, where $\mathcal{F}$ is the formal group of height 1 for $E$ and $\overline{\mathfrak{m}}$ is the maximal ideal of the integers of $\overline{F}_v$.) A characterization in terms of $E[p^\infty]$ is that $C_v$ is $G_{F_v}$-invariant and $E[p^\infty]/C_v$ is the maximal unramified quotient of $E[p^\infty]$. Let $M$ be a finite extension of $F$. If $\eta$ is a residue of $M$ lying above $v$, then we can consider $M_\eta$ as a subfield of $\overline{F}_v$ containing $F_v$. (The identification will not matter.) We then have a natural map

$$\lambda_\eta : H^1(M_\eta, C_v) \rightarrow H^1(M_\eta, E[p^\infty]).$$

Here is a description of $\mathrm{Im}(\kappa_\eta)$.

**Proposition 2.2.** $\mathrm{Im}(\kappa_\eta) = \mathrm{Im}(\lambda_\eta)_{\mathrm{div}}$.

*Proof.* The idea is quite simple. We know that $\mathrm{Im}(\kappa_\eta)$ and $\mathrm{Im}(\lambda_\eta)$ are $p$-primary groups, that $\mathrm{Im}(\kappa_\eta)$ is divisible, and has $\mathbb{Z}_p$-corank $[M_\eta : \mathbb{Q}_p]$. It suffices to prove two things: (i) $\mathrm{Im}(\kappa_\eta) \subseteq \mathrm{Im}(\lambda_\eta)$ and (ii) $\mathrm{Im}(\lambda_\eta)$ has $\mathbb{Z}_p$-corank equal to $[M_\eta : \mathbb{Q}_p]$. To prove (i), let $c \in \mathrm{Im}(\kappa_\eta)$. We show that $c \in \ker(H^1(M_\eta, E[p^\infty]) \rightarrow H^1(M_\eta, \widetilde{E}[p^\infty]))$, which coincides with $\mathrm{Im}(\lambda_\eta)$. Let $f_v$ denote the residue field of $F_v$, $\overline{f}_v$ its algebraic closure—the residue field of $\overline{F}_v$. If $b \in E(\overline{F}_v)$, we let $\widetilde{b} \in \widetilde{E}(\overline{f}_v)$ denote its reduction. Let $\phi$ be a cocycle representing $c$. Then $\phi(g) = g(b) - b$ for all $g \in G_{M_\eta}$, where $b \in E(\overline{F}_v)$. The 1-cocycle induced by $E[p^\infty] \rightarrow \widetilde{E}[p^\infty]$ is $\widetilde{\phi}$, given by $\widetilde{\phi}(g) = g(\widetilde{b}) - \widetilde{b}$ for all $g \in G_{M_\eta}$. But $\widetilde{\phi}$ represents a class $\widetilde{c}$ in $H^1(M_\eta, \widetilde{E}[p^\infty])$ which becomes trivial in $H^1(M_\eta, \widetilde{E}(\overline{f}_v))$, i.e. $\widetilde{\phi}$ is a 1-coboundary. Finally, the key point is that $\widetilde{E}(\overline{f}_v)$ is a torsion group, $\widetilde{E}[p^\infty]$ is its $p$-primary subgroup, and hence the map $H^1(M_\eta, \widetilde{E}[p^\infty]) \rightarrow H^1(M_\eta, \widetilde{E}(\overline{f}_v))$ must be injective. Thus, $\widetilde{c}$ is trivial, and therefore $c \in \mathrm{Im}(\lambda_\eta)$.

Now we calculate the $\mathbb{Z}_p$-corank of $\mathrm{Im}(\lambda_\eta)$. We have the exact sequence

$$E[p^\infty]^{G_{M_\eta}} \rightarrow \widetilde{E}[p^\infty]^{G_{M_\eta}} \rightarrow H^1(M_\eta, C_v) \xrightarrow{\lambda_\eta} H^1(M_\eta, E[p^\infty]).$$

If $m_\eta$ denotes the residue field of $M_\eta$, then $\widetilde{E}[p^\infty]^{G_{M_\eta}}$ is just the $p$-primary subgroup of $\widetilde{E}(m_\eta)$, a finite group. Thus, $\ker(\lambda_\eta)$ is finite. The following lemma then suffices to prove (ii). If $\psi : G_{F_v} \to \mathbb{Z}_p^\times$ is a continuous homomorphism, we will let $(\mathbb{Q}_p/\mathbb{Z}_p)(\psi)$ denote the group $\mathbb{Q}_p/\mathbb{Z}_p$ together with the action of $G_{F_v}$ given by $\psi$.

**Lemma 2.3.** $H^1(M_\eta, (\mathbb{Q}_p/\mathbb{Z}_p)(\psi))$ has $\mathbb{Z}_p$-corank equal to $[M_\eta : \mathbb{Q}_p] + \delta$, where $\delta = 1$ if $\psi|_{G_{M_\eta}}$ is either the trivial character or the cyclotomic character of $G_{M_\eta}$ and $\delta = 0$ otherwise.

*Remark.* Because of the importance of this lemma, we will give a fairly self-contained proof using local class field theory and techniques of Iwasawa Theory. But we then show how to obtain the same result as a simple application of the Duality theorems of Poitou and Tate.

*Proof.* The case where $\psi$ is trivial follows from local class field theory. Then $H^1(M_\eta, (\mathbb{Q}_p/\mathbb{Z}_p)(\psi)) = \mathrm{Hom}(\mathrm{Gal}(M_\eta^{ab}/M_\eta), \mathbb{Q}_p/\mathbb{Z}_p)$. The well-known structure of $M_\eta^\times$ implies that $\mathrm{Gal}(M_\eta^{ab}/M_\eta) \cong \mathbb{Z}_p^{[M_\eta : \mathbb{Q}_p]} \times \widehat{\mathbb{Z}} \times (M_\eta^\times)_{\mathrm{tors}}$, where $\widehat{\mathbb{Z}}$ is the profinite completion of $\mathbb{Z}$. The lemma is clear in this case. If $\psi|_{G_{M_\eta}}$ is the cyclotomic character, then $(\mathbb{Q}_p/\mathbb{Z}_p)(\psi) \cong \mu_{p^\infty}$ as $G_{M_\eta}$-modules. Then $H^1(M_\eta, \mu_{p^\infty}) \cong (M_\eta^\times) \otimes (\mathbb{Q}_p/\mathbb{Z}_p)$, which indeed has the stated $\mathbb{Z}_p$-corank.

Now suppose we are not in one of the above two cases. For brevity, we will write $M$ for $M_\eta$. Let $M_\infty$ be the extension of $M$ cut out by $\psi|_{G_M}$. Thus, $G = \mathrm{Gal}(M_\infty/M) \cong \mathrm{Im}(\psi|_{G_M})$. If $\psi$ has finite order, one can reduce to studying the action of $G$ on $\mathrm{Gal}(M_\infty^{ab}/M_\infty)$ since $M_\infty$ would just be a finite extension of $\mathbb{Q}_p$. We will do something similar if $\psi$ has infinite order. Then, $G \cong \Delta \times H$, where $\Delta$ is finite and $H \cong \mathbb{Z}_p$. If $p$ is odd, $|\Delta|$ divides $p - 1$. If $p = 2$, $|\Delta| = 1$ or $2$. Let $C = (\mathbb{Q}_p/\mathbb{Z}_p)(\psi)$. The inflation-restriction sequence gives

$$0 \to H^1(G, C) \to H^1(M, C) \to H^1(M_\infty, C)^G \to H^2(G, C).$$

Now let $h$ be a topological generator of $H$. Then $H^1(H, C) = C/(h-1)C = 0$ because, considering $h - 1$ as an endomorphism of $C$, $\ker(h-1)$ is finite and $\mathrm{Im}(h - 1)$ is divisible. Thus, $H^1(G, C) = 0$ if $p$ is odd, and has order $\leq 2$ if $p = 2$. On the other hand, $H^2(H, C) = 0$ since $H$ has $p$-cohomological dimension 1. Then $H^2(G, C) = 0$ if $p$ is odd, and again has order $\leq 2$ if $p = 2$. Thus, it is enough to study

$$H^1(M_\infty, C)^G = \mathrm{Hom}_G(\mathrm{Gal}(M_\infty^{ab}/M_\infty), C).$$

Let $X = \mathrm{Gal}(L_\infty/M_\infty)$, where $L_\infty$ is the maximal abelian pro-$p$ extension of $M_\infty$. We will prove the rest of lemma 2.3 by studying the structure of $X$ as a module for $\mathbb{Z}_p[[\Delta \times H]] = \Lambda[\Delta]$, where $\Lambda = \mathbb{Z}_p[[H]] \cong \mathbb{Z}_p[[T]]$, with $T = h - 1$. The results are due to Iwasawa.

For any $n \geq 0$, let $H_n = H^{p^n}$. Let $M_n = M_\infty^{H_n}$. The commutator subgroup of $\mathrm{Gal}(L_\infty/M_n)$ is $(h^{p^n}-1)X$ and so, if $L_n$ is the maximal abelian extension of $M_n$ contained in $L_\infty$, then $\mathrm{Gal}(L_n/M_n) \cong H_n \times (X/(h^{p^n}-1)X)$. But $L_n$ is the maximal abelian pro-$p$ extension of $M_n$ and, by local class field theory, this Galois group is isomorphic to $\mathbb{Z}_p^{[M_n:\mathbb{Q}_p]+1} \times W_n$, where $W_n$ denotes the group of $p$-power roots of unity contained in $M_n$. Consequently, if we put $t = [M_0 : \mathbb{Q}_p] = |\Delta| \cdot [M : \mathbb{Q}_p]$, we have

$$X/(h^{p^n}-1)X \cong \mathbb{Z}_p^{tp^n} \times W_n.$$

Now, the structure theory for $\Lambda$-modules states that $X/X_{\Lambda\text{-tors}}$ is isomorphic to a submodule of $\Lambda^r$, with finite index, where $r = \mathrm{rank}_\Lambda(X)$. Also, we have $\Lambda/(h^{p^n}-1)\Lambda \cong \mathbb{Z}_p^{p^n}$ for $n \geq 0$. It follows that $r = t$. One can also see that $X_{\Lambda\text{-tors}} \cong \varprojlim W_n$, where this inverse limit is defined by the norm maps $M_m^\times \to M_n^\times$ for $m \geq n$. If $W_n$ has bounded order (i.e., if $\mu_{p^\infty} \not\subseteq M_\infty$), then $X_{\Lambda\text{-tors}} = 0$. Thus, $X \subseteq \Lambda^t$. To get more precise information about the structure of $X$, choose $n$ large enough so that $h^{p^n}-1$ annihilates $\Lambda^t/X$. We then have

$$(h^{p^n}-1)X \subseteq (h^{p^n}-1)\Lambda^t \subseteq X \subseteq \Lambda^t.$$

We can see easily from this that $\Lambda^t/X$ is isomorphic to the torsion subgroup of $X/(h^{p^n}-1)X$. That is, $\Lambda^t/X \cong W$, where $W = M_\infty^\times \cap \mu_{p^\infty}$. On the other hand, if $\mu_{p^\infty} \subseteq M_\infty$, then $X_{\Lambda\text{-tors}} \cong \mathbb{Z}_p(1)$, the Tate module for $\mu_{p^\infty}$. In this case, $X/X_{\Lambda\text{-tors}}$ is free and hence $X \cong \Lambda^t \times \mathbb{Z}_p(1)$.

In the preceding discussion, the $\Lambda$-module $\Lambda^t$ is in fact canonical. It is the reflexive hull of $X/X_{\Lambda\text{-tors}}$. Thus, the action of $\Delta$ on $X$ gives an action on $\Lambda^t$. Examining the above arguments more carefully, one finds that, for $p$ odd, $\Lambda^t$ is isomorphic to $\Lambda[\Delta]^{[M:\mathbb{Q}_p]}$. (One just studies the $\Lambda$-module $X^\phi$ for each character $\phi$ of $\Delta$. Recall that $|\Delta|$ divides $p-1$ and hence each character $\phi$ has values in $\mathbb{Z}_p^\times$.) For $p = 2$, we can at least make such an identification up to a group of exponent 2. For the proof of lemma 2.3, it suffices to point out that $\mathrm{Hom}_{\Delta \times H}(\Lambda[\Delta], C)$ is isomorphic to $\mathbb{Q}_p/\mathbb{Z}_p$ and that $\mathrm{Hom}_{\Delta \times H}(\mathbb{Z}_p(1), C)$ is finite. (We are assuming now that $C \not\cong \mu_{p^\infty}$ as $G_M$-modules.) This completes the proof of lemma 2.3 and consequently proposition 2.2, since one sees easily that $\delta = 0$ when $C = C_v$. ∎

The above discussion of the $\Lambda[\Delta]$-module structure of $X$ gives a more precise result concerning $H^1(M_\eta, (\mathbb{Q}_p/\mathbb{Z}_p)(\psi))$. Assume that $p$ is odd and that $\psi$ has infinite order. If the extension of $M_\eta$ cut out by the character $\psi$ of $G_{M_\eta}$ contains $\mu_{p^\infty}$, then we see that

$$H^1(M_\eta, C) \cong (\mathbb{Q}_p/\mathbb{Z}_p)^{[M_\eta:\mathbb{Q}_p]} \times \mathrm{Hom}_{G_{M_\eta}}(\mathbb{Z}_p(1), C), \qquad (1)$$

where as above $C = (\mathbb{Q}_p/\mathbb{Z}_p)(\psi)$. The factor $\mathrm{Hom}_{G_{M_\eta}}(\mathbb{Z}_p(1), C)$ is just $H^0(M_\eta, C \otimes \chi^{-1})$, where $\chi$ denotes the cyclotomic character. Even if $W$ is finite, we can prove (1). For if $g_0$ is a topological generator of $\Delta \times H$, then the torsion subgroup of $X/(g_0 - \psi(g_0))X$ is isomorphic to the kernel of $g_0 - \psi(g_0)$ acting on $\Lambda^t/X \cong W$. (It is seen to be $((g_0 - \psi(g_0))\Lambda^t \cap X)/(g_0 - \psi(g_0))X$.) But this in turn is isomorphic to $W/(g_0 - \psi(g_0))W$, whose dual is easily identified with $\mathrm{Hom}_{G_{M_\eta}}(\mathbf{Z}_p(1), (\mathbf{Q}_p/\mathbf{Z}_p)(\psi))$.

We have attempted to give a rather self-contained "Iwasawa-theoretic" approach to studying the above local Galois cohomology group. This suffices for the proof of proposition 2.2. But using results of Poitou and Tate is often easier and more effective. Let us illustrate this. Let $C = (\mathbf{Q}_p/\mathbf{Z}_p)(\psi)$. Let $T$ denote its Tate module and $V = T \otimes_{\mathbf{Z}_p} \mathbf{Q}_p$. The $\mathbf{Z}_p$-corank of $H^1(G_{M_\eta}, C)$ is just $\dim_{\mathbf{Q}_p}(H^1(M_\eta, V))$. (Cocycles are required to be continuous. $V$ has its $\mathbf{Q}_p$-vector space topology. Similarly, $T$ has its natural topology and is compact.) Letting $h_i$ denote $\dim_{\mathbf{Q}_p}(H^i(M_\eta, V))$, then the Euler characteristic for $V$ over $M_\eta$ is given by

$$h_0 - h_1 + h_2 = -[M_\eta : \mathbf{Q}_p] \dim_{\mathbf{Q}_p}(V)$$

for any $G_{M_\eta}$-representation space $V$. We have $\dim_{\mathbf{Q}_p}(V) = 1$ and so the $\mathbf{Z}_p$-corank of $H^1(M_\eta, (\mathbf{Q}_p/\mathbf{Z}_p)(\psi))$ is $[M_\eta : \mathbf{Q}_p] + h_0 + h_2$. Poitou-Tate Duality implies that $H^2(M_\eta, V)$ is dual to $H^0(M_\eta, V^*)$, where $V^* = \mathrm{Hom}(V, \mathbf{Q}_p(1))$. It is easy to see from this that $\delta = h_0 + h_2$, proving lemma 2.3 again.

The exact sequence $0 \to T \to V \to C \to 0$ induces the exact sequence

$$H^1(M_\eta, V) \xrightarrow{\alpha} H^1(M_\eta, C) \xrightarrow{\beta} H^2(M_\eta, T) \xrightarrow{\gamma} H^2(M_\eta, V).$$

The image of $\alpha$ is the maximal divisible subgroup of $H^1(G_{M_\eta}, C)$. The kernel of $\gamma$ is the torsion subgroup of $H^2(M_\eta, T)$. Of course, $\mathrm{coker}(\alpha) \cong \mathrm{Im}(\beta) \cong \ker(\gamma)$. Poitou-Tate Duality implies that $H^2(M_\eta, T)$ is dual to $H^0(M_\eta, \mathrm{Hom}(T, \mu_{p^\infty})) = \mathrm{Hom}_{G_{M_\eta}}(T, \mu_{p^\infty})$. The action of $G_{M_\eta}$ on $T$ is by $\psi$; the action on $\mu_{p^\infty}$ is by $\chi$. Thus, $\mathrm{Hom}_{G_{M_\eta}}(T, \mu_{p^\infty})$ can be identified with the dual of $H^0(M_\eta, (\mathbf{Q}_p/\mathbf{Z}_p)(\chi\psi^{-1}))$. If $\psi|_{G_{M_\eta}} = \chi|_{G_{M_\eta}}$, then we find that $H^2(M_\eta, T) \cong \mathbf{Z}_p$, $\mathrm{Im}(\beta) = 0$, and therefore $H^1(M_\eta, C)$ is divisible. Otherwise, we find that $H^2(M_\eta, T)$ is finite and that

$$H^1(M_\eta, C)/H^1(M_\eta, T)_{\mathrm{div}} \cong (\mathbf{Q}_p/\mathbf{Z}_p)(\psi\chi^{-1})^{G_{M_\eta}}, \qquad (2)$$

which is a finite cyclic group, indeed isomorphic to $\mathrm{Hom}_{G_{M_\eta}}(\mathbf{Z}_p(1), C)$. This argument works even for $p = 2$.

We want to mention here one useful consequence of the above discussion. Again we let $C = (\mathbf{Q}_p/\mathbf{Z}_p)(\psi)$, where $\psi : G_{F_v} \to \mathbf{Z}_p^\times$ is a continuous homomorphism, $v$ is any prime of $F$ lying over $p$. If $\eta$ is a prime of $F_\infty$ lying over $v$, then $(F_\infty)_\eta$ is the cyclotomic $\mathbf{Z}_p$-extension of $F_v$. By lemma 2.3, the $\mathbf{Z}_p$-corank of $H^1((F_n)_\eta, C)$ differs from $[(F_n)_\eta : F_v]$ by at most 1. Thus, if we let $\Gamma_v = \mathrm{Gal}((F_\infty)_\eta/F_v)$, then it follows that as $n \to \infty$

$$\mathrm{corank}_{\mathbf{Z}_p}(H^1((F_\infty)_\eta, C)^{\Gamma_v^{p^n}}) = p^n[F_v : \mathbf{Q}_p] + O(1).$$

The structure theory of $\Lambda$-modules then implies that $H^1((F_\infty)_\eta, C)$ has corank equal to $[F_v : \mathbb{Q}_p]$ as a $\mathbb{Z}_p[[\Gamma_v]]$-module. Assume that $\psi$ is unramified and that the maximal unramified extension of $F_v$ contains no $p$-th roots of unity. (If the ramification index $e_v$ for $v$ over $p$ is $\le p-2$, then this will be true. If $F = \mathbb{Q}$, this is true for all $p \ge 3$.) Then by (2) we see that $H^1(F_v, C)$ is divisible. The $\mathbb{Z}_p$-corank of $H^1(F_v, C)$ is $[F_v : \mathbb{Q}_p] + \delta$, where $\delta = 0$ if $\psi$ is nontrivial, $\delta = 1$ if $\psi$ is trivial. By the inflation-restriction sequence we see that $H^1((F_\infty)_\eta, C)^{\Gamma_v} \cong (\mathbb{Q}_p/\mathbb{Z}_p)^{[F_v:\mathbb{Q}_p]}$. It follows that $H^1((F_\infty)_\eta, C)$ is $\mathbb{Z}_p[[\Gamma_v]]$-cofree of corank $[F_v : \mathbb{Q}_p]$, under the hypotheses that $\psi$ is unramified and $e_v \le p-2$. These remarks are a special case of results proved in [Gr2].

Now we return to the case where $C_v = \ker(E[p^\infty] \to \widetilde{E}[p^\infty])$. The action of $G_{F_v}$ on $C_v$ is by a character $\psi$, the action on $\widetilde{E}[p^\infty]$ is by a character $\phi$, and we have $\psi\phi = \chi$ since the Weil pairing $T_p(E) \wedge T_p(E) \cong \mathbb{Z}_p(1)$ means that $\chi$ is the determinant of the representation of $G_{F_v}$ on $T_p(E)$. Note that $\phi$ has infinite order. The same is true for $\psi$ since $\psi$ and $\chi$ become equal after restriction to the inertia subgroup $G_{F_v^{\text{unr}}}$. This explains why $\delta = 0$ for $\psi|_{G_{M_\eta}}$, as used to prove proposition 2.2. In this case, $\chi\psi^{-1} = \phi$ and hence $H^0(G_{M_\eta}, (\mathbb{Q}_p/\mathbb{Z}_p)(\chi\psi^{-1}))$ is isomorphic to $\widetilde{E}(m_\eta)_p$, where $m_\eta$ is the residue field for $M_\eta$. These facts lead to a version of proposition 2.2 for some infinite extensions of $F_v$.

**Proposition 2.4.** *Assume that $K$ is a Galois extension of $F_v$, that $\mathrm{Gal}(K/F_v)$ contains an infinite pro-$p$ subgroup, and that the inertia subgroup of $\mathrm{Gal}(K/F_v)$ is of finite index. Then $\mathrm{Im}(\kappa_K) = \mathrm{Im}(\lambda_K)$, where $\kappa_K$ is the Kummer homomorphism for $E$ over $K$ and $\lambda_K$ is the canonical homomorphism*

$$H^1(K, C_v) \to H^1(K, E[p^\infty]).$$

*Proof.* Let $M$ run over all finite extensions of $F_v$ contained in $K$. Then $\mathrm{Im}(\kappa_K) = \varinjlim \mathrm{Im}(\kappa_M)$, $\mathrm{Im}(\lambda_K) = \varinjlim \mathrm{Im}(\lambda_M)$, and $\mathrm{Im}(\kappa_M) = \mathrm{Im}(\lambda_M)_{\text{div}}$ by proposition 2.2. But $\mathrm{Im}(\lambda_M)/\mathrm{Im}(\lambda_M)_{\text{div}}$ has order bounded by $|\widetilde{E}(m)_p|$, where $m$ is the residue field of $M$. Now $|m|$ is bounded by assumption. Hence it follows that $\mathrm{Im}(\lambda_K)/\mathrm{Im}(\kappa_K)$ is a finite group. On the other hand, $G_K$ has $p$-cohomological dimension 1 because of the hypothesis that $\mathrm{Gal}(K/F_v)$ contains an infinite pro-$p$ subgroup. (See Serre, Cohomologie Galoisienne, Chapitre II, §3.) Thus if $C$ is a divisible, $p$-primary $G_K$-module, then the exact sequence $0 \to C[p] \to C \xrightarrow{p} C \to 0$ induces the cohomology exact sequence $H^1(K, C) \xrightarrow{p} H^1(K, C) \to H^2(K, C[p])$. The last group is zero and hence $H^1(K, C)$ is divisible. Applying this to $C = C_v$, we see that $\mathrm{Im}(\lambda_K)$ is divisible and so $\mathrm{Im}(\kappa_K) = \mathrm{Im}(\lambda_K)$. ∎

If $F_\infty$ denotes the cyclotomic $\mathbb{Z}_p$-extension of $F$, then every prime $v$ of $F$ lying over $p$ is ramified in $F_\infty/F$. If $\eta$ is a prime of $F_\infty$ over $v$, then $K = (F_\infty)_\eta$ satisfies the hypothesis of proposition 2.4 since the inertia subgroup of $\Gamma = \mathrm{Gal}(F_\infty/F)$ for $\eta$ is infinite, pro-$p$, and has finite index in $\Gamma$. Propositions 2.1, 2.2, and 2.4 will allow us to give a fairly straightforward proof of theorem 1.2, which we will do in section 3. However, in section 4 it will be useful to have more precise information about $\mathrm{Im}(\lambda_\eta)/\mathrm{Im}(\kappa_\eta)$, where $\eta$ is a prime for a finite extension $M$ of $F$ lying over $p$. What we will need is the following.

**Proposition 2.5.** *Let $M_\eta$ be a finite extension of $F_v$, where $v|p$. Let $m_\eta$ be the residue field for $M_\eta$. Then*

$$\mathrm{Im}(\lambda_\eta)/\mathrm{Im}(\kappa_\eta) \cong \widetilde{E}(m_\eta)_p.$$

*Proof.* The proof comes out of the following diagram:

$$
\begin{array}{ccccccccc}
0 & \longrightarrow & \mathcal{F}(\mathfrak{m}) \otimes (\mathbb{Q}_p/\mathbb{Z}_p) & \xrightarrow{\kappa_{\mathcal{F}}} & H^1(M_\eta, C_v) & \longrightarrow & H^1(M_\eta, \mathcal{F}(\overline{\mathfrak{m}}))_p & \longrightarrow & 0 \\
& & \downarrow & & \downarrow{\scriptstyle \lambda_\eta} & & \downarrow{\scriptstyle \epsilon} & & \\
0 & \longrightarrow & E(M_\eta) \otimes (\mathbb{Q}_p/\mathbb{Z}_p) & \xrightarrow{\kappa_\eta} & H^1(M_\eta, E[p^\infty]) & \longrightarrow & H^1(M_\eta, \widetilde{E}(\overline{F}_v))_p & \longrightarrow & 0
\end{array}
$$

Here $\mathcal{F}$ is the formal group for $E$ (which has height 1), $\mathfrak{m}$ is the maximal ideal of $M_\eta$. The upper row is the Kummer sequence for $\mathcal{F}(\mathfrak{m})$, based on the fact that $\mathcal{F}(\overline{\mathfrak{m}})$ is divisible. The first vertical arrow is surjective since $\mathcal{F}(\mathfrak{m})$ has finite index in $E(M_\eta)$. Comparing $\mathbb{Z}_p$-coranks, one sees that $\mathrm{Im}(\kappa_{\mathcal{F}}) = H^1(G_{M_\eta}, C_v)_{\mathrm{div}}$. A simple diagram chase shows that the map

$$H^1(M_\eta, C_v)/H^1(M_\eta, C_v)_{\mathrm{div}} \longrightarrow \mathrm{Im}(\lambda_\eta)/\mathrm{Im}(\lambda_\eta)_{\mathrm{div}} \qquad (3)$$

is surjective and has kernel isomorphic to $\ker(\epsilon)$. The exact sequence

$$0 \to \mathcal{F}(\overline{\mathfrak{m}}) \to E(\overline{F}_v) \to \widetilde{E}(\overline{f}_v) \to 0,$$

together with the fact that the reduction map $E(M_\eta) \to \widetilde{E}(m_\eta)$ is surjective implies that $\epsilon$ is injective. (For the surjectivity of the reduction map, see proposition 2.1 of [Si].) Therefore, the map (3) is an isomorphism. Combining this with the observation preceding proposition 2.4, we get the stated conclusion. ∎

Assume now that $E$ has split, multiplicative reduction at $v$. Then one has an exact sequence

$$0 \to C_v \to E[p^\infty] \to \mathbb{Q}_p/\mathbb{Z}_p \to 0$$

where $C_v \cong \mu_{p^\infty}$. The proof of proposition 2.2 can be made to work and gives the following result. *For any algebraic extension $K$ of $F_v$, we have $Im(\kappa_K) = Im(\lambda_K)$.* It is enough to prove this when $[K:F_v] < \infty$. Then $\mathrm{Im}(\kappa_K)$ is divisible and has $\mathbb{Z}_p$-corank $[K:\mathbb{Q}_p]$. $H^1(K, C_v)$ is divisible and has $\mathbb{Z}_p$-corank $[K:\mathbb{Q}_p]+1$. But the kernel of $\lambda_K: H^1(K, C_v) \to H^1(K, E[p^\infty])$ is isomorphic to $\mathbb{Q}_p/\mathbb{Z}_p$. Thus, $\mathrm{Im}(\lambda_K)$ and $\mathrm{Im}(\kappa_K)$ are both divisible and have the same $\mathbb{Z}_p$-corank. The inclusion $\mathrm{Im}(\kappa_K) \subseteq \mathrm{Im}(\lambda_K)$ can be seen by noting that in defining $\kappa_K$, one can assume that $\alpha \in E(K) \otimes (\mathbb{Q}_p/\mathbb{Z}_p)$ has been written as $\alpha = a \otimes (1/p^t)$, where $a \in \mathcal{F}(\mathfrak{m})$. Here $\mathcal{F}$ is the formal group for $E$, $\mathfrak{m}$ is the maximal ideal for $K$. Then, since $\mathcal{F}(\overline{\mathfrak{m}})$ is divisible, one can choose $b \in \mathcal{F}(\overline{\mathfrak{m}})$ so that $p^t b = a$. The 1-cocycle $\phi_\alpha$ then has values in $C_v = \mathcal{F}(\overline{\mathfrak{m}})[p^\infty]$. Alternatively, the equality $\mathrm{Im}(\kappa_K) = \mathrm{Im}(\lambda_K)$ can be verified quite directly by using the Tate parametrization for $E$.

If $E$ has nonsplit, multiplicative reduction, then the above assertion still holds for $p$ odd. That is, $\mathrm{Im}(\kappa_K) = \mathrm{Im}(\lambda_K)$ for every algebraic extension $K$ of $F_v$. We can again assume that $[K:F_v] < \infty$. If $E$ becomes split over $K$, then the argument in the preceding paragraph applies. If not, then lemma 2.3 and (2) imply that $H^1(K, C_v)$ is divisible and has $\mathbb{Z}_p$-corank $[K:\mathbb{Q}_p]$. Just as in the case of good, ordinary reduction, we see that $\mathrm{Im}(\kappa_K) = \mathrm{Im}(\lambda_K)$. (It is analogous to the case where $\widetilde{E}(k)_p = 0$, where $k$ is the residue field of $K$.) Now assume that $p = 2$. If $[K:F_v] < \infty$ and $E$ is nonsplit over $K$, then we have $H^1(K, C_v)/H^1(K, C_v)_{\mathrm{div}} \cong \mathbb{Z}/2\mathbb{Z}$ by (2), since $\psi\chi^{-1}$ will be the unramified character of $G_K$ of order 2. Thus, we obtain that $\mathrm{Im}(\kappa_K) = \mathrm{Im}(\lambda_K)_{\mathrm{div}}$ and that $[\mathrm{Im}(\lambda_K) : \mathrm{Im}(\kappa_K)] \leq 2$. Using the same argument as in the proof of proposition 2.5, we find that this index is equal to the Tamagawa factor $[E(K) : \mathcal{F}(\mathfrak{m}_K)]$ for $E$ over $K$. This equals 1 or 2 depending on whether $\mathrm{ord}_K(j_E)$ is odd or even. Finally, we remark that proposition 2.4 holds when $E$ has multiplicative reduction. The proof given there works because the index $[\mathrm{Im}(\lambda_M):\mathrm{Im}(\kappa_M)]$ is bounded.

For completeness, we will state a result of Bloch and Kato describing $\mathrm{Im}(\kappa_K)$ when $E$ has good, supersingular reduction and $[K:F_v] < \infty$. It involves the ring $B_{\mathrm{cris}}$ of Fontaine. Define

$$H_f^1\left(K, V_p(E)\right) = \ker\left(H^1(K, V_p(E)) \to H^1(K, V_p(E) \otimes B_{\mathrm{cris}})\right).$$

The result is that $\mathrm{Im}(\kappa_K)$ is the image of $H_f^1(K, V_p(E))$ under the canonical map $H^1(K, V_p(E)) \to H^1(K, V_p(E)/T_p(E))$, noting that $V_p(E)/T_p(E)$ is isomorphic to $E[p^\infty]$. This description is also correct if $E$ has good, ordinary reduction.

If $E$ has supersingular reduction at $v$, where $v|p$, and if $K$ is any ramified $\mathbb{Z}_p$-extension of $F_v$, then the analogue of proposition 2.4 is true. In this case, $C_v = E[p^\infty]$ since $\widetilde{E}[p^\infty] = 0$. Thus, the result is that $\mathrm{Im}(\kappa_K) = H^1(K, E[p^\infty])$. Perhaps the easiest way to prove this is to use the analogue of Hilbert's theorem 90 for formal groups proved in [CoGr]. If $\mathcal{F}$ denotes the formal group (of height 2) associated to $E$, then $H^1(K, \mathcal{F}(\overline{\mathfrak{m}})) = 0$. (This is a special case of Corollary 3.2 in [CoGr].) Just as in the case of Kummer theory for the multiplicative group, we then obtain an isomorphism

$$\kappa_K^{\mathcal{F}} : \mathcal{F}(\mathfrak{m}_K) \otimes (\mathbb{Q}_p/\mathbb{Z}_p) \xrightarrow{\sim} H^1(K, C_v)$$

because $C_v = \mathcal{F}(\overline{\mathfrak{m}})[p^\infty]$. We get the result stated above immediately, since

$$E(K) \otimes (\mathbb{Q}_p/\mathbb{Z}_p) = \mathcal{F}(\mathfrak{m}_K) \otimes (\mathbb{Q}_p/\mathbb{Z}_p).$$

The assertion that $\mathrm{Im}(\kappa_K) = H^1(K, E[p^\infty])$ is proved in [CoGr] under the hypotheses that $E$ has potentially supersingular reduction at $v$ and that $K/F_v$ is a "deeply ramified extension" (which means that $K/F_v$ has infinite conductor, i.e., $K \not\subset F_v^{(t)}$ for any $t \geq 1$, where $F_v^{(t)}$ denotes the fixed field for the $t$-th ramification subgroup of $\mathrm{Gal}(\overline{F}_v/F_v)$). A ramified $\mathbb{Z}_p$-extension $K$ of $F_v$ is the simplest example of a deeply ramified extension. As an illustration of how this result affects the structure of Selmer groups, consider the definition of $\mathrm{Sel}_E(M)_p$ given near the beginning of this section. If $E$ has potentially supersingular reduction at a prime $v$ of $F$ lying over $p$ and if $M_\eta/F_v$ is deeply ramified for all $\eta|v$, then the groups $H^1(M_\eta, E[p^\infty])/\mathrm{Im}(\kappa_\eta)$ occurring in the definition of $\mathrm{Sel}_E(M)_p$ are simply zero. In particular, if $M = F_\infty$, the cyclotomic $\mathbb{Z}_p$-extension of $F$, then the primes $\eta$ of $F_\infty$ lying over primes of $F$ where $E$ has potentially supersingular reduction can be omitted in the local conditions defining $\mathrm{Sel}_E(F_\infty)_p$. This is the key to proving theorem 1.7.

One extremely important consequence of the fact that the Selmer group for an elliptic curve $E$ has a description involving just the Galois representations attached to the torsion points on $E$ is that one can then attempt to introduce analogously-defined "Selmer groups" and to study all the natural questions associated to such objects in a far more general context. We will illustrate this idea by considering $\Delta$, the normalized cusp form of level 1, weight 12. Its $q$-expansion is $\Delta = \sum_{n=1}^{\infty} \tau(n)q^n$, where $\tau(n)$ is Ramanujan's tau function. Deligne attached to $\Delta$ a compatible system $\{V_l(\Delta)\}$ of $l$-adic representations of $G_{\mathbb{Q}}$. Consider a prime $p$ such that $p \nmid \tau(p)$. For such a prime $p$, Mazur and Wiles have proved that the action of $G_{\mathbb{Q}_p}$ on $V_p(\Delta)$ is reducible (where one fixes an embedding $\overline{\mathbb{Q}} \to \overline{\mathbb{Q}}_p$, identifying $G_{\mathbb{Q}_p}$ with a subgroup of $G_{\mathbb{Q}}$). More precisely, there is an exact sequence

$$0 \to W_p(\Delta) \to V_p(\Delta) \to U_p(\Delta) \to 0$$

where $W_p(\Delta)$ is 1-dimensional and $G_{\mathbb{Q}_p}$-invariant, the action of $G_{\mathbb{Q}_p}$ on $U_p(\Delta)$ is unramified, and the action of $\mathrm{Frob}_p$ on $U_p(\Delta)$ is multiplication by $\alpha_p$ (where $\alpha_p$ is the $p$-adic unit root of $t^2 - \tau(p)t + p^{11}$). Let $T_p(\Delta)$ be any $G_{\mathbb{Q}}$-invariant $\mathbb{Z}_p$-lattice in $V_p(\Delta)$. (It turns out to be unique up to homothety for $p \nmid \tau(p)$, except for $p = 691$, when there are two possible choices up to homothety.) Let $A = V_p(\Delta)/T_p(\Delta)$. As a group, $A \cong (\mathbb{Q}_p/\mathbb{Z}_p)^2$. Let $C$ denote the image of $W_p(\Delta)$ in $A$. Then $C \cong \mathbb{Q}_p/\mathbb{Z}_p$ as a group. Here then is a definition of the $p$-Selmer group $S_A(\mathbb{Q})_p$ for $A$ over $\mathbb{Q}$.

$$S_A(\mathbb{Q})_p = \ker\!\left(H^1(\mathbb{Q}, A) \to \prod_v H^1(\mathbb{Q}_v, A)/L_v\right)$$

where $v$ runs over all primes of $\mathbb{Q}$. Here we take $L_v = 0$ for $v \neq p$, analogously to the elliptic curve case. One defines $L_p = \mathrm{Im}(\lambda_p)_{\mathrm{div}}$, where

$$\lambda_p : H^1(\mathbb{Q}_p, C) \to H^1(\mathbb{Q}_p, A)$$

is the natural map. In [Gr3], one can find a calculation of $S_A(\mathbb{Q})_p$, and also $S_A(\mathbb{Q}_\infty)_p$, for $p = 11, 23$, and $691$. One can make similar definitions whenever one has a $p$-adic Galois representation with suitable properties.

## 3.   Control Theorems.

We will now give a proof of theorem 1.2. It is based on the description of the images of the local Kummer homomorphisms presented in section 2, specifically propositions 2.1, 2.2, and 2.4. We will also prove a special case of conjecture 1.6. Let $E$ be any elliptic curve defined over $F$. Let $M$ be an algebraic extension of $F$. For every prime $\eta$ of $M$, we let

$$\mathcal{H}_E(M_\eta) = H^1(M_\eta, E[p^\infty])/\mathrm{Im}(\kappa_\eta).$$

Let $\mathcal{P}_E(M) = \prod_\eta \mathcal{H}_E(M_\eta)$, where $\eta$ runs over all primes of $M$. Thus,

$$\mathrm{Sel}_E(M)_p = \ker\left(H^1(M, E[p^\infty]) \to \mathcal{P}_E(M)\right),$$

where the map is induced by restricting cocycles to decomposition groups. Also, we put

$$\mathcal{G}_E(M) = \mathrm{Im}\left(H^1(M, E[p^\infty]) \to \mathcal{P}_E(M)\right).$$

Let $F_\infty = \bigcup_n F_n$ be the cyclotomic $\mathbb{Z}_p$-extension. Consider the following commutative diagram with exact rows.

$$
\begin{array}{ccccccccc}
0 & \longrightarrow & \mathrm{Sel}_E(F_n)_p & \longrightarrow & H^1(F_n, E[p^\infty]) & \longrightarrow & \mathcal{G}_E(F_n) & \longrightarrow & 0 \\
 & & \downarrow{s_n} & & \downarrow{h_n} & & \downarrow{g_n} & & \\
0 & \longrightarrow & \mathrm{Sel}_E(F_\infty)_p^{\Gamma_n} & \longrightarrow & H^1(F_\infty, E[p^\infty])^{\Gamma_n} & \longrightarrow & \mathcal{G}_E(F_\infty)^{\Gamma_n}. & &
\end{array}
$$

Here $\Gamma_n = \mathrm{Gal}(F_\infty/F_n) = \Gamma^{p^n}$. The maps $s_n$, $h_n$, and $g_n$ are the natural restriction maps. The snake lemma then gives the exact sequence

$$0 \to \ker(s_n) \to \ker(h_n) \to \ker(g_n) \to \mathrm{coker}(s_n) \to \mathrm{coker}(h_n).$$

Therefore, we must study $\ker(h_n)$, $\mathrm{coker}(h_n)$, and $\ker(g_n)$, which we do in a sequence of lemmas.

**Lemma 3.1.**  *The kernel of $h_n$ is finite and has bounded order as $n$ varies.*

*Proof.* By the inflation-restriction sequence, $\ker(h_n) \cong H^1(\Gamma_n, B)$, where $B$ is the $p$-primary subgroup of $E(F_\infty)$. This group $B$ is in fact finite and hence $H^1(\Gamma_n, B) = \mathrm{Hom}(\Gamma_n, B)$ for $n \gg 0$. Lemma 3.1 follows immediately. But it is not necessary to know the finiteness of $B$. If $\gamma$ denotes a topological generator of $\Gamma$, then $H^1(\Gamma_n, B) = B/(\gamma^{p^n} - 1)B$. Since $E(F_n)$ is finitely generated, the kernel of $\gamma^{p^n} - 1$ acting on $B$ is finite. Now $B_{\mathrm{div}}$ has finite $\mathbb{Z}_p$-corank. It is clear that

$$B_{\mathrm{div}} \subseteq (\gamma^{p^n} - 1)B \subseteq B.$$

Thus, $H^1(\Gamma_n, B)$ has order bounded by $[B : B_{\mathrm{div}}]$, which is independent of $n$. If we use the fact that $B$ is finite, then $\ker(h_n)$ has the same order as $H^0(\Gamma_n, B)$, namely $|E(F_n)_p|$. ∎

**Lemma 3.2.** $\mathrm{Coker}(h_n) = 0$.

*Proof.* The sequence $H^1(F_n, E[p^\infty]) \to H^1(F_\infty, E[p^\infty])^{\Gamma_n} \to H^2(\Gamma_n, B)$ is exact, where $B = H^0(F_\infty, E[p^\infty])$ again. But $\Gamma_n \cong \mathbb{Z}_p$ is a free pro-$p$ group. Hence $H^2(\Gamma_n, B) = 0$. Thus, $h_n$ is surjective as claimed. ∎

Let $v$ be any prime of $F$. We will let $v_n$ denote any prime of $F_n$ lying over $v$. To study $\ker(g_n)$, we focus on each factor in $\mathcal{P}_E(F_n)$ by considering

$$r_{v_n} : \mathcal{H}_E((F_n)_{v_n}) \to \mathcal{H}_E((F_\infty)_\eta)$$

where $\eta$ is any prime of $F_\infty$ lying above $v_n$. ($\mathcal{P}_E(F_\infty)$ has a factor for all such $\eta$'s, but the kernels will be the same.) If $v$ is archimedean, then $v$ splits completely in $F_\infty/F$, i.e., $F_v = K_\eta$. Thus, $\ker(r_{v_n}) = 0$. For nonarchimedean $v$, we consider separately $v \nmid p$ and $v|p$.

**Lemma 3.3.** *Suppose $v$ is a nonarchimedean prime not dividing $p$. Then $\ker(r_{v_n})$ is finite and has bounded order as $n$ varies. If $E$ has good reduction at $v$, then $\ker(r_{v_n}) = 0$ for all $n$.*

*Proof.* By proposition 2.1, $\mathcal{H}_E(M_\eta) = H^1(M_\eta, E[p^\infty])$ for every algebraic extension $M_\eta$ of $F_v$. Let $B_v = H^0(K, E[p^\infty])$, where $K = (F_\infty)_\eta$. Since $v$ is unramified and finitely decomposed in $F_\infty/F$, $K$ is the unramified $\mathbb{Z}_p$-extension of $F_v$ (in fact, the only $\mathbb{Z}_p$-extension of $F_v$). The group $B_v$ is isomorphic to $(\mathbb{Q}_p/\mathbb{Z}_p)^e \times$ (a finite group), where $0 \le e \le 2$. Let $\Gamma_{v_n} = \mathrm{Gal}(K/(F_n)_{v_n})$, which is isomorphic to $\mathbb{Z}_p$, topologically generated by $\gamma_{v_n}$, say. Then

$$\ker(r_{v_n}) \cong H^1(\Gamma_{v_n}, B_v) \cong B_v/(\gamma_{v_n} - 1)B_v.$$

Since $E((F_n)_{v_n})$ has a finite $p$-primary subgroup, it is clear that $(\gamma_{v_n} - 1)B_v$ contains $(B_v)_{\mathrm{div}}$ (just as in the proof of lemma 3.1) and hence

$$|\ker(r_{v_n})| \le |B_v/(B_v)_{\mathrm{div}}|. \tag{4}$$

This bound is independent of $n$ and $v_n$. We have equality if $n \gg 0$. Now assume that $E$ has good reduction at $v$. Then, since $v \nmid p$, $F_v(E[p^\infty])/F_v$ is unramified. It is clear that $K \subseteq F_v(E[p^\infty])$ and that $\Delta = \mathrm{Gal}(F_v(E[p^\infty])/K)$ is a finite, cyclic group of order prime to $p$. It then follows that $B_v = E[p^\infty]^\Delta$ is divisible. Therefore, $\ker(r_{v_n}) = 0$ as stated. ∎

One can determine the precise order of $\ker(r_{v_n})$, where $v_n | v$ and $v$ is any nonarchimedean prime of $F$ not dividing $p$ where $E$ has bad reduction. This will be especially useful in section 4, where we will need $|\ker(r_v)|$. The result is: $|\ker(r_v)| = c_v^{(p)}$, *where $c_v^{(p)}$ is the highest power of $p$ dividing the Tamagawa factor $c_v$ for $E$ at $v$.* Recall that $c_v = [E(F_v) : E_0(F_v)]$, where $E_0(F_v)$ is the subgroup of local points which have nonsingular reduction at $v$. First we consider the case where $E$ has additive reduction at $v$. Then $H^0(I_v, E[p^\infty])$ is finite, where $I_v$ denotes the inertia subgroup of $G_{F_v}$. Hence $B_v$ is finite because $I_v \subseteq G_K$. Also, $E_0(F_v)$ is a pro-$l$ group, where $l$ is the characteristic of the residue field for $v$, i.e., $v|l$. (Note: Using the notation in [Si], chapter 5, we have $|\widetilde{E}_{ns}(f_v)| = |f_v| =$ a power of $l$ and $E_1(F_v)$ is pro-$l$.) Since $l \neq p$, we have $c_v^{(p)} = |E(F_v)_p|$, which in turn equals $|B_v/(\gamma_v - 1)B_v|$. Hence $|\ker(r_v)| = c_v^{(p)}$ when $E$ has additive reduction at $v$. (It is known that $c_v \leq 4$ when $E$ has additive reduction at $v$. Thus, for such $v$, $\ker(r_v) = 0$ if $p \geq 5$.) Now assume that $E$ has split, multiplicative reduction at $v$. Then $c_v = \mathrm{ord}_v(q_E^{(v)}) = -\mathrm{ord}_v(j_E)$, where $q_E^{(v)}$ denotes the Tate period for $E$ at $v$. Thus, $q_E^{(v)} = \pi_v^{c_v} \cdot u$, where $u$ is a unit of $F_v$ and $\pi_v$ is a uniformizing parameter. One can verify easily that the group of units in $K$ is divisible by $p$. By using the Tate parametrization one can show that $B_v/(B_v)_{\mathrm{div}}$ is cyclic of order $c_v^{(p)}$ and that $\Gamma_v$ acts trivially on this group. Thus, $|\ker(r_{v_n})| = c_v^{(p)}$ for all $n \geq 0$. $B_v$ might be infinite. In fact, $(B_v)_{\mathrm{div}} = \mu_{p^\infty}$ if $\mu_p \subseteq F_v$; $(B_v)_{\mathrm{div}} = 0$ if $\mu_p \not\subseteq F_v$. Finally, assume that $E$ has nonsplit, multiplicative reduction at $v$. Then $c_v = 1$ or 2, depending on whether $\mathrm{ord}_v(j_E)$ is odd or even. Using the Tate parametrization, one can see that $B_v$ is divisible when $p$ is odd (and then $\ker(r_v) = 0$). If $p = 2$, $E$ will have split, multiplicative reduction over $K$ and so again $B_v/(B_v)_{\mathrm{div}}$ has order related to $\mathrm{ord}_v(q_E^{(v)})$. But $\gamma_v$ acts by $-1$ on this quotient. Hence $H^1(\Gamma_v, B_v)$ has order 1 or 2, depending on the parity of $\mathrm{ord}_v(q_E^{(v)})$. Hence, in all cases, $|\ker(r_v)| = c_v^{(p)}$.

Now assume that $v|p$. For each $n$, we let $f_{v_n}$ denote the residue field for $(F_n)_{v_n}$. It doesn't depend on the choice of $v_n$. Also, since $v_n$ is totally ramified in $F_\infty/F_n$ for $n \gg 0$, the finite field $f_{v_n}$ stabilizes to $f_\eta$, the residue field of $(F_\infty)_\eta$. We let $\widetilde{E}$ denote the reduction of $E$ at $v$. Then we have

**Lemma 3.4.** *Assume that $E$ has good, ordinary reduction at $v$. Then*

$$|\ker(r_{v_n})| = |\widetilde{E}(f_{v_n})_p|^2.$$

*It is finite and has bounded order as $n$ varies.*

*Proof.* Let $C_v = \ker(E[p^\infty] \to \widetilde{E}[p^\infty])$, where we regard $E[p^\infty]$ as a subgroup of $E(\overline{F}_v)$. Considering $(F_n)_{v_n}$ as a subfield of $\overline{F}_v$, we have $\mathrm{Im}(\kappa_{v_n}) = \mathrm{Im}(\lambda_{v_n})_{\mathrm{div}}$ by proposition 2.2. By proposition 2.4, we have $\mathrm{Im}(\kappa_\eta) = \mathrm{Im}(\lambda_\eta)$, since the inertia subgroup of $\mathrm{Gal}(F_\infty/F)$ for $v$ has finite index. Thus, we can factor $r_{v_n}$ as follows.

$$
\begin{array}{ccc}
H^1((F_n)_{v_n}, E[p^\infty])/\mathrm{Im}(\lambda_{v_n})_{\mathrm{div}} & \xrightarrow{a_{v_n}} & H^1((F_n)_{v_n}, E[p^\infty])/\mathrm{Im}(\lambda_{v_n}) \\
& \searrow{\scriptstyle r_{v_n}} & \downarrow{\scriptstyle b_{v_n}} \\
& & H^1((F_\infty)_\eta, E[p^\infty])/\mathrm{Im}(\lambda_\eta)
\end{array}
$$

Now $a_{v_n}$ is clearly surjective. Hence $|\ker(r_{v_n})| = |\ker(a_{v_n})| \cdot |\ker(b_{v_n})|$. By proposition 2.5, we have $|\ker(a_{v_n})| = |\widetilde{E}(f_{v_n})_p|$. For the proof of proposition 1.2, just the boundedness of $|\ker(a_{v_n})|$ (and of $|\ker(b_{v_n})|$) suffices. To study $\ker(b_{v_n})$ we use the following commutative diagram.

$$
\begin{array}{ccccccc}
H^1((F_n)_{v_n}, C_v) & \xrightarrow{\lambda_{v_n}} & H^1((F_n)_{v_n}, E[p^\infty]) & \xrightarrow{\pi_{v_n}} & H^1((F_n)_{v_n}, \widetilde{E}[p^\infty]) & \longrightarrow & 0 \\
\downarrow & & \downarrow & & \downarrow{\scriptstyle d_{v_n}} & & \\
H^1((F_\infty)_\eta, C_v) & \xrightarrow{\lambda_\eta} & H^1((F_\infty)_\eta, E[p^\infty]) & \xrightarrow{\pi_\eta} & H^1((F_\infty)_\eta, \widetilde{E}[p^\infty]) & \longrightarrow & 0
\end{array}
\tag{5}
$$

The surjectivity of the first row follows from Poitou-Tate Duality, which gives $H^2(M, C_v) = 0$ for any finite extension $M$ of $F_v$. (Note that $C_v \not\cong \mu_{p^\infty}$ for the action of $G_M$.) Thus, $\ker(b_{v_n}) \cong \ker(d_{v_n})$. But

$$
\ker(d_{v_n}) \cong H^1((F_\infty)_\eta/(F_n)_{v_n}, \widetilde{E}(f_\eta)_p) \cong \widetilde{E}(f_\eta)_p/(\gamma_{v_n} - 1)\widetilde{E}(f_\eta)_p
$$

where $\gamma_{v_n}$ is a topological generator of $\mathrm{Gal}((F_\infty)_\eta/(F_n)_{v_n})$. Now $\widetilde{E}(f_\eta)_p$ is finite and the kernel and cokernel of $\gamma_{v_n} - 1$ have the same order, namely $|\widetilde{E}(f_{v_n})_p|$. This is the order of $\ker(d_{v_n})$. Lemma 3.4. follows. ∎

Let $\Sigma_0$ denote the finite set of nonarchimedean primes of $F$ which either lie over $p$ or where $E$ has bad reduction. If $v \notin \Sigma_0$ and $v_n$ is a prime of $F_n$ lying over $v$, then $\ker(r_{v_n}) = 0$. For each $v \in \Sigma_0$, lemmas 3.3 and 3.4 show that $|\ker(r_{v_n})|$ is bounded as $n$ varies. The number of primes $v_n$ of $F_n$ lying over any nonarchimedean prime $v$ is also bounded. Consequently, we have proved the following lemma.

**Lemma 3.5.** *The order of* $\ker(g_n)$ *is bounded as* $n$ *varies.*

Lemma 3.1 implies that $\ker(s_n)$ is finite and has bounded order no matter what type of reduction $E$ has at $v|p$. Lemmas 3.2 and 3.5 show that $\mathrm{coker}(s_n)$ is finite and of bounded order, assuming that $E$ has good, ordinary reduction at all $v|p$. Thus, theorem 1.2 is proved.

It is possible for $s_n$ to be injective for all $n$. A simple sufficient condition for this is: $E(F)$ *has no element of order* $p$. For then $E(F_\infty)$ will have no $p$-torsion, since $\Gamma = \mathrm{Gal}(F_\infty/F)$ is a pro-$p$ group. Thus $\ker(h_n)$ and hence $\ker(s_n)$ would be trivial for all $n$. A somewhat more subtle result will be proved later, in proposition 3.9.

It is also possible for $s_n$ to be surjective for all $n$. Still assuming that $E$ has good, ordinary reduction at all primes of $F$ lying over $v$, here is a sufficient condition for this: *For each* $v|p$, $\widetilde{E}_v(f_v)$ *has no element of order* $p$ *and, for each* $v$ *where* $E$ *has bad reduction,* $E[p^\infty]^{I_v}$ *is divisible.* The first part of this condition implies that $\widetilde{E}_v(f_{v_n})_p = 0$ for all $v|p$ and all $n$, again using the fact that $\Gamma$ is pro-$p$. Thus, $\ker(r_{v_n}) = 0$ by lemma 3.4. In the second part of this condition, $I_v$ denotes the inertia subgroup of $G_{F_v}$. Note that $v \nmid p$. It is easy to see that if $E[p^\infty]^{I_v}$ is divisible, the same is true of $B_v = H^0((F_\infty)_\eta, E[p^\infty])$ for $\eta|v$. Thus, $\ker(r_{v_n}) = 0$ for $v_n|v$, because of (4). The second part of this condition is equivalent to $p \nmid c_v$.

We want to now discuss the case where $E$ has multiplicative reduction at some $v|p$. In this case, one can attempt to imitate the proof of lemma 3.4, taking $C_v = \mathcal{F}(\overline{\mathfrak{m}})[p^\infty]$. We first assume that $E$ has split, multiplicative reduction. Then $C_v \cong \mu_{p^\infty}$ and we have an exact sequence

$$0 \to \mu_{p^\infty} \to E[p^\infty] \to \mathbb{Q}_p/\mathbb{Z}_p \to 0$$

of $G_{F_v}$-modules, where the action on $\mathbb{Q}_p/\mathbb{Z}_p$ is trivial. Then $H^1((F_n)_{v_n}, \mu_{p^\infty})$ and hence $\mathrm{Im}(\lambda_{v_n})$ are divisible. We have $\mathrm{Im}(\kappa_{v_n}) = \mathrm{Im}(\lambda_{v_n})$ as well as $\mathrm{Im}(\kappa_\eta) = \mathrm{Im}(\lambda_\eta)$. Thus, $\ker(r_{v_n}) = \ker(b_{v_n})$, where $b_{v_n}$ is the map

$$b_{v_n} : H^1((F_n)_{v_n}, E[p^\infty])/\mathrm{Im}(\lambda_{v_n}) \to H^1((F_\infty)_\eta, E[p^\infty])/\mathrm{Im}(\lambda_\eta).$$

For any algebraic extension $M$ of $F_v$, we have an exact sequence

$$H^1(M, \mu_{p^\infty}) \xrightarrow{\lambda_M} H^1(M, E[p^\infty]) \xrightarrow{\pi_M} H^1(M, \mathbb{Q}_p/\mathbb{Z}_p) \xrightarrow{\delta_M} H^2(M, \mu_{p^\infty}) \to 0.$$

If $[M{:}F_v] < \infty$, then Poitou-Tate Duality shows that $H^2(M, \mu_{p^\infty}) \cong \mathbb{Q}_p/\mathbb{Z}_p$, whereas $H^2(M, E[p^\infty]) = 0$. Thus, $\pi_M$ gives the surjectivity of $\delta_M$. Thus, $\pi_M$ is not surjective in contrast to the case where $E$ has good, ordinary reduction at $v$. We let $\pi_{v_n} = \pi_{(F_n)_{v_n}}$, $\pi_\eta = \pi_{(F_\infty)_\eta}$. Thus, $\ker(b_{v_n})$ can be identified with $\mathrm{Im}(\pi_{v_n}) \cap \ker(d_{v_n})$, where $d_{v_n}$ is the map

$$d_{v_n} : H^1((F_n)_{v_n}, \mathbb{Q}_p/\mathbb{Z}_p) \to H^1((F_\infty)_\eta, \mathbb{Q}_p/\mathbb{Z}_p).$$

The kernel of $d_{v_n}$ is quite easy to describe. We have

$$\ker(d_{v_n}) = \mathrm{Hom}(\mathrm{Gal}((F_\infty)_\eta/(F_n)_{v_n}), \mathbb{Q}_p/\mathbb{Z}_p)$$

which is isomorphic to $\mathbb{Q}_p/\mathbb{Z}_p$ as a group. The image of $\pi_{v_n}$ is more interesting to describe. It depends on the Tate period $q_E$ for $E$, which is defined by the equation $j(q_E) = j_E$, solving this equation for $q_E \in F_v^\times$. Here $j(q) = $ $q^{-1} + 744 + 196884q + \cdots$ for $|q|_v < 1$ and $j_E$ is the $j$-invariant for $E$. Since $j_E \in F$ is algebraic, the theorem of [B-D-G-P] referred to in section 1 implies that $q_E$ is transcendental. Also, we have $|q_E|_v = |j_E|_v^{-1}$. Let

$$\text{rec}_M : M^\times \to \text{Gal}(M^{ab}/M)$$

denote the reciprocity map of local class field theory. We will prove the following result.

**Proposition 3.6.** *Let $M$ be a finite extension of $F_v$. Then*

$$\text{Im}(\pi_M) = \{\psi \in \text{Hom}(\text{Gal}(M^{ab}/M), \mathbb{Q}_p/\mathbb{Z}_p) \mid \psi(\text{rec}_M(q_E)) = 0\}.$$

*If $M$ is a $\mathbb{Z}_p$-extension of $F_v$, then $\pi_M$ is surjective.*

*Proof.* The last statement is clear since $G_M$ has $p$-cohomological dimension 1 if $M/F_v$ has profinite degree divisible by $p^\infty$. For the first statement, the exact sequence

$$0 \to \mu_{p^n} \to E[p^n] \to \mathbb{Z}/p^n\mathbb{Z} \to 0$$

induces a map $\pi_M^{(n)} : H^1(M, E[p^n]) \to H^1(M, \mathbb{Z}/p^n\mathbb{Z})$ for every $n \geq 1$. Because of the Weil pairing, we have $\text{Hom}(E[p^n], \mu_{p^n}) \cong E[p^n]$. Thus, by Poitou-Tate Duality, $\pi_M^{(n)}$ is adjoint to the natural map

$$H^1(M, \mu_{p^n}) \to H^1(M, E[p^n])$$

whose kernel is easy to describe. It is generated by the class of the 1-cocycle $\phi : G_M \to \mu_{p^n}$ given by $\phi(g) = g(\sqrt[p^n]{q_E}) / \sqrt[p^n]{q_E}$ for all $g \in G_M$. The pairing

$$H^1(M, \mu_{p^n}) \times H^1(M, \mathbb{Z}/p^n\mathbb{Z}) \to \mathbb{Z}/p^n\mathbb{Z}$$

is just $(\phi_q, \psi) \to \psi(\text{rec}_M(q))$ for $q \in M^\times$, where $\phi_q$ is the 1-cocycle associated to $q$ as above, i.e., the image of $q$ under the Kummer homomorphism $M^\times/(M^\times)^{p^n} \to H^1(M, \mu_{p^n})$. This implies that

$$\text{Im}(\pi_M^{(n)}) = \{\psi \in \text{Hom}(\text{Gal}(M^{ab}/M), \mathbb{Z}/p^n\mathbb{Z}) \mid \psi(\text{rec}_M(q_E)) = 0\},$$

from which the first part of proposition 3.6 follows by just taking a direct limit. ∎

Still assuming that $E$ has split, multiplicative reduction at $v$, the statement that $\ker(r_{v_n})$ is finite is equivalent to the assertion that $\ker(d_{v_n}) \not\subseteq \text{Im}(\pi_{v_n})$. In this case, we show that $|\ker(r_{v_n})|$ is bounded as $n$ varies. For let $\sigma = \text{rec}_{F_v}(q_E)|_{(F_\infty)_\eta} \in \text{Gal}((F_\infty)_\eta/F_v)$. Let $e_n = [(F_n)_{v_n} : F_v]$. Then we have $\text{rec}_{(F_n)_{v_n}}(q_E)|_{(F_\infty)_\eta} = \sigma^{e_n}$. It is clear that

$$\ker(r_{v_n}) = \{\psi \in \text{Hom}(\text{Gal}((F_\infty)_\eta/(F_n)_{v_n}), \mathbb{Q}_p/\mathbb{Z}_p) \mid \psi(\sigma^{e_n}) = 0\}$$

has order equal to $[\mathrm{Gal}((F_\infty)_\eta/(F_n)_{v_n}) : \overline{\langle \sigma^{e_n} \rangle}]$. But $\mathrm{Gal}((F_\infty)_\eta/F_v) \cong \mathbb{Z}_p$. This index is constant for $n \geq 0$. Thus, $\ker(r_{v_n})$ is finite and of constant order as $n$ varies provided that $\sigma \neq id$. Let $\mathbb{Q}_p^{\mathrm{cyc}}$ denote the cyclotomic $\mathbb{Z}_p$-extension of $\mathbb{Q}_p$. Then $(F_\infty)_\eta = F_v \mathbb{Q}_p^{\mathrm{cyc}}$. We have the following diagram

$$
\begin{array}{ccc}
F_v^\times & \xrightarrow{\ \mathrm{rec}\ } & \mathrm{Gal}((F_\infty)_\eta/F_v) \\
\Big\downarrow{\scriptstyle N_{F_v/\mathbb{Q}_p}} & & \Big\downarrow{\scriptstyle \mathrm{rest}} \\
\mathbb{Q}_p^\times & \xrightarrow{\ \mathrm{rec}\ } & \mathrm{Gal}(\mathbb{Q}_p^{\mathrm{cyc}}/\mathbb{Q}_p)
\end{array}
$$

where the horizontal arrows are the reciprocity maps. It is known that the group of universal norms for $\mathbb{Q}_p^{\mathrm{cyc}}/\mathbb{Q}_p$ is precisely $\mu \cdot \langle p \rangle$, where $\mu$ denotes the roots of unity in $\mathbb{Q}_p$. This of course coincides with the kernel of the reciprocity map $\mathbb{Q}_p^\times \to \mathrm{Gal}(\mathbb{Q}_p^{\mathrm{cyc}}/\mathbb{Q}_p)$ and also coincides with the kernel of $\log_p$ (where we take Iwasawa's normalization $\log_p(p) = 0$.) Also, it is clear that $\sigma \neq \mathrm{id} \Leftrightarrow \sigma|_{\mathbb{Q}_p^{\mathrm{cyc}}} \neq \mathrm{id}$. Thus we have shown that $\ker(r_{v_n})$ is finite if and only if $\log_p(N_{F_v/\mathbb{Q}_p}(q_E)) \neq 0$. The order will then be constant and is determined by the projection of $N_{F_v/\mathbb{Q}_p}(q_E)$ to $\mathbb{Z}_p^\times$ in the decomposition $\mathbb{Q}_p^\times = \langle p \rangle \times \mathbb{Z}_p^\times$. One finds that

$$
|\ker(r_{v_n})| \sim \log_p(N_{F_v/\mathbb{Q}_p}(q_E))/2p[F_v \cap \mathbb{Q}_p^{\mathrm{cyc}} : \mathbb{Q}_p]
$$

where $\sim$ indicates that the two sides have the same $p$-adic valuation.

Assume now that $p$ is odd and that $E$ has nonsplit, multiplicative reduction. We then show that $\ker(r_{v_n}) = 0$. We have an exact sequence

$$
0 \to \mu_{p^\infty} \otimes \phi \to E[p^\infty] \to (\mathbb{Q}_p/\mathbb{Z}_p) \otimes \phi \to 0
$$

where $\phi$ is the unramified character of $G_{F_v}$ of order 2. As discussed in section 2, we have $\mathrm{Im}(\kappa_{v_n}) = \mathrm{Im}(\lambda_{v_n})$. Also $\pi_{v_n}$ is surjective. We can identify $\ker(r_{v_n})$ with $\ker(d_{v_n})$, where $d_{v_n}$ is the map

$$
H^1((F_n)_{v_n}, (\mathbb{Q}_p/\mathbb{Z}_p)(\phi)) \to H^1((F_\infty)_\eta, (\mathbb{Q}_p/\mathbb{Z}_p)(\phi))
$$

whose kernel is clearly zero. Thus, as stated, $\ker(r_{v_n}) = 0$. (The value of $N_{F_v/\mathbb{Q}_p}(q_E)$ is not relevant in this case.) If $p = 2$, then $|\mathrm{Im}(\lambda_{v_n})/\mathrm{Im}(\lambda_{v_n})_{\mathrm{div}}|$ is easily seen to be at most 2. Hence, if $E$ has nonsplit, multiplicative reduction over $(F_n)_{v_n}$, we have $|\ker(r_{v_n})| \leq 2$. (Note: It can happen that $(F_\infty)_\eta$ contains the unramified quadratic extension of $F_v$. Thus $E$ can become split over $(F_n)_{v_n}$ for $n > 0$.) We will give the order of $\ker(r_v)$ when $E$ has nonsplit, multiplicative reduction at $v|2$. The kernel of $a_v$ has order $[\mathrm{Im}(\lambda_v) : \mathrm{Im}(\kappa_v)]$, which is just the Tamagawa factor for $E$ at $v$. (See the discussion following the proof of proposition 2.5.) On the other hand, $\ker(b_v) \cong \ker(d_v)$ and this group has order 2. Thus, $|\ker(r_v)| \sim 2c_v$, where $c_v$ denotes the Tamagawa factor for $E$ at $v$.

The above observations together with lemmas 3.1–3.3 provide a proof of the following result in the direction of conjecture 1.6.

**Proposition 3.7.** *Assume that $E$ is an elliptic curve defined over $F$ which has good, ordinary reduction or multiplicative reduction at all primes $v$ of $F$ lying over $p$. Assume also that $\log_p(N_{F_v/\mathbb{Q}_p}(q_E^{(v)})) \neq 0$ for all $v$ where $E$ has multiplicative reduction. Then the maps*

$$s_n \colon \operatorname{Sel}_E(F_n)_p \to \operatorname{Sel}_E(F_\infty)_p^{\Gamma_n}$$

*have finite kernel and cokernel, of bounded order as $n$ varies.*

In the above result, $q_E^{(v)}$ denotes the Tate period for $E$ over $F_v$. If $j_E \in \mathbb{Q}_p$, then so is $q_E^{(v)}$. Thus, $N_{F_v/\mathbb{Q}_p}(q_E^{(v)}) = (q_E^{(v)})^{[F_v:\mathbb{Q}_p]}$ is transcendental according to the theorem of Barré-Sirieix, Diaz, Gramain, and Philibert. Perhaps, it is reasonable to conjecture in general that $N_{F_v/\mathbb{Q}_p}(q_E^{(v)})$ is transcendental whenever $j_E \in F_v \cap \overline{\mathbb{Q}}$. Then the hypothesis $\log_p(N_{F_v/\mathbb{Q}_p}(q_E^{(v)})) \neq 0$ obviously holds. This hypothesis is unnecessary in proposition 3.7, if $p$ is odd, for those $v$'s where $E$ has nonsplit, multiplicative reduction. (For $p = 2$, one needs the hypothesis when $E$ has split reduction over $(F_\infty)_\eta$.)

Let $X$ be a profinite $\Lambda$-module, where $\Lambda = \mathbb{Z}_p[[T]]$, $T = \gamma - 1$, as in section 1. Here are some facts which are easily proved or can be found in [Wa2].

(1) $X = TX \Rightarrow X = 0$.
(2) $X/TX$ *finite* $\Rightarrow X$ *is a finitely generated, torsion $\Lambda$-module.*
(3) $X/TX$ *finitely generated over $\mathbb{Z}_p \Rightarrow X$ is finitely generated over $\Lambda$.*
(4) *Assume that $X$ is a finitely generated, torsion $\Lambda$-module. Let $\theta_n$ denote $\gamma^{p^n} - 1 \in \Lambda$ for $n \geq 0$. Then there exists integers $a$, $b$, and $c$ such that the $\mathbb{Z}_p$-torsion subgroup of $X/\theta_n X$ has order $p^{e_n}$, where $e_n = an + bp^n + c$ for $n \gg 0$.*

We sketch an argument for (4). Let $f(T)$ be a generator for the characteristic ideal of $X$, assuming that $X$ is finitely generated and torsion over $\Lambda$. If we have $f(\zeta - 1) \neq 0$ for all $p$-power roots of unity, then $X/\theta_n X$ is finite for all $n \geq 0$ and one estimates its order by studying $\prod f(\zeta - 1)$, where $\zeta$ runs over the $p^n$-th roots of unity. One then could take $a = \lambda(f)$, $b = \mu(f)$ in (4). Suppose $X = \Lambda/(h(T)^e)$, where $h(T)$ is an irreducible element of $\Lambda$. If $h(T) \nmid \theta_n$ for all $n$, then we are in the case just discussed. This is true for $(h(T)) = p\Lambda$ for example. If $h(T) | \theta_{n_0}$ for some $n_0 \geq 0$, then write $\theta_n = h(T)\phi_n$, for $n \geq n_0$, where $\phi_n \in \Lambda$. Since $\theta_n = (1 + T)^{p^n} - 1$ has no multiple factors, we have $h(T) \nmid \phi_n$. Then we get an exact sequence

$$0 \to Y/\phi_n Y \to X/\theta_n X \to \Lambda/h(T)\Lambda \to 0$$

for $n \geq n_0$. Here $Y = (h(T))/(h(T)^e) \cong \Lambda/(h(T)^{e-1})$. Then $Y/\phi_n Y$ is finite and one estimates its growth essentially as mentioned above. Now $\Lambda/h(T)\Lambda$ is a free $\mathbb{Z}_p$-module of rank $= \lambda(h)$. Thus the $\mathbb{Z}_p$-torsion subgroup of $X/\theta_n X$ is $Y/\phi_n Y$ whose order is given by a formula as above. In general, $X$ is pseudo-isomorphic to a direct sum of $\Lambda$-modules of the form $\lambda/(h(T)^e)$ and one can reduce to that case. One sees that $b = \mu(f)$, where $f = f(T)$ generates the characteristic ideal of $X$. Also, $a = \lambda(f) - \lambda_0$, where $\lambda_0 = \max(\mathrm{rank}_{\mathbb{Z}_p}(X/\theta_n X))$. The $\mathbb{Z}_p$-rank of $X/\theta_n X$ clearly stabilizes, equal to $\lambda_0$ for $n \gg 0$.

These facts together with the results of this section have some immediate consequences, some of which we state here without trying to be as general as possible. For simplicity, we take $\mathbb{Q}$ as the base field.

**Proposition 3.8.** *Let $E$ be an elliptic curve with good, ordinary reduction at $p$. We make the following assumptions:*

*(i) $p$ does not divide $|\widetilde{E}(\mathbb{F}_p)|$, where $\widetilde{E}$ denotes the reduction of $E$ at $p$.*

*(ii) If $E$ has split, multiplicative reduction at $l$, where $l \neq p$, then $p \nmid \mathrm{ord}_l(j_E)$. If $E$ has nonsplit, multiplicative reduction at $l$, then either $p$ is odd or $\mathrm{ord}_l(j_E)$ is odd.*

*(iii) If $E$ has additive reduction at $l$, then $E(\mathbb{Q}_l)$ has no point of order $p$.*

*Then the map $\mathrm{Sel}_E(\mathbb{Q})_p \to \mathrm{Sel}_E(\mathbb{Q}_\infty)_p^\Gamma$ is surjective. If $\mathrm{Sel}_E(\mathbb{Q})_p = 0$, then $\mathrm{Sel}_E(\mathbb{Q}_\infty)_p = 0$ also.*

*Remark.* The comments in the paragraph following the proof of lemma 3.3 allow us to restate hypotheses (ii) and (iii) in the following way: $p \nmid c_l$ *for all* $l \neq p$. Here $c_l$ is the Tamagawa factor for $E$ at $l$. If $E$ has good reduction at $l$, then $c_l = 1$. If $E$ has additive reduction at $l$, then $c_l \leq 4$. Thus, hypothesis (iii) is automatically satisfied for any $p \geq 5$. If $E$ has nonsplit, multiplicative reduction at $l$, then hypothesis (ii) holds for any $p \geq 3$. On the other hand, if $E$ has split, multiplicative reduction at $l$, then there is no restriction on the primes which could possibly divide $c_l$. Hypothesis (i) is equivalent to $a_p \not\equiv 1 \pmod{p}$, where $a_p = 1 + p - |\widetilde{E}(\mathbb{F}_p)|$.

*Proof.* We refer back to the sequence at the beginning of this section. We have $\mathrm{coker}(h_n) = 0$ by lemma 3.2. The surjectivity of the map $s_0$ would follow from the assertion $\ker(g_0) = 0$. But the above assumptions simply guarantee that the map $\mathcal{P}_E(\mathbb{Q}) \to \mathcal{P}_E(\mathbb{Q}_\infty)$ is injective and hence that $\ker(g_0) = 0$. For by lemma 3.4, (i) implies that $\ker(r_p) = 0$. If $E$ has multiplicative reduction at $l \neq p$ then (ii) implies that $\mathrm{ord}_l(q_E^{(l)})$ is not divisible by $p$. This means $\mathbb{Q}_l(\sqrt[p]{q_E^{(l)}})/\mathbb{Q}_l$ is ramified. Thus $H^0(L, E[p^\infty])$ is a divisible group, where $L$ denotes the maximal unramified extension of $\mathbb{Q}_l$. Now $\mathrm{Gal}(L/\mathbb{Q}_l) \cong \widehat{\mathbb{Z}}$. The cyclotomic $\mathbb{Z}_p$-extension of $\mathbb{Q}_l$ is $(\mathbb{Q}_\infty)_\eta$, where $\eta | l$. Thus, $(\mathbb{Q}_\infty)_\eta \subseteq L$. Let $H = \mathrm{Gal}(L/(\mathbb{Q}_\infty)_\eta)$. Then $H$ acts on $H^0(L, E[p^\infty])$ through a finite cyclic group of order prime to $p$. Thus, it is easy to see that $H^0((\mathbb{Q}_\infty)_\eta, E[p^\infty])$ is divisible and hence, from (4), we have $\ker(r_l) = 0$. Assume now that $E$ has additive reduction at $l$ (where, of course, $l \neq p$). Then $E[p^\infty]^{I_l}$ is finite, where $I_l$ denotes $G_L$, the inertia subgroup of $G_{\mathbb{Q}_l}$. We know that if $E$ has potentially good reduction at $l$, then $I_l$ acts on $E[p^\infty]$ through a quotient of order $2^a 3^b$. Thus $E[p^\infty]^{I_l} = 0$ if $p \geq 5$, and (iii) is then not important. If $p = 2$ or 3, (iii) suffices to conclude that $H^0((\mathbb{Q}_\infty)_\eta, E[p^\infty]) = 0$ since $\mathrm{Gal}((\mathbb{Q}_\infty)_\eta / \mathbb{Q}_l)$ is pro-$p$. Thus, again, $\ker(r_l) = 0$. (We are essentially repeating some previous observations.) Finally, if $\mathrm{Sel}_E(\mathbb{Q})_p$ is trivial, then so is $\mathrm{Sel}_E(\mathbb{Q}_\infty)_p^\Gamma$. Let $X = X_E(\mathbb{Q}_\infty)$. Then $X/TX = 0$, which implies that $X = 0$. Hence $\mathrm{Sel}_E(\mathbb{Q}_\infty) = 0$ as stated.    ∎

If we continue to take $F = \mathbb{Q}$, then we now know that the restriction map $\mathrm{Sel}_E(\mathbb{Q})_p \to \mathrm{Sel}_E(\mathbb{Q}_\infty)_p^\Gamma$ has finite cokernel if $E$ has good, ordinary or multiplicative reduction at $p$. (In fact, potentially ordinary or potentially multiplicative reduction would suffice.) Thus, if $\mathrm{Sel}_E(\mathbb{Q})_p$ is finite, then so is $\mathrm{Sel}_E(\mathbb{Q}_\infty)_p^\Gamma$. Hence, for $X = X_E(\mathbb{Q}_\infty)$, we would have that $X/TX$ is finite. Thus, $X$ would be a $\Lambda$-torsion module. In addition, we would have $T \nmid f_E(T)$.

Assume that $E$ has good, ordinary reduction at $p$. If $p$ is odd, then the map $\mathrm{Sel}_E(\mathbb{Q}_n)_p \to \mathrm{Sel}_E(\mathbb{Q}_\infty)_p^{\Gamma_n}$ is actually injective for all $n \geq 0$. To see this, let $B = H^0(\mathbb{Q}_\infty, E[p^\infty])$. Then $\ker(h_n) = H^1(\Gamma_n, B)$. The inertia subgroup $I_p$ of $G_{\mathbb{Q}_p}$ acts on $\ker(E[p] \to \widetilde{E}[p])$ by the Teichmüller character $\omega$. That is,

$$\ker(E[p] \to \widetilde{E}[p]) \cong \mu_p$$

for the action of $I_p$. On the other hand, $I_p$ acts on $B$ through $\mathrm{Gal}((\mathbb{Q}_\infty)_\eta / \mathbb{Q}_p)$, where $\eta$ denotes the unique prime of $\mathbb{Q}_\infty$ lying over $p$. This Galois group is pro-$p$, being isomorphic to $\mathbb{Z}_p$. Since $p > 2$, $\omega$ has nontrivial order and this order is relatively prime to $p$. It follows that

$$B \cap \ker(E[p^\infty] \to \widetilde{E}[p^\infty]) = \{O_E\}$$

and therefore $B$ maps injectively into $\widetilde{E}[p^\infty]$. Thus, $I_p$ acts trivially on $B$ Since $p$ is totally ramified in $\mathbb{Q}_\infty / \mathbb{Q}$, it is clear that $\Gamma = \mathrm{Gal}(\mathbb{Q}_\infty / \mathbb{Q})$ also acts trivially on $B$. That is,

$$B = E(\mathbb{Q}_\infty)_p = E(\mathbb{Q})_p.$$

Hence $\ker(h_n) = \mathrm{Hom}(\Gamma_n, B)$ for all $n \geq 0$. Now suppose that $\phi$ is a nontrivial element of $\mathrm{Hom}(\Gamma_n, B)$. Let $I_p^{(n)}$ denote the inertia subgroup of $G_{(\mathbb{Q}_n)_\eta}$. Then $\phi$ clearly remains nontrivial when restricted to

$$H^1(I_p^{(n)}, \widetilde{E}[p^\infty]) = \mathrm{Hom}(I_p^{(n)}, \widetilde{E}[p^\infty]).$$

But this implies that $[\phi] \notin \mathrm{Sel}_E(\mathbb{Q}_n)_p$. Hence $\ker(s_n) = \ker(h_n) \cap \mathrm{Sel}_E(\mathbb{Q}_n)_p$ is trivial as claimed. This argument also applies if $E$ has multiplicative reduction at $p$. More generally, the argument gives the following result. We let $F$ be any number field. For any prime $v$ of $F$ lying over $p$, we let $e(v/p)$ denote the ramification index for $F_v / \mathbb{Q}_p$.

**Proposition 3.9.** *Let $E$ be an elliptic curve defined over $F$. Assume that there is at least one prime $v$ of $F$ lying over $p$ with the following properties:*

*(i) $E$ has good, ordinary reduction or multiplicative reduction at $v$,*
*(ii) $e(v/p) \leq p - 2$.*

*Then the map $\operatorname{Sel}_E(F_n)_p \to \operatorname{Sel}_E(F_\infty)_p$ is injective for all $n \geq 0$.*

Theorem 1.10 is also an application of the results described in this section. One applies the general fact (4) about torsion $\Lambda$-modules to $X = X_E(F_\infty)$. Then, $X/\theta_n X$ is the Pontryagin dual of $\operatorname{Sel}_E(F_\infty)_p^{\Gamma_n}$. The torsion subgroup of $X/\theta_n X$ is then dual to $\operatorname{Sel}_E(F_\infty)_p^{\Gamma_n}/(\operatorname{Sel}_E(F_\infty)_p^{\Gamma_n})_{\mathrm{div}}$. One compares this to $\operatorname{Sel}_E(F_n)_p/(\operatorname{Sel}_E(F_n)_p)_{\mathrm{div}}$, which is precisely $\text{Ш}_E(F_n)_p$ under the assumption of finiteness. One must show that the orders of the relevant kernels and cokernels stabilize, which we leave for the reader. One then obtains the formula for the growth of $|\text{Ш}_E(F_n)_p|$, with the stated $\lambda$ and $\mu$.

We want to mention one other useful result. It plays a role in Li Guo's proof of a parity conjecture for elliptic curves with complex multiplication. (See [Gu2].)

**Proposition 3.10.** *Assume that $E$ is an elliptic curve/$F$ and that $\operatorname{Sel}_E(F_\infty)_p$ is $\Lambda$-cotorsion. Let $\lambda_E = \operatorname{corank}_{\mathbb{Z}_p}(\operatorname{Sel}_E(F_\infty)_p)$. Assume also that $p$ is odd. Then*

$$\operatorname{corank}_{\mathbb{Z}_p}(\operatorname{Sel}_E(F)_p) \equiv \lambda_E \pmod{2}.$$

*Proof.* The maps $H^1(F_n, E[p^\infty]) \to H^1(F_\infty, E[p^\infty])$ have finite kernels of bounded order as $n$ varies, by lemma 3.1. Thus, $\operatorname{corank}_{\mathbb{Z}_p}(\operatorname{Sel}_E(F_n)_p)$ is bounded above by $\lambda_E$. Let $\lambda'_E$ denote the maximum of these $\mathbb{Z}_p$-coranks. Then $\operatorname{corank}_{\mathbb{Z}_p}(\operatorname{Sel}_E(F_n)_p) = \lambda'_E$ for all $n \geq n_0$, say. For brevity, we let $S_n = \operatorname{Sel}_E(F_n)_p$, $T_n = (S_n)_{\mathrm{div}}$, and $U_n = S_n/T_n$, which is finite. The restriction map $S_0 \to S_n^{\operatorname{Gal}(F_n/F)}$, and hence the map $T_0 \to T_n^{\operatorname{Gal}(F_n/F)}$, have finite kernel and cokernel. Since the nontrivial $\mathbb{Q}_p$-irreducible representations of $\operatorname{Gal}(F_n/F)$ have degree divisible by $p - 1$, it follows easily that $\operatorname{corank}_{\mathbb{Z}_p}(T_n) \equiv \operatorname{corank}_{\mathbb{Z}_p}(T_0) \pmod{p-1}$. Hence

$$\operatorname{corank}_{\mathbb{Z}_p}(\operatorname{Sel}_E(F)_p) \equiv \lambda'_E \pmod{p-1}.$$

Since $p$ is odd, this gives a congruence modulo 2. Let $S_\infty = \operatorname{Sel}_E(F_\infty)_p$ and let $T_\infty = \varprojlim T_n$, which is a $\Lambda$-submodule of $S_\infty$. Also, $T_\infty \cong (\mathbb{Q}_p/\mathbb{Z}_p)^{\lambda'_E}$. Let $U_\infty = \overrightarrow{S_\infty/T_\infty} = \varinjlim U_n$. The map $T_n \to T_\infty$ is obviously surjective for all $n \geq n_0$ (since the kernel is finite). This implies that

$$|\ker(U_n \to U_\infty)| \leq |\ker(S_n \to S_\infty)|$$

for $n \geq n_0$, which is of bounded order as $n$ varies. Now a well-known theorem of Cassels states that there exists a nondegenerate, skew-symmetric pairing

$$U_n \times U_n \to \mathbb{Q}_p/\mathbb{Z}_p.$$

This forces $|U_n|$ to be a perfect square. More precisely, if the abelian group $U_n$ is decomposed as a direct product of cyclic groups of orders $p^{e_n^{(i)}}$, $1 \le i \le g_n$, say, then $g_n$ is even and one can arrange the terms so that $e_n^{(1)} = e_n^{(2)} \ge \cdots \ge e_n^{(g_n-1)} = e_n^{(g_n)}$. We refer to [Gu1] for a proof of this elementary result. (See lemma 3, page 157 there.) Since the kernels of the maps $U_n \to U_\infty$ have bounded order, the $\mathbf{Z}_p$-corank $u$ of $U_\infty$ can be determined from the behavior of the $e_n^{(i)}$'s as $n \to \infty$, namely, the first $u$ of the $e_n^{(i)}$'s will be unbounded, the rest bounded as $n \to \infty$. Thus $u$ is even. Since $u = \lambda_E - \lambda'_E$, it follows that

$$\lambda_E \equiv \lambda'_E \pmod{2}.$$

Combining that with the previous congruence, we get proposition 3.10. ∎

**Appendix to Section 3.** We would like to give a different and rather novel proof of a slightly weaker form of proposition 3.6, which is in fact adequate for proving proposition 3.7. We let $\widetilde{M}$ denote the composition of all $\mathbf{Z}_p$-extensions of $M$. For any $q \in M^\times$, we let $\widetilde{M}_q$ denote the composition of all $\mathbf{Z}_p$-extensions $M_\infty$ of $M$ such that $\mathrm{rec}_M(q)|_{M_\infty}$ is trivial, i.e., the image of $q$ under the reciprocity map $M^\times \to \mathrm{Gal}(M_\infty/M)$ is trivial. This means that $q \in N_{M_n/M}(M_n^\times)$ for all $n \ge 0$, where $M_n$ denotes the $n$-th layer in $M_\infty/M$. We then say that $q$ is a universal norm for the $\mathbf{Z}_p$-extension $M_\infty/M$. We will show that

$$\mathrm{Im}(\pi_M)_{\mathrm{div}} = \mathrm{Hom}(\mathrm{Gal}(\widetilde{M}_{q_E}/M), \mathbf{Q}_p/\mathbf{Z}_p). \tag{6}$$

The proof is based on the following observation:

**Proposition 3.11.** *Assume that $q \in M^\times$ is a universal norm for the $\mathbf{Z}_p$-extension $M_\infty/M$. Then the image of $\langle q \rangle \otimes (\mathbf{Q}_p/\mathbf{Z}_p)$ under the composite map*

$$\langle q \rangle \otimes (\mathbf{Q}_p/\mathbf{Z}_p) \to M^\times \otimes (\mathbf{Q}_p/\mathbf{Z}_p) \xrightarrow{\sim} H^1(M, \mu_{p^\infty}) \to H^1(M_\infty, \mu_{p^\infty})$$

*is contained in $H^1(M_\infty, \mu_{p^\infty})_{\Lambda\text{-div}}$, where $\Lambda = \mathbf{Z}_p[[\mathrm{Gal}(M_\infty/M)]]$.*

*Proof.* To justify this, note that the inflation-restriction sequence shows that the natural map

$$H^1(M_n, \mu_{p^\infty}) \to H^1(M_\infty, \mu_{p^\infty})^{\Gamma_n}$$

is surjective and has finite kernel. Here $\Gamma = \mathrm{Gal}(M_\infty/M)$, $\Gamma_n = \Gamma^{p^n} = \mathrm{Gal}(M_\infty/M_n)$. But $H^1(M_n, \mu_{p^n})$ is isomorphic to $(\mathbf{Q}_p/\mathbf{Z}_p)^{tp^n+1}$ as a group, where $t = [M : \mathbf{Q}_p]$. Thus, $H^1(M_\infty, \mu_{p^\infty})^{\Gamma_n}$ is divisible and has $\mathbf{Z}_p$-corank $tp^n + 1$. If $X = H^1(M_\infty, \mu_{p^\infty})\widehat{\ }$, then $X$ is a finitely generated $\Lambda$-module with the property that $X/\theta_n X \cong \mathbb{Z}_p^{tp^n+1}$ for all $n \geq 0$, where $\theta_n = \gamma^{p^n} - 1 \in \Lambda$, and $\gamma$ is some topological generator of $\Gamma$. It is not hard to deduce from this that $X \cong \Lambda^t \times \mathbb{Z}_p$, where $\mathbb{Z}_p = X_{\Lambda\text{-tors}}$ is just $\Lambda/\theta_0 \Lambda$. Letting $\widehat{\Lambda}$ denote the Pontryagin dual of $\Lambda$, regarded as a discrete $\Lambda$-module, we have

$$H^1(M_\infty, \mu_{p^\infty}) \cong \widehat{\Lambda}^t \times (\mathbb{Q}_p/\mathbb{Z}_p),$$

where the action of $\Gamma$ on $\mathbb{Q}_p/\mathbb{Z}_p$ is trivial. Thus, $H^1(M_\infty, \mu_{p^\infty})_{\Lambda\text{-div}} \cong \widehat{\Lambda}^t$, noting that the Pontryagin dual of a torsion-free $\Lambda$-module is $\Lambda$-divisible. Hence $(H^1(M_\infty, \mu_{p^\infty})_{\Lambda\text{-div}})^\Gamma$ has $\mathbb{Z}_p$-corank $t$. The maximal divisible subgroup of its inverse image in $M^\times \otimes (\mathbb{Q}_p/\mathbb{Z}_p)$ is isomorphic to $(\mathbb{Q}_p/\mathbb{Z}_p)^t$. We must show that this "canonical subgroup" of $M^\times \otimes (\mathbb{Q}_p/\mathbb{Z}_p)$, which the $\mathbb{Z}_p$-extension $M_\infty/M$ determines, contains $\langle q \rangle \otimes (\mathbb{Q}_p/\mathbb{Z}_p)$ whenever $q$ is a universal norm for $M_\infty/M$. Since $\mathrm{Gal}(M_\infty/M)$ is torsion-free, we may assume that $q \notin (M^\times)^p$. For every $n \geq 0$, choose $q_n \in M_n^\times$ so that $N_{M_n/M}(q_n) = q$. Fix $m \geq 1$. Consider $\alpha = q \otimes (1/p^m)$. In $M_n^\times \otimes (\mathbb{Q}_p/\mathbb{Z}_p)$, we have $N_{M_n/M}(\alpha_n) = \alpha$, where $\alpha_n = q_n \otimes (1/p^m)$. Let $\widetilde{\alpha}$, $\widetilde{\alpha}_n$ denote the images of $\alpha$, $\alpha_n$ in $M_\infty^\times \otimes (\mathbb{Q}_p/\mathbb{Z}_p)/(M_\infty^\times \otimes (\mathbb{Q}_p/\mathbb{Z}_p))_{\Lambda\text{-div}}$. The action of $\Gamma$ on this group is trivial. Hence $p^n \widetilde{\alpha}_n = \widetilde{\alpha}$. But $\widetilde{\alpha}_n$ has order dividing $p^m$. Since $n$ is arbitrary, we have $\widetilde{\alpha} = 0$, which of course means that the image of $q \otimes (1/p^m)$ is in $H^1(M_\infty, \mu_{p^\infty})_{\Lambda\text{-div}}$. This is true for any $m \geq 1$, as claimed. ∎

We now will prove (6). We know that $H^1(M, \mathbb{Q}_p/\mathbb{Z}_p)$ has $\mathbb{Z}_p$-corank $t+1$. Thus, $\mathrm{Im}(\pi_M)$ has $\mathbb{Z}_p$-corank $t$, which is also the $\mathbb{Z}_p$-corank of $\mathrm{Gal}(\widehat{M}_{q_E}/M)$. To justify (6), it therefore suffices to prove that $\mathrm{Hom}(\mathrm{Gal}(M_\infty/M), \mathbb{Q}_p/\mathbb{Z}_p)$ is contained in $\mathrm{Im}(\pi_M)$ for all $\mathbb{Z}_p$-extensions $M_\infty$ of $M$ contained in $\widehat{M}_{q_E}$. We do this by studying the following diagram

$$
\begin{array}{ccccc}
0 \longrightarrow H^1(M, \mu_{p^\infty})/B & \longrightarrow & H^1(M, E[p^\infty]) & \xrightarrow{\;\pi_M\;} & H^1(M, \mathbb{Q}_p/\mathbb{Z}_p) \\
\downarrow{\scriptstyle a} & & \downarrow{\scriptstyle b} & & \downarrow{\scriptstyle c} \\
0 \longrightarrow (H^1(M_\infty, \mu_{p^\infty})/B_\infty)^\Gamma & \longrightarrow & (H^1(M_\infty, E[p^\infty]))^\Gamma & \xrightarrow{\;e\;} & (H^1(M_\infty, \mathbb{Q}_p/\mathbb{Z}_p)^\Gamma
\end{array}
$$

where $B$ is the image of $\langle q_E \rangle \otimes (\mathbb{Q}_p/\mathbb{Z}_p)$ in $H^1(M, \mu_{p^\infty})$, which is the kernel of the map $H^1(M, \mu_{p^\infty}) \to H^1(M, E[p^\infty])$. Thus the first row is exact. We define $B_\infty$ as the image of $B$ under the restriction map. The exactness of the second row follows similarly, noting that $B_\infty$ is the image of $\langle q_E \rangle \otimes (\mathbb{Q}_p/\mathbb{Z}_p)$ in $H^1(M_\infty, \mu_{p^\infty})$. Now $\ker(c) = \mathrm{Hom}(\mathrm{Gal}(M_\infty/M), \mathbb{Q}_p/\mathbb{Z}_p)$ is isomorphic to $\mathbb{Q}_p/\mathbb{Z}_p$. We prove that $\ker(c) \subseteq \mathrm{Im}(\pi_M)$ by showing that $\mathrm{Im}(c \circ \pi_M) = \mathrm{Im}(e \circ b)$ has $\mathbb{Z}_p$-corank $t-1$. The first row shows that $H^1(M, E[p^\infty])$ has $\mathbb{Z}_p$-corank $2t$. Since $b$ is surjective and has finite kernel, the $\mathbb{Z}_p$-corank of $H^1(M_\infty, E[p^\infty])^\Gamma$ is also $2t$. But $H^1(M_\infty, \mu_{p^\infty}) \cong \widehat{\Lambda}^t \times (\mathbb{Q}_p/\mathbb{Z}_p)$ and $B_\infty$ is contained in the $\Lambda$-divisible submodule corresponding to $\widehat{\Lambda}^t$ by proposition 3.11. One can see from this that $H^1(M_\infty, \mu_{p^\infty})/B_\infty$ is also isomorphic to $\widehat{\Lambda}^t \times (\mathbb{Q}_p/\mathbb{Z}_p)$. (This is an exercise on $\Lambda$-modules: If $X$ is a free $\Lambda$-module of finite rank and $Y$ is a $\Lambda$-submodule such that $X/Y$ has no $\mathbb{Z}_p$-torsion, then $Y$ is a free $\Lambda$-module too.) It now follows that $(H^1(M_\infty, \mu_{p^\infty})/B_\infty)^\Gamma$ has $\mathbb{Z}_p$-corank $t + 1$. Therefore, $\mathrm{Im}(e)$ indeed has $\mathbb{Z}_p$-corank $t - 1$.

# 4.    Calculation of an Euler Characteristic.

This section will concern the evaluation of $f_E(0)$. We will assume that $E$ has good, ordinary reduction at all primes of $F$ lying over $p$. We will also assume that $\mathrm{Sel}_E(F)_p$ is finite. By theorem 1.4, $\mathrm{Sel}_E(F_\infty)_p$ is then $\Lambda$-cotorsion. By definition, $f_E(T)$ is a generator of the characteristic ideal of the $\Lambda$-module $X_E(F_\infty) = \mathrm{Hom}(\mathrm{Sel}_E(F_\infty)_p, \mathbb{Q}_p/\mathbb{Z}_p)$. Since $\mathrm{Sel}_E(F_\infty)_p^\Gamma$ is finite by theorem 1.2, it follows that $X_E(F_\infty)/TX_E(F_\infty)$ is finite. Hence $T \nmid f_E(T)$ and so $f_E(0) \neq 0$. The following theorem is a special case of a result of B. Perrin-Riou (if $E$ has complex multiplication) and of P. Schneider (in general). (See [Pe1] and [Sch1].) For every prime $v$ of $F$ lying over $p$, we let $\widetilde{E}_v$ denote the reduction of $E$ modulo $v$, which is defined over the residue field $f_v$. For primes $v$ where $E$ has bad reduction, we let $c_v = [E(F_v):E_0(F_v)]$ as before, where $E_0(F_v)$ denotes the subgroup of points with nonsingular reduction modulo $v$. The highest power of $p$ dividing $c_v$ is denoted by $c_v^{(p)}$. Also, if $a, b \in \mathbb{Q}_p^\times$, we write $a \sim b$ to indicate that $a$ and $b$ have the same $p$-adic valuation.

**Theorem 4.1.** *Assume that $E$ is an elliptic curve defined over $F$ with good, ordinary reduction at all primes of $F$ lying over $p$. Assume also that $\mathrm{Sel}_E(F)_p$ is finite. Then*

$$f_E(0) \sim (\prod_{v \text{ bad}} c_v^{(p)})(\prod_{v|p} |\widetilde{E}_v(f_v)_p|^2)|\mathrm{Sel}_E(F)_p|/|E(F)_p|^2.$$

Note that under the above hypotheses, $\mathrm{Sel}_E(F)_p = \text{Ш}_E(F)_p$. Also, we have $|\widetilde{E}_v(f_v)| = (1 - \alpha_v)(1 - \beta_v)$, where $\alpha_v \beta_v = N(v)$, $\alpha_v + \beta_v = a_v \in \mathbb{Z}$, and $p \nmid a_v$. It follows that $\alpha_v, \beta_v \in \mathbb{Q}_p$. We can assume that $\alpha_v \in \mathbb{Z}_p^\times$. Hence $p \mid |\widetilde{E}_v(f_v)|$ if and only if $a_v \equiv 1 \pmod{p}$. We say in this case that $v$ is an anomalous prime for $E$, a terminology introduced by Mazur who first pointed out the interest of such primes for the Iwasawa theory of $E$. In [Maz1], one finds an extensive discussion of them.

We will prove theorem 4.1 by a series of lemmas. We begin with a general fact about $\Lambda$-modules.

**Lemma 4.2.** *Assume that $S$ is a cofinitely generated, cotorsion $\Lambda$-module. Let $f(T)$ be a generator of the characteristic ideal of $X = \mathrm{Hom}(S, \mathbb{Q}_p/\mathbb{Z}_p)$. Assume that $S^\Gamma$ is finite. Then $S_\Gamma$ is finite, $f(0) \neq 0$, and $f(0) \sim |S^\Gamma|/|S_\Gamma|$.*

*Remark.* Note that $H^i(\Gamma, S) = 0$ for $i > 1$. Hence the quantity $|S^\Gamma|/|S_\Gamma|$ is the Euler characteristic $|H^0(\Gamma, S)|/|H^1(\Gamma, S)|$. Also, the assumption that $S^\Gamma$ is finite in fact implies that $S$ is cofinitely generated and cotorsion as a $\Lambda$-module.

*Proof of lemma 4.2.* By assumption, we have that $X/TX$ is finite. Now $X$ is pseudo-isomorphic to a direct sum of $\Lambda$-modules of the form $Y = \Lambda/(g(T))$. For each such $Y$, we have $Y/TY = \Lambda/(T, g(T)) = \mathbb{Z}_p/(g(0))$. Thus, $Y/TY$ is finite if and only if $g(0) \neq 0$. In this case, we have $\ker(T: Y \to Y) = 0$. ¿From this, one sees that $X/TX$ is finite if and only if $f(0) \neq 0$, and then obviously $\ker(T: X \to X)$ would be finite. Thus, $S_\Gamma$ is finite. Since both Euler characteristics and the characteristic power series of $\Lambda$-modules behave multiplicatively in exact sequences, it is enough to verify the final statement when $S$ is finite and when $\mathrm{Hom}(S, \mathbb{Q}_p/\mathbb{Z}_p) = \Lambda/(g(T))$. In the first case, the Euler characteristic is 1 and the characteristic ideal is $\Lambda$. The second case is clear from the above remarks about $Y$. ∎

Referring to the diagram at the beginning of section 3, we will denote $s_0$, $h_0$, and $g_0$ simply by $s$, $h$, and $g$.

**Lemma 4.3.** *Under the assumptions of theorem 4.1, we have*

$$|\mathrm{Sel}_E(F_\infty)_p^\Gamma| = |\mathrm{Sel}_E(F)_p||\ker(g)|/|E(F)_p|.$$

*Proof.* We have $|(\mathrm{Sel}_E(F_\infty)_p^\Gamma|/|\mathrm{Sel}_E(F)_p| = |\mathrm{coker}(s)|/|\ker(s)|$, where all the groups occurring are finite. By lemma 3.2, $\mathrm{coker}(h) = 0$. Thus, we have an exact sequence: $0 \to \ker(s) \to \ker(h) \to \ker(g) \to \mathrm{coker}(s) \to 0$. It follows that $|\mathrm{coker}(s)|/|\ker(s)| = |\ker(g)|/|\ker(h)|$. Now we use the fact that $E(F_\infty)_p$ is finite. Then

$$\ker(h) = H^1(\Gamma, E(F_\infty)_p) = (E(F_\infty)_p)_\Gamma$$

has the same order as $H^0(\Gamma, E(F_\infty)_p) = E(F)_p$. These facts give the formula in lemma 4.3. ∎

The proof of theorem 4.1 clearly rests now on studying $|\ker(g)|$. The results of section 3 allow us to study $\ker(r)$, factor by factor, where $r$ is the natural map

$$r : \mathcal{P}_E(F) \to \mathcal{P}_E(F_\infty).$$

It will be necessary for us to replace $\mathcal{P}_E(*)$ by a much smaller group. Let $\Sigma$ denote the set of primes of $F$ where $E$ has bad reduction or which divide $p$ or $\infty$. By lemma 3.3, we have $\ker(r_v) = 0$ if $v \notin \Sigma$. Let $\mathcal{P}_E^\Sigma(F) = \prod_v \mathcal{H}_E(F_v)$, where the product is over all primes of $F$ in $\Sigma$. We consider $\mathcal{P}_E^\Sigma(F)$ as a subgroup of $\mathcal{P}_E(F)$. Clearly, $\ker(r) \subseteq \mathcal{P}_E^\Sigma(F)$. Thus $|\ker(r)| = \prod_v |\ker(r_v)|$, where $v$ again varies over all primes in $\Sigma$. For $v|p$, the order of $\ker(r_v)$ is given in lemma 3.4. For $v \nmid p$, the remarks after the proof of lemma 3.3 show that $|\ker(r_v)| \sim c_v^{(p)}$. We then obtain the following result.

**Lemma 4.4.** *Assume that $E/F$ has good, ordinary reduction at all $v|p$. Then*
$$|\ker(r)| \sim \Big( \prod_{v \text{ bad}} c_v^{(p)} \Big) \Big( \prod_{v|p} |\widetilde{E}_v(f_v)_p|^2 \Big).$$

Now let $\mathcal{G}_E^\Sigma(F) = \operatorname{Im}\big( H^1(F_\Sigma/F, E[p^\infty]) \to \mathcal{P}_E^\Sigma(F) \big)$, where $F_\Sigma$ denotes the maximal extension of $F$ unramified outside of $\Sigma$. Then

$$\ker(g) = \ker(r) \cap \mathcal{G}_E^\Sigma(F).$$

We now recall a theorem of Cassels which states that $\mathcal{P}_E^\Sigma(F)/\mathcal{G}_E^\Sigma(F) \cong E(F)_p$. (We will sketch a proof of this later, using the Duality Theorem of Poitou and Tate.) It is interesting to consider theorem 4.1 in the case where $E(F)_p = 0$, which is of course true for all but finitely many primes $p$. Then, by Cassels' theorem, $\ker(g) = \ker(r)$. Lemmas 4.3, 4.4 then show that the right side of $\sim$ in theorem 4.1 is precisely $|\mathrm{Sel}_E(F_\infty)_p^\Gamma|$. Therefore, in this special case, by lemma 4.2, theorem 4.1 is equivalent to asserting that $(\mathrm{Sel}_E(F_\infty)_p)_\Gamma = 0$. It is an easy exercise to see that this in turn is equivalent to asserting that the $\Lambda$-module $X_E(F_\infty)$ has no finite, nonzero $\Lambda$-submodules. In section 5 we will give an example where $X_E(F_\infty)$ does have a finite, nonzero $\Lambda$-submodule. All the hypotheses of this section will hold, but of course $E(F)$ will have an element of order $p$.

The following general fact will be useful in the rest of the proof of theorem 4.1. We will assume that $G$ is a profinite group and that $A$ is a discrete, $p$-primary abelian group on which $G$ acts continuously.

**Lemma 4.5.** *Assume that $G$ has $p$-cohomological dimension $n \geq 1$ and that $A$ is a divisible group. Then $H^n(G, A)$ is a divisible group.*

*Proof.* Consider the exact sequence $0 \to A[p] \to A \xrightarrow{p} A \to 0$, where the map $A \xrightarrow{p} A$ is of course multiplication by $p$. This induces an exact sequence

$$H^n(G, A) \xrightarrow{p} H^n(G, A) \to H^{n+1}(G, A).$$

Since the last group is zero, $H^n(G, A)$ is divisible by $p$. The lemma follows because $H^n(G, A)$ is a $p$-primary group. ∎

We have actually already applied this lemma once, namely in the proof of proposition 2.4. We will apply it to some other cases. A good reference for the facts we use is [Se2]. Let $v$ be a nonarchimedean prime of $F$, $\eta$ a prime of $F_\infty$ lying above $v$. Then $\mathrm{Gal}((F_\infty)_\eta/F_v) \cong \mathbb{Z}_p$, as mentioned earlier. Thus, $G_{(F_\infty)_\eta}$ has $p$-cohomological dimension 1. Hence $H^1((F_\infty)_\eta, E[p^\infty])$ must be divisible, and consequently the same is true for $\mathcal{H}_E((F_\infty)_\eta)$. As another example, $\mathrm{Gal}(F_\Sigma/F)$ has $p$-cohomological dimension 2 if $p$ is any odd prime. Let $A_s = E[p^\infty] \otimes (\kappa^s)$, where $\kappa\colon \Gamma \to 1 + 2p\mathbb{Z}_p$ is an isomorphism and $s \in \mathbb{Z}$. ($A_s$ is something like a Tate twist of the $G_F$-module $E[p^\infty]$. One could even take $s \in \mathbb{Z}_p$.) It then follows that $H^2(F_\Sigma/F, A_s)$ is a divisible group.

**Lemma 4.6.** *Assume that* $\mathrm{Sel}_E(F_\infty)_p$ *is $\Lambda$-cotorsion. Then the map*

$$H^1(F_\Sigma/F_\infty, E[p^\infty]) \to \mathcal{P}_E^\Sigma(F_\infty)$$

*is surjective.*

*Remark.* We must define $\mathcal{P}_E^\Sigma(F_\infty)$ carefully. For any prime $v$ in $\Sigma$, we define

$$\mathcal{P}_E^{(v)}(F_\infty) = \varinjlim_n \mathcal{P}_E^{(v)}(F_n)$$

where $\mathcal{P}_E^{(v)}(F_n) = \prod_{v_n|v} \mathcal{H}_E((F_n)_{v_n})$ and $\mathcal{H}_E(*)$ is as defined at the beginning of section 3. The maps $\mathcal{P}_E^{(v)}(F_n) \to \mathcal{P}_E^{(v)}(F_{n+1})$ are easily defined, considering separately the case where $v_n$ is inert or ramified in $F_{n+1}/F_n$ (where one uses a restriction map) or where $v_n$ splits completely in $F_{n+1}/F_n$ (where one uses a "diagonal" map). If $v$ is nonarchimedean, then $v$ is finitely decomposed in $F_\infty/F$ and one can more simply define $\mathcal{P}_E^{(v)}(F_\infty) = \prod_{\eta|v} \mathcal{H}_E((F_\infty)_\eta)$, where $\eta$ runs over the finite set of primes of $F_\infty$ lying over $v$. If $v$ is archimedean, then $v$ splits completely in $F_\infty/F$. We know that $\mathrm{Im}(\kappa_{v_n}) = 0$ for $v_n|v$. Thus, $\mathcal{H}_E((F_n)_{v_n}) = \mathcal{H}_E(F_v) = H^1(F_v, E[p^\infty])$. Usually, this group is zero. But it can be nonzero if $p = 2$ and $F_v = \mathbb{R}$. In fact,

$$H^1(F_v, E[2^\infty]) \cong E(F_v)/E(F_v)_{\mathrm{con}},$$

where $E(F_v)_{\mathrm{con}}$ denotes the connected component of the identity of $E(F_v)$. Therefore, obviously $H^1(F_v, E[2^\infty])$ has order 1 or 2. The order is 2 if $E[2]$ is contained in $E(F_v)$. We have

$$\mathcal{P}_E^{(v)}(F_n) \cong H^1(F_v, E[2^\infty]) \otimes \mathbb{Z}_2[\mathrm{Gal}(F_n/F)],$$

which is either zero or isomorphic to $(\mathbb{Z}/2\mathbb{Z})[\mathrm{Gal}(F_n/F)]$. In each of the above cases, $\mathcal{P}_E^{(v)}(F_\infty)$ can be regarded naturally as a $\Lambda$-module. If $v$ is nonarchimedean then the remarks following lemma 4.5 show that, as a group, $\mathcal{P}_E^{(v)}(F_\infty)$ is divisible. If $v$ is archimedean, then usually $\mathcal{P}_E^{(v)}(F_\infty) = 0$. But, if $p = 2$, $F_v = \mathbb{R}$, and $E[2]$ is contained in $E(F_v)$, then one sees that $\mathcal{P}_E^{(v)}(F_\infty) \cong \mathrm{Hom}(\Lambda/2\Lambda, \mathbb{Z}/2\mathbb{Z})$ as a $\Lambda$-module. (One uses the fact that $\mathcal{P}_E^{(v)}(F_\infty)^{\Gamma_n} \cong \mathcal{P}_E^{(v)}(F_n)$ for all $n \geq 0$ and the structure of $\mathcal{P}_E^{(v)}(F_n)$ mentioned above.) Finally, we define $\mathcal{P}_E^\Sigma(F_\infty) = \prod_{v \in \Sigma} \mathcal{P}_E^{(v)}(F_\infty)$.

*Proof of Lemma 4.6.* We can regard $\mathcal{P}_E^\Sigma(F_\infty)$ as a $\Lambda$-module. The idea of the proof is to show that the image of the above map is a $\Lambda$-submodule of $\mathcal{P}_E^\Sigma(F_\infty)$ with finite index and that any such $\Lambda$-submodule must be $\mathcal{P}_E^\Sigma(F_\infty)$. We will explain the last point first. If $p$ is odd, the remarks above show that each factor in $\mathcal{P}_E^\Sigma(F_\infty)$ is divisible. Hence $\mathcal{P}_E^\Sigma(F_\infty)$ is divisible and therefore has no proper subgroups of finite index. If $p = 2$, one has to observe that the factor $\mathcal{P}_E^{(v)}(F_\infty)$ of $\mathcal{P}_E^\Sigma(F_\infty)$ coming from an archimedean prime $v$ of $F$ is a $\Lambda$-module whose Pontryagin dual is either zero or isomorphic to $(\Lambda/2\Lambda)$. Since $\Lambda/2\Lambda$ has no nonzero, finite $\Lambda$-submodules, we see that $\mathcal{P}_E^{(v)}(F_\infty)$ has no proper $\Lambda$-submodules of finite index. Since the factors $\mathcal{P}_E^{(v)}(F_\infty)$ for nonarchimedean $v$ are still divisible, it follows again that $\mathcal{P}_E^\Sigma(F_\infty)$ has no proper $\Lambda$-submodules of finite index.

Now we will prove that the image of the map in the lemma has finite index. (It is clearly a $\Lambda$-submodule.) To give the idea of the proof, assume first that $\mathrm{Sel}_E(F_n)_p$ is finite for all $n \geq 0$. Then the cokernel of the map $H^1(F_\Sigma/F_n, E[p^\infty]) \to \mathcal{P}_E^\Sigma(F_n)$ is isomorphic to $E(F_n)_p$ by a theorem of Cassels. But $|E(F_n)_p|$ is bounded since it is known that $E(F_\infty)_p$ is finite. It clearly follows that the cokernel of the corresponding map over $F_\infty$ is also finite. To give the proof in general, we use a trick of twisting the Galois module $E[p^\infty]$. We let $A_s$ be defined as above, where $s \in \mathbb{Z}$. As $G_{F_\infty}$-modules, $A_s = E[p^\infty]$. Thus, $H^1(F_\infty, A_s) = H^1(F_\infty, E[p^\infty])$. But the action of $\Gamma$ changes in a simple way, namely $H^1(F_\infty, A_s) = H^1(F_\infty, E[p^\infty]) \otimes (\kappa^s)$. Now we can define Selmer groups for $A_s$ as suggested at the end of section 2. One just imitates the description of the $p$-Selmer group for $E$. For the local condition at $v$ dividing $p$, one uses $C_v \otimes (\kappa^s)$. For $v$ not dividing $p$, we require 1-cocycles to be locally trivial. We let $S_{A_s}(F_n)$, $S_{A_s}(F_\infty)$ denote the Selmer groups defined in this way. Then $S_{A_s}(F_\infty) = \mathrm{Sel}_E(F_\infty)_p \otimes (\kappa^s)$ as $\Lambda$-modules. Now we are assuming that $\mathrm{Sel}_E(F_\infty)_p$ is $\Lambda$-cotorsion. It is not hard to show from this that for all but finitely many values of $s$, $S_{A_s}(F_\infty)^{\Gamma_n}$ will be finite for all $n \geq 0$. Since there is a map $S_{A_s}(F_n) \to S_{A_s}(F_\infty)^{\Gamma_n}$ with finite kernel, it follows that $S_{A_s}(F_n)$ is finite for all $n \geq 0$. There is also a variant of Cassels' theorem for $A_s$: the cokernel of the global-to-local map for the $G_{F_n}$-module $A_s$ is isomorphic to $H^0(F_n, A_{-s})$. But this last group is finite and has order bounded by $|E(F_\infty)_p|$. The surjectivity of the global-to-local map for $A_s$ over $F_\infty$ follows just as before. Lemma 4.6 follows since $A_s \cong E[p^\infty]$ as $G_{F_\infty}$-modules. (Note: the variant of Cassels' theorem is a consequence of proposition 4.13. It may be necessary to exclude one more value of $s$.) ∎

The following lemma, together with lemmas 4.2–4.4 implies theorem 4.1.

**Lemma 4.7.** *Under the assumptions of theorem 4.1, we have*

$$|\ker(g)| = |\ker(r)| \, |(\mathrm{Sel}_E(F_\infty)_p)_\Gamma| / |E(F)_p|.$$

*Proof.* By lemma 4.6, the following sequence is exact:

$$0 \to \mathrm{Sel}_E(F_\infty)_p \to H^1(F_\Sigma/F_\infty, E[p^\infty]) \to \mathcal{P}_E^\Sigma(F_\infty) \to 0.$$

Now $\Gamma$ acts on these groups. We can take the corresponding cohomology sequence obtaining

$$H^1(F_\Sigma/F_\infty, E[p^\infty])^\Gamma \to \mathcal{P}_E^\Sigma(F_\infty)^\Gamma \to ((\mathrm{Sel}_E(F_\infty)_p)_\Gamma \to H^1(F_\Sigma/F_\infty, E[p^\infty])_\Gamma.$$

In the appendix, we will give a proof that the last term is zero. Thus we get the following commutative diagram with exact rows and columns.

$$\begin{array}{ccccccc}
H^1(F_\Sigma/F, E[p^\infty]) & \xrightarrow{a} & \mathcal{P}_E^\Sigma(F) & \longrightarrow & \mathcal{P}_E^\Sigma(F)/\mathcal{G}_E^\Sigma(F) & \longrightarrow & 0 \\
\downarrow & & \downarrow & & \downarrow & & \\
H^1(F_\Sigma/F_\infty, E[p^\infty])^\Gamma & \xrightarrow{b} & \mathcal{P}_E^\Sigma(F_\infty)^\Gamma & \longrightarrow & (\mathrm{Sel}_E(F_\infty)_p)_\Gamma & \longrightarrow & 0 \\
\downarrow & & \downarrow & & \downarrow & & \\
0 & & 0 & & 0 & &
\end{array}$$

The exactness of the first row is clear. The remark above gives the exactness of the second row. The surjectivity of the first vertical arrow is because $\Gamma$ has $p$-cohomological dimension 1. The surjectivity of the second vertical arrow can be verified similarly. One must consider each $v \in \Sigma$ separately, showing that $\mathcal{P}_E^{(v)}(F) \to \mathcal{P}_E^{(v)}(F_\infty)^\Gamma$ is surjective. One must take into account the fact that $v$ can split completely in $F_n/F$ for some $n$. But then it is easy to see that $\mathcal{P}_E^{(v)}(F) \xrightarrow{\sim} \mathcal{P}_E^{(v)}(F_n)^{\mathrm{Gal}(F_n/F)}$. One then uses the fact that $\mathrm{Gal}((F_\infty)_\eta/(F_n)_{v_n})$ has $p$-cohomological dimension 1, looking at the maps $r_{v_n}$ for $v \nmid p$ or $d_{v_n}$ for $v|p$. For archimedean $v$, one easily verifies that $\mathcal{P}_E^{(v)}(F) \xrightarrow{\sim} \mathcal{P}_E^{(v)}(F_\infty)^\Gamma$. The surjectivity of the third vertical arrow follows. It is also clear that $\mathrm{Im}(a)$ is mapped surjectively to $\mathrm{Im}(b)$. We then obtain the following commutative diagram

$$\begin{array}{ccccccccc}
0 & \longrightarrow & \mathcal{G}_E^\Sigma(F) & \longrightarrow & \mathcal{P}_E^\Sigma(F) & \longrightarrow & \mathcal{P}_E^\Sigma(F)/\mathcal{G}_E^\Sigma(F) & \longrightarrow & 0 \\
& & \downarrow{g} & & \downarrow{r} & & \downarrow{t} & & \\
0 & \longrightarrow & \mathrm{Im}(b) & \longrightarrow & \mathcal{P}_E^\Sigma(F_\infty)^\Gamma & \longrightarrow & (\mathrm{Sel}_E(F_\infty)_p)_\Gamma & \longrightarrow & 0 \\
& & \downarrow & & \downarrow & & \downarrow & & \\
& & 0 & & 0 & & 0 & &
\end{array}$$

¿From the snake lemma, we then obtain $0 \to \ker(g) \to \ker(r) \to \ker(t) \to 0$. Thus, $|\ker(g)| = |\ker(r)|/|\ker(t)|$. Combining this with Cassels' theorem and the obvious value of $|\ker(t)|$ proves lemma 4.7. ∎

The last commutative diagram, together with Cassels' theorem, gives the following consequence which will be quite useful in the discussion of various examples in section 5. A more general result will be proved in the appendix.

**Proposition 4.8.** *Assume that $E$ is an elliptic curve defined over $F$ with good, ordinary reduction at all primes of $F$ lying over $p$. Assume that $\mathrm{Sel}_E(F)_p$ is finite and that $E(F)_p = 0$. Then $\mathrm{Sel}_E(F_\infty)_p$ has no proper $\Lambda$-submodules of finite index. In particular, if $\mathrm{Sel}_E(F_\infty)_p$ is nonzero, then it must be infinite.*

*Proof.* We have the map $t\colon E(F)_p \to \mathrm{Sel}_E(F_\infty)_\Gamma$, which is surjective. Since $E(F)_p = 0$, it follows that $(\mathrm{Sel}_E(F_\infty)_p)_\Gamma = 0$ too. Suppose that $\mathrm{Sel}_E(F_\infty)_p$ has a finite, nonzero $\Lambda$-module quotient $M$. Then $M$ is just a nonzero, finite, abelian $p$-group on which $\Gamma$ acts. Obviously, $M_\Gamma \neq 0$. But $M_\Gamma$ is a homomorphic image of $(\mathrm{Sel}_E(F_\infty)_p)_\Gamma$, which is impossible. ∎

Theorem 4.1 gives a conjectural relationship of $f_E(0)$ to the value of the Hasse-Weil $L$-function $L(E/F, s)$ at $s = 1$. This is based on the Birch and Swinnerton-Dyer conjecture for $E$ over $F$, for the case where $E(F)$ is assumed to be finite. We assume of course that $\mathrm{III}_E(F)_p$ is finite and hence so is $\mathrm{Sel}_E(F)_p = \mathrm{III}_E(F)_p$. We also assume that $L(E/F, s)$ has an analytic continuation to $s = 1$. The conjecture then asserts that $L(E/F, 1) \neq 0$ and that for a suitably defined period $\Omega(E/F)$, the value $L(E/F, 1)/\Omega(E/F)$ is rational and

$$L(E/F, 1)/\Omega(E/F) \sim (\prod_{v \text{ bad}} c_v^{(p)}) |\mathrm{Sel}_E(F)_p| / |E(F)_p|^2 \, .$$

As before, $\sim$ means that the two sides have the same $p$-adic valuation. If $\mathcal{O}$ denotes the ring of integers in $F$, then one must choose a minimal Weierstrass equation for $E$ over $\mathcal{O}_{(p)}$, the localization of $\mathcal{O}$ at $p$, to define $\Omega(E/F)$ (as a product of periods over the archimedean primes of $F$). For $v|p$, the Euler factor for $v$ in $L(E/F, s)$ is

$$(1 - \alpha_v N(v)^{-s})(1 - \beta_v N(v)^{-s}),$$

where $\alpha_v, \beta_v$ are as defined just after theorem 4.1. Recall that $\alpha_v \in \mathbb{Z}_p^\times$. (We are assuming that $E$ has good, ordinary reduction at all $v|p$.) Then we have

$$|\widetilde{E}(f_v)_p| \sim (1 - \alpha_v) \sim (1 - \alpha_v^{-1}) = (1 - \beta_v N(v)^{-1}).$$

The last quantity is one factor in the Euler factor for $v$, evaluated at $s = 1$. Thus, theorem 4.1 conjecturally states that

$$f_E(0) \sim (\prod_{v|p} (1 - \beta_v N(v)^{-1})^2) L(E/F, 1)/\Omega(E/F).$$

For $F = \mathbb{Q}$, one should compare this with conjecture 1.13.

As we mentioned in the introduction, there is a result of P. Schneider (generalizing a result of B. Perrin-Riou for elliptic curves with complex multiplication) which concerns the behavior of $f_E(T)$ at $T = 0$. We assume that $E$ is an elliptic curve/$F$ with good, ordinary reduction at all primes of $F$ lying over $p$, that $p$ is odd and that $F \cap \mathbb{Q}_\infty = \mathbb{Q}$ (to slightly simplify the statement). Let $r = \mathrm{rank}(E(F))$. We will state the result for the case where $r = 1$ and $\mathrm{III}_E(F)_p$ is finite. (Then $\mathrm{Sel}_E(F)_p$ has $\mathbb{Z}_p$-corank 1.) Since then $T|f_E(T)$, one can write $f_E(T) = T g_E(T)$, where $g_E(T) \in \Lambda$. The result is that

$$g_E(0) \sim \frac{h_p(P)}{p} (\prod_{v \text{ bad}} c_v^{(p)})(\prod_{v|p} |\widetilde{E}_v(f_v)_p|^2) |\mathrm{III}_E(F)_p|/|E(F)_p^2|.$$

Here $P \in E(F)$ is a generator of $E(F)/E(F)_{\mathrm{tors}}$ and $h_p(P)$ is its analytic $p$-adic height. (See [Sch2] and the references there for the definition of $h_p(P)$.) The other factors are as in theorem 4.1. Conjecturally, one should have $h_p(P) \neq 0$. This would mean that $f_E(T)$ has a simple zero at $T = 0$. But if $h_p(P) = 0$, the result means that $g_E(0) = 0$, i.e., $T^2 | f_E(T)$. If $F = \mathbb{Q}$ and $E$ is modular, then B. Perrin-Riou [Pe3] has proven an analogue of a theorem of Gross and Zagier for the $p$-adic $L$-function $L_p(E/\mathbb{Q}, s)$. Assume that $L(E/\mathbb{Q}, s)$ has a simple zero at $s = 1$. Then a result of Kolyvagin shows that rank$(E(\mathbb{Q})) = 1$ and $\text{III}_E(\mathbb{Q})$ is finite. Assume that $P$ generates $E(\mathbb{Q})/E(\mathbb{Q})_{\text{tors}}$. Assume that $h_p(P) \neq 0$. Perrin-Riou's result asserts that $L_p(E/\mathbb{Q}, s)$ also has a simple zero at $s = 1$ and that

$$L_p'(E/\mathbb{Q}, 1)/h_p(P) = (1 - \beta_p p^{-1})^2 L'(E/\mathbb{Q}, 1)/\Omega_E h_\infty(P)$$

where $h_\infty(P)$ is the canonical height of $P$. If one assumes the validity of the Birch and Swinnerton-Dyer conjecture, then this result and Schneider's result are compatible with conjecture 1.13.

The proof of theorem 4.1 can easily be adapted to the case where $E$ has multiplicative reduction at some primes of $F$ lying over $p$. One then obtains a special case of a theorem of J. Jones [Jo]. Jones determines the $p$-adic valuation of $(f_E(T)/T^r)|_{T=0}$, where $r = \text{rank}(E(F))$, generalizing the results of P. Schneider. He studies certain natural $\Lambda$-modules which can be larger, in some sense, than $\text{Sel}_E(F_\infty)_p$. Their characteristic ideal will contain $T^e f_E(T)$, where $e$ is the number of primes of $F$ where $E$ has split, multiplicative reduction. This is an example of the phenomenon of "trivial zeros". Another example of this phenomenon is the $\Lambda$-module $S_\infty$ in the case where $p$ splits in an imaginary quadratic field $F$. As we explained in the introduction, $S_\infty^\Gamma$ is infinite. That is, a generator of its characteristic ideal will vanish at $T = 0$. For a general discussion of this phenomenon, we refer the reader to [Gr4].

To state the analogue of theorem 4.1, we assume that $\text{Sel}_E(F)_p$ is finite, that $E$ has either good, ordinary or multiplicative reduction at all primes of $F$ over $p$, and that $\log_p(N_{F_v/\mathbb{Q}_p}(q_E^{(v)})) \neq 0$ for all $v$ lying over $p$ where $E$ has split, multiplicative reduction. (As in section 3, $q_E^{(v)}$ denotes the Tate period for $E$ over $F_v$.) Under these assumptions, $\ker(r_v)$ will be finite for all $v | p$. It follows from proposition 3.7 that $\text{Sel}_E(F_\infty)_p^\Gamma$ will be finite and hence $\text{Sel}_E(F_\infty)_p$ will be $\Lambda$-cotorsion. In theorem 4.1, the only necessary change is to replace the factor $|\widetilde{E}_v(f_v)_p|^2$ for those $v | p$ where $E$ has multiplicative reduction by the factor $|\ker(r_v)|/c_v^{(p)}$. (Note that the factor $c_v^{(p)}$ for such $v$ will occur in $\prod_{v \text{ bad}} c_v^{(p)}$.) The analogue of theorem 4.1 can be expressed as

$$f_E(0) \sim (\prod_{v|p} l_v)(\prod_{v \text{ bad}} c_v^{(p)}) |\text{Sel}_E(F)_p|/|E(F)_p|^2.$$

If $E$ has good, ordinary reduction at $v$, then $l_v = |\widetilde{E}_v(f_v)_p|^2$. Assume that $E$ has nonsplit, multiplicative reduction at $v$. If $p$ is odd, then both $|\ker(r_v)|$ and $c_v^{(p)}$ are equal to 1. If $p = 2$, then $|\ker(r_v)| = 2c_v^{(p)}$. (Recalling the discussion concerning $\ker(r_v)$ after the proof of proposition 3.6, the 2 corresponds to $|\ker(b_v)|$, and the $c_v^{(p)}$ corresponds to $|\ker(a_v)| = [\mathrm{Im}(\lambda_v) : \mathrm{Im}(\kappa_v)]$. In the case of good, ordinary reduction at $v$, both $\ker(a_v)$ and $\ker(b_v)$ have order $|\widetilde{E}_v(f_v)_p|$.) Thus, if $E$ has nonsplit, multiplicative reduction at $v$, one can take $l_v = 2$ (for any prime $p$). We remark that the Euler factor for $v$ in $L(E/F, s)$ is $(1 + N(v)^{-s})^{-1}$. One should take $\alpha_v = -1$, $\beta_v = 0$. Perhaps this factor $l_v = 2$ should be thought of as $(1 - \alpha_v^{-1})$. (This is suggested by the fact that, for a modular elliptic curve $E$ defined over $F = \mathbb{Q}$, the $p$-adic $L$-function constructed in [M-T-T] has a factor $(1 - \alpha_p^{-1})$ in its interpolation property when $E$ has multiplicative reduction at $p$. This is in place of $(1 - \alpha_p^{-1})^2 = (1 - \beta_p p^{-1})^2$ when $E$ has good, ordinary reduction at $p$.)

Finally, assume that $E$ has split, multiplicative reduction at $v$. (Then $(1 - \alpha_v^{-1})$ would be zero.) We have $c_v^{(p)} = \mathrm{ord}_v(q_E^{(v)})$. If we let $\mathbb{Q}_p^{\mathrm{unr}}$ denote the unramified $\mathbb{Z}_p$-extension of $\mathbb{Q}_p$ and $\mathbb{Q}_p^{\mathrm{cyc}}$ denote the cyclotomic $\mathbb{Z}_p$-extension of $\mathbb{Q}_p$, then we should take

$$l_v = \frac{\log_p(N_{F_v/\mathbb{Q}_p}(q_E^{(v)}))}{\mathrm{ord}_p(N_{F_v/\mathbb{Q}_p}(q_E^{(v)}))} \cdot \frac{[F_v \cap \mathbb{Q}_p^{\mathrm{unr}} : \mathbb{Q}_p]}{2p[F_v \cap \mathbb{Q}_p^{\mathrm{cyc}} : \mathbb{Q}_p]}.$$

(Again, we refer to the discussion of $\ker(r_v)$ following proposition 3.6. This time, $\ker(a_v) = 0$ and $\ker(r_v) \cong \ker(b_v)$.) We will give another way to define $l_v$, at least up to a $p$-adic unit, which comes directly from the earlier discussion of $\ker(r_v)$. Let $F_v^{\mathrm{cyc}}$ and $F_v^{\mathrm{unr}}$ denote the cyclotomic and the unramified $\mathbb{Z}_p$-extensions of $F_v$. Fix isomorphisms

$$\theta_{F_v}^{\mathrm{cyc}} : \mathrm{Gal}(F_v^{\mathrm{cyc}}/F_v) \xrightarrow{\sim} \mathbb{Z}_p, \qquad \theta_{F_v}^{\mathrm{unr}} : \mathrm{Gal}(F_v^{\mathrm{unr}}/F_v) \xrightarrow{\sim} \mathbb{Z}_p.$$

Then $l_v \sim \theta_{F_v}^{\mathrm{cyc}}(\mathrm{rec}_{F_v}(q_E^{(v)})|_{F_v^{\mathrm{cyc}}})/\theta_{F_v}^{\mathrm{unr}}(\mathrm{rec}_{F_v}(q^{(v)})|_{F_v^{\mathrm{unr}}})$. The value of $l_v$ given above comes from choosing specific isomorphisms.

**Appendix to Section 4.** We will give a proof of the following important result, which will allow us to justify the assertion used in the proof of lemma 4.7 that, under the hypotheses of theorem 4.1, $H^1(F_\Sigma/F_\infty, E[p^\infty])_\Gamma = 0$. Later, we will prove a rather general form of Cassels' theorem as well as a generalization of proposition 4.8.

**Proposition 4.9.** *Assume that* $\mathrm{Sel}_E(F_\infty)_p$ *is $\Lambda$-cotorsion. Then the $\Lambda$-module* $H^1(F_\Sigma/F_\infty, E[p^\infty])$ *has no proper $\Lambda$-submodules of finite index.*

In the course of the proof, we will show that $H^1(F_\Sigma/F_\infty, E[p^\infty])$ has $\Lambda$-corank $[F : \mathbb{Q}]$ and also that $H^2(F_\Sigma/F_\infty, E[p^\infty])$ is $\Lambda$-cotorsion. For odd $p$, these results are contained in [Gr2]. (See section 7 there.) For $p = 2$, one can modify the arguments given in that article. However, we will present a rather different approach here which has the advantage of avoiding the use of a spectral sequence. In either approach, the crucial point is that the group

$$R^2(F_\Sigma/F_\infty, E[p^\infty]) = \ker\left(H^2(F_\Sigma/F_\infty, E[p^\infty]) \to \prod_{\eta|\infty} H^2((F_\infty)_\eta, E[p^\infty])\right)$$

is zero, under the assumption that $\mathrm{Sel}_E(F_\infty)_p$ is $\Lambda$-cotorsion. (Note: It probably seems more natural to take the product over all $\eta$ lying over primes in $\Sigma$. However, if $\eta$ is nonarchimedean, then $G_{(F_\infty)_\eta}$ has $p$-cohomological dimension 1 and hence $H^2((F_\infty)_\eta, E[p^\infty]) = 0$.)

First of all, we will determine the $\Lambda$-corank of $\mathcal{P}_E^{(v)}(F_\infty)$. Now $\mathcal{P}_E^{(v)}(F_\infty)$ is $\Lambda$-cotorsion if $v \nmid p$. This is clear if $v$ is archimedean because $\mathcal{P}_E^{(v)}(F_\infty)$ then has exponent 2. (It is zero if $p$ is odd.) If $v$ is nonarchimedean, then $\mathcal{P}_E^{(v)}(F) = H^1(F_v, E[p^\infty])$ is finite. The map $\mathcal{P}_E^{(v)}(F_v) \to \mathcal{P}_E^{(v)}(F_\infty)^\Gamma$ is surjective. Hence $\mathcal{P}_E^{(v)}(F_\infty)^\Gamma$ is finite, which suffices to prove that $\mathcal{P}_E^{(v)}(F_\infty)$ is $\Lambda$-cotorsion, using Fact (2) about $\Lambda$-modules mentioned in section 3. Alternatively, one can refer to proposition 2 of [Gr2], which gives a more precise result concerning the structure of $\mathcal{P}_E^{(v)}(F_\infty)$. Assume $v|p$. Let $\Gamma_v \subseteq \Gamma$ be the decomposition group for any prime $\eta$ of $F_\infty$ lying over $v$. Then by proposition 1 of [Gr2], $H^1((F_\infty)_\eta, E[p^\infty])$ has corank equal to $2[F_v{:}\mathbb{Q}_p]$ over the ring $\mathbb{Z}_p[[\Gamma_v]]$. Also, $H^1((F_\infty)_\eta, C_v)$ has corank $[F_v{:}\mathbb{Q}_p]$. Both of these facts could be easily proved using lemma 2.3, applied to the layers in the $\mathbb{Z}_p$-extension $(F_\infty)_\eta/F_v$. Consequently, $\mathcal{H}_E((F_\infty)_\eta)$ has $\mathbb{Z}_p[[\Gamma_v]]$-corank equal to $[F_v{:}\mathbb{Q}_p]$. It follows that $\mathcal{P}_E^{(v)}(F_\infty)$ has $\Lambda$-corank equal to $[F_v{:}\mathbb{Q}_p]$. Combining these results, we find that

$$\mathrm{corank}_\Lambda(\mathcal{P}_E^\Sigma(F_\infty)) = [F{:}\mathbb{Q}],$$

using the fact that $\sum\limits_{v|p}[F_v{:}\mathbb{Q}_p] = [F{:}\mathbb{Q}]$.

Secondly, we consider the coranks of the $\Lambda$-modules $H^1(F_\Sigma/F_\infty, E[p^\infty])$ and $H^2(F_\Sigma/F_\infty, E[p^\infty])$. These are related by the equation

$$\mathrm{corank}_\Lambda(H^1(F_\Sigma/F_\infty, E[p^\infty]) = \mathrm{corank}_\Lambda(H^2(F_\Sigma/F_\infty, E[p^\infty]) + \delta,$$

where $\delta = \sum\limits_{v|\infty}[F_v{:}\mathbb{R}] = [F{:}\mathbb{Q}]$. As a consequence, we have the inequalities

$$\mathrm{corank}_\Lambda(H^1(F_\Sigma/F_\infty, E[p^\infty])) \geq [F{:}\mathbb{Q}].$$

(For more discussion of this relationship, see [Gr2], section 4. It is essentially the fact that $-\delta$ is the Euler characteristic for the $\mathrm{Gal}(F_\Sigma/F_\infty)$-module $E[p^\infty]$ together with the fact that $H^0(F_\Sigma/F_\infty, E[p^\infty])$ is clearly $\Lambda$-cotorsion. This Euler characteristic of $\Lambda$-coranks is in turn derived from the fact that

$$\sum_{i=0}^{2}(-1)^i \mathrm{corank}_{\mathbb{Z}_p}(H^i(F_\Sigma/F_n, E[p^\infty])) = -\delta p^n$$

for all $n \geq 0$. That is, $-\delta p^n$ is the Euler characteristic for the $\mathrm{Gal}(F_\Sigma/F_n)$-module $E[p^\infty]$.) Recalling the exact sequence

$$0 \to \mathrm{Sel}_E(F_\infty)_p \to H^1(F_\Sigma/F_\infty, E[p^\infty]) \to \mathcal{G}_E^\Sigma(F_\infty) \to 0,$$

we see that $\mathrm{Sel}_E(F_\infty)_p$ is $\Lambda$-cotorsion if and only if $H^1(F_\Sigma/F_\infty, E[p^\infty])$ and $\mathcal{G}_E^\Sigma(F_\infty)$ have the same $\Lambda$-corank, both equal to $[F:\mathbb{Q}]$. (The last equality is because $[F:\mathbb{Q}]$ is a lower bound for the $\Lambda$-corank of $H^1(F_\Sigma/F_\infty, E[p^\infty])$ and an upper bound for the $\Lambda$-corank of $\mathcal{G}_E^\Sigma(F_\infty)$ (which is a $\Lambda$-submodule of $\mathcal{P}_E^\Sigma(F_\infty)$). Thus, if we assume that $\mathrm{Sel}_E(F_\infty)_p$ is $\Lambda$-cotorsion, then it follows that $H^1(F_\Sigma/F_\infty, E[p^\infty])$ has $\Lambda$-corank $[F:\mathbb{Q}]$ and that $H^2(F_\Sigma/F_\infty, E[p^\infty])$ has $\Lambda$-corank 0 (and hence is $\Lambda$-cotorsion). By lemma 4.6, we already would know that $\mathcal{G}_E^\Sigma(F_\infty)$ has $\Lambda$-corank $[F:\mathbb{Q}]$.

We will use a version of Shapiro's Lemma. Let $\mathcal{A} = \mathrm{Hom}(\Lambda, E[p^\infty])$. We consider $\mathcal{A}$ as a $\Lambda$-module as follows: if $\phi \in \mathcal{A}$ and $\theta \in \Lambda$, then $\theta\phi$ is defined by $(\theta\phi)(\lambda) = \phi(\theta\lambda)$ for all $\lambda \in \Lambda$. The Pontryagin dual of $\mathcal{A}$ is $\Lambda^2$ and so $\mathcal{A}$ has $\Lambda$-corank 2. We define a $\Lambda$-linear action of $\mathrm{Gal}(F_\Sigma/F)$ on $\mathcal{A}$ as follows: if $\phi \in \mathcal{A}$ and $g \in \mathrm{Gal}(F_\Sigma/F)$, then $g(\phi)$ is defined by $g(\phi)(\lambda) = g(\phi(\widetilde{\kappa}(g)^{-1}\lambda))$ for all $\lambda \in \Lambda$. Here $\widetilde{\kappa}$ is defined as the composite

$$\mathrm{Gal}(F_\Sigma/F) \to \Gamma \to \Lambda^\times$$

where the second map is just the natural inclusion of $\Gamma$ in its completed group ring $\Lambda$. The above definition is just the usual way to define the action of a group on $\mathrm{Hom}(*, *)$, where we let $\mathrm{Gal}(F_\Sigma/F)$ act on $\Lambda$ by $\widetilde{\kappa}$ and on $E[p^\infty]$ as usual. The $\Lambda$-linearity is easily verified, using the fact that $\Lambda$ is a commutative ring. For any $\theta \in \Lambda$, we will let $\mathcal{A}[\theta]$ denote the kernel of the map $\mathcal{A} \xrightarrow{\theta} \mathcal{A}$, which is just multiplication by $\theta$. Then clearly

$$\mathcal{A}[\theta] \cong \mathrm{Hom}(\Lambda/\Lambda\theta, E[p^\infty]).$$

Let $\kappa : \Gamma \to 1 + 2p\mathbb{Z}_p$ be a fixed isomorphism. If $s \in \mathbb{Z}$ (or in $\mathbb{Z}_p$), then the homomorphism $\kappa^s : \Gamma \to 1 + 2p\mathbb{Z}_p$ induces a homomorphism $\sigma_s : \Lambda \to \mathbb{Z}_p$ of $\mathbb{Z}_p$-algebras. If we write $\Lambda = \mathbb{Z}_p[[T]]$, where $T = \gamma - 1$ as before, then $\sigma_s$ is defined by $\sigma_s(T) = \kappa^s(\gamma) - 1 \in p\mathbb{Z}_p$. We have $\ker(\sigma_s) = (\theta_s)$, where we have let $\theta_s = (T - (\kappa^s(\gamma) - 1))$. Then $\Lambda/\Lambda\theta_s \cong \mathbb{Z}_p(\kappa^s)$, a $\mathbb{Z}_p$-module of rank 1 on which $\mathrm{Gal}(F_\Sigma/F)$ acts by $\kappa^s$. Then

$$\mathcal{A}[\theta_s] \cong \mathrm{Hom}(\mathbb{Z}_p(\kappa^s), E[p^\infty]) \cong E[p^\infty] \otimes (\kappa^{-s}) = A_{-s}$$

as $\mathrm{Gal}(F_\Sigma/F)$-modules.

The version of Shapiro's Lemma that we will use is the following.

**Proposition 4.10.** *For all $i \geq 0$, $H^i(F_\Sigma/F_\infty, E[p^\infty]) \cong H^i(F_\Sigma/F, \mathcal{A})$ as $\Lambda$-modules.*

*Remark.* The first cohomology group is a $\Lambda$-module by virtue of the natural action of $\Gamma$ on $H^i(F_\Sigma/F_\infty, E[p^\infty])$; the second cohomology group is a $\Lambda$-module by virtue of the $\Lambda$-module structure on $\mathcal{A}$.

*Proof.* We let $A$ denote $E[p^\infty]$. The map $\phi \to \phi(1)$, for each $\phi \in \mathcal{A}$, defines a $\mathrm{Gal}(F_\Sigma/F_\infty)$-equivariant homomorphism $\mathcal{A} \to A$. The isomorphism in the proposition is defined by

$$H^i(F_\Sigma/F, \mathcal{A}) \xrightarrow{\mathrm{rest.}} H^i(F_\Sigma/F_\infty, \mathcal{A}) \to H^i(F_\Sigma/F_\infty, A).$$

One can verify that this composite map is a $\Lambda$-homomorphism as follows. $\mathrm{Gal}(F_\Sigma/F_\infty)$ acts trivially on $\Lambda$. We therefore have a canonical isomorphism

$$H^i(F_\Sigma/F_\infty, \mathcal{A}) \cong \mathrm{Hom}(\Lambda, H^i(F_\Sigma/F_\infty, A)) \tag{7}$$

The image of the restriction map in (7) is contained in $H^i(F_\Sigma/F_\infty, \mathcal{A})^\Gamma$, which corresponds under (7) to $\mathrm{Hom}_\Gamma(\Lambda, H^i(F_\Sigma/F_\infty, A))$. The action of $\Gamma$ on $\Lambda$ is given by $\widetilde{\kappa}$. But this is simply the usual structure of $\Lambda$ as a $\Lambda$-module, restricted to $\Gamma \subseteq \Lambda$. Thus, by continuity, we have

$$H^i(F_\Sigma/F_\infty, \mathcal{A})^\Gamma \cong \mathrm{Hom}_\Lambda(\Lambda, H^i(F_\Sigma/F_\infty, A))$$

under (7). Now $\mathrm{Hom}_\Lambda(\Lambda, H^i(F_\Sigma/F_\infty, A)) \cong H^i(F_\Sigma/F_\infty, A)$ as $\Lambda$-modules, under the map defined by evaluating a homomorphism at $\lambda = 1$.

To verify that the map $H^i(F_\Sigma/F, \mathcal{A}) \to H^i(F_\Sigma/F_\infty, A)$ is bijective, note that both groups and the map are direct limits:

$$H^i(F_\Sigma/F, \mathcal{A}) = \varinjlim_n H^i(F_\Sigma/F, \mathcal{A}[\theta^{(n)}]),$$
$$H^i(F_\Sigma/F_\infty, A) = \varinjlim_n H^i(F_\Sigma/F_n, A).$$

Here $\theta^{(n)} = (1+T)^{p^n} - 1$ and so $\mathcal{A}[\theta^{(n)}] = \mathrm{Hom}(\mathbb{Z}_p[\mathrm{Gal}(F_n/F)], A)$. On each term the composite map

$$H^i(F_\Sigma/F, \mathcal{A}[\theta^{(n)}]) \to H^i(F_\Sigma/F_n, \mathcal{A}[\theta^{(n)}]) \to H^i(F_\Sigma/F_n, A)$$

defined analogously to (7) is known to be bijective by the usual version of Shapiro's Lemma. The map (7) is the direct limit of these maps (which are compatible) and so is bijective too. ∎

For the proof of proposition 4.9, we may assume that $H^1(F_\Sigma/F, \mathcal{A})$ has $\Lambda$-corank $[F:\mathbb{Q}]$ and that $H^2(F_\Sigma/F, \mathcal{A})$ is $\Lambda$-cotorsion. Let $s \in \mathbb{Z}$. The exact sequence

$$0 \longrightarrow \mathcal{A}[\theta_s] \longrightarrow \mathcal{A} \xrightarrow{\theta_s} \mathcal{A} \longrightarrow 0$$

induces an exact sequence

$$H^1(F_\Sigma/F, \mathcal{A})/\theta_s H^1(F_\Sigma/F, \mathcal{A}) \xrightarrow{a} H^2(F_\Sigma/F, \mathcal{A}[\theta_s]) \xrightarrow{b} H^2(F_\Sigma/F, \mathcal{A})[\theta_s]$$

where of course $a$ is injective and $b$ is surjective. Let $X$ denote the Pontryagin dual of $H^1(F_\Sigma/F, \mathcal{A})$. Since $X$ is a finitely generated $\Lambda$-module, it is clear that $\ker(X \xrightarrow{\theta_s} X)$ will be finite for all but finitely many values of $s$. (Just choose $s$ so that $\theta_s \nmid f(T)$, where $f(T)$ is a generator of the characteristic ideal of $X_{\Lambda\text{-tors}}$. The $\theta_s$'s are irreducible and relatively prime.) Now $\text{Im}(a) = \ker(b)$ is the Pontryagin dual of $\ker(X \xrightarrow{\theta_s} X)$. We will show that $\ker(b)$ is always a divisible group. Hence, for suitable $s$, $\ker(X \xrightarrow{\theta_s} X) = 0$. Now if $Z$ is a nonzero, finite $\Lambda$-module, then $\ker(Z \xrightarrow{\theta_s} Z)$ is also clearly nonzero, since $\theta_s \notin \Lambda^\times$. Therefore, $X$ cannot contain a nonzero, finite $\Lambda$-submodule, which is equivalent to the assertion in proposition 4.9.

Assume that $p$ is odd. Then $\text{Gal}(F_\Sigma/F)$ has $p$-cohomological dimension 2. Since $\mathcal{A}[\theta_s] = A_{-s}$ is divisible, it follows from lemma 4.5 that $H^2(F_\Sigma/F, \mathcal{A}[\theta_s])$ is also divisible. Hence the same is true for $H^2(F_\Sigma/F, \mathcal{A})[\theta_s]$. But since $H^2(F_\Sigma/F, \mathcal{A})$ is $\Lambda$-cotorsion, $H^2(F_\Sigma/F, \mathcal{A})[\theta_s]$ will be finite for some value of $s$. Hence it must be zero. But this implies that $H^2(F_\Sigma/F, \mathcal{A}) = 0$, using Fact 1 about $\Lambda$-modules. Thus $\ker(b) = H^2(F_\Sigma/F, \mathcal{A}[\theta_s])$ for all $s$ and this is indeed divisible, proving proposition 4.9 if $p$ is odd.

The difficulty with the prime $p = 2$ is that $\text{Gal}(F_\Sigma/F)$ doesn't have finite $p$-cohomological dimension (unless $F$ is totally complex, in which case the argument in the preceding paragraph works). But we use the following fact: *the map*

$$\beta_n : H^n(\text{Gal}(F_\Sigma/F), M) \to \prod_{v|\infty} H^n(F_v, M)$$

*is an isomorphism for all $n \geq 3$.* Here $M$ can be any $p$-primary $\text{Gal}(F_\Sigma/F)$-module. (This is proved in [Mi], theorem 4.10(c) for the case where $M$ is finite. The general case follows from this.) The groups $H^n(F_v, M)$ have exponent $\leq 2$ for all $n \geq 1$. The following lemma is the key to dealing with the prime 2.

**Lemma 4.11.** *Assume that $M$ is divisible. Then the kernel of the map*

$$\beta_2 : H^2(F_\Sigma/F, M) \to \prod_{v|\infty} H^2(F_v, M)$$

*is a divisible group.*

*Proof.* Of course, if $p$ is odd, then $H^2(F_v, M) = 0$ for $v|\infty$. We already know that $H^2(F_\Sigma/F, M)$ is divisible in this case. Let $p = 2$. For any $m \geq 1$, consider the following commutative diagram with exact rows

$$
\begin{array}{ccccc}
H^2(F_\Sigma/F, M) & \xrightarrow{\;\;2^m\;\;} & H^2(F_\Sigma/F, M) & \xrightarrow{\;\alpha\;} & H^3(F_\Sigma/F, M[2^m]) \\
\downarrow{\scriptstyle\beta_2} & & \downarrow{\scriptstyle\beta_2} & & \downarrow{\scriptstyle\beta} \\
\prod_{v|\infty} H^2(F_v, M) & \xrightarrow{(.25)2^m} & \prod_{v|\infty} H^2(F_v, M) & \xrightarrow{(.15)\gamma} & \prod_{v|\infty} H^3(F_v, M[2^m])
\end{array}
$$

induced from the exact sequence $0 \to M[2^m] \to M \xrightarrow{2^m} M \to 0$. Since the group $H^2(F_v, M)$ is of exponent $\leq 2$, the map $\gamma$ is injective. Since $\beta$ is injective too, it follows that $\ker(\alpha) = \ker(\beta_2)$. Thus $\ker(\beta_2) = 2^m H^2(F_\Sigma/F, M)$ for any $m \geq 1$. Using this for $m = 1, 2$, we see that $\ker(\beta_2) = 2\ker(\beta_2)$, which implies that $\ker(\beta_2)$ is indeed divisible. ∎

Now we can prove that $b : H^2(F_\Sigma/F, \mathcal{A}[\theta_s]) \to H^2(F_\Sigma/F, \mathcal{A})$ has a divisible kernel even when $p = 2$. We use the following commutative diagram:

$$
\begin{array}{ccccccc}
0 & \longrightarrow & R^2(F_\Sigma/F, \mathcal{A}[\theta_s]) & \longrightarrow & H^2(F_\Sigma/F, \mathcal{A}[\theta_s]) & \longrightarrow & \prod_{v|\infty} H^2(F_v, \mathcal{A}[\theta_s]) \\
& & \downarrow{\scriptstyle d} & & \downarrow{\scriptstyle b} & & \downarrow{\scriptstyle e} \\
0 & \longrightarrow & R^2(F_\Sigma/F, \mathcal{A})[\theta_s] & \longrightarrow & H^2(F_\Sigma/F, \mathcal{A})[\theta_s] & \longrightarrow & \prod_{v|\infty} H^2(F_v, \mathcal{A})[\theta_s]
\end{array}
$$

The rows are exact by definition. (We define $R^2(F_\Sigma/F, M)$ as the kernel of the map $H^2(F_\Sigma/F, M) \to \prod_{v|\infty} H^2(F_v, M)$.) The map $b$ is surjective. Now $\mathcal{A}[\theta_s] \cong A_{-s}$ is divisible and hence, by lemma 4.11, $R^2(F_\Sigma/F, \mathcal{A}[\theta_s])$ is divisible. Under the assumption that $H^2(F_\Sigma/F, \mathcal{A})$ is $\Lambda$-cotorsion, we will show that $\ker(b)$ coincides with the divisible group $R^2(F_\Sigma/F, \mathcal{A}[\theta_s])$, completing the proof of proposition 4.9 for all $p$. Suppose that $v|\infty$. Since $v$ splits completely in $F_\infty/F$, we have $H^1(F_v, \mathcal{A}) = \operatorname{Hom}(\Lambda, H^1(F_v, E[p^\infty]))$. Of course, this group is zero unless $p = 2$ and $H^1(F_v, E[2^\infty]) \cong \mathbb{Z}/2\mathbb{Z}$, in which case $H^1(F_v, \mathcal{A}) = \operatorname{Hom}(\Lambda, \mathbb{Z}/2\mathbb{Z}) \cong (\Lambda/2\Lambda)\hat{\ }$. This last group is divisible by $\theta_s$ for any $s$, which implies that the map $e$ must be injective. The snake lemma then implies that the map $d$ is surjective. Thus $R^2(F_\Sigma/F, \mathcal{A})[\theta_s]$ is divisible for all $s \in \mathbb{Z}$. But this group is finite for all but finitely many $s$, since $H^2(F_\Sigma/F, \mathcal{A})$ is $\Lambda$-cotorsion. Hence, for some $s$, $R^2(F_\Sigma/F, \mathcal{A})[\theta_s] = 0$. This implies that the $\Lambda$-module $R^2(F_\Sigma/F, \mathcal{A})$ is zero. Therefore, since $e$ is injective, $\ker(b) = \ker(d) = R^2(F_\Sigma/F, \mathcal{A}[\theta_s])$ for all $s$, as claimed. ∎

The following proposition summarizes several consequences of the above arguments, which we translate back to the traditional form.

**Proposition 4.12.** *$H^1(F_\Sigma/F_\infty, E[p^\infty])$ has $\Lambda$-corank $[F:\mathbb{Q}]$ if and only if $H^2(F_\Sigma/F_\infty, E[p^\infty])$ is $\Lambda$-cotorsion. If this is so, then $H^1(F_\Sigma/F_\infty, E[p^\infty])$ has no proper $\Lambda$-submodule of finite index. Also, $H^2(F_\Sigma/F_\infty, E[p^\infty])$ will be zero if $p$ is odd and $(\Lambda/2\Lambda)$-cofree if $p = 2$.*

In this form, proposition 4.12 should apply to all primes $p$, since one conjectures that $H^2(F_\Sigma/F_\infty, E[p^\infty])$ is always $\Lambda$-cotorsion. (See conjecture 3 in [Gr2].) If $E$ has potentially good or multiplicative reduction at all primes over $p$, then, as mentioned in section 1, one expects that $\operatorname{Sel}_E(F_\infty)_p$ is $\Lambda$-cotorsion, which suffices to prove that $H^2(F_\Sigma/F_\infty, E[p^\infty])$ is indeed $\Lambda$-cotorsion. For any prime $p$, the conjecture that $\mathrm{Sel}_E(F_n)_p$ has bounded $\mathbb{Z}_p$-corank as $n$ varies can also be shown to suffice.

We must now explain why $H^1(F_\Sigma/F_\infty, E[p^\infty])_\Gamma$ is zero, under the hypotheses of theorem 4.1. We can assume that $\mathrm{Sel}_E(F_\infty)_p$ is $\Lambda$-cotorsion and that $H^1(F_\Sigma/F_\infty, E[p^\infty])$ has $\Lambda$-corank equal to $[F:\mathbb{Q}]$. By proposition 4.9, it is enough to prove that $H^1(F_\Sigma/F_\infty, E[p^\infty])_\Gamma$ is finite. Let

$$Q = H^1(F_\Sigma/F_\infty, E[p^\infty])/H^1(F_\Sigma/F_\infty, E[p^\infty])_{\Lambda\text{-div}}.$$

Thus, $Q$ is cofinitely generated and cotorsion as a $\Lambda$-module. Its Pontryagin dual is the torsion $\Lambda$-submodule of the Pontryagin dual of $H^1(F_\Sigma/F_\infty, E[p^\infty])$. We have

$$H^1(F_\Sigma/F_\infty, E[p^\infty])_\Gamma = Q_\Gamma.$$

But $Q_\Gamma$ and $Q^\Gamma$ have the same $\mathbb{Z}_p$-corank. Also, $\widehat{\Lambda}^\Gamma \cong \mathbb{Q}_p/\mathbb{Z}_p$ has $\mathbb{Z}_p$-corank 1. Since the map $H^1(F_\Sigma/F, E[p^\infty]) \to H^1(F_\Sigma/F_\infty, E[p^\infty])^\Gamma$ is surjective and has finite kernel, we see that

$$\mathrm{corank}_{\mathbb{Z}_p}(H^1(F_\Sigma/F, E[p^\infty])) = [F:\mathbb{Q}] + \mathrm{corank}_{\mathbb{Z}_p}(Q_\Gamma).$$

Now

$$\mathrm{Sel}_E(F)_p = \ker(H^1(F_\Sigma/F, E[p^\infty]) \to \mathcal{P}_E^\Sigma(F)).$$

The $\mathbb{Z}_p$-corank of $\mathcal{P}_E^\Sigma(F)$ is equal to $[F:\mathbb{Q}]$. Since we are assuming that $\mathrm{Sel}_E(F)_p$ is finite, it follows that $H^1(F_\Sigma/F, E[p^\infty])$ has $\mathbb{Z}_p$-corank $[F:\mathbb{Q}]$ and hence that, indeed, $Q_\Gamma$ is finite which completes the argument. We should point out that sometimes $H^1(F_\Sigma/F_\infty, E[p^\infty])_\Gamma$ is nonzero. This clearly happens for example when $\mathrm{rank}_{\mathbb{Z}}(E(F)) > [F:\mathbb{Q}]$. For then $H^1(F_\Sigma/F, E[p^\infty])$ must have $\mathbb{Z}_p$-corank at least $[F:\mathbb{Q}] + 1$, which implies that $Q_\Gamma$ is nonzero.

We will now prove a rather general version of Cassels' theorem. Let $\Sigma$ be a finite set of primes of a number field $F$, containing at least all primes of $F$ lying above $p$ and $\infty$. We suppose that $M$ is a $\mathrm{Gal}(F_\Sigma/F)$-module isomorphic to $(\mathbb{Q}_p/\mathbb{Z}_p)^d$ as a group (for any $d \geq 1$). For each $v \in \Sigma$, we assume that $L_v$ is a divisible subgroup of $H^1(F_v, M)$. Then we define a "Selmer group"

$$S_M(F) = \ker\left(H^1(F_\Sigma/F, M) \to \prod_{v \in \Sigma} H^1(F_v, M)/L_v\right).$$

This is a discrete, $p$-primary group which is cofinitely generated over $\mathbb{Z}_p$. Let

$$T^* = \mathrm{Hom}(M, \mu_{p^\infty})$$

which is a free $\mathbb{Z}_p$-module of rank $d$. For each $v \in \Sigma$, we define a subgroup $U_v^*$ of $H^1(F_v, T^*)$ as the orthogonal complement of $L_v$ under the perfect pairing (from Tate's local duality theorems)

$$H^1(F_v, M) \times H^1(F_v, T^*) \to \mathbb{Q}_p/\mathbb{Z}_p. \tag{8}$$

Since $L_v$ is divisible, it follows that $H^1(F_v, T^*)/U_v^*$ is $\mathbb{Z}_p$-torsion free. Thus $U_v^*$ contains $H^1(F_v, T^*)_{\text{tors}}$. We define the Selmer group

$$S_{T^*}(F) = \ker\big(H^1(F_\Sigma/F, T^*) \to \prod_{v \in \Sigma} H^1(F_v, T^*)/U_v^*\big)$$

which will be a finitely generated $\mathbb{Z}_p$-module. Let $V^* = T^* \otimes \mathbb{Q}_p$. Let $M^* = V^*/T^* = T^* \otimes (\mathbb{Q}_p/\mathbb{Z}_p)$. For each $v \in \Sigma$, we can define a divisible subgroup $L_v^*$ of $H^1(F_v, M^*)$ as follows: Under the map $H^1(F_v, T^*) \to H^1(F_v, V^*)$, the image of $U_v^*$ generates a $\mathbb{Q}_p$-subspace of $H^1(F_v, V^*)$. We define $L_v^*$ as the image of this subspace under the map $H^1(F_v, V^*) \to H^1(F_v, M^*)$. Thus, we can define a Selmer group

$$S_{M^*}(F) = \ker\big(H^1(F_\Sigma/F, M^*) \to \prod_{v \in \Sigma} H^1(F_v, M^*)/L_v^*\big).$$

One can verify that the $\mathbb{Z}_p$-corank of $S_{M^*}(F)$ is equal to the $\mathbb{Z}_p$-rank of $S_{T^*}(F)$.

We will use the following notation. Let

$$P = \prod_{v \in \Sigma} H^1(F_v, M), \qquad\qquad P^* = \prod_{v \in \Sigma} H^1(F_v, T^*)$$

$$L = \prod_{v \in \Sigma} L_v, \qquad\qquad U^* = \prod_{v \in \Sigma} U_v^*.$$

Then (8) induces a perfect pairing $P \times P^* \to \mathbb{Q}_p/\mathbb{Z}_p$, under which $L$ and $U^*$ are orthogonal complements. Furthermore, we let

$$G = \text{Im}\big(H^1(F_\Sigma/F, M) \to P\big), \qquad G^* = \text{Im}\big(H^1(F_\Sigma/F, T^*) \to P^*\big).$$

The duality theorems of Poitou and Tate imply that $G$ and $G^*$ are also orthogonal complements under the above perfect pairing. Consider the map

$$\gamma : H^1(F_\Sigma/F, M) \to P/L,$$

whose kernel is, by definition, $S_M(F)$. The cokernel of $\gamma$ is clearly $P/GL$. But the orthogonal complement of $GL$ under the pairing $P \times P^* \to \mathbb{Q}_p/\mathbb{Z}_p$ must be $G^* \cap U^*$. Thus $\text{coker}(\gamma) \cong (G^* \cap U^*)\widehat{\phantom{x}}$. Again by definition, $S_{T^*}(F)$ is the inverse image of $U^*$ under the map $H^1(F_\Sigma/F, T^*) \to P^*$. Thus clearly $G^* \cap U^*$ is a homomorphic image of $S_{T^*}(F)$. As we mentioned above, $\text{rank}_{\mathbb{Z}_p}(S_{T^*}(F))$ is equal to $\text{corank}_{\mathbb{Z}_p}(S_{M^*}(F))$. On the other hand, since $H^1(F_v, T^*)_{\text{tors}}$ is contained in $U_v^*$ for all $v \in \Sigma$, it follows that

$$S_{T^*}(F)_{\text{tors}} = H^1(F_\Sigma/F, T^*)_{\text{tors}},$$

which in turn is isomorphic to $H^0(F_\Sigma/F, M^*)/H^0(F_\Sigma/F, M^*)_{\text{div}}$. (This last assertion follows from the cohomology sequence induced from the exact sequence $0 \to T^* \to V^* \to M^* \to 0$.) We denote $H^0(F_\Sigma/F, M^*) = (M^*)^{G_F}$ by $M^*(F)$ as usual. Then, as a $\mathbb{Z}_p$-module, we have

$$S_{T^*}(F) \cong (M^*(F)/M^*(F)_{\text{div}}) \times \mathbb{Z}_p^{\text{corank}_{\mathbb{Z}_p}(S_{M^*}(F))}.$$

The preceding discussion proves that the Pontryagin dual of the cokernel of the map $\gamma$ is a homomorphic image of $S_{T^*}(F)$. In particular, one important special case is: *if $S_{M^*}(F)$ is finite and $M^*(F) = 0$, then $coker(\gamma) = 0$.*

We now make the following slightly restrictive hypothesis: $M^*(F_v) = H^0(F_v, M^*)$ *is finite for at least one $v \in \Sigma$.* This implies that $M^*(F)$ is also finite. Consider the following commutative diagram.

$$\begin{array}{ccc} H^0(F_\Sigma/F, M^*) & \xrightarrow{\sim} & H^1(F_\Sigma/F, T^*)_{\mathrm{tors}} \\ \downarrow & & \downarrow \\ H^0(F_v, M^*) & \xrightarrow{\sim} & H^1(F_v, T^*)_{\mathrm{tors}} \end{array} \qquad (9)$$

Since the first vertical arrow is obviously injective, so is the second. Hence the map $H^1(F_\Sigma/F, T^*)_{\mathrm{tors}} \to P^*$ is injective. It follows that if $S_{T^*}(F)$ is finite, then

$$\mathrm{coker}(\gamma) \cong (G^* \cap U^*)\,\widehat{} \cong S_{T^*}(F)\,\widehat{} = H^1(F_\Sigma/F, T^*)_{\mathrm{tors}}.$$

This last group is isomorphic to $M^*(F)$. We obtain the following general version of Cassels' theorem.

**Proposition 4.13.** *Assume that $m^* = \mathrm{corank}_{\mathbb{Z}_p}(\mathrm{Sel}_{M^*}(F))$. Assume also that $H^0(F_v, M^*)$ is finite for at least one $v \in \Sigma$. Then the cokernel of the map*

$$\gamma : H^1(F_\Sigma/F, M) \to \prod_{v \in \Sigma} H^1(F_v, M)/L_v$$

*has $\mathbb{Z}_p$-corank $\leq m^*$. Also,*

$$\dim_{\mathbb{Z}/p\mathbb{Z}}(\mathrm{coker}(\gamma)[p]) \leq m^* + \dim_{\mathbb{Z}/p\mathbb{Z}}(H^0(F, M^*[p])).$$

*If $m^* = 0$, then $\mathrm{coker}(\gamma) \cong H^0(F, M^*)\,\widehat{}$.*

It is sometimes useful to know how $\mathrm{Im}(\gamma)$ sits inside of $P/L$. We can make the following remark. Let $v_0$ be any prime in $\Sigma$ for which $H^0(F_{v_0}, M^*)$ is finite. Assume that $S_{M^*}(F)$ is finite. Then

$$\mathrm{Im}(\gamma)(H^1(F_{v_0}, M)/L_{v_0}) = P/L.$$

Here $H^1(F_{v_0}, M)/L_v$ is a direct factor in $P/L$. To justify this, one must just show that the map

$$\gamma' : H^1(F_\Sigma/F, M) \to \prod_{\substack{v \in \Sigma \\ v \neq v_0}} H^1(F_v, M)/L_v$$

is surjective under the above assumptions about $S_{M^*}(F)$ and $v_0$. In the above arguments, one can study $\mathrm{coker}(\gamma')$ by changing $L_{v_0}$ to $L'_{v_0} = H^1(F_{v_0}, M)$ and leaving $L_v$ for $v \neq v_0$ unchanged. Now $L'_{v_0}$ may not be divisible, but we still have $\operatorname{coker}(\gamma') \cong (G^* \cap U'^*)^\wedge$, where now $U^*_{v_0}$ has been replaced by $U'^*_{v_0} = 0$. Since $U'^* \subseteq U^*$, the corresponding Selmer group $S'_{T^*}(F) = 0$. Thus an element $\sigma$ in $S'_{T^*}(F)$ is in $H^1(F_\Sigma/F, T^*)_{\mathrm{tors}}$ and has the property that $\sigma|_{G_{F_{v_0}}}$ is trivial. But the diagram (9) shows that

$$H^1(F_\Sigma/F, T^*)_{\mathrm{tors}} \to H^1(F_{v_0}, T^*)_{\mathrm{tors}}$$

is injective. Hence $\sigma$ is trivial. Thus, $S'_{T^*}(F)$ is trivial and hence so is $\operatorname{coker}(\gamma')$.

Cassels' theorem is the following special case of proposition 4.13: $M = E[p^\infty]$, $\Sigma =$ any finite set of primes of $F$ containing the primes lying over $p$ or $\infty$ and the primes where $E$ has bad reduction, and $L_v = \operatorname{Im}(\kappa_v)$ for all $v \in \Sigma$. Then $T^* = T_p(E)$ by the Weil pairing. Thus $M^* = E[p^\infty]$, $L^*_v = \operatorname{Im}(\kappa_v)$, and $S_{M^*}(F) = S_M(F) = \operatorname{Sel}_E(F)_p$. It is clear that $H^0(F_v, M^*)$ is finite for any nonarchimedean $v \in \Sigma$. Thus, proposition 4.13 implies that

$$\operatorname{coker}\big(H^1(F_\Sigma/F, E[p^\infty]) \to \prod_{v \in \Sigma} H^1(F_v, E[p^\infty])/\operatorname{Im}(\kappa_v)\big) \cong E(F)\hat{_p}$$

if $\operatorname{Sel}_E(F)_p$ is finite. (Of course, as a group, $E(F)\hat{_p} \cong E(F)_p$.) In the proof of lemma 4.6 we need the following case: $E$ is an elliptic curve which we assume has (potentially) good, ordinary or multiplicative reduction at all $v|p$, $M = E[p^\infty] \otimes \kappa^s$ where $s \in \mathbb{Z}$, $L_v = \operatorname{Im}\big(H^1(F_v, C_v \otimes \kappa^s) \to H^1(F_v, M)\big)_{\mathrm{div}}$ if $v|p$, $L_v = 0$ if $v \nmid p$. Then $T^* = T_p(E) \otimes \kappa^{-s}$, $M^* = E[p^\infty] \otimes \kappa^{-s}$, and $L^*_v$ is defined just as $L_v$. Assuming that $\operatorname{Sel}_E(F_\infty)_p$ is $\Lambda$-cotorsion, we can choose $s \in \mathbb{Z}$ so that $S_{M^*}(F)$ is finite. The hypothesis that $H^0(F_v, M^*)$ is finite for some $v \in \Sigma$ is also easily satisfied (possibly avoiding one value of $s$). Then the cokernel of the map $\gamma$ will be isomorphic to the finite group $H^0(F, M^*)\hat{\ }$.

We can now prove the following generalization of proposition 4.8.

**Proposition 4.14.** *Assume that $E$ is an elliptic curve defined over $F$ and that $\operatorname{Sel}_E(F_\infty)_p$ is $\Lambda$-cotorsion. Assume that $E(F)_p = 0$. Then $\operatorname{Sel}_E(F_\infty)_p$ has no proper $\Lambda$-submodules of finite index.*

*Proof.* As in the proof of lemma 4.6, we will use the twisted Galois modules $A_s = E[p^\infty] \otimes (\kappa^s)$, where $s \in \mathbb{Z}$. Since $E(F)_p = 0$, it follows that $E(F_\infty)_p = 0$ too. (One uses the fact that $\Gamma$ is pro-$p$.) Since $A_s \cong E[p^\infty]$ as $G_{F_\infty}$-modules, it is clear that $H^0(F, A_s) = 0$ for all $s$. Now $E$ must have potentially ordinary or multiplicative reduction at all $v|p$, since we are assuming that $\operatorname{Sel}_E(F_\infty)_p$ is $\Lambda$-cotorsion. So we can define a Selmer group $S_{A_s}(K)$ for any algebraic extension $K$ of $F$. If we take $K$ to be a subfield of $F_\Sigma$, then $S_{A_s}(K)$ is the kernel of a map $H^1(F_\Sigma/K, A_s) \to \mathcal{P}^\Sigma(A_s, K)$, where this last group is defined in a way analogous to $\mathcal{P}^\Sigma_E(K)$. As we pointed out in the proof of lemma 4.6, we have

$$S_{A_s}(F_\infty) \cong \operatorname{Sel}_E(F_\infty)_p \otimes (\kappa^s)$$

as $\Lambda$-modules. We also have $\mathcal{P}^\Sigma(A_s, F_\infty) \cong \mathcal{P}_E^\Sigma(F_\infty) \otimes (\kappa^s)$ as $\Lambda$-modules. The hypothesis that $\mathrm{Sel}_E(F_\infty)_p$ is $\Lambda$-cotorsion implies that $S_{A_s}(F_\infty)^\Gamma$, and hence $S_{A_s}(F)$, will be finite for all but finitely many values of $s$. (We will add another requirement on $s$ below.) We let $M = A_s$, where $s \in \mathbb{Z}$ has been chosen so that $S_{A_{-s}}(F)$ is finite. Note that $M^* = A_{-s}$. Since $S_{M^*}(F)$ is finite and $M^*(F) = 0$, we can conclude that the map

$$\gamma : H^1(F_\Sigma/F, M) \to \mathcal{P}^\Sigma(M, F)$$

is surjective. Since $\Gamma$ has cohomological dimension 1, the restriction maps $H^1(F_\Sigma/F, M) \to H^1(F_\Sigma/F_\infty, M)^\Gamma$ and $\mathcal{P}^\Sigma(M, F) \to \mathcal{P}^\Sigma(M, F_\infty)^\Gamma$ are both surjective. Hence it follows that the map

$$H^1(F_\Sigma/F_\infty, M)^\Gamma \to \mathcal{P}^\Sigma(M, F_\infty)^\Gamma$$

must be surjective. We have the exact sequence defining $S_M(F_\infty)$:

$$0 \to S_M(F_\infty) \to H^1(F_\Sigma/F_\infty, M) \to \mathcal{P}^\Sigma(M, F_\infty) \to 0.$$

This is just the exact sequence defining $\mathrm{Sel}_E(F_\infty)_p$, twisted by $\kappa^s$. The corresponding cohomology sequence induces an injective map

$$S_M(F_\infty)_\Gamma \to H^1(F_\Sigma/F_\infty, M)_\Gamma.$$

If we let $Q = H^1(F_\Sigma/F_\infty, E[p^\infty])/H^1(F_\Sigma/F_\infty, E[p^\infty])_{\Lambda\text{-div}}$, as before, then

$$H^1(F_\Sigma/F_\infty, M)_\Gamma \cong (Q \otimes (\kappa^s))_\Gamma$$

and, since $Q$ is $\Lambda$-cotorsion, we can choose $s$ so that $(Q \otimes (\kappa^s))_\Gamma$ is finite. (This will be true for all but finitely many values of $s$.) But since $H^1(F_\Sigma/F_\infty, E[p^\infty])$ has no proper $\Lambda$-submodules of finite index, neither does $H^1(F_\Sigma/F_\infty, M)$. It follows that, for suitably chosen $s$, $H^1(F_\Sigma/F_\infty, M)_\Gamma = 0$. Hence $S_M(F_\infty)_\Gamma = 0$. This implies that $S_M(F_\infty)$ has no proper $\Lambda$-submodules of finite index, from which proposition 4.14 follows. ∎

We will give two other sufficient conditions for the nonexistence of proper $\Lambda$-submodules of finite index in $\mathrm{Sel}_E(F_\infty)_p$. We want to mention that a rather different proof of proposition 4.14 and part of the following proposition has been found by Hashimori and Matsuno [HaMa]. This proof is based on the Cassels-Tate pairing for $\mathrm{III}_E(F_n)_p$. This topic will be pursued much more generally in [Gr6].

**Proposition 4.15.** *Assume that $E$ is an elliptic curve defined over $F$ and that $\mathrm{Sel}_E(F_\infty)_p$ is $\Lambda$-cotorsion. Assume that at least one of the following two hypotheses holds:*

*(i) there is a prime $v_0$ of $F$, $v_0 \nmid p$, where $E$ has additive reduction.*

*(ii) there exists a prime $v_0$ of $F$, $v_0|p$, such that the ramification index $e_{v_0}$ of $F_{v_0}/\mathbb{Q}_p$ satisfies $e_{v_0} \leq p - 2$ and such that $E$ has good, ordinary or multiplicative reduction at $v_0$.*

*Then $\mathrm{Sel}_E(F_\infty)_p$ has no proper $\Lambda$-submodules of finite index.*

*Remark.* If condition (i) holds, then $H^0(I_{\mathfrak{Q}_{v_0}}, E[p^\infty])$ is finite. This group will be zero if $p \geq 5$. Then $E(F)_p = 0$ and we are in the situation of proposition 4.14.

*Proof.* We will modify the proof of proposition 4.14. In addition to the requirements on $M = A_s$ occurring in that proof, we also require that $H^0(F_{v_0}, M^*)$ be finite, which is true for all but finitely many values of $s \in \mathbb{Z}$. Here $v_0$ is the prime of $F$ satisfying (i) or (ii). (If $E$ has additive reduction at $v_0$, $v_0 \nmid p$, then this holds for all $s$.) Assume first that (i) holds. In this case, let $S'_M(F)$, $S'_M(F_\infty)$ denote the Selmer groups where one omits the local condition at $v_0$ (or the primes above $v_0$). If $\eta$ is a prime of $F_\infty$ lying over $v_0$, then $H^0((F_\infty)_\eta, M^*)$ is finite. This implies that $H^1((F_\infty)_\eta, M) = 0$. Thus, $S'_M(F_\infty) = S_M(F_\infty)$. The remark following proposition 4.13 shows that the map

$$\gamma' : H^1(F_\Sigma/F, M) \to \mathcal{P}^{\Sigma'}(M, F)$$

is surjective, where $\Sigma' = \Sigma - \{v_0\}$ and $\mathcal{P}^{\Sigma'}(M, F)$ is the product over all primes of $\Sigma'$. The proof then shows that $S'_M(F_\infty)$ has no proper $\Lambda$-submodules of finite index. This obviously gives the same statement for $\mathrm{Sel}_E(F_\infty)_p$.

Now assume (ii). We again define $S'_M(F_\infty)$ by omitting the local condition at all primes $\eta$ of $F_\infty$ lying over $v_0$. Just as above, we see that $S'_M(F_\infty)$ has no proper $\Lambda$-submodules of finite index. Thus, the same is true for $\mathrm{Sel}'_E(F_\infty)_p$. By lemma 4.6, we see that

$$\mathrm{Sel}'_E(F_\infty)_p / \mathrm{Sel}_E(F_\infty)_p \cong \prod_{\eta|v_0} \mathcal{H}_E((F_\infty)_\eta).$$

But $\mathcal{H}_E((F_\infty)_\eta) \cong H^1((F_\infty)_\eta, E[p^\infty])/\mathrm{Im}(\kappa_\eta) \cong H^1((F_\infty)_\eta, D_{v_0})$ by proposition 2.4 and the analogous statement proved in section 2 for the case where $E$ has multiplicative reduction at $v$. Here $D_{v_0} = E[p^\infty]/C_{v_0}$ is an unramified $G_{F_{v_0}}$-module isomorphic to $\mathbb{Q}_p/\mathbb{Z}_p$. We can use a remark made in section 2 (preceding proposition 2.4) to conclude that $H^1((F_\infty)_\eta, D_{v_0})$ is $\mathbb{Z}_p[[\mathrm{Gal}((F_\infty)_\eta/F_{v_0})]]$-cofree. Proposition 4.15 in case (ii) is then a consequence of the following fact about finitely generated $\Lambda$-modules: *Suppose that $X'$ is a finitely generated $\Lambda$-module which has no nonzero, finite $\Lambda$-submodules. Assume that $Y$ is a free $\Lambda$-submodule of $X'$. Then $X = X'/Y$ has no nonzero, finite $\Lambda$-submodules.* The proof is quite easy. By induction, one can assume that $Y \cong \Lambda$. Suppose that $X$ does have a nonzero, finite $\Lambda$-submodule. Then $Y \subseteq Y_0$, where $[Y_0 : Y] < \infty$, $Y \neq Y_0$, and $Y_0$ is a $\Lambda$-submodule of $X'$. Then $Y_0$ is pseudo-isomorphic to $\Lambda$ and has no nonzero, finite $\Lambda$-submodules. Hence $Y_0$ would be isomorphic to a submodule of $\Lambda$ of finite index. It would follow that $\Lambda$ contains a proper ideal of finite index which is isomorphic to $\Lambda$, i.e., a principal ideal. But if $f \in \Lambda$, then $(f)$ can't have finite index unless $f \in \Lambda^\times$, in which case $(f) = \Lambda$. Hence in fact $X$ has no nonzero, finite $\Lambda$-submodules.

## 5.  Conclusion.

In this final section we will discuss the structure of $\mathrm{Sel}_E(F_\infty)_p$ in various special cases, making full use of the results of sections 3 and 4. In particular, we will see that each of the invariants $\mu_E$, $\lambda_E^{M-W}$, and $\lambda_E^{\mathrm{III}}$ can be positive. We will assume (usually) that the base field $F$ is $\mathbb{Q}$ and that $E/\mathbb{Q}$ has good, ordinary reduction at $p$. Our examples will be based on the predicted order of the Shafarevich-Tate groups given in Cremona's tables. In principle, these orders can be verified by using results of Kolyvagin.

We start with the following corollary to proposition 3.8.

**Proposition 5.1.** *Assume that $E$ is an elliptic curve$/\mathbb{Q}$ and that both $E(\mathbb{Q})$ and $\mathrm{III}_E(\mathbb{Q})$ are finite. Let $p$ vary over the primes where $E$ has good, ordinary reduction. Then $\mathrm{Sel}_E(\mathbb{Q}_\infty)_p = 0$ except for $p$ in a set of primes of zero density. This set of primes is finite if $E$ is $\mathbb{Q}$-isogenous to an elliptic curve $E'$ such that $|E'(\mathbb{Q})| > 1$.*

*Remark.* Recall that if $p$ is a prime where $E$ has supersingular (or potentially supersingular) reduction, then $\mathrm{Sel}_E(\mathbb{Q}_\infty)_p$ has positive $\Lambda$-corank. Under the hypothesis that $E(\mathbb{Q})$ and $\mathrm{III}_E(\mathbb{Q})$ are finite, this $\Lambda$-corank can be shown to equal 1, agreeing with the conjecture stated after theorem 1.7. If $E$ doesn't have complex multiplication, the set of supersingular primes for $E$ also has zero density.

*Proof.* We are assuming that $\mathrm{Sel}_E(\mathbb{Q})$ is finite. Thus, excluding finitely many primes, we can assume that $\mathrm{Sel}_E(\mathbb{Q})_p = 0$. If we also exclude the finite set of primes dividing $\prod_l c_l$, where $l$ varies over the primes where $E$ has bad reduction and $c_l$ is the corresponding Tamagawa factor, then hypotheses (ii) and (iii) in proposition 3.8 are satisfied. As for hypothesis (i), it is equivalent to $a_p \equiv 1 \pmod{p}$, where $a_p = 1 + p - |\widetilde{E}(\mathbb{F}_p)|$. Now we have Hasse's result that $|a_p| < 2\sqrt{p}$ and hence $a_p \equiv 1 \pmod{p} \Rightarrow a_p = 1$ if $p > 5$. By using the Chebotarev Density Theorem, one can show that $\{p \mid a_p = 1\}$ has zero density. (That is, the cardinality of $\{p \mid a_p = 1, p < x\}$ is $o(x/\log(x))$ as $x \to \infty$.) The argument is a standard one, using the $l$-adic representation attached to $E$ for any fixed prime $l$. The trace of a Frobenius element for $p \ (\neq l)$ is $a_p$. One considers separately the cases where $E$ does or does not have complex multiplication. For the non-CM case, see [Se1], IV–13, exercise 1. These remarks show that the hypotheses in proposition 3.8 hold if $p$ is outside a set of primes of zero density. For such $p$, $\mathrm{Sel}_E(\mathbb{Q}_\infty)_p = 0$. The final part of proposition 5.1 follows from the next lemma.

**Lemma 5.2.** *Suppose that $E$ is an elliptic curve defined over $\mathbb{Q}$ and that $p$ is a prime where $E$ has good reduction. If $E(\mathbb{Q})$ has a point of order 2 and $p > 5$, then $a_p \not\equiv 1 \pmod{p}$. If $E$ is $\mathbb{Q}$-isogenous to an elliptic curve $E'$ such that $E'(\mathbb{Q})_{\mathrm{tors}}$ has a subgroup of order $q > 2$ and if $p \nmid q$, then $a_p \not\equiv 1 \pmod{p}$.*

*Proof.* $\mathbb{Q}$-isogenous elliptic curves have the same set of primes of bad reduction. If $E$ has good reduction at $p$, then the prime-to-$p$ part of $E'(\mathbb{Q})_{\mathrm{tors}}$ maps injectively into $\widetilde{E}'(\mathbb{F}_p)$, which has the same order as $\widetilde{E}(\mathbb{F}_p)$. For the first part, $a_p \equiv 1 \pmod{p}$ implies that $2p$ divides $|\widetilde{E}(\mathbb{F}_p)|$. Hence $2p < 1 + p + 2\sqrt{p}$, which is impossible for $p > 5$. For the second part, if $a_p \equiv 1 \pmod{p}$ and $p \nmid q$, then $qp$ divides $|\widetilde{E}(\mathbb{F}_p)|$. Hence $qp < 1 + p + 2\sqrt{p}$, which again is impossible since $q > 2$.

Here are several specific examples.

$\mathbf{E = X_0(11)}$. The equation $y^2 + y = x^3 - x^2 - 10x - 20$ defines this curve, which is 11(A1) in [Cre]. $E$ has split, multiplicative reduction at $p = 11$ and good reduction at all other primes. We have $\mathrm{ord}_{11}(j_E) = -5$, $E(\mathbb{Q}) \cong \mathbb{Z}/5\mathbb{Z}$, and we will assume that $\mathrm{Sel}_E(\mathbb{Q}) = 0$ as predicted. If $p \neq 11$, then $a_p \equiv 1 \pmod{p}$ happens only for $p = 5$. Therefore, if $E$ has good, ordinary reduction at $p \neq 5$, then $\mathrm{Sel}_E(\mathbb{Q}_\infty)_p = 0$ according to proposition 3.8. We will discuss the case $p = 5$ later, showing that $\mathrm{Sel}_E(\mathbb{Q}_\infty)_p \cong \mathrm{Hom}(\Lambda/p\Lambda, \mathbb{Z}/p\mathbb{Z})$ and hence that $\mu_E = 1$, $\lambda_E = 0$. We just mention now that, by theorem 4.1, $f_E(0) \sim 5$. We will also discuss quite completely the other two elliptic curves/$\mathbb{Q}$ of conductor 11 for the case $p = 5$. If $p = 11$, then $\mathrm{Sel}_E(\mathbb{Q}_\infty)_p = 0$. This is verified in [Gr3], example 3.

$\mathbf{E = X_0(32)}$. This curve is defined by $y^2 = x^3 - 4x$ and is 32(A1) in [Cre]. It has complex multiplication by $\mathbb{Z}[i]$. $E$ has potentially supersingular reduction at 2. For an odd prime $p$, $E$ has good, ordinary reduction at $p$ if and only if $p \equiv 1 \pmod 4$. We have $E(\mathbb{Q}) \cong \mathbb{Z}/4\mathbb{Z}$, $\mathrm{III}_E(\mathbb{Q}) = 0$ (as verified in Rubin's article in this volume), and $c_2 = 4$. By lemma 5.2, there are no anomalous primes for $E$. Therefore, $\mathrm{Sel}_E(\mathbb{Q}_\infty)_p = 0$ for all primes $p$ where $E$ has good, ordinary reduction.

$\mathbf{E_1 : y^2 = x^3 + x^2 - 7x + 5}$ and $\mathbf{E_2 : y^2 = x^3 + x^2 - 647x - 6555}$. Both of these curves have conductor 768. They are 768(D1) and 768(D3) in [Cre]. They are related by a 5-isogeny defined over $\mathbb{Q}$. We will assume that $\mathrm{Sel}_{E_1}(\mathbb{Q})$ is trivial as predicted by the Birch and Swinnerton-Dyer conjecture. This implies that $\mathrm{Sel}_{E_2}(\mathbb{Q})_p = 0$ for all primes $p \neq 5$. We will verify later that this is true for $p = 5$ too. Both curves have additive reduction at $p = 2$, and split, multiplicative reduction at $p = 3$. For $E_1$, the Tamagawa factors are $c_2 = 2$, $c_3 = 1$. For $E_2$, they are $c_2 = 2$, $c_3 = 5$. We have $E_1(\mathbb{Q}) \cong \mathbb{Z}/2\mathbb{Z} \cong E_2(\mathbb{Q})$. By lemma 5.2, no prime $p > 5$ is anomalous for $E_1$ or $E_2$. If $E_1$ (and hence $E_2$) have ordinary reduction at a prime $p > 5$, then proposition 3.8 implies that $\mathrm{Sel}_{E_1}(\mathbb{Q}_\infty)_p = 0 = \mathrm{Sel}_{E_2}(\mathbb{Q}_\infty)_p$. Both of these curves have good, ordinary reduction at $p = 5$. (In fact, $\widetilde{E}_1 = \widetilde{E}_2 : y^2 = x^3 + x^2 + \widetilde{3}x$ and one finds 4 points. That is, $a_5 = 2$ and so $p = 5$ is not anomalous for $E_1$ or $E_2$.) The hypotheses of proposition 3.8 are satisfied for $E_1$ and $p = 5$. Hence $\mathrm{Sel}_{E_1}(\mathbb{Q}_\infty)_5 = 0$. But, by using either the results of section 3 or theorem 4.1, one sees that $\mathrm{Sel}_{E_2}(\mathbb{Q}_\infty)_5 \neq 0$. (One can either point out that $\mathrm{coker}(\mathrm{Sel}_{E_2}(\mathbb{Q})_5 \to \mathrm{Sel}_{E_2}(\mathbb{Q}_\infty)_5^{\Gamma})$ is nonzero or that $f_{E_2}(0) \sim 5$. We remark that proposition 4.8 tells us that $\mathrm{Sel}_{E_2}(\mathbb{Q}_\infty)$ cannot just be finite if it is nonzero.) Now if $\phi : E_1 \to E_2$ is a 5-isogeny defined over $\mathbb{Q}$, the induced map $\Phi : \mathrm{Sel}_{E_1}(\mathbb{Q}_\infty)_5 \to \mathrm{Sel}_{E_2}(\mathbb{Q}_\infty)_5$ will have kernel and cokernel of exponent 5. Hence $\lambda_{E_2} = \lambda_{E_1} = 0$ (for $p = 5$). Since $f_{E_2}(0) \sim 5$, it is clear that $\mu_{E_2} = 1$. Below we will verify directly that $\mathrm{Sel}_{E_2}(\mathbb{Q}_\infty)_5 \cong \mathrm{Hom}(\Lambda/5\Lambda, \mathbb{Z}/5\mathbb{Z})$. Note that this example illustrates conjecture 1.11.

**E : $y^2 + y = x^3 + x^2 - 12x - 21$**. This is 67(A1) in [Cre]. It has split, multiplicative reduction at $p = 67$, good reduction at all other primes. We have $E(\mathbb{Q}) = 0$ and $c_{67} = 1$. It should be true that $\mathrm{Sel}_E(\mathbb{Q}) = 0$, which we will assume. According to proposition 3.8, $\mathrm{Sel}_E(\mathbb{Q}_\infty)_p = 0$ for any prime $p \neq 67$ where $a_p \not\equiv 0, 1 \pmod p$. If $a_p \equiv 0 \pmod p$, then $E$ has supersingular reduction at $p$, and hence $\mathrm{Sel}_E(\mathbb{Q}_\infty)_p$ is not even $\Lambda$-cotorsion. (In fact, the $\Lambda$-corank will be 1.) If $a_p \equiv 1 \pmod p$, then $\mathrm{Sel}_E(\mathbb{Q}_\infty)_p$ must be nonzero and hence infinite. (Proposition 4.8 applies.) By proposition 4.1, we in fact have $f_E(0) \sim |\widetilde{E}(\mathbb{F}_p)|^2 \sim p^2$ for any such prime $p$. (Here we use Hasse's estimate on $|\widetilde{E}(\mathbb{F}_p)|$, noting that $1 + p + 2\sqrt{p} < p^2$ for $p \geq 3$. The prime $p = 2$ is supersingular for this elliptic curve.) Now it seems reasonable to expect that $E$ has infinitely many anomalous primes. The first such $p$ is $p = 3$ (and the only such $p < 100$). Conjecture 1.11 implies that $\mu_E = 0$. Assuming this, we will later see that $\lambda_E^{M-W} = 0$ and $\lambda_E^{III} = 2$.

**E : $y^2 + y = x^3 - x^2 - 460x - 11577$**. This curve has conductor 915. It is 915(A1) in [Cre]. It has split, multiplicative reduction at 5 and 61, nonsplit at 3. We have $c_3 = c_{61} = 1$ and $c_5 = 7$. $\mathrm{Sel}_E(\mathbb{Q}) = 0$, conjecturally. $E(\mathbb{Q}) = 0$. Proposition 3.8 implies that $\mathrm{Sel}_E(\mathbb{Q}_\infty)_p = 0$ for any prime $p$ where $E$ has good, ordinary reduction, unless either $p = 7$ or $a_p \equiv 1 \pmod p$. In these two cases, $\mathrm{Sel}_E(\mathbb{Q}_\infty)_p$ must be infinite by proposition 4.8. More precisely, theorem 4.1 implies the following: Let $p = 7$. Then $f_E(0) \sim 7$. (One must note that $a_7 = 3 \not\equiv 1 \pmod 7$.) This implies that $f_E(T)$ is an irreducible element of $\Lambda$. On the other hand, suppose $a_p \equiv 1 \pmod p$ but $p \neq 5$ or 61. Then $f_E(0) \sim p^2$. The only such anomalous prime $p \leq 100$ is $p = 43$. Assuming the validity of conjecture 1.11 for $E$, we will see later that $\lambda_E^{M-W} = 0$ and $\lambda_E^{III} = 2$ for $p = 43$.

$\mathbf{E : y^2 + xy = x^3 - 3x + 1}$. This is 34(A1) in [Cre]. $\mathrm{Sel}_E(\mathbb{Q})$ should be trivial. $E$ has multiplicative reduction at 2 and 17, $c_2 = 6$, $c_{17} = 1$. Also, $E(\mathbb{Q}) \cong \mathbb{Z}/6\mathbb{Z}$. The prime $p = 3$ is anomalous: $a_3 = -2$ and so $|\widetilde{E}(\mathbb{F}_3)| = 6$. If $p$ is any other prime where $E$ has good, ordinary reduction, then $a_p \not\equiv 1(\mathrm{mod}\ p)$ and we clearly have $\mathrm{Sel}_E(\mathbb{Q}_\infty)_p = 0$. For $p = 3$, proposition 4.1 gives $f_E(0) \sim 3$. Thus, $f_E(T)$ is irreducible. Let $F$ be the first layer of the cyclotomic $\mathbb{Z}_3$-extension of $\mathbb{Q}$. Then $F = \mathbb{Q}(\beta)$, where $\beta = \zeta + \zeta^{-1}$, $\zeta$ denoting a primitive 9-th root of unity. Notice that $\beta$ is a root of $x^3 - 3x + 1$. Thus $(\beta, -\beta)$ is a point in $E(F)$, which is not in $E(\mathbb{Q})$. Now the residue field for $\mathbb{Q}_\infty$ at the unique prime $\eta$ above 3 is $\mathbb{F}_3$. The prime-to-3 torsion of $E(\mathbb{Q}_\infty)$ is mapped by reduction modulo $\eta$ injectively into $\widetilde{E}(\mathbb{F}_3)$, and thus is $\mathbb{Z}/2\mathbb{Z}$. It is defined over $\mathbb{Q}$. The discussion preceding proposition 3.9 shows that $E(\mathbb{Q}_\infty)_3 = E(\mathbb{Q})_3$. Thus, $E(\mathbb{Q}_\infty)_{\mathrm{tors}} = E(\mathbb{Q})_{\mathrm{tors}}$. It follows that $(\beta, -\beta)$ has infinite order. Now $\mathrm{Gal}(F/\mathbb{Q})$ acts faithfully on $E(F) \otimes \mathbb{Q}_3$. It is clear that this $\mathbb{Q}_3$-representation must be isomorphic to $\rho^t$ where $\rho$ is the unique 2-dimensional irreducible $\mathbb{Q}_3$-representation of $\mathrm{Gal}(F/\mathbb{Q})$ and $t \geq 1$. If $\gamma$ generates $\Gamma = \mathrm{Gal}(\mathbb{Q}_\infty/\mathbb{Q})$ topologically, then $\rho(\gamma|_F)$ is given by a matrix with trace $-1$, determinant 1. Regarding $\rho$ as a representation of $\Gamma$ and letting $T = \gamma - 1$, $\rho(\gamma - 1)$ has characteristic polynomial

$$\theta_1 = (1 + T)^2 + (1 + T) + 1 = T^2 + 3T + 3.$$

Since $E(F) \otimes (\mathbb{Q}_3/\mathbb{Z}_3)$ is a $\Lambda$-submodule of $\mathrm{Sel}_E(\mathbb{Q}_\infty)_3$, it follows that $\theta_1^t$ divides $f_E(T)$. Comparing the valuation of $\theta_1(0)$ and $f_E(0)$, we clearly have $t = 1$ and $f_E(T) = \theta_1$, up to a factor in $\Lambda^\times$. Therefore $\lambda_E^{M-W} = 2$, $\lambda_E^{\mathrm{III}} = 0$, and $\mu_E = 0$.

When is $\mathrm{Sel}_E(\mathbb{Q}_\infty)_p$ infinite? A fairly complete answer is given by the following partial converse to proposition 3.8.

**Proposition 5.3.** *Assume that $E$ has good, ordinary reduction at $p$ and that $E(\mathbb{Q})$ has no element of order $p$. Assume also that at least one of the following statements is true:*

(i) $\mathrm{Sel}_E(\mathbb{Q})_p \neq 0$.
(ii) $a_p \equiv 1(\mathrm{mod}\ p)$.
(iii) *there exists at least one prime $l$ where $E$ has multiplicative reduction such that $a_l \equiv 1(\mathrm{mod}\ p)$ and $\mathrm{ord}_l(j_E) \equiv 0(\mathrm{mod}\ p)$.*
(iv) *there exists at least one prime $l$ where $E$ has additive reduction such that $E(\mathbb{Q}_l)$ has a point of order $p$.*

*Then $\mathrm{Sel}_E(\mathbb{Q}_\infty)_p$ is infinite.*

*Remark.* If $E$ has multiplicative reduction at $l$, then $a_l = \pm 1$. Thus, in (iii), $a_l \equiv 1(\mathrm{mod}\ p)$ is always true if $p = 2$ and is equivalent to $a_l = 1$ if $p$ is odd. Also, (iii) and (iv) simply state that there exists an $l$ such that $p|c_l$. If $E$ has additive reduction at $l$, then the only possible prime factors of $c_l$ are 2, 3, or $l$. Since $E$ has good reduction at $p$, (iv) can only occur for $p = 2$ or 3.

*Proof.* If $\mathrm{Sel}_E(\mathbb{Q})_p$ is infinite, the conclusion follows from theorem 1.2, or more simply from lemma 3.1. If $\mathrm{Sel}_E(\mathbb{Q})_p$ is finite, then we can apply proposition 4.1 to say that $f_E(0)$ is not in $\mathbb{Z}_p^\times$. Hence $f_E(T)$ is not invertible and so $X_E(\mathbb{Q}_\infty)$ must indeed be infinite. (The characteristic ideal of a finite $\Lambda$-module is $\Lambda$.) Alternatively, one can point out that since $E(\mathbb{Q})_p = 0$, it follows that $\ker(h_0) = 0$ and $\mathcal{G}_E^\Sigma(\mathbb{Q}) = \mathcal{P}_E^\Sigma(\mathbb{Q})$, where $\Sigma$ consists of $p, \infty$, and all primes of bad reduction. Hence, if (i) holds, then $\mathrm{Sel}_E(\mathbb{Q}_\infty)_p \neq 0$. If (ii), (iii), or (iv) holds, then $\ker(g_0) \neq 0$. Therefore, since $\ker(h_0)$ and $\mathrm{coker}(h_0)$ are both zero, we have $\mathrm{coker}(s_0) \neq 0$. This implies again that $\mathrm{Sel}_E(\mathbb{Q}_\infty)_p \neq 0$. Finally, proposition 4.8 then shows that $\mathrm{Sel}_E(\mathbb{Q}_\infty)_p$ must be infinite.    ■

As our examples show, quite a variety of possibilities for the data going into theorem 4.1 can arise. This is made even more clear from the following observation, where is a variant on lemma 8.19 of [Maz1].

**Proposition 5.4.** *Let $P$ and $L$ be disjoint, finite sets of primes. Let $Q$ be any finite set of primes. For each $p \in P$, let $a_p^*$ be any integer satisfying $|a_p^*| < 2\sqrt{p}$. For each $l \in L$, let $a_l^* = +1$ or $-1$. If $a_l^* = +1$, let $c_l^*$ be any positive integer. If $a_l^* = -1$, let $c_l^* = 1$ or 2. Then there exist infinitely many non-isomorphic elliptic curves $E$ defined over $\mathbb{Q}$ such that*

(i) *For each $p \in P$, $E$ has good reduction at $p$ and $a_p = a_p^*$.*

(ii) *For each $l \in L$, $E$ has multiplicative reduction at $l$, $a_l = a_l^*$, and $c_l = c_l^*$.*

(iii) *For each $q \in Q$, $E[q]$ is irreducible as a $\mathbb{F}_q$-representation space of $G_{\mathbb{Q}}$.*

*Proof.* This is an application of the Chinese Remainder Theorem. For each $p \in P$, a theorem of Deuring states that an elliptic curve $\widetilde{E}_p$ defined over $\mathbb{F}_p$ exists such that $|\widetilde{E}_p(\mathbb{F}_p)| = 1 + p - a_p^*$. One can then choose arbitrarily a lifting $E_p^*$ of $\widetilde{E}_p$ defined by a Weierstrass equation (as described in Tate's article [Ta]). We write this equation as $f_p^*(x, y) = 0$ where $f_p^*(x, y) \in \mathbb{Z}_p[x, y]$. Let $l \in L$. If $a_l^* = +1$, we let $E_l^*$ denote the Tate curve over $\mathbb{Q}_l$ with $j_{E_l^*} = l^{-c_l^*}$. Then $E_l^*$ has split, multiplicative reduction at $l$ and $\mathrm{ord}_l(j_E) = -c_l^*$. If $a_l^* = -1$, then we instead take $E_l^*$ as the unramified quadratic twist of this Tate curve, so that $E_l^*$ has non-split, multiplicative reduction. The index $[E_l^*(\mathbb{Q}_l) : E_{l,0}^*(\mathbb{Q}_l)]$ is then 1 or 2, depending on the parity of $c_l^*$. In either case, we let $f_l^*(x, y) = 0$ be a Weierstrass equation for $E_l^*$, where $f_l^*(x, y) \in \mathbb{Z}_l[x, y]$. Let $q \in Q$. Then we can choose a prime $r = r_q \neq q$ and an elliptic curve $\widetilde{E}_r$ defined over $\mathbb{F}_r$ such that $\widetilde{E}_r[q]$ is irreducible for the action of $G_{\mathbb{F}_r}$. If $q = 2$, this is easy. We take $r$ to be an odd prime and define $\widetilde{E}_r$ by $y^2 = g(x)$, where $g(x) \in \mathbb{F}_r[x]$ is an irreducible cubic polynomial. Then $\widetilde{E}_r(\mathbb{F}_r)$ has no element of order of 2, which suffices. If $q$ is odd, we choose $r$ to be an odd prime such that $-r$ is a quadratic nonresidue mod $q$. We can choose $\widetilde{E}/\mathbb{F}_r$ to be supersingular. Then the action of $\mathrm{Frob}_r \in G_{\mathbb{F}_r}$ on $\widetilde{E}_r[q]$ has characteristic polynomial $t^2 + r$. Since this has no roots in $\mathbb{F}_q$, $\widetilde{E}_r[q]$ indeed has no proper invariant subspaces under the action of $\mathrm{Frob}_r$. We can choose a lifting $E_r^*$ of $\widetilde{E}_r$ defined over $\mathbb{Q}_r$ by a Weierstrass equation $f_r^*(x, y) = 0$, where $f_r^*(x, y) \in \mathbb{Z}_r[x, y]$. For each $q \in Q$, infinitely many suitable $r_q$'s exist. Hence we can also require that the $r_q$'s are distinct and outside of $P \cap L$. We let $R = \{r_q\}_{q \in Q}$, choosing one $r_q$ for each $q \in Q$. We then choose an equation $f(x, y) = 0$ in Weierstrass form, where $f(x, y) \in \mathbb{Z}[x, y]$ and satisfies $f(x, y) \equiv f_m^*(x, y) \pmod{m^{t_m}}$ for all $m \in P \cup L \cup R$, where $t_m$ is chosen sufficiently large. The equation $f(x, y) = 0$ determines an elliptic curve defined over $\mathbb{Q}$. If $p \in P$, we just take $t_p = 1$. Then $E$ has good reduction at $p$, $\widetilde{E} = \widetilde{E}_p$, and hence $a_p = a_p^*$, as desired. If $r \in R$, then $r = r_q$ for some $q \in Q$. We take $t_r = 1$ again. $E$ has good reduction at $r$, $\widetilde{E} = \widetilde{E}_r$, and hence the action of a Frobenius automorphism in $\mathrm{Gal}(\mathbb{Q}(E[q])/\mathbb{Q})$ (for any prime above $r$) on $E[q]$ has no invariant subspaces. Hence obviously $E[q]$ is irreducible as an $\mathbb{F}_q$-representation space of $G_{\mathbb{Q}}$. Finally, suppose $l \in L$. If we take $t_l$ sufficiently large, then clearly $j_E$ will be close enough to $j_{E_l^*}$ to guarantee that $\mathrm{ord}_l(j_E) = \mathrm{ord}_l(j_{E_l^*}) = -c_l^*$. In terms of the coefficients of a Weierstrass equation over $\mathbb{Z}_l$, there is a simple criterion for an elliptic curve to have split or nonsplit reduction at $l$. (It involves the coset in $\mathbb{Q}_l^\times/(\mathbb{Q}_l^\times)^2$ containing the quantity $-c_4/c_6$ in the notation of Tate.) Hence it is clear that $E$ will have multiplicative reduction at $l$ and that $a_l = a_l^*$ if $t_l$ is taken large enough. Thus $E$ will have the required properties. The fact that infinitely many non-isomorphic $E$'s exist is clear, since we can vary $L$ and thus the set of primes where $E$ has bad reduction. ∎

*Remark.* We can assume that $P \cup L$ contains 3 and 5. Any elliptic curve $E$ defined over $\mathbb{Q}$ and satisfying (i) and (ii) will be semistable at 3 and 5 and therefore will be modular. This follows from a theorem of Diamond [D]. Furthermore, let $E_d$ denote the quadratic twist of $E$ by some square-free integer $d$. If we assume that all primes in $P \cup L$ split in $\mathbb{Q}(\sqrt{d})$, then $E_d$ also satisfies (i), (ii), and (iii). One can choose such $d$ so that $L(E_d/\mathbb{Q}, 1) \neq 0$. (See [B-F-H] for a discussion of this result which was first proved by Waldspurger.) A theorem of Kolyvagin then would imply that $E_d(\mathbb{Q})$ and $\mathrm{III}_{E_d}(\mathbb{Q})$ are finite. Thus, there in fact exist infinitely many non-isomorphic modular elliptic curves $E$ satisfying (i), (ii), and (iii) and such that $\mathrm{Sel}_E(\mathbb{Q})$ is finite.

**Corollary 5.5.** *Let $P$ be any finite set of primes. Then there exist infinitely many elliptic curves $E/\mathbb{Q}$ such that $E$ has good, ordinary reduction at $p$, $a_p = 1$, and $E[p]$ is an irreducible $\mathbb{F}_p$-representation space for $G_{\mathbb{Q}}$, for all $p \in P$.*

*Proof.* This follows immediately from proposition 5.9. One takes $P = Q$, $a_p^* = 1$ for all $p \in P$, and $L = \phi$.

**Corollary 5.6.** *Let $p$ be any prime. Assume that conjecture 1.11 is true when $F = \mathbb{Q}$. Then $\lambda_E$ is unbounded as $E$ varies over elliptic curves defined over $\mathbb{Q}$ with good, ordinary reduction at $p$.*

*Proof.* Take $P = \{p\} = Q$. Let $a_p^*$ be such that $p \nmid a_p^*$. Take $L$ to be a large finite set of primes. For each $l \in L$, let $a_l^* = +1$, $c_l^* = p$. Let $E/\mathbb{Q}$ be any elliptic curve satisfying the statements in proposition 5.4. As remarked above, we can assume $E$ is modular. Now $E$ has good, ordinary reduction at $p$. According to Theorem 1.5, $\mathrm{Sel}_E(\mathbb{Q}_\infty)_p$ is $\Lambda$-cotorsion. (Alternatively, we could assume that $\mathrm{Sel}_E(\mathbb{Q})_p$ is finite and then use theorem 1.4. The rest of this proof becomes somewhat easier if we make this assumption on $E$.) Also $\mu_E = 0$ by conjecture 1.11. We will show that $\lambda_E \geq |L|$, which certainly implies the corollary. Let $t = |L|$. Let $n = \mathrm{corank}_{\mathbb{Z}_p}(\mathrm{Sel}_E(\mathbb{Q})_p)$. Of course, $\lambda_E \geq n$ by theorem 1.2. So we can assume now that $n \leq t$. Let $\Sigma$ be the set of primes $p, \infty$, and all primes where $E$ has bad reduction. Then, by proposition 4.13, there are at most $p^n$ elements of order $p$ in $\mathcal{P}_E^\Sigma(\mathbb{Q})/\mathcal{G}_E^\Sigma(\mathbb{Q})$. Also, for each $l \in L$, we have $|\ker(r_l)| = c_l = p$. Thus, the kernel of the restriction map $\mathcal{P}_E^\Sigma(\mathbb{Q}) \to \mathcal{P}_E^\Sigma(\mathbb{Q}_\infty)$ contains a subgroup isomorphic to $(\mathbb{Z}/p\mathbb{Z})^t$. It follows that $\ker(g_0)$ contains a subgroup isomorphic to $(\mathbb{Z}/p\mathbb{Z})^{t-n}$. Now $\ker(h_0) = \mathrm{coker}(h_0) = 0$. Thus it follows that $\mathrm{coker}(s_0)$ contains a subgroup isomorphic to $(\mathbb{Z}/p\mathbb{Z})^{t-n}$. By proposition 4.14, and the assumption that $\mu_E = 0$, we have

$$\mathrm{Sel}_E(\mathbb{Q}_\infty)_p \cong (\mathbb{Q}_p/\mathbb{Z}_p)^{\lambda_E}.$$

But $\mathrm{Sel}_E(\mathbb{Q}_\infty)_p$ contains a subgroup $\mathrm{Im}(s_0)_{\mathrm{div}}$ isomorphic to $(\mathbb{Q}_p/\mathbb{Z}_p)^n$ and the corresponding quotient has a subgroup isomorphic to $(\mathbb{Z}/p\mathbb{Z})^{t-n}$. It follows that $\lambda_E \geq t$, as we claimed. ■

*Remark.* If we don't assume conjecture 1.11, then one still gets the weaker result that $\lambda_E + \mu_E$ is unbounded as $E$ varies over modular elliptic curves with good, ordinary reduction at a fixed prime $p$. For the above argument shows that $\dim_{\Lambda/\mathbf{m}\Lambda}(X_E(\mathbb{Q}_\infty)/\mathbf{m}X_E(\mathbb{Q}_\infty))$ is unbounded, where $\mathbf{m}$ denotes the maximal ideal of $\Lambda$. We then use the following result about $\Lambda$-modules: *Suppose $X$ is a finitely generated, torsion $\Lambda$-module and that $X$ has no nonzero, finite $\Lambda$-submodules. Let $\lambda$ and $\mu$ denote the corresponding invariants. Then*

$$\lambda + \mu \geq \dim_{\Lambda/\mathbf{m}\Lambda}(X/\mathbf{m}X).$$

The proof is not difficult. One first notes that the right-hand side, which is just the minimal number of generators of $X$ as a $\Lambda$-module, is "sub-additive" in an exact sequence $0 \to X_1 \to X_2 \to X_3 \to 0$ of $\Lambda$-modules. Both $\lambda$ and $\mu$ are additive. One then reduces to the special cases where either (a) $X$ has exponent $p$ and has no finite, nonzero $\Lambda$-submodule or (b) $X$ has no $\mathbb{Z}_p$-torsion. In the first case, $X$ is a $(\Lambda/p\Lambda)$-module. One then uses the fact that $\Lambda/p\Lambda$ is a PID. In the second case, $\lambda$ is the minimal number of generators of $X$ as a $\mathbb{Z}_p$-module. The inequality is clear.

We will now discuss the $\mu$-invariant $\mu_E$ of $\text{Sel}_E(\mathbb{Q}_\infty)_p$. We always assume that $E$ is defined over $\mathbb{Q}$ and has either good, ordinary or multiplicative reduction at $p$. According to conjecture 1.11, we should have $\mu_E = 0$ if $E[p]$ is irreducible as a $G_\mathbb{Q}$-module. Unfortunately, it seems very difficult to verify this even for specific examples. In this discussion we will assume that $E[p]$ is reducible as a $G_\mathbb{Q}$-module, i.e., that $E$ admits a cyclic $\mathbb{Q}$-isogeny of degree $p$. In [Maz2], Mazur proves that this can happen only for a certain small set of primes $p$. With the above restriction on the reduction type of $E$ at $p$, then $p$ is limited to the set $\{2, 3, 5, 7, 13, 37\}$. For $p = 2, 3, 5, 7,$ or $13$, there are infinitely many possible $E$'s, even up to quadratic twists. For $p = 37$, $E$ must be the elliptic curve defined by $y^2 + xy + y = x^3 + x^2 - 8x + 6$ (which has conductor $35^2$) or another elliptic curve related to this by a $\mathbb{Q}$-isogeny of degree 37, up to a quadratic twist.

Assume at first that $E[p]$ contains a $G_\mathbb{Q}$-invariant subgroup $\Phi$ isomorphic to $\mu_p$. We will let $\Sigma$ be the finite set consisting of $p$, $\infty$, and all primes where $E$ has bad reduction. Then we have a natural map

$$\epsilon : H^1(\mathbb{Q}_\Sigma/\mathbb{Q}_\infty, \Phi) \to H^1(\mathbb{Q}_\Sigma/\mathbb{Q}_\infty, E[p^\infty]).$$

It is easy to verify that $\ker(\epsilon)$ is finite. We also have the Kummer homomorphism

$$\beta : \mathcal{U}_\infty/\mathcal{U}_\infty^p \to H^1(\mathbb{Q}_\Sigma/\mathbb{Q}_\infty, \mu_p),$$

where $\mathcal{U}_\infty$ denotes the unit group of $\mathbb{Q}_\infty$. The map $\beta$ is injective. Dirichlet's unit theorem implies that the $(\Lambda/p\Lambda)$-module $\mathcal{U}_\infty/\mathcal{U}_\infty^p$ has corank 1. Consider a prime $l \neq p$. Let $\eta$ be a prime of $\mathbb{Q}_\infty$ lying over $l$. Then $(\mathbb{Q}_\infty)_\eta$ is the unramified $\mathbb{Z}_p$-extension of $\mathbb{Q}_l$ (which is the only $\mathbb{Z}_p$-extension of $\mathbb{Q}_l$). All units of $(\mathbb{Q}_\infty)_\eta$ are $p$-th powers. Thus, if $u \in \mathcal{U}_\infty$, then $u$ is a $p$-th power in $(\mathbb{Q}_\infty)_\eta$. Therefore, if $\varphi \in \text{Im}(\beta)$, then $\varphi|_{G_{(\mathbb{Q}_\infty)_\eta}}$ is trivial. If we fix an isomorphism $\Phi \cong \mu_p$, then it follows that the elements of $\text{Im}(\epsilon \circ \beta)$ satisfy the local conditions defining $\text{Sel}_E(\mathbb{Q}_\infty)_p$ at all primes $\eta$ of $\mathbb{Q}_\infty$ not lying over $p$ or $\infty$. Now assume that $p$ is odd. We can then ignore the archimedean primes of $\mathbb{Q}_\infty$. Since the inertia subgroup $I_{\mathbb{Q}_p}$ acts nontrivially on $\mu_p$ (because $p$ is odd) and acts trivially on $E[p^\infty]/C_p$, it follows that $\Phi \subseteq C_p$. If $\pi$ denotes the unique prime of $\mathbb{Q}_\infty$ lying over $p$, recall that $\text{Im}(\kappa_\pi) = \text{Im}(\lambda_\pi)$, where $\lambda_\pi$ is the map

$$H^1((\mathbb{Q}_\infty)_\pi, C_p) \to H^1((\mathbb{Q}_\infty)_\pi, E[p^\infty]).$$

Therefore, it is obvious that if $\varphi \in \text{Im}(\epsilon)$, then $\varphi|_{G_{(\mathbb{Q}_\infty)_\pi}} \in \text{Im}(\kappa_\pi)$. Combining the above observations, it follows that $\text{Im}(\epsilon \circ \beta) \subseteq \text{Sel}_E(\mathbb{Q}_\infty)_p$ if $p$ is odd. Thus, $\text{Sel}_E(\mathbb{Q}_\infty)_p$ contains a $\Lambda$-submodule of exponent $p$ with $(\Lambda/p\Lambda)$-corank equal to 1, which implies that either $\mu_E \geq 1$ or $\text{Sel}_E(\mathbb{Q}_\infty)_p$ is not $\Lambda$-cotorsion.

We will prove a more general result. Suppose that $E[p^\infty]$ has a $G_\mathbb{Q}$-invariant subgroup $\Phi$ which is cyclic of order $p^m$, with $m \geq 1$. If $E$ has semistable reduction at $p$, then it actually follows that $E$ has either good, ordinary reduction or multiplicative reduction at $p$. $\Phi$ has a $G_{\mathbb{Q}}$-composition series with composition factors isomorphic to $\Phi[p]$. We assume again that $p$ is an odd prime. Then the action of $I_{\mathbb{Q}_p}$ on $\Phi[p]$ is either trivial or given by the Teichmüller character $\omega$. In the first case, $\Phi$ is isomorphic as a $G_{\mathbb{Q}_p}$-module to a subgroup of $E[p^\infty]/C_p$. The action of $I_{\mathbb{Q}_p}$ on $\Phi$ is trivial and so we say that $\Phi$ is unramified at $p$. In the second case, we have $\Phi \subseteq C_p$ and we say that $\Phi$ is ramified at $p$. The action of $\mathrm{Gal}(\mathbb{C}/\mathbb{R})$ on $\Phi[p]$ determines its action on $\Phi$. We say that $\Phi$ is even or odd, depending on whether the action of $\mathrm{Gal}(\mathbb{C}/\mathbb{R})$ is trivial or nontrivial. With this terminology, we can state the following result.

**Proposition 5.7.** *Assume that $p$ is odd and that $E$ is an elliptic curve/$\mathbb{Q}$ with good, ordinary or multiplicative reduction at $p$. Assume that $\mathrm{Sel}_E(\mathbb{Q}_\infty)_p$ is $\Lambda$-cotorsion. Assume also that $E[p^\infty]$ contains a cyclic $G_{\mathbb{Q}}$-invariant subgroup $\Phi$ of order $p^m$ which is ramified at $p$ and odd. Then $\mu_E \geq m$.*

*Proof.* We will show that $\mathrm{Sel}_E(\mathbb{Q}_\infty)_p$ contains a $\Lambda$-submodule pseudo-isomorphic to $\widehat{\Lambda}[p^m]$. Consider the map

$$\epsilon : H^1(\mathbb{Q}_\Sigma/\mathbb{Q}_\infty, \Phi) \to H^1(\mathbb{Q}_\Sigma/\mathbb{Q}_\infty, E[p^\infty]).$$

The kernel is finite. Let $l \in \Sigma$, $l \neq p$ or $\infty$. There are just finitely many primes $\eta$ of $\mathbb{Q}_\infty$ lying over $l$. For each $\eta$, $H^1((\mathbb{Q}_\infty)_\eta, \Phi)$ is finite. (An easy way to verify this is to note that any Sylow pro-$p$ subgroup $V$ of $G_{(\mathbb{Q}_\infty)_\eta}$ is isomorphic to $\mathbb{Z}_p$ and that the restriction map $H^1((\mathbb{Q}_\infty)_\eta, \Phi) \to H^1(V, \Phi)$ is injective.) Therefore

$$\ker\big(H^1(\mathbb{Q}_\Sigma/\mathbb{Q}_\infty, \Phi) \to \prod_{\substack{l \in \Sigma \\ l \neq p, \infty}} \prod_{\eta \mid l} H^1((\mathbb{Q}_\infty)_\eta, \Phi)$$

has finite index in $H^1(\mathbb{Q}_\Sigma/\mathbb{Q}_\infty, \Phi)$. On the other hand, $\Phi \subseteq C_p$. Hence, elements in $\mathrm{Im}(\epsilon)$ automatically satisfy the local condition at $\pi$ occurring in the definition of $\mathrm{Sel}_E(\mathbb{Q}_\infty)_p$. These remarks imply that $\mathrm{Im}(\epsilon) \cap \mathrm{Sel}_E(\mathbb{Q}_\infty)_p$ has finite index in $\mathrm{Im}(\epsilon)$ and therefore $\mathrm{Sel}_E(\mathbb{Q}_\infty)_p$ contains a $\Lambda$-submodule pseudo-isomorphic to $H^1(\mathbb{Q}_\Sigma/\mathbb{Q}_\infty, \Phi)$.

One can study the structure of $H^1(\mathbb{Q}_\Sigma/\mathbb{Q}_\infty), \Phi)$ either by restriction to a subgroup of finite index in $\mathrm{Gal}(\mathbb{Q}_\Sigma/\mathbb{Q}_\infty)$ which acts trivially on $\Phi$ or by using Euler characteristics. We will sketch the second approach. The restriction map

$$H^1(\mathbb{Q}_\Sigma/\mathbb{Q}_n, \Phi) \to H^1(\mathbb{Q}_\Sigma/\mathbb{Q}_\infty, \Phi)^{\Gamma_n}$$

is surjective and its kernel is finite and has bounded order as $n \to \infty$. The Euler characteristic of the $\mathrm{Gal}(\mathbb{Q}_\Sigma/\mathbb{Q}_n)$-module $\Phi$ is $\prod_{v \mid \infty} |\Phi/\Phi^{D_v}|^{-1}$, where $v$ runs over the infinite primes of the totally real field $\mathbb{Q}_n$ and $D_v = G_{(\mathbb{Q}_n)_v}$. By assumption, $\Phi^{D_v} = 0$ and hence this Euler characteristic is $p^{-mp^n}$ for all $n \geq 0$. Therefore, $H^1(\mathbb{Q}_{\Sigma}/\mathbb{Q}_n, \Phi)$ has order divisible by $p^{mp^n}$. It follows that the $\Lambda$-module $H^1(\mathbb{Q}_{\Sigma}/\mathbb{Q}_{\infty}, \Phi)$, which is of exponent $p^m$ and hence certainly $\Lambda$-cotorsion, must have $\mu$-invariant $\geq m$. This suffices to prove that $\mu_E \geq m$. Under the assumptions that $E[p]$ is reducible as a $G_{\mathbb{Q}}$-module and that $\mathrm{Sel}_E(\mathbb{Q}_{\infty})_p$ is $\Lambda$-cotorsion, it follows from the next proposition that $\mathrm{Sel}_E(\mathbb{Q}_{\infty})_p$ contains a $\Lambda$-submodule pseudo-isomorphic to $\widehat{\Lambda}[p^{\mu_E}]$ and that the corresponding quotient has finite $\mathbb{Z}_p$-corank. Also, $\mathrm{Im}(\varepsilon)$ must almost coincide with $H^1(\mathbb{Q}_{\Sigma}/\mathbb{Q}_{\infty}, E[p^\infty])[p^m]$. (That is, the intersection of the two groups must have finite index in both.) This last $\Lambda$-module is pseudo-isomorphic to $\widehat{\Lambda}[p^m]$ according to the proposition below.

If $E$ is any elliptic curve/$\mathbb{Q}$ and $p$ is any prime, the weak Leopoldt conjecture would imply that $H^1(\mathbb{Q}_{\Sigma}/\mathbb{Q}_{\infty}, E[p^\infty])$ has $\Lambda$-corank equal to 1. That is, $H^1(\mathbb{Q}_{\Sigma}/\mathbb{Q}_{\infty}, E[p^\infty])_{\Lambda\text{-div}}$ should be pseudo-isomorphic to $\widehat{\Lambda}$. (This has been proven by Kato if $E$ is modular.) Here we will prove a somewhat more precise statement under the assumption that $E[p]$ is reducible as a $G_{\mathbb{Q}}$-module. It will be a rather simple consequence of the Ferrero-Washington theorem mentioned in the introduction. As usual, $\Sigma$ is a finite set of primes of $\mathbb{Q}$ containing $p, \infty$, and all primes where $E$ has bad reduction.

**Proposition 5.8.** *Assume that $E$ is an elliptic curve defined over $\mathbb{Q}$ and that $E$ admits a $\mathbb{Q}$-isogeny of degree $p$ for some prime $p$. Then $H^1(\mathbb{Q}_{\Sigma}/\mathbb{Q}_{\infty}, E[p^\infty])$ has $\Lambda$-corank 1. Furthermore, $H^1(\mathbb{Q}_{\Sigma}/\mathbb{Q}_{\infty}, E[p^\infty])/H^1(\mathbb{Q}_{\Sigma}/\mathbb{Q}_{\infty}, E[p^\infty])_{\Lambda\text{-div}}$ has $\mu$-invariant equal to 0 if $p$ is odd. If $p = 2$, this quotient has $\mu$-invariant equal to 0 or 1, depending on whether $E(\mathbb{R})$ has 1 or 2 connected components.*

*Proof.* First assume that $p$ is odd. Then we have an exact sequence

$$0 \to \Phi \to E[p] \to \Psi \to 0$$

where $\mathrm{Gal}(\mathbb{Q}_{\Sigma}/\mathbb{Q})$ acts on the cyclic groups $\Phi$ and $\Psi$ of order $p$ by characters $\varphi, \psi : \mathrm{Gal}(\mathbb{Q}_{\Sigma}/\mathbb{Q}) \to (\mathbb{Z}/p\mathbb{Z})^\times$. We know that $H^1(\mathbb{Q}_{\Sigma}/\mathbb{Q}_{\infty}, E[p^\infty])$ has $\Lambda$-corank $\geq 1$. Also, the exact sequence

$$0 \to E[p] \to E[p^\infty] \xrightarrow{p} E[p^\infty] \to 0$$

induces a surjective map $H^1(\mathbb{Q}_{\Sigma}/\mathbb{Q}_{\infty}, E[p]) \to H^1(\mathbb{Q}_{\Sigma}/\mathbb{Q}_{\infty}, E[p^\infty])[p]$ with finite kernel. Thus, it clearly is sufficient to prove that $H^1(\mathbb{Q}_{\Sigma}/\mathbb{Q}_{\infty}, E[p])$ has $(\Lambda/p\Lambda)$-corank 1. Now the determinant of the action of $G_{\mathbb{Q}}$ on $E[p]$ is the Teichmüller character $\omega$. Hence, $\varphi\psi = \omega$. Since $\omega$ is an odd character, one of the characters $\varphi$ or $\psi$ is odd, the other even. We have the following exact sequence:

$$H^1(\mathbb{Q}_{\Sigma}/\mathbb{Q}_{\infty}, \Phi) \to H^1(\mathbb{Q}_{\Sigma}/\mathbb{Q}_{\infty}, E[p]) \to H^1(\mathbb{Q}_{\Sigma}/\mathbb{Q}_{\infty}, \Psi)$$

and hence proposition 5.8 (for odd primes $p$) is a consequence of the following lemma.

**Lemma 5.9.** *Let $p$ be any prime. Let $\Theta$ be a $\mathrm{Gal}(\mathbb{Q}_{\Sigma}/\mathbb{Q})$-module which is cyclic of order $p$. Then $H^1(\mathbb{Q}_{\Sigma}/\mathbb{Q}_{\infty}, \Theta)$ has $(\Lambda/p\Lambda)$-corank 1 if $\Theta$ is odd or if $p = 2$. Otherwise, $H^1(\mathbb{Q}_{\Sigma}/\mathbb{Q}_{\infty}, \Theta)$ is finite. If $p = 2$, then the map*

$$\alpha : H^1(\mathbb{Q}_{\Sigma}/\mathbb{Q}_{\infty}, \Theta) \to \mathcal{P}_{\Theta}^{(\infty)}(\mathbb{Q}_{\infty})$$

*is surjective and has finite kernel. Here $\mathcal{P}_{\Theta}^{(\infty)}(\mathbb{Q}_{\infty}) = \varinjlim_{n} \prod_{v_n \mid \infty} H^1((\mathbb{Q}_n)_{v_n}, \Theta)$.*

*Remark.* We will use a similar notation to that introduced in the remark following lemma 4.6. For example, $\mathcal{P}_C^{(\ell)}(\mathbb{Q}_{\infty})$, which occurs in the following proof, is defined as $\varinjlim_{n} \prod_{v_n \mid \ell} H^1((\mathbb{Q}_n)_{v_n}, C)$. If $\ell$ is a nonarchimedean prime, then $\mathcal{P}_C^{(\ell)}(\mathbb{Q}_{\infty}) = \prod_{\eta \mid \ell} H^1((\mathbb{Q}_{\infty})_{\eta}, C)$, since there are only finitely many primes $\eta$ of $\mathbb{Q}_{\infty}$ lying over $\ell$.

*Proof.* Let $\theta$ be the character (with values in $(\mathbb{Z}/p\mathbb{Z})^{\times}$) which gives the action of $\mathrm{Gal}(\mathbb{Q}_{\Sigma}/\mathbb{Q})$ on $\Theta$. Let $C = (\mathbb{Q}_p/\mathbb{Z}_p)(\theta)$, where we now regard $\theta$ as a character of $\mathrm{Gal}(\mathbb{Q}_{\Sigma}/\mathbb{Q})$ with values in $\mathbb{Z}_p^{\times}$. Then $\Theta = C[p]$. We have an isomorphism

$$H^1(\mathbb{Q}_{\Sigma}/\mathbb{Q}_{\infty}, \Theta) \xrightarrow{\sim} H^1(\mathbb{Q}_{\Sigma}/\mathbb{Q}_{\infty}, C)[p].$$

(The surjectivity is clear. The injectivity follows from the fact that $H^0(\mathbb{Q}_{\Sigma}/\mathbb{Q}_{\infty}, C)$ is either $C$ or $0$, depending on whether $\theta$ is trivial or nontrivial.) We will relate the structure of $H^1(\mathbb{Q}_{\Sigma}/\mathbb{Q}_{\infty}, C)$ to various classical Iwasawa modules. Let $\Sigma' = \Sigma - \{p\}$. Consider

$$S_C'(\mathbb{Q}_{\infty}) = \ker\big(H^1(\mathbb{Q}_{\Sigma}/\mathbb{Q}_{\infty}, C) \to \prod_{\ell \in \Sigma'} \mathcal{P}_C^{(\ell)}(\mathbb{Q}_{\infty})\big).$$

If $\ell \in \Sigma'$ is nonarchimedean, then $H^1((\mathbb{Q}_{\infty})_{\eta}, C)$ is either trivial or isomorphic to $\mathbb{Q}_p/\mathbb{Z}_p$, for any prime $\eta$ of $\mathbb{Q}_{\infty}$ lying over $\ell$. $\mathcal{P}_C^{(\ell)}(\mathbb{Q}_{\infty})$ is then a cotorsion $\Lambda$-module with $\mu$-invariant 0. If $\ell = \infty$, then $(\mathbb{Q}_{\infty})_{\eta} = \mathbb{R}$ for any $\eta \mid \ell$. $H^1(\mathbb{R}, C)$ is, of course, trivial if $p$ is odd. But if $p = 2$, then $\theta$ is trivial and $H^1(\mathbb{R}, C) \cong \mathbb{Z}/2\mathbb{Z}$. Thus, in this case, $\mathcal{P}_C^{(\infty)}(\mathbb{Q}_{\infty})$ is isomorphic to $\mathrm{Hom}(\Lambda/2\Lambda, \mathbb{Z}/2\mathbb{Z}) = \widehat{\Lambda}[2]$, which is $\Lambda$-cotorsion and has $\mu$-invariant 1. It follows that $H^1(\mathbb{Q}_{\Sigma}/\mathbb{Q}_{\infty}, C)/S_C'(\mathbb{Q}_{\infty})$ is $\Lambda$-cotorsion and has $\mu$-invariant 0 if $p$ is odd. If $p = 2$, then the $\mu$-invariant is $\leq 1$.

Assume that $p$ is odd. Let $F$ be the cyclic extension of $\mathbb{Q}$ corresponding to $\theta$. (Thus, $F \subseteq \mathbb{Q}_{\Sigma}$ and $\theta$ is a faithful character of $\mathrm{Gal}(F/\mathbb{Q})$.) Then $F_{\infty} = F\mathbb{Q}_{\infty}$ is the cyclotomic $\mathbb{Z}_p$-extension of $F$. We let $\Delta = \mathrm{Gal}(F_{\infty}/\mathbb{Q}_{\infty}) \cong \mathrm{Gal}(F/\mathbb{Q})$. Let

$$X = \mathrm{Gal}(L_{\infty}/F_{\infty}), \qquad Y = \mathrm{Gal}(M_{\infty}/F_{\infty})$$

where $M_\infty$ is the maximal abelian pro-$p$ extension of $F_\infty$ unramified at all primes of $F_\infty$ not lying over $p$ and $L_\infty$ is the maximal subfield of $M_\infty$ unramified at the primes of $F_\infty$ over $p$ too. Now $\mathrm{Gal}(F_\infty/\mathbb{Q}) \cong \Delta \times \Gamma$ acts on both $X$ and $Y$ by inner automorphisms. Thus, they are both $\Lambda$-modules on which $\Delta$ acts $\Lambda$-linearly. That is, $X$ and $Y$ are $\Lambda[\Delta]$-modules.

The restriction map $H^1(\mathbb{Q}_\Sigma/\mathbb{Q}_\infty, C) \to H^1(\mathbb{Q}_\Sigma/F_\infty, C)^\Delta$ is an isomorphism. Also, $\mathrm{Gal}(\mathbb{Q}_\Sigma/F_\infty)$ acts trivially on $C$. Hence the elements of $H^1(\mathbb{Q}_\Sigma/F_\infty, C)$ are homomorphisms. Taking into account the local conditions, the restriction map induces an isomorphism

$$S_C'(\mathbb{Q}_\infty) \xrightarrow{\sim} \mathrm{Hom}_\Delta(Y, C) = \mathrm{Hom}(Y^\theta, C)$$

as $\Lambda$-modules, where $Y^\theta = e_\theta Y$, the $\theta$-component of the $\Delta$-module $Y$. (Here $e_\theta$ denotes the idempotent for $\theta$ in $\mathbb{Z}_p[\Delta]$.) Iwasawa proved that $Y^\theta$ is $\Lambda$-torsion if $\theta$ is even and has $\Lambda$-rank 1 if $\theta$ is odd. One version of the Ferrero-Washington theorem states that the $\mu$-invariant of $Y^\theta$ vanishes if $\theta$ is even. Thus, in this case, $H^1(\mathbb{Q}_\Sigma/\mathbb{Q}_\infty, C)$ must be $\Lambda$-cotorsion and have $\mu$-invariant 0. It then follows that $H^1(\mathbb{Q}_\Sigma/\mathbb{Q}_\infty, \Theta)$ must be finite. On the other hand, if $\theta$ is odd, then $S_C'(\mathbb{Q}_\infty)$ will have $\Lambda$-corank 1. Hence, the same is true of $H^1(\mathbb{Q}_\Sigma/\mathbb{Q}_\infty, C)$ and so $H^1(\mathbb{Q}_\Sigma/\mathbb{Q}_\infty, C)[p]$ will have $(\Lambda/p\Lambda)$-corank $\geq 1$. We will prove that equality holds and, therefore, $H^1(\mathbb{Q}_\Sigma/\mathbb{Q}_\infty, \Theta)$ indeed has $(\Lambda/p\Lambda)$-corank 1. It is sufficient to prove that $S_C'(\mathbb{Q}_\infty)[p]$ has $(\Lambda/p\Lambda)$-corank 1. We will deduce this from another version of the Ferrero-Washington theorem—the assertion that the torsion $\Lambda$-module $X$ has $\mu$-invariant 0. Let $\pi$ be the unique prime of $\mathbb{Q}_\infty$ lying over $p$. Consider

$$S_C(\mathbb{Q}_\infty) = \ker(S_C'(\mathbb{Q}_\infty) \to H^1((\mathbb{Q}_\infty)_\pi, C)).$$

In the course of proving lemma 2.3, we actually determined the structure of $H^1((\mathbb{Q}_\infty)_\pi, C)$. (See also section 3 of [Gr2].) It has $\Lambda$-corank 1 and the quotient $H^1((\mathbb{Q}_\infty)_\pi, C)/H^1((\mathbb{Q}_\infty)_\pi, C)_{\Lambda\text{-div}}$ is either trivial or isomorphic to $\mathbb{Q}_p/\mathbb{Z}_p$ as a group. To show that $S_C'(\mathbb{Q}_\infty)[p]$ has $(\Lambda/p\Lambda)$-corank 1, it suffices to prove that $S_C(\mathbb{Q}_\infty)[p]$ is finite. Now the restriction map identifies $S_C(\mathbb{Q}_\infty)$ with the subgroup of $\mathrm{Hom}_\Delta(Y, C)$ which is trivial on all the decomposition subgroups of $Y$ corresponding to primes of $F_\infty$ lying over $p$. Thus, $S_C(\mathbb{Q}_\infty)$ is isomorphic to a $\Lambda$-submodule of $\mathrm{Hom}_\Delta(X, C) = \mathrm{Hom}(X^\theta, C)$. Since the $\mu$-invariant of $X$ vanishes, it is clear that $S_C(\mathbb{Q}_\infty)[p]$ is indeed finite. This completes the proof of lemma 5.9 when $p$ is odd.

Now assume that $p = 2$. Thus, $\theta$ is trivial. We let $F_\infty = \mathbb{Q}_\infty$. Let $M_\infty$ be as defined above. Then it is easy to see that $M_\infty = \mathbb{Q}_\infty$. For let $M_0$ be the maximal abelian extension of $\mathbb{Q}$ contained in $M_\infty$. Thus, $\mathrm{Gal}(M_0/\mathbb{Q}_\infty) \cong Y/TY$. We must have $M_0 \subseteq \mathbb{Q}(\mu_{2^\infty})$. But $M_0$ is totally real and so clearly $M_0 = \mathbb{Q}_\infty$. Hence $Y/TY = 0$ and hence that $M_\infty = \mathbb{Q}_\infty$. Therefore, $S_C'(\mathbb{Q}_\infty) = 0$. It follows that $H^1(\mathbb{Q}_\Sigma/\mathbb{Q}_\infty, C)$ is $\Lambda$-cotorsion and has $\mu$-invariant $\leq 1$. In fact, the $\mu$-invariant is 1 and arises in the following way. Let $\mathcal{U}_\infty$ denote the unit group of $\mathbb{Q}_\infty$. Let $K_\infty =$ $\mathbb{Q}_\infty(\{\sqrt{u}|u \in \mathcal{U}_\infty\})$. Then $K_\infty \subseteq M_\infty^*$, the maximal abelian pro-2 extension of $\mathbb{Q}_\infty$ unramified outside of the primes over $p$ and $\infty$. Also, one can see that $\mathrm{Gal}(K_\infty/\mathbb{Q}_\infty) \cong \Lambda/2\Lambda$. Thus, clearly $H^1(\mathbb{Q}_\Sigma/\mathbb{Q}_\infty, C)[2] = H^1(\mathbb{Q}_\Sigma/\mathbb{Q}_\infty, \Theta)$ contains the $\Lambda$-submodule $\mathrm{Hom}(\mathrm{Gal}(K_\infty/\mathbb{Q}_\infty), \Theta)$ which has $\mu$-invariant 1. To complete the proof of lemma 5.9, we point out that $K_\infty$ can't contain any totally real subfield larger than $\mathbb{Q}_\infty$, since $M_\infty = \mathbb{Q}_\infty$. That is, $\ker(\alpha) \cap \mathrm{Hom}(\mathrm{Gal}(K_\infty/\mathbb{Q}_\infty), \Theta)$ is trivial. It follows that $\ker(\alpha)$ is finite. We also see that $\alpha$ must be surjective because $\mathcal{P}_\Theta^{(\infty)}(\mathbb{Q}_\infty)$ is isomorphic to $\widehat{\Lambda}[2]$. ∎

We must complete the proof of proposition 5.8 for $p = 2$. Consider the following commutative diagram with exact rows:

$$
\begin{array}{ccccc}
H^1(\mathbb{Q}_\Sigma/\mathbb{Q}_\infty, \Phi) & \xrightarrow{a} & H^1(\mathbb{Q}_\Sigma/\mathbb{Q}_\infty, E[2]) & \xrightarrow{b} & H^1(\mathbb{Q}_\Sigma/\mathbb{Q}_\infty, \Psi) \\
\downarrow{\alpha_1} & & \downarrow{\alpha_2} & & \downarrow{\alpha_3} \\
\mathcal{P}_\Phi^{(\infty)}(\mathbb{Q}_\infty) & \xrightarrow{c} & \mathcal{P}_{E[2]}^{(\infty)}(\mathbb{Q}_\infty) & \xrightarrow{d} & \mathcal{P}_\Psi^{(\infty)}(\mathbb{Q}_\infty)
\end{array}
$$

By lemma 5.9, both $\alpha_1$ and $\alpha_3$ are surjective and have finite kernel. Also, $\ker(a)$ is finite. We see that $H^1(\mathbb{Q}_\Sigma/\mathbb{Q}_\infty, E[2])$ has $(\Lambda/2\Lambda)$-corank equal to 1 or 2. First assume that $E(\mathbb{R})$ is connected, i.e., that the discriminant of a Weierstrass equation for $E$ is negative. Then $H^1(\mathbb{R}, E[2]) = 0$, and so $\mathcal{P}_{E[2]}^{(\infty)}(\mathbb{Q}_\infty) = 0$. It follows that $d \circ \alpha_2$ is the zero map and hence $\mathrm{Im}(b)$ is finite. Thus, $H^1(\mathbb{Q}_\Sigma/\mathbb{Q}_\infty, E[2])$ is pseudo-isomorphic to $H^1(\mathbb{Q}_\Sigma/\mathbb{Q}_\infty, \Phi)$ and so has $(\Lambda/2\Lambda)$-corank 1. In this case, $H^1(\mathbb{Q}_\Sigma/\mathbb{Q}_\infty, E[2^\infty])$ must have $\Lambda$-corank 1 and its maximal $\Lambda$-cotorsion quotient must have $\mu$-invariant 0. This proves proposition 5.8 in the case that $E(\mathbb{R})$ is connected.

Now assume that $E(\mathbb{R})$ has two components, i.e., that a Weierstrass equation for $E$ has positive discriminant. Then $E[2] \subseteq E(\mathbb{R})$ and $H^1(\mathbb{R}, E[2]) \cong (\mathbb{Z}/2\mathbb{Z})^2$. The $(\Lambda/2\Lambda)$-module $\mathcal{P}_{E[2]}^{(\infty)}(\mathbb{Q}_\infty)$ is isomorphic to $\widehat{\Lambda}[2]^2$. In this case, we will see that $H^1(\mathbb{Q}_\Sigma/\mathbb{Q}_\infty, E[2])$ has $(\Lambda/2\Lambda)$-corank 2. This is clear if $E[2] \cong \Phi \times \Psi$ as a $G_\mathbb{Q}$-module. If $E[2]$ is a nonsplit extension of $\Psi$ by $\Phi$, then $F = \mathbb{Q}(E[2])$ is a real quadratic field contained in $\mathbb{Q}_\Sigma$. Let $F_\infty = F\mathbb{Q}_\infty$. Considering the field $K_\infty = F_\infty(\{\sqrt{u}|u \in \mathcal{U}_{F_\infty}\})$, where $\mathcal{U}_{F_\infty}$ is the group of units of $F_\infty$, one finds that $H^1(\mathbb{Q}_\Sigma/F_\infty, \Phi)$ and $H^1(\mathbb{Q}_\Sigma/F_\infty, \Psi)$ have $(\Lambda/2\Lambda)$-corank 2. Now $E[2] \cong \Phi \times \Psi$ as a $G_F$-module and so $H^1(\mathbb{Q}_\Sigma/F_\infty, E[2])$ has $(\Lambda/2\Lambda)$-corank 4. The inflation-restriction sequence then will show that $H^1(\mathbb{Q}_\Sigma/\mathbb{Q}_\infty, E[2])$ is pseudo-isomorphic to $H^1(\mathbb{Q}_\Sigma/F_\infty, E[2])^\Delta$, where $\Delta = \mathrm{Gal}(F_\infty/\mathbb{Q}_\infty)$. One then sees that $H^1(\mathbb{Q}_\Sigma/\mathbb{Q}_\infty, E[2])$ must have $(\Lambda/2\Lambda)$-corank 2. The facts that $c$ is injective and that both $\alpha_1$ and $\alpha_3$ have finite kernel implies that $\alpha_2$ has finite kernel too. The map $\alpha_2$ must therefore be surjective. Now consider the commutative diagram

$$
\begin{array}{ccc}
H^1(\mathbb{Q}_\Sigma/\mathbb{Q}_\infty, E[2]) & \xrightarrow{\alpha_2} & \mathcal{P}^{(\infty)}_{E[2]}(\mathbb{Q}_\infty) \\
\downarrow & & \downarrow{\scriptstyle e} \\
H^1(\mathbb{Q}_\Sigma/\mathbb{Q}_\infty, E[2^\infty]) & \xrightarrow{\alpha_E} & \mathcal{P}^{(\infty)}_{E[2^\infty]}(\mathbb{Q}_\infty)
\end{array}
$$

Note that $\mathcal{P}^{(\infty)}_{E[2^\infty]}(\mathbb{Q}_\infty)$ is what we denoted by $\mathcal{P}^{(\infty)}_E(\mathbb{Q}_\infty)$ in section 4. The map $H^1(\mathbb{R}, E[2]) \to H^1(\mathbb{R}, E[2^\infty])$ is surjective. (But it's not injective since $H^1(\mathbb{R}, E[2^\infty]) \cong \mathbb{Z}/2\mathbb{Z}$ when $E(\mathbb{R})$ has two components.) Hence the map $e$ is surjective. Thus $e \circ \alpha_2$ is surjective and this implies that $\alpha_E$ is surjective. In fact, more precisely, the above diagram shows that the restriction of $\alpha_E$ to $H^1(\mathbb{Q}_\Sigma/\mathbb{Q}_\infty, E[2^\infty])[2]$ is surjective.

We can now easily finish the proof of proposition 5.8. Clearly

$$
H^1(\mathbb{Q}_\Sigma/\mathbb{Q}_\infty, E[2^\infty])_{\Lambda\text{-div}} \subseteq \ker(\alpha_E).
$$

Since $\ker(\alpha_E)[2]$ has $(\Lambda/2\Lambda)$-corank 1, it is clear that $H^1(\mathbb{Q}_\Sigma/\mathbb{Q}_\infty, E[2^\infty])$ has $\Lambda$-corank 1 and that $\ker(\alpha_E)/H^1(\mathbb{Q}_\Sigma/\mathbb{Q}_\infty, E[2^\infty])_{\Lambda\text{-div}}$ has $\mu$-invariant 0. Hence the maximal $\Lambda$-cotorsion quotient of $H^1(\mathbb{Q}_\Sigma/\mathbb{Q}_\infty, E[2^\infty])$ has $\mu$-invariant 1. ∎

*Remark.* Assume that $E$ is an elliptic curve/$\mathbb{Q}$ which has a $\mathbb{Q}$-isogeny of degree $p$. Assuming that $\mathrm{Sel}_E(\mathbb{Q}_\infty)_p$ is $\Lambda$-cotorsion, the above results show that $\mathrm{Sel}_E(\mathbb{Q}_\infty)_p$ contains a $\Lambda$-submodule pseudo-isomorphic to $\widehat{\Lambda}[p^{\mu_E}]$. Thus the $\mu$-invariant of $\mathrm{Sel}_E(\mathbb{Q}_\infty)_p$ arises "non-semisimply" if $\mu_E > 1$. For odd $p$, we already noted this before. For $p = 2$, it follows from the above discussion of $\ker(\alpha_E)$ and the fact that $\mathrm{Sel}_E(\mathbb{Q}_\infty)_p \subseteq \ker(\alpha_E)$. If $E$ has no $\mathbb{Q}$-isogeny of degree $p$, then $\mu_E$ is conjecturally 0, although there has been no progress on proving this.

Before describing various examples where $\mu_E$ is positive, we will prove another consequence of lemma 5.9 (and its proof).

**Proposition 5.10.** *Assume that $p$ is odd and that $E$ is an elliptic curve/$\mathbb{Q}$ with good, ordinary or multiplicative reduction at $p$. Assume also that $E[p^\infty]$ contains a $G_\mathbb{Q}$-invariant subgroup $\Phi$ of order $p$ which is either ramified at $p$ and even or unramified at $p$ and odd. Then $\mathrm{Sel}_E(\mathbb{Q}_\infty)_p$ is $\Lambda$-cotorsion and $\mu_E = 0$.*

*Proof.* We will show that $\mathrm{Sel}_E(\mathbb{Q}_\infty)[p]$ is finite. This obviously implies the conclusion. We have the exact sequence

$$
H^1(\mathbb{Q}_\Sigma/\mathbb{Q}_\infty, \Phi) \xrightarrow{a} H^1(\mathbb{Q}_\Sigma/\mathbb{Q}_\infty, E[p]) \xrightarrow{b} H^1(\mathbb{Q}_\Sigma/\mathbb{Q}_\infty, \Psi)
$$

as before. Under the above hypotheses, both $\Phi^{G_{\mathbb{Q}_\infty}}$ and $\Psi^{G_{\mathbb{Q}_\infty}}$ are trivial. Hence $H^0(\mathbb{Q}_\Sigma/\mathbb{Q}_\infty, E[p]) = 0$. This implies that

$$H^1(\mathbb{Q}_\Sigma/\mathbb{Q}_\infty, E[p]) \xrightarrow{\sim} H^1(\mathbb{Q}_\Sigma/\mathbb{Q}_\infty, E[p^\infty])[p]$$

under the natural map. Thus we can regard $\mathrm{Sel}_E(\mathbb{Q}_\infty)[p]$ as a subgroup of $H^1(\mathbb{Q}_\Sigma/\mathbb{Q}_\infty, E[p])$. Assume that $\mathrm{Sel}_E(\mathbb{Q}_\infty)[p]$ is infinite. Hence either $B = b(\mathrm{Sel}_E(\mathbb{Q}_\infty)[p])$ or $A = \mathrm{Im}(a) \cap \mathrm{Sel}_E(\mathbb{Q}_\infty)[p]$ is infinite. Assume first that $B$ is infinite. Then, by lemma 5.9, $\Psi$ must be odd. Hence $\Psi$ is unramified, $\Phi$ is ramified at $p$. Let $\bar{\pi}$ be any prime of $\mathbb{Q}_\Sigma$ lying over the prime $\pi$ of $\mathbb{Q}_\infty$ over $p$. Then $\Phi = C_{\bar{\pi}}[p]$, where $C_{\bar{\pi}}$ is the subgroup of $E[p^\infty]$ occurring in propositions 2.2, 2.4. (For example, if $E$ has good reduction at $p$, then $C_{\bar{\pi}}$ is the kernel of reduction modulo $\bar{\pi} : E[p^\infty] \to \widetilde{E}[p^\infty]$.) The inertia subgroup $I_{\bar{\pi}}$ of $\mathrm{Gal}(\mathbb{Q}_\Sigma/\mathbb{Q}_\infty)$ for $\bar{\pi}$ acts trivially on $D_{\bar{\pi}} = E[p^\infty]/C_{\bar{\pi}}$. Thus, $\Psi$ can be identified with $D_{\bar{\pi}}[p]$. Let $\sigma$ be a 1-cocycle with values in $E[p]$ representing a class in $\mathrm{Sel}_E(\mathbb{Q}_\infty)[p]$. Let $\tilde{\sigma}$ be the induced 1-cocycle with values in $\Psi$. Since $H^1(I_{\bar{\pi}}, D_{\bar{\pi}}) = \mathrm{Hom}(I_{\bar{\pi}}, D_{\bar{\pi}})$, it is clear that $\tilde{\sigma}|_{I_{\bar{\pi}}} = 0$. Thus, $\tilde{\sigma} \in H^1(\mathbb{Q}_\Sigma/\mathbb{Q}_\infty, \Psi)$ is unramified at $\bar{\pi}$. Now for each of the finite number of primes $\eta$ of $\mathbb{Q}_\infty$ lying over some $\ell \in \Sigma, \ell \neq p$, $H^1((\mathbb{Q}_\infty)_\eta, \Psi)$ is finite. Thus, it is clear that $B \cap H^1_{\mathrm{unr}}(\mathbb{Q}_\Sigma/\mathbb{Q}_\infty, \Psi)$ is of finite index in $B$ and is therefore infinite, where $H^1_{\mathrm{unr}}(\mathbb{Q}_\Sigma/\mathbb{Q}_\infty, \Psi)$ denotes the group of everywhere unramified cocycle classes. However, if we let $F$ denote the extension of $\mathbb{Q}$ corresponding to $\psi$, then we see that

$$H^1_{\mathrm{unr}}(\mathbb{Q}_\Sigma/\mathbb{Q}_\infty, \Psi) = \mathrm{Hom}(X^\psi, \Psi),$$

where we are using the same notation as in the proof of proposition 5.9. The Ferrero-Washington theorem implies that $H^1_{\mathrm{unr}}(\mathbb{Q}_\Sigma/\mathbb{Q}_\infty, \Psi)$ is finite. Hence in fact $B$ must be finite. Similarly, if $A$ is infinite, then $\Phi$ must be odd and hence unramified. Thus, $\Phi \cap C_{\bar{\pi}} = 0$. If $\sigma$ is as above, then $\sigma|_{I_{\bar{\pi}}}$ must have values in $C_{\bar{\pi}}$. But if $\sigma$ represents a class in $A$, then we can assume that its values are in $\Phi$. Thus $\sigma|_{I_{\bar{\pi}}} = 0$. Now the map $H^1(I_{\bar{\pi}}, \Phi) \to H^1(I_{\bar{\pi}}, E[p])$ is injective. Thus, we see just as above, that $H^1_{\mathrm{unr}}(\mathbb{Q}_\Sigma/\mathbb{Q}_\infty, \Phi)$ is infinite, again contradicting the Ferrero-Washington theorem. ∎

Later we will prove analogues of propositions 5.7 and 5.10 for $p = 2$. One can pursue the situation of proposition 5.10 much further, obtaining for example a simple formula for $\lambda_E$ in terms of the $\lambda$-invariant of $X^\theta$, where $\theta$ is the odd character in the pair $\varphi, \psi$. (Remark: Obviously, $\varphi\psi = \omega$. It is known that $X^\theta$ and $Y^{\omega\theta^{-1}}$ have the same $\lambda$-invariants, when $\theta$ is odd. Both $\Lambda$-modules occur in the arguments.) As mentioned in the introduction, one can prove conjecture 1.13 when $E/\mathbb{Q}$ has good, ordinary reduction at $p$ and satisfies the other hypotheses in proposition 5.10. The key ingredients are Kato's theorem and a comparison of $\lambda$-invariants based on a congruence between $p$-adic $L$-functions. We will pursue these ideas fully in [GrVa].

Another interesting idea, which we will pursue more completely elsewhere, is to study the relationship between $\mathrm{Sel}_E(\mathbb{Q}_\infty)_p$ and $\mathrm{Sel}_{E'}(\mathbb{Q}_\infty)_p$ when $E$ and $E'$ are elliptic curves/$\mathbb{Q}$ such that $E[p] \cong E'[p]$ as $G_\mathbb{Q}$-modules. If $E$ and $E'$ have good, ordinary or multiplicative reduction at $p$ and if $p$ is odd, then it is not difficult to prove the following result: *if* $\mathrm{Sel}_E(\mathbb{Q}_\infty)_p[p]$ *is finite, then so is* $\mathrm{Sel}_{E'}(\mathbb{Q}_\infty)_p[p]$. It follows that if $\mathrm{Sel}_E(\mathbb{Q}_\infty)_p$ is $\Lambda$-cotorsion and if $\mu_E = 0$, then $\mathrm{Sel}_{E'}(\mathbb{Q}_\infty)_p$ is also $\Lambda$-cotorsion and $\mu_{E'} = 0$. Furthermore, it is then possible to relate the $\lambda$-invariants $\lambda_E$ and $\lambda_{E'}$ to each other. (They usually will not be equal. The relationship involves the sets of primes of bad reduction and the Euler factors at those primes.)

A theorem of Washington [Wa1] as well as a generalization due to E. Friedman [F], which are somewhat analogous to the Ferrero-Washington theorem, can also be used to obtain nontrivial results. This idea was first exploited in [R-W] to prove that $E(K)$ is finitely generated for certain elliptic curves $E$ and certain infinite abelian extensions $K$ of $\mathbb{Q}$. The proof of proposition 5.10 can be easily modified to prove some results of this kind. Here is one.

**Proposition 5.11.** *Assume that $E$ and $p$ satisfy the hypotheses of proposition 5.10. Let $K$ denote the cyclotomic $\mathbb{Z}_q$-extension of $\mathbb{Q}$, where $q$ is any prime different than $p$. Then* $\mathrm{Sel}_E(K)[p]$ *is finite. Hence*

$$\mathrm{Sel}_E(K)_p \cong (\mathbb{Q}_p/\mathbb{Z}_p)^t \times \ (\text{a finite group})$$

*for some $t \geq 0$.*

Washington's theorem would state that the power of $p$ dividing the class number of the finite layers in the $\mathbb{Z}_q$-extension $FK/F$ is bounded. To adapt the proof of proposition 5.10, one can replace $\mathrm{Im}(\kappa_\eta)$ by $\mathrm{Im}(\lambda_\eta)$ for each prime $\eta$ of $K$ lying over $p$, obtaining a possibly larger subgroup of $H^1(K, E[p^\infty])$. The arguments also work if $\mathrm{Gal}(K/\mathbb{Q}) \cong \prod_{i=1}^n \mathbb{Z}_q$, where $q_1, q_2, \ldots, q_n$ are distinct primes, possibly including $p$. Then one uses the main result of [F]. One consequence is that $E(K)$ is finitely generated. If $E$ is any modular elliptic curve/$\mathbb{Q}$, this same statement is a consequence of the work of Kato and Rohrlich.

We will now discuss various examples where $\mu_E > 0$. We will take the base field to be $\mathbb{Q}$ and assume always that $E$ is an elliptic curve/$\mathbb{Q}$ with good, ordinary or multiplicative reduction at $p$. We assume first that $p$ is odd. Since $V_p(E)$ is irreducible as a representation space for $G_\mathbb{Q}$, there is a maximal subgroup $\Phi$ of $E[p^\infty]$ such that $\Phi$ is cyclic, $G_\mathbb{Q}$-invariant, ramified at $p$, and odd. Define $m_E$ by $|\Phi| = p^{m_E}$. Thus, $m_E \geq 0$. Proposition 5.8 states that $\mu_E \geq m_E$. It is not hard to see that conjecture 1.11 is equivalent to the assertion that $\mu_E = m_E$. For $p = 2$, $m_E$ can be $0, 1, 2, 3$, or $4$. For $p = 3$ or $5$, $m_E$ can be $0$, $1$, or $2$. For $p = 7$, $13$, or $37$, $m_E$ can be $0$ or $1$. For other odd primes (where $E$ has the above reduction type), there are no $\mathbb{Q}$-isogenies of degree $p$ and so $m_E = 0$. In [Maz1], there is a complete discussion of conductor 11 and numerous other examples having non-trivial $p$-isogenies.

**Conductor = 11**. If $E$ has conductor 11, then $E[p]$ is irreducible except for $p = 5$. Let $E_1$, $E_2$, and $E_3$ denote the curves 11A1, 11A2, and 11A3 in Cremona's tables. Thus $E_1 = X_0(11)$ and one has $E_1[5] \cong \mu_5 \times \mathbb{Z}/5\mathbb{Z}$ as a $G_\mathbb{Q}$-module. For $E_2$ (which is $E_1/(\mathbb{Z}/5\mathbb{Z})$), one has the nonsplit exact sequence

$$0 \to \mu_5 \to E_2[5] \to \mathbb{Z}/5\mathbb{Z} \to 0.$$

Now $E_2/\mu_5 = E_1$ and so one sees that $E_2[5^\infty]$ contains a subgroup $\Phi$ which is cyclic of order 25, $G_\mathbb{Q}$-invariant, ramified at 5, and odd. ($\Phi$ is an extension of $\mu_5$ by $\mu_5$.) For $E_3$ (which is $E_1/\mu_5$), one has a nonsplit exact sequence

$$0 \to \mathbb{Z}/5\mathbb{Z} \to E_3[5] \to \mu_5 \to 0.$$

All of these statements follow from the data about isogenies and torsion subgroups given in [Cre]. One then sees easily that $m_{E_1} = 1$, $m_{E_2} = 2$, and $m_{E_3} = 0$. We will show that $\mathrm{Sel}_{E_1}(\mathbb{Q}_\infty)_5 \cong \widehat{\Lambda}[5]$, $\mathrm{Sel}_{E_2}(\mathbb{Q}_\infty)_5 \cong \widehat{\Lambda}[5^2]$, and $\mathrm{Sel}_{E_3}(\mathbb{Q}_\infty)_5 = 0$. Thus, $\lambda_{E_i} = 0$ and $\mu_{E_i} = m_{E_i}$ for $1 \le i \le 3$.

We will let $\Phi_i = \mu_5$ and $\Psi_i = \mathbb{Z}/5\mathbb{Z}$ as $G_\mathbb{Q}$-modules for $1 \le i \le 3$. Then we have the following exact sequences of $G_\mathbb{Q}$-modules

$$0 \longrightarrow \Phi_2 \longrightarrow E_2[5] \longrightarrow \Psi_2 \longrightarrow 0$$

$$0 \longrightarrow \Psi_3 \longrightarrow E_3[5] \longrightarrow \Phi_3 \longrightarrow 0.$$

These exact sequences are nonsplit. For $E_1$, we have $E_1[5] = \Phi_1 \times \Psi_1$. As $G_{\mathbb{Q}_5}$-modules, we have exact sequences

$$0 \longrightarrow C_5 \longrightarrow E_i[5^\infty] \longrightarrow D_5 \longrightarrow 0$$

where $D_5$ is unramified and $C_5 \cong \mu_{5^\infty}$ for the action of $I_{\mathbb{Q}_5}$, the inertia subgroup of $G_{\mathbb{Q}_5}$. There will be no need to index $C_5$ and $D_5$ by $i$. As $G_{\mathbb{Q}_{11}}$-modules, we have exact sequences

$$0 \longrightarrow C_{11} \longrightarrow E_i[5^\infty] \longrightarrow D_{11} \longrightarrow 0$$

where $C_{11} \cong \mu_{5^\infty}$ and $D_{11} \cong \mathbb{Q}_5/\mathbb{Z}_5$ for the action of $G_{\mathbb{Q}_{11}}$. It will again not be necessary to include an index $i$ on these groups. The homomorphisms $E_i[5^\infty] \to D_5$ and $E_i[5^\infty] \to D_{11}$ induce natural identifications. As $G_{\mathbb{Q}_5}$-modules, $\Psi_1$, $\Psi_2$, $\Psi_3$ are all identified with $D_5[5]$. This is clear from the action of $G_{\mathbb{Q}_5}$ on these groups (which is trivial). But, as $G_{\mathbb{Q}_{11}}$-modules, $\Phi_1$, $\Psi_1$, $\Psi_2$, and $\Phi_3$ are all identified with $D_{11}[5]$. One verifies this by using the isogeny data and the fact that the Tate periods for the $E_i$'s in $\mathbb{Q}_{11}^\times$ have valuations 5, 1, 1, respectively. For example, if $\Phi_1$ or $\Psi_1$ were contained in $C_{11}$, then the

Tate period for $E_2$ or $E_3$ would have valuation divisible by 5. We will use the fact that the maps

$$H^1(\mathbb{Q}_{11}, D_{11}[5]) \to H^1(\mathbb{Q}_{11}, D_{11}), \qquad H^1(I_{\mathbb{Q}_5}, D_5[5]) \to H^1(I_{\mathbb{Q}_5}, D_5)$$

are both injective. This is so because $G_{\mathbb{Q}_{11}}$ acts trivially on $D_{11}$ and $I_{\mathbb{Q}_5}$ acts trivially on $D_5$. Our calculations of the Selmer groups will be in several steps and depend mostly on the results of section 2 and 3. We take $\Sigma = \{\infty, 5, 11\}$.

$\mathrm{Sel}_{E_3}(\mathbb{Q})_5 = 0$. Suppose $[\sigma] \in \mathrm{Sel}_{E_3}(\mathbb{Q})[5]$. It is enough to prove that $[\sigma] = 0$. We can assume that $\sigma$ has values in $E_3[5]$. (But note that in this case the map $H^1(\mathbb{Q}_\Sigma/\mathbb{Q}, E_3[5]) \to H^1(\mathbb{Q}_\Sigma/\mathbb{Q}_1 E_3[5^\infty])$ has a nontrivial kernel.) The image of $\sigma$ in $H^1(\mathbb{Q}_{11}, \Phi_3) = H^1(\mathbb{Q}_{11}, D_{11}[5])$ must become trivial in $H^1(\mathbb{Q}_{11}, D_{11})$. Thus this image must be trivial. Now $\Phi_3 = \mu_5$ and $H^1(\mathbb{Q}_\Sigma/\mathbb{Q}, \mu_5) \cong (\mathbb{Z}/5\mathbb{Z})^2$, the classes for the 1-cocycles associated to $\sqrt[5]{5^i 11^j}$, $0 \le i, j \le 4$. The restriction of such a 1-cocycle to $G_{\mathbb{Q}_{11}}$ is trivial when $5^i 11^j \in (\mathbb{Q}_{11}^\times)^5$, which happens only when $i = j = 0$. Thus, the image of $[\sigma]$ in $H^1(\mathbb{Q}_\Sigma/\mathbb{Q}, \Phi_3)$ must be trivial. Hence we can assume that $\sigma$ has values in $\Psi_3 = \mathbb{Z}/5\mathbb{Z}$.

Now $H^1(\mathbb{Q}_\Sigma/\mathbb{Q}, \Psi_3) \cong (\mathbb{Z}/5\mathbb{Z})^2$ by class field theory, but its image in $H^1(\mathbb{Q}_\Sigma/\mathbb{Q}, E[5^\infty])$ is of order 5. Since $[\sigma] \in \mathrm{Sel}_{E_3}(\mathbb{Q})_5$, it has a trivial image in $H^1(I_{\mathbb{Q}_5}, D_5)$. Hence, regarding $\sigma$ as an element of $\mathrm{Hom}(\mathrm{Gal}(\mathbb{Q}_\Sigma/\mathbb{Q}), \Psi_3)$, it must be unramified at 5 and hence factor through $\mathrm{Gal}(K/\mathbb{Q})$, where $K$ is the cyclic extension of $\mathbb{Q}$ of conductor 11. But this implies that $[\sigma] = 0$ in $H^1(\mathbb{Q}_\Sigma/\mathbb{Q}, E[5^\infty])$ because $\mathrm{Hom}(\mathrm{Gal}(K/\mathbb{Q}), \Psi_3)$ is the kernel of the map $H^1(\mathbb{Q}_\Sigma/\mathbb{Q}, \Psi_3) \to H^1(\mathbb{Q}_\Sigma/\mathbb{Q}, E[5^\infty])$. To see this, note that this kernel has order 5 and that the map $H^1(I_{\mathbb{Q}_5}, \Psi_3) \to H^1(I_{\mathbb{Q}_5}, D_5)$ is injective. Hence $\mathrm{Sel}_{E_3}(\mathbb{Q})_5 = 0$.

$\mathrm{Sel}_{E_2}(\mathbb{Q})_5 = 0$. We have $H^1(\mathbb{Q}_\Sigma/\mathbb{Q}, E_2[5]) \cong H^1(\mathbb{Q}_\Sigma/\mathbb{Q}, E_2[5^\infty])[5]$. Let $[\sigma] \in \mathrm{Sel}_{E_2}(\mathbb{Q})[5]$. We can assume that $\sigma$ has values in $E_2[5]$. The image of $\sigma$ in $H^1(\mathbb{Q}_\Sigma/\mathbb{Q}, \Psi_2)$ must have a trivial restriction to $G_{\mathbb{Q}_{11}}$. But

$$H^1(\mathbb{Q}_\Sigma/\mathbb{Q}, \Psi_2) = \mathrm{Hom}(\mathrm{Gal}(KL/\mathbb{Q}), \mathbb{Z}/5\mathbb{Z}),$$

where $K$ is as above and $L$ is the first layer of the cyclotomic $\mathbb{Z}_5$-extension of $\mathbb{Q}$. Now 11 is inert in $L/\mathbb{Q}$ and ramified in $KL/L$. Thus it is clear that $\sigma$ has trivial image in $H^1(\mathbb{Q}_\Sigma/\mathbb{Q}, \Psi_2)$ and hence has values in $\Phi_2 = \mu_5$.

Now $H^1(\mathbb{Q}_\Sigma/\mathbb{Q}, \mu_5) \cong (\mathbb{Z}/5\mathbb{Z})^2$, but the map

$$\epsilon_0 : H^1(\mathbb{Q}_\Sigma/\mathbb{Q}, \mu_5) \to H^1(\mathbb{Q}_\Sigma/\mathbb{Q}, E_2[5]) = H^1(\mathbb{Q}_\Sigma/\mathbb{Q}, E_2[5^\infty])[5]$$

has $\ker(\epsilon_0) \cong \mathbb{Z}/5\mathbb{Z}$. Now $[\sigma] \in \mathrm{Im}(\epsilon_0)$, which we will show is not contained in $\mathrm{Sel}_{E_2}(\mathbb{Q})_5$. This will imply that $\mathrm{Sel}_{E_2}(\mathbb{Q})_5 = 0$. Consider the commutative

diagram

$$\begin{array}{ccccc}
H^1(\mathbb{Q}_\Sigma/\mathbb{Q}, \mu_5) & \xrightarrow{\ a\ } & H^1(\mathbb{Q}_5, \mu_5) & \xrightarrow{\ b\ } & H^1(\mathbb{Q}_5, C_5) \\
\downarrow{\scriptstyle \epsilon_0} & & & & \downarrow{\scriptstyle \lambda_5} \\
H^1(\mathbb{Q}_\Sigma/\mathbb{Q}, E_2[5^\infty]) & \xrightarrow{\hspace{3.5cm} c \hspace{3.5cm}} & & & H^1(\mathbb{Q}_5, E_2[5^\infty])
\end{array}$$

One sees easily that $a$ is an isomorphism. Also, $H^1(\mathbb{Q}_5, \mu_5) \cong (\mathbb{Z}/5\mathbb{Z})^2$ and $b$ induces an isomorphism $H^1(\mathbb{Q}_5, \mu_5) \cong H^1(\mathbb{Q}_5, C_5)[5]$. Referring to (2) following the proof of lemma 2.3, one sees that $H^1(\mathbb{Q}_5, C_5) \cong (\mathbb{Q}_5/\mathbb{Z}_5) \times \mathbb{Z}/5\mathbb{Z}$. (One needs the fact that $|\widetilde{E}_2(\mathbb{Z}/5\mathbb{Z})|$ is divisible by 5, but not by $5^2$.) In section 2, one also finds a proof that the map

$$H^1(\mathbb{Q}_5, C_5)/H^1(\mathbb{Q}_5, C_5)_{\mathrm{div}} \to \mathrm{Im}(\lambda_5)/\mathrm{Im}(\lambda_5)_{\mathrm{div}}$$

is an isomorphism. (See (3) in the proof of proposition 2.5.) If we had $\mathrm{Im}(\epsilon_0) \subseteq \mathrm{Sel}_{E_2}(\mathbb{Q})_5$, then we must have $\mathrm{Im}(c \circ \epsilon_0) \subseteq \mathrm{Im}(\lambda_5)_{\mathrm{div}}$, which is the image of the local Kummer homomorphism $\kappa_5$. But this can't be so because clearly $\mathrm{Im}(b \circ a) \not\subset H^1(\mathbb{Q}_5, C_5)_{\mathrm{div}}$. It follows that $\mathrm{Sel}_{E_2}(\mathbb{Q})_5 = 0$.

Although we don't need it, we will determine $\ker(\epsilon_0)$. The discussion in the previous paragraph shows that $\ker(\epsilon_0)$ is the inverse image under $b \circ a$ of $H^1(\mathbb{Q}_5, C_5)_{\mathrm{div}}[5]$. One can use proposition 3.11 to determine this. Let $\varphi$ be the unramified character of $G_{\mathbb{Q}_5}$ giving the action in $D_5 = \widetilde{E}_2[5^\infty]$. Since 5 is an anomalous prime for $E_2$, one gets an isomorphism

$$\varphi : \mathrm{Gal}(M_\infty/\mathbb{Q}_5) \to 1 + 5\mathbb{Z}_5$$

where $M_\infty$ denotes the unramified $\mathbb{Z}_5$-extension of $\mathbb{Q}_5$. One has

$$H^1(M_\infty, \mu_{5^\infty}) \cong \widehat{R} \times \mathbb{Q}_5/\mathbb{Z}_5,$$

where now $R = \mathbb{Z}_p[[G]]$, $G = \mathrm{Gal}(M_\infty/\mathbb{Q}_5)$. We have $C_5 \cong \mu_{5^\infty} \otimes \varphi^{-1}$ and $H^1(M_\infty, C_5) = H^1(M_\infty, \mu_{5^\infty}) \otimes \varphi^{-1}$. Now $H^1(\mathbb{Q}_5, C_5) \xrightarrow{\sim} H^1(M_\infty, C_5)^G$, by the inflation-restriction sequence. The image of $H^1(\mathbb{Q}_5, C_5)_{\mathrm{div}}$ under the restriction map is $H^1(M_\infty, C_5)^G_{R\text{-div}}$. But $H^1(M_\infty, C_5)_{R\text{-div}}$ coincides with $H^1(M_\infty, \mu_{5^\infty})_{R\text{-div}}$, with the action of $G$ twisted by $\varphi^{-1}$. Let $q \in \mathbb{Q}_5^\times$ and let $\sigma_q$ be the 1-cocycle with values in $\mu_5$ associated to $\sqrt[5]{q}$. Then $\sigma_q \in H^1(\mathbb{Q}_5, C_5)_{\mathrm{div}}$ if and only if $\sigma_q|_{G_{M_\infty}} \in H^1(M_\infty, \mu_{5^\infty})_{R\text{-div}}$. By proposition 3.11, this means that $q$ is a universal norm for $M_\infty/\mathbb{Q}_5$, i.e., $q \in \mathbb{Z}_5^\times$. Now $H^1(\mathbb{Q}_\Sigma/\mathbb{Q}, \mu_5)$ consists of the classes of 1-cocycles associated to $\sqrt[5]{u}$, where $u = 5^i 11^j$, $0 \le i, j \le 4$. It follows that $\ker(\epsilon_0)$ is generated by the 1-cocycle corresponding to $\sqrt[5]{11}$. There are other ways to interpret this result. The extension class of $\mathbb{Z}/5\mathbb{Z}$ by $\mu_5$ given by $E_2[5]$ corresponds to the 1-cocycle associated to $\sqrt[5]{11}$. The field $\mathbb{Q}(E_2[5])$ is $\mathbb{Q}(\mu_5, \sqrt[5]{11})$. The Galois module $E_2[5]$ is "peu ramifiée" at 5, in the sense of Serre. (This of course must be so because $E_2$ has good reduction at 5.)

$\mathrm{Sel}_{E_1}(\mathbb{Q})_5 = 0$. We have an exact sequence

$$0 \to H^1(\mathbb{Q}_\Sigma/\mathbb{Q}, \Psi_1) \to H^1(\mathbb{Q}_\Sigma/\mathbb{Q}, E_1[5^\infty]) \to H^1(\mathbb{Q}_\Sigma/\mathbb{Q}, E_2[5^\infty]).$$

Since $\mathrm{Sel}_{E_2}(\mathbb{Q})_5 = 0$, it is clear that $\mathrm{Sel}_{E_1}(\mathbb{Q})_5 \subseteq \mathrm{Im}(H^1(\mathbb{Q}_\Sigma/\mathbb{Q}, \Psi_1))$. But $\Psi_1 = \mathbb{Z}/5\mathbb{Z}$ and $H^1(\mathbb{Q}_\Sigma/\mathbb{Q}, \Psi_1) = \mathrm{Hom}(\mathrm{Gal}(KL/\mathbb{Q}), \mathbb{Z}/5\mathbb{Z})$, where $K$ and $L$ are as defined before. Since the decomposition group for 11 in $\mathrm{Gal}(KL/\mathbb{Q})$ is the entire group and since $\Psi_1$ is mapped to $D_{11}[5]$, we see as before that $H^1(\mathbb{Q}_\Sigma/\mathbb{Q}, \Psi_1) \to H^1(\mathbb{Q}_{11}, E_1[5^\infty])$ is injective. Hence $\mathrm{Sel}_{E_1}(\mathbb{Q})_5 = 0$.

$f_{E_i}(T) = 5^{m_{E_i}}$. We can now apply theorem 4.1 to see that $f_{E_1}(0) \sim 5$, $f_{E_2}(0) \sim 5^2$, and $f_{E_3}(0) \sim 1$, using the fact that $\widetilde{E}_i(\mathbb{Z}/5\mathbb{Z})$ has order 5. But we know that $5^{m_{E_i}}$ divides $f_{E_i}(T)$. Hence it follows that, after multiplication by a factor in $\Lambda^\times$, we can take $f_{E_1}(T) = 5$, $f_{E_2}(T) = 5^2$, and $f_{E_3}(T) = 1$. We now determine directly the precise structure of the Selmer groups $\mathrm{Sel}_{E_i}(\mathbb{Q}_\infty)_5$ as $\Lambda$-modules.

$\mathrm{Sel}_{E_3}(\mathbb{Q}_\infty)_5 = 0$. The fact that $f_{E_3}(T) = 1$ shows that $\mathrm{Sel}_{E_3}(\mathbb{Q}_\infty)_5$ is finite. Proposition 4.15 then implies that $\mathrm{Sel}_{E_3}(\mathbb{Q}_\infty)_5 = 0$. However, it is interesting to give a more direct argument. We will show that the restriction map $s_0^{(3)} : \mathrm{Sel}_{E_3}(\mathbb{Q})_5 \to \mathrm{Sel}_{E_3}(\mathbb{Q}_\infty)_5^\Gamma$ is surjective, which then implies that $\mathrm{Sel}_{E_3}(\mathbb{Q}_\infty)_5^\Gamma$ and hence $\mathrm{Sel}_{E_3}(\mathbb{Q}_\infty)_5$ are both zero. Here and in the following discussions, we will let $s_0^{(i)}$, $h_0^{(i)}$, $g_0^{(i)}$, and $r_v^{(i)}$ for $v \in \{5, 11\}$ denote the maps considered in sections 3 and 4 for the elliptic curve $E_i$, $1 \le i \le 3$. Thus, $\ker(s_0^{(i)}) = 0$ for $1 \le i \le 3$, by proposition 3.9. But $|\ker(h_0^{(3)})| = 5$. We have the exact sequence

$$0 \to \ker(h_0^{(3)}) \to \ker(g_0^{(3)}) \to \mathrm{coker}(s_0^{(3)}) \to 0.$$

Thus it suffices to show that $|\ker(g_0^{(3)})| = 5$. We let

$$A = \ker(\mathcal{P}_{E_3}^\Sigma(\mathbb{Q}) \to \mathcal{P}_{E_3}^\Sigma(\mathbb{Q}_\infty)).$$

Now $\mathcal{P}_{E_3}^\Sigma(\mathbb{Q}) = \mathcal{H}_{E_3}(\mathbb{Q}_5) \times \mathcal{H}_{E_3}(\mathbb{Q}_{11})$, $\mathcal{P}_{E_3}^\Sigma(\mathbb{Q}_\infty) = \mathcal{H}_{E_3}(\mathbb{Q}_5^{\mathrm{cyc}}) \times \mathcal{H}_{E_3}(\mathbb{Q}_{11}^{\mathrm{cyc}})$. The local duality theorems easily imply that

$$\mathcal{H}_{E_3}(\mathbb{Q}_5) = H^1(\mathbb{Q}_5, E_3[5^\infty])/\mathrm{Im}(\kappa_5) \cong (\mathbb{Q}_5/\mathbb{Z}_5) \times \mathbb{Z}/5\mathbb{Z}$$

$$\mathcal{H}_{E_3}(\mathbb{Q}_{11}) = H^1(\mathbb{Q}_{11}, E_3[5^\infty]) \cong \mathbb{Z}/5\mathbb{Z}.$$

The kernels of the maps $r_v^{(i)} : \mathcal{H}_{E_i}(\mathbb{Q}_v) \to \mathcal{H}_{E_i}(\mathbb{Q}_v^{\mathrm{cyc}})$ can be determined by the results in section 3. In particular, one finds that $|\ker(r_5^{(3)})| = 5^2$, while $r_{11}^{(3)}$ is injective. Also, $\mathcal{H}_{E_i}(\mathbb{Q}_5^{\mathrm{cyc}}) \cong H^1(\mathbb{Q}_5^{\mathrm{cyc}}, D_5) \cong \widehat{\Lambda}$ for $1 \le i \le 3$. Thus, $\mathrm{Im}(r_5^{(i)})$ is obviously isomorphic to $\mathbb{Q}_5/\mathbb{Z}_5$ for each $i$. It follows that $A$ has order $5^2$, $A \subseteq \mathcal{H}_{E_3}(\mathbb{Q}_5)$, $A \cap \mathcal{H}_{E_3}(\mathbb{Q}_5)_{\mathrm{div}}$ has order 5, and $A\mathcal{H}_{E_3}(\mathbb{Q}_5)_{\mathrm{div}} = \mathcal{H}_{E_3}(\mathbb{Q}_5)$. Now $\mathcal{G}_{E_3}^\Sigma(\mathbb{Q})$ has index 5 in $\mathcal{P}_{E_3}^\Sigma(\mathbb{Q})$ and projects onto $\mathcal{H}_{E_3}(\mathbb{Q}_{11})$. It follows easily that

$$|\ker(g_0^{(3)})| = |A \cap \mathcal{G}_{E_3}^\Sigma(\mathbb{Q})| = 5.$$

As we said, this implies that $\mathrm{Sel}_{E_3}(\mathbb{Q}_\infty)_5 = 0$.

$\mathrm{Sel}_{E_2}(\mathbb{Q}_\infty)_5 \cong \widehat{\Lambda}[5^2]$. Let $\Phi$ be the $G_{\mathbb{Q}}$-invariant subgroup of $E_2[5^\infty]$ which is cyclic of order $5^2$. (This $\Phi$ is an extension of $\Phi_1$ by $\Phi_2$.) We have $E_2/\Phi \cong E_3$. Since $\mathrm{Sel}_{E_3}(\mathbb{Q}_\infty)_5 = 0$, it follows that

$$\mathrm{Sel}_{E_2}(\mathbb{Q}_\infty)_5 \subseteq \mathrm{Im}(H^1(\mathbb{Q}_\Sigma/\mathbb{Q}_\infty, \Phi) \to H^1(\mathbb{Q}_\Sigma/\mathbb{Q}_\infty, E_2[5^\infty])).$$

The index is finite by proposition 5.7. Thus it is clear that $\mathrm{Sel}_{E_2}(\mathbb{Q}_\infty)_5$ is pseudo-isomorphic to $\widehat{\Lambda}[5^2]$ and has exponent $5^2$. Since $E_2(\mathbb{Q}) = 0$, we have $\mathcal{G}_{E_2}^\Sigma(\mathbb{Q}) = \mathcal{P}_{E_2}^\Sigma(\mathbb{Q})$ and $\ker(h_0^{(2)}) = 0$. Hence

$$\mathrm{coker}(s_0^{(2)}) \cong \ker(g_0^{(2)}) \cong \ker(r_5^{(2)}) \times \ker(r_{11}^{(2)}).$$

Now $\ker(r_{11}^{(2)}) = 0$ because $5 \nmid \mathrm{ord}_{11}(q_{E_2}^{(11)})$, where $q_{E_2}^{(11)}$ denotes the Tate period for $E_2$ in $\mathbb{Q}_{11}^\times$. Also, $|\ker(r_5^{(2)})| = 5^2$. We pointed out earlier that the $G_{\mathbb{Q}}$-module $E_2[5]$ is the nonsplit extension of $\mathbb{Z}/5\mathbb{Z}$ by $\mu_5$ corresponding to $\sqrt[5]{11}$. Since $11 \notin (\mathbb{Q}_5^\times)^5$, this extension remains nonsplit as a $G_{\mathbb{Q}_5}$-module. Thus, $H^0(\mathbb{Q}_5, E_2[5^\infty]) = 0$. One deduces from this that $H^1(\mathbb{Q}_5, E_2[5^\infty]) \cong \mathbb{Q}_5/\mathbb{Z}_5$ and $\mathcal{H}_{E_2}(\mathbb{Q}_5) \cong \mathbb{Q}_5/\mathbb{Z}_5$. This implies that $\ker(r_5^{(2)}) \cong \mathbb{Z}/5^2\mathbb{Z}$. Hence $\ker(g_0^{(2)})$, $\mathrm{coker}(s_0^{(2)})$ and hence $\mathrm{Sel}_{E_2}(\mathbb{Q}_\infty)_5^\Gamma$ are all cyclic of order $5^2$. Therefore, $X_{E_2}(\mathbb{Q}_\infty) = \mathrm{Sel}_{E_2}(\mathbb{Q}_\infty)_5^\wedge$ is a cyclic $\Lambda$-module of exponent $5^2$. That is, $X_{E_2}(\mathbb{Q}_\infty)$ is a quotient of $\Lambda/5^2\Lambda$ and, since the two are pseudo-isomorphic, it follows easily that $X_{E_2}(\mathbb{Q}_\infty) \cong \Lambda/5^2\Lambda$. This gives the stated result about the structure of $\mathrm{Sel}_{E_2}(\mathbb{Q}_\infty)_5$.

$\mathrm{Sel}_{E_1}(\mathbb{Q}_\infty)_5 \cong \widehat{\Lambda}[5]$. Since $E_1/\Phi_1 \cong E_3$, it follows

$$\mathrm{Sel}_{E_1}(\mathbb{Q}_\infty)_5 \subseteq \mathrm{Im}(H^1(\mathbb{Q}_\Sigma/\mathbb{Q}, \Phi_1) \to H^1(\mathbb{Q}_\Sigma/\mathbb{Q}_\infty, E_1[5^\infty])).$$

Hence $\mathrm{Sel}_{E_1}(\mathbb{Q}_\infty)_5$ has exponent 5 and is pseudo-isomorphic to $\widehat{\Lambda}[5]$. Also, by proposition 4.15, $\mathrm{Sel}_{E_1}(\mathbb{Q}_\infty)_5$ has no proper $\Lambda$-submodules of finite index. Thus, $X_{E_1}(\mathbb{Q}_\infty)$ is a $(\Lambda/5\Lambda)$-module pseudo-isomorphic to $(\Lambda/5\Lambda)$ and with no nonzero, finite $\Lambda$-submodules. Since $\Lambda/5\Lambda$ is a PID, it follows that $X_{E_1}(\mathbb{Q}_\infty) \cong \Lambda/5\Lambda$, which gives the stated result concerning the structure of $\mathrm{Sel}_{E_1}(\mathbb{Q}_\infty)_5$.

*Twists.* Let $\xi$ be a quadratic character for $\mathbb{Q}$. Then $\xi$ corresponds to a quadratic field $\mathbb{Q}(\sqrt{d})$, where $d = d_\xi \in \mathbb{Z}$ and $|d|$ is the conductor of $\xi$. We consider separately the cases where $\xi$ is even or odd. For even $\xi$, the following conjecture seems reasonable. It can be deduced from conjecture 1.11, but may be more approachable. We let $E^\xi$ denote the quadratic twist of $E$ by $d$.

**Conjecture 5.12.** *Let $E$ be an elliptic curve/$\mathbb{Q}$ with potentially ordinary or multiplicative reduction at $p$, where $p$ is an odd prime. Let $\xi$ be an even quadratic character. Then $\mathrm{Sel}_{E^\xi}(\mathbb{Q}_\infty)_p$ and $\mathrm{Sel}_E(\mathbb{Q}_\infty)_p$ have the same $\mu$-invariants.*

We remark that the $\lambda$-invariants can certainly be different. For example, if $E$ is any one of the three elliptic curves of conductor 11, then $\lambda_E = 0$ for any prime $p$ satisfying the hypothesis in the above conjecture. But if $\xi$ is the quadratic character corresponding to $\mathbb{Q}(\sqrt{2})$ (of conductor 8), then rank$(E^\xi(\mathbb{Q})) = 1$. (In fact, $E^\xi$ is 704(A1, 2, or 3) in [Cre].) Then of course $\lambda_{E^\xi} \geq 1$ for all such $p$.

Assume now that $\xi$ is an odd character and that $\xi(5) \neq 0$. Let $E$ be any one of the elliptic curves of conductor 11. Let $p = 5$. Then $E^\xi[5]$ is $G_\mathbb{Q}$-reducible with composition factors $\mu_5 \otimes \xi$ and $(\mathbb{Z}/5\mathbb{Z}) \otimes \xi$. The hypotheses in proposition 5.10 are satisfied and so the $\mu$-invariant of $\mathrm{Sel}_{E^\xi}(\mathbb{Q}_\infty)_5$ is zero. The $\lambda$-invariant $\lambda_{E^\xi}$ is unchanged by isogeny and so doesn't depend on the choice of $E$. It follows from proposition 4.14 that $\mathrm{Sel}_{E^\xi}(\mathbb{Q}_\infty)_5$ is divisible. Hence $\lambda_{E^\xi}$, which is the $\mathbb{Z}_5$-corank of $\mathrm{Sel}_{E^\xi}(\mathbb{Q}_\infty)_5$, is obviously equal to the $(\mathbb{Z}/5\mathbb{Z})$-dimension of $\mathrm{Sel}_{E^\xi}(\mathbb{Q}_\infty)_5[5]$. We will not give the verification (which we will discuss more generally elsewhere), but one finds the following formula:

$$\lambda_{E^\xi} = 2\lambda_\xi + \epsilon_\xi,$$

where $d = d_\xi$ and $\lambda_\xi$ denotes the classical $\lambda$-invariant $\lambda(F_\infty/F)$ for the imaginary quadratic field $F = \mathbb{Q}(\sqrt{d})$ and for the prime $p = 5$ and where $\epsilon_\xi = 1$ if 11 splits in $\mathbb{Q}(\sqrt{d})/\mathbb{Q}$, $\epsilon_\xi = 0$ if 11 is inert or ramified. By proposition 3.10, it follows that corank$_{\mathbb{Z}_5}(\mathrm{Sel}_{E^\xi}(\mathbb{Q})_5)) \equiv \epsilon_\xi \pmod 2$, which is in agreement with the Birch and Swinnerton-Dyer conjecture since the sign in the functional equation for the Hasse-Weil $L$-function $L(E^\xi/\mathbb{Q}, 5) = L(E/\mathbb{Q}, \xi, 5)$ is $(-1)^{\epsilon_\xi}$. As an example, consider the case where $\xi$ corresponds to $\mathbb{Q}(\sqrt{-2})$. Then $E^\xi$ is 704(K1, 2, or 3). The class number of $\mathbb{Q}(\sqrt{-2})$ is 1. The prime $p = 5$ is inert in $F = \mathbb{Q}(\sqrt{-2})$. Hence the discussion of Iwasawa's theorem in the introduction shows that the $\lambda$-invariant for this quadratic field is 0. But 11 splits in $F$. Therefore, $\lambda_{E^\xi} = 1$. Since rank$(E^\xi(\mathbb{Q})) = 1$, it is clear that $\mathrm{Sel}_{E^\xi}(\mathbb{Q}_\infty)_5 = E^\xi(\mathbb{Q}) \otimes (\mathbb{Q}_p/\mathbb{Z}_p)$. As another example, suppose that $\xi$ corresponds to $F = \mathbb{Q}(\sqrt{-1})$. Then $E^\xi$ is 176(B1, 2, or 3) in [Cre]. The prime 5 splits in $F/\mathbb{Q}$ and so $\lambda(F_\infty/F) \geq 1$. In fact, $\lambda(F_\infty/F) = 1$. Since 11 is inert in $F$, we have $\lambda_{E^\xi} = 2$. But $E^\xi(\mathbb{Q})$ is trivial. If $E^\xi(\mathbb{Q}_1)$ had positive rank, one would have rank$(E^\xi(\mathbb{Q}_1)) \geq 4$ (because the nontrivial irreducible $\mathbb{Q}$-representation of Gal$(\mathbb{Q}_1/\mathbb{Q})$ has degree 4). Hence it is clear that $\lambda_{E^\xi}^\mathrm{III} = 2$, $\lambda_{E^\xi}^{M-W} = 0$. T. Fukuda has done extensive calculations of $\lambda(F_\infty/F)$ when $F$ is an imaginary quadratic field and $p = 3, 5$, or 7. Some of these $\lambda$-invariants are quite large. Presumably they are unbounded as $F$ varies. For $p = 5$, he finds that $\lambda_\xi = 10$ if $\xi$ corresponds to $F = \mathbb{Q}(\sqrt{-3,624,233})$. Since 11 splits in $F/\mathbb{Q}$, we have $\lambda_{E^\xi} = 21$ in this case. However, we don't know the values of $\lambda_{E^\xi}^{M-W}$ and $\lambda_{E^\xi}^\mathrm{III}$.

We will briefly explain in the case of $E^\xi$ (where $E$ and $\xi$ are as in the previous paragraph and $p = 5$) how to prove conjecture 1.13. Kato's theorem states that $f_E^\xi(T)$ divides $f_{E^\xi}^\mathrm{anal}(T)$, up to multiplication by a power of $p$. Thus, $\lambda(f_{E^\xi}^\mathrm{anal}) \geq \lambda_{E^\xi}$. Now it is known that $\lambda_\xi$ is equal to the $\lambda$-invariant of the Kubota-Leopoldt 5-adic $L$-function $L_5(\omega\xi, s)$. The $\mu$-invariant is zero (by [Fe-Wa]). In [Maz3], Mazur proves the following congruence formula

$$L_5(E^\xi/\mathbb{Q}, s) \equiv (1 - \xi(11)11^{1-s})L_5(\omega\xi, s-1)L_5(\omega\xi, 1-s) \pmod{5\mathbb{Z}_5}$$

for all $s \in \mathbb{Z}_5$. More precisely, one can interpret this as a congruence in the Iwasawa algebra $\Lambda$ modulo the ideal $5\Lambda$. The left side corresponds to $f_{E^\xi}^{\mathrm{anal}}(T)$, and each factor on the right side corresponds to an element of $\Lambda$. The two sides are congruent modulo $5\Lambda$. Now, if $f(T) \in \Lambda$ is any power series with $\mu(f) = 0$, then one has $f(T) \equiv uT^{\lambda(f)} \pmod{p\Lambda}$, where $u \in \Lambda^\times$. Applying this, we obtain $\lambda(f_{E^\xi}^{\mathrm{anal}}) = 2\lambda_\xi + \epsilon_\xi$ and therefore $\lambda(f_{E^\xi}^{\mathrm{anal}}) = \lambda_{E^\xi}$. Since both $f_{E^\xi}^{\mathrm{anal}}(T)$ and $f_E(T)$ have $\mu$-invariant equal to zero, it follows that indeed $(f_E(T)) = (f_{E^\xi}^{\mathrm{anal}}(T))$.

Theorems 4.8, 4.14, and 4.15 give sufficient conditions for the nonexistence of proper $\Lambda$-submodules of finite index in $\mathrm{Sel}_E(F_\infty)_p$. In particular, if $F = \mathbb{Q}$ and if $E$ has good, ordinary or multiplicative reduction at $p$, where $p$ is any odd prime, then no such $\Lambda$-submodule of $\mathrm{Sel}_E(\mathbb{Q}_\infty)_p$ can exist. (This is also true for $p = 2$, although the above results don't cover this case completely.) The following example shows that in general some restrictive hypotheses are needed. We let $F = \mathbb{Q}(\mu_5)$, $F_\infty = \mathbb{Q}(\mu_{5^\infty})$. Let $E = E_2$, the elliptic curve of conductor 11 with $E(\mathbb{Q}) = 0$. We shall show that $\mathrm{Sel}_E(F_\infty)_5$ has a $\Lambda$-submodule of index 5. To be more precise, note that $\mathrm{Gal}(F_\infty/\mathbb{Q}) = \Delta \times \Gamma$, where $\Delta = \mathrm{Gal}(F_\infty/\mathbb{Q}_\infty)$ and $\Gamma = \mathrm{Gal}(F_\infty/F)$. Now $\Delta$ has order 4 and its characters are $\omega^i$, $0 \le i \le 3$. We can decompose $\mathrm{Sel}_E(F_\infty)_5$ as a $\Lambda$-module by the action of $\Delta$:

$$\mathrm{Sel}_E(F_\infty)_5 = \bigoplus_{i=0}^{3} \mathrm{Sel}_E(F_\infty)_5^{\omega^i}.$$

As we will see, it turns out that $\mathrm{Sel}_E(F_\infty)_5^{\omega^3} \cong \mathbb{Z}/5\mathbb{Z}$, which of course is a $\Lambda$-module quotient of $\mathrm{Sel}_E(F_\infty)_5$. This component is $(\mathrm{Sel}_E(F_\infty)_5 \otimes \omega^{-3})^\Delta$, which can be identified with a subgroup of $H^1(\mathbb{Q}_\Sigma/\mathbb{Q}_\infty, E[5^\infty] \otimes \omega^{-3})$, where $\Sigma = \{\infty, 5, 11\}$. For brevity, we let $A = E[5^\infty] \otimes \omega^{-3}$. We let $S_A(\mathbb{Q}_\infty)$ denote the subgroup of $H^1(\mathbb{Q}_\Sigma/\mathbb{Q}_\infty, A)$ which is identified with $\mathrm{Sel}_E(F_\infty)_5^{\omega^3}$ by the restriction map. Noting that $\omega^{-3} = \omega$, we have a nonsplit exact sequence of $G_\mathbb{Q}$-modules

$$0 \to \mu_5^{\otimes 2} \to A[5] \to \mu_5 \to 0.$$

This is even nonsplit as a sequence of $G_{\mathbb{Q}_5}$-modules or $G_{\mathbb{Q}_{11}}$-modules. The $G_\mathbb{Q}$-submodule $\mu_5^{\otimes 2}$ of $A[5]$ is just $\Phi_2 \otimes \omega$, which we will denote simply by $\Phi$. We let $\Psi = A[5]/\Phi$. We will show that

$$S_A(\mathbb{Q}_\infty) \cong H^1(\mathbb{Q}_\Sigma/\mathbb{Q}_\infty, \Phi)$$

where the isomorphism is by the map $\epsilon : H^1(\mathbb{Q}_\Sigma/\mathbb{Q}_\infty, \Phi) \to H^1(\mathbb{Q}_\Sigma/\mathbb{Q}_\infty, A)$. This map is clearly injective. Since $\Phi \subseteq C_5 \otimes \omega$, it follows that the local condition defining $S_A(\mathbb{Q}_\infty)$ at the prime of $\mathbb{Q}_\infty$ lying over 5 is satisfied by the elements of $\mathrm{Im}(\epsilon)$. We now verify that the local condition at the prime of $\mathbb{Q}_\infty$ over 11 is also satisfied. This is because $\Phi \subseteq C_{11} \otimes \omega$, which is true because the above exact sequence is nonsplit over $\mathbb{Q}_{11}$. Since 11 splits completely in $F/\mathbb{Q}$, $\omega|_{G_{\mathbb{Q}_{11}}}$ is trivial. Thus $A = E[5^\infty]$ as $G_{\mathbb{Q}_{11}}$-modules. One then easily sees that the map

$$H^1(\mathbb{Q}_{11}^{\mathrm{cyc}}, C_{11}) \to H^1(\mathbb{Q}_{11}^{\mathrm{cyc}}, E[5^\infty])$$

is the zero map. That is, the map $H^1(\mathbb{Q}_{11}^{\mathrm{cyc}}, A) \to H^1(\mathbb{Q}_{11}^{\mathrm{cyc}}, D_{11})$ is an isomorphism. Elements of $\mathrm{Im}(\epsilon)$ are mapped to 0 and hence are trivial already in $H^1(\mathbb{Q}_{11}^{\mathrm{cyc}}, A)$, therefore satisfying the local condition defining $S_A(\mathbb{Q}_\infty)$ at the prime over 11.

So it is clear that $\mathrm{Im}(\epsilon) \subseteq S_A(\mathbb{Q}_\infty)$. We will prove that equality holds and that $H^1(\mathbb{Q}_\Sigma/\mathbb{Q}_\infty, \Phi) \cong \mathbb{Z}/5\mathbb{Z}$. This last assertion is rather easy to verify. Let $F^+ = \mathbb{Q}(\sqrt{5})$, the maximal real subfield of $F$. By class field theory, one finds that there is a unique cyclic extension $K/F^+$ of degree 5 such that $K/\mathbb{Q}$ is dihedral and $K \subseteq \mathbb{Q}_\Sigma$. Thus, $H^1(\mathbb{Q}_\Sigma/\mathbb{Q}, \Phi) = \mathrm{Hom}(\mathrm{Gal}(K/F^+), \Phi)$ has order 5. It follows that $H^1(\mathbb{Q}_\Sigma/\mathbb{Q}_\infty, \Phi)$ is nontrivial. Also, one can see that $K/F^+$ is ramified at the primes of $F^+$ lying over 5 and 11. Let $\Sigma' = \{\infty, 5\}$. Then $\Phi$ is a $\mathrm{Gal}(\mathbb{Q}_{\Sigma'}/\mathbb{Q})$-module and one can verify that $H^1(\mathbb{Q}_{\Sigma'}/\mathbb{Q}_\infty, \Phi) = 0$. (It is enough to show that $H^1(\mathbb{Q}_{\Sigma'}/\mathbb{Q}_\infty, \Phi)^\Gamma = H^1(\mathbb{Q}_{\Sigma'}/\mathbb{Q}, \Phi)$ vanishes. This is clear since $K/F^+$ is ramified at 11.) Therefore, the restriction map $H^1(\mathbb{Q}_\Sigma/\mathbb{Q}_\infty, \Phi) \to H^1(\mathbb{Q}_{11}^{\mathrm{cyc}}, \Phi)$ must be injective. But $H^1(\mathbb{Q}_{11}^{\mathrm{cyc}}, \Phi) \cong \mathbb{Z}/5\mathbb{Z}$, from which it follows that $H^1(\mathbb{Q}_\Sigma/\mathbb{Q}_\infty, \Phi)$ indeed has order 5.

It remains to show that $S_A(\mathbb{Q}_\infty) \subseteq \mathrm{Im}(\epsilon)$. Let $B = A/\Phi$. Then $B \cong E_1[5^\infty] \otimes \omega$ and $B[5] \cong \mu_5^{\otimes 2} \times \mu_5$ as $G_\mathbb{Q}$-modules. We will prove that $S_B(\mathbb{Q}_\infty) = 0$, from which it follows that $S_A(\mathbb{Q}_\infty) \subseteq \mathrm{Im}(\epsilon)$. Consider $S_B(\mathbb{Q}_\infty)[5]$, any element of which is represented by a 1-cocycle $\sigma$ with values in $B[5]$. The map $B[5] \to \mu_5$ sends $\sigma$ to a 1-cocycle $\widetilde{\sigma}$ such that $[\widetilde{\sigma}|_{G_5^{\mathrm{cyc}}}]$ is trivial as an element of $H^1(\mathbb{Q}_5^{\mathrm{cyc}}, D_5 \otimes \omega)$. Thus, $[\widetilde{\sigma}]$ is in the kernel of the composite map

$$H^1(\mathbb{Q}_\Sigma/\mathbb{Q}_\infty, \mu_5) \to H^1(\mathbb{Q}_5^{\mathrm{cyc}}, \mu_5) \to H^1(\mathbb{Q}_5^{\mathrm{cyc}}, D_5 \otimes \omega).$$

The second map is clearly injective. If the kernel of the first map were nontrivial, ,then it would have a nonzero intersection with $H^1(\mathbb{Q}_\Sigma/\mathbb{Q}_\infty, \mu_5)^\Gamma = H^1(\mathbb{Q}_\Sigma/\mathbb{Q}, \mu_5)$. One then sees that the map $a : H^1(\mathbb{Q}_\Sigma/\mathbb{Q}, \mu_5) \to H^1(\mathbb{Q}_5, \mu_5)$ would have a nonzero kernel. But, as we already used before, the map $a$ is injective. (The elements of $H^1(\mathbb{Q}_\Sigma/\mathbb{Q}, \mu_5)$ are represented by the 1-cocycles associated to $\sqrt[5]{5^i 11^j}$, $0 \le i, j \le 3$. But $5^i 11^j \in (\mathbb{Q}_5^\times)^5 \Leftrightarrow i = j = 0$.) Thus, the first map is injective too. Thus $[\widetilde{\sigma}] = 0$. Hence we may assume that $\sigma$ has values in $\mu_5^{\otimes 2}$. Now, in contrast to $A$, we have $\mu_5^{\otimes 2} \not\subseteq C_{11} \otimes \omega$. That is, the map $B \to D_{11}$ induces an isomorphism $\mu_5^{\otimes 2} \xrightarrow{\sim} D_{11}[5]$. The composite map

$$H^1(\mathbb{Q}_{11}^{\mathrm{cyc}}, \mu_5^{\otimes 2}) \to H^1(\mathbb{Q}_{11}^{\mathrm{cyc}}, B) \to H^1(\mathbb{Q}_{11}^{\mathrm{cyc}}, D_{11})$$

is clearly injective. Since $[\sigma]$ becomes trivial in $H^1(\mathbb{Q}_{11}^{\text{cyc}}, B)$, it follows that

$$[\sigma] \in \ker\left(H^1(\mathbb{Q}_\Sigma/\mathbb{Q}_\infty, \mu_5^{\otimes 2}) \to H^1(\mathbb{Q}_{11}^{\text{cyc}}, \mu_5^{\otimes 2})\right)$$

But we already showed that this kernel is trivial. (Recall that $\mu_5^{\otimes 2} \cong \Phi$.) Hence $[\sigma] = 0$, proving that $S_B(\mathbb{Q}_\infty) = 0$ as claimed.

**Conductor = 768**. We return now to the elliptic curves 768(D1, D3) which we denoted previously by $E_1$ and $E_2$. We take $p = 5$. As we mentioned earlier, $\text{Sel}_{E_1}(\mathbb{Q})_5 = 0$ and $\text{Sel}_{E_1}(\mathbb{Q}_\infty)_5 = 0$. Also, $E_1$ and $E_2$ are related by an isogeny of degree 5. Let $\Phi$ denote the $G_\mathbb{Q}$-invariant subgroup of $E_2[5^\infty]$ such that $E_2/\Phi \cong E_1$, $|\Phi| = 5$. Let $\Psi = E_2[5]/\Phi$. Then $G_\mathbb{Q}$ acts on $\Phi$ and $\Psi$ by characters $\varphi, \psi$ with values in $(\mathbb{Z}/5\mathbb{Z})^\times$ which factor through $\text{Gal}(\mathbb{Q}_\Sigma/\mathbb{Q})$, where now $\sigma = \{\infty, 2, 3, 5\}$. Since $E_2$ has good, ordinary reduction at 5, one of the characters $\varphi, \psi$ will be unramified at 5. Denote this character by $\theta$. By looking at the Fourier coefficients for the modular form associated to $E_2$ (which are given in [Cre]), one finds that $\theta$ is the even character of conductor 16 determined by $\theta(5) = 2 + 5\mathbb{Z}$. Then $\theta(3) = 3 + 5\mathbb{Z}$. Now $E_1$ and $E_2$ have split, multiplicative reduction at 3. One has a nonsplit exact sequence of $G_\mathbb{Q}$-modules

$$0 \to \Psi \to E_1[5] \to \Phi \to 0$$

which remains nonsplit for the action of $G_{\mathbb{Q}_3}$ since the Tate period for $E_1$ over $\mathbb{Q}_3$ has valuation not divisible by 5. Thus, $\Psi \cong \mu_5$ as $G_{\mathbb{Q}_3}$-modules. Thus, $\psi(3) = 3 + 5\mathbb{Z}$, $\varphi(3) = 1 + 5\mathbb{Z}$. Hence we have $\theta = \psi$. Therefore, $\Psi$ is even and unramified at 5, $\Phi$ is odd and ramified at 5. By theorem 5.7, we see that $\text{Sel}_{E_2}(\mathbb{Q}_\infty)_5$ has positive $\mu$-invariant. But since $\text{Sel}_{E_1}(\mathbb{Q}_\infty)_5 = 0$, it is clear that $\text{Sel}_{E_2}(\mathbb{Q}_\infty)_5 \subseteq H^1(\mathbb{Q}_\Sigma/\mathbb{Q}_\infty, \Phi)$. Thus, $\mu_{E_2} = 1$ and $\text{Sel}_{E_2}(\mathbb{Q}_\infty)_5$ has exponent 5. One then sees easily (using proposition 4.8 and the fact that $\Lambda/5\Lambda$ is a PID) that $\text{Sel}_{E_2}(\mathbb{Q}_\infty)_5 \cong \widehat{\Lambda}[5]$, as we stated earlier. Theorem 4.1 then implies that $\text{Sel}_{E_2}(\mathbb{Q})_5 = 0$.

**Conductor = 14**. Let $p = 3$. The situation is quite analogous to that for elliptic curves of conductor 11 and for $p = 5$. The $\mu$-invariants of $\text{Sel}_E(\mathbb{Q}_\infty)_3$ if $E$ has conductor 14 can be 0, 1, or 2. The $\lambda$-invariant is 0.

**Conductor = 34**. Let $p = 3$. There are four isogenous curves of conductor 34. We considered earlier the curve $E = 34(\text{A1})$, showing that $f_E(T) = \theta_1$, up to a factor in $\Lambda^\times$, where $\theta_1 = T^2 + 3T + 3$. The curve $34(\text{A2})$ is related to $E$ by a $\mathbb{Q}$-isogeny of degree 2 and so again has $\mu$-invariant 0 and $\lambda$-invariant equal to 2. The two other curves of conductor 34 have $\mathbb{Q}$-isogenies of degree 3 with kernel isomorphic to $\mu_3$ as a $G_\mathbb{Q}$-module. It then follows that they have $\mu$-invariant 1. Denoting either of them by $E'$, the characteristic ideal of the Pontryagin dual of $\text{Sel}_{E'}(\mathbb{Q}_\infty)_3$ is generated by $3\theta_1$.

**Conductor = 306**. Take $p = 3$. We will consider just the elliptic curve $E$ defined by $y^2 + xy = x^3 - x^2 - 927x + 11097$. This is 306(B3) in [Cre].

It is the quadratic twist of 34(A3) by the character $\omega$ of conductor 3. The Mordell-Weil group $E(\mathbb{Q})$ is of rank 1, isomorphic to $\mathbb{Z} \times (\mathbb{Z}/6\mathbb{Z})$. $E$ has potentially ordinary reduction at 3, and has good ordinary reduction over $K = \mathbb{Q}(\mu_3)$ at the prime $\mathfrak{p}$ lying over 3. The unique subgroup $\Phi$ of $E(\mathbb{Q})$ of order 3 is contained in the kernel of reduction modulo $\mathfrak{p}$ for $E(K)$. Although the hypotheses of proposition 5.10 are not satisfied by $E$, the proof can still be followed to show that the $\mu$-invariant of $\mathrm{Sel}_E(\mathbb{Q}_\infty)_3$ is 0. Let $F$ denote the first layer of $\mathbb{Q}_\infty$, $F = \mathbb{Q}(\beta)$ where $\beta = \zeta + \zeta^{-1}$, $\zeta$ being a primitive 9-th root of unity. The prime 17 splits completely in $F/\mathbb{Q}$. Using this fact, it is easy to verify directly that

$$\mathrm{Hom}(\mathrm{Gal}(F/\mathbb{Q}), \Phi) \subseteq \mathrm{Sel}_E(\mathbb{Q})_3.$$

This clearly implies that $\ker(\mathrm{Sel}_E(\mathbb{Q})_3 \to \mathrm{Sel}_E(\mathbb{Q}_\infty)_3)$ is nontrivial. (Contrast this with proposition 3.9.) However, we can explain this in the following more concrete way, using the results of a calculation carried out by Karl Rubin. The point $P = (9, 54)$ on $E(\mathbb{Q})$ is a generator of $E(\mathbb{Q})/E(\mathbb{Q})_{\mathrm{tors}}$. But $P = 3Q$, where $Q = (-6\beta^2 + 9\beta + 15, 15\beta^2 - 48\beta + 9)$ is in $E(F)$. This implies that the map $E(\mathbb{Q}) \otimes (\mathbb{Q}_3/\mathbb{Z}_3) \to E(F) \otimes (\mathbb{Q}_3/\mathbb{Z}_3)$ has a nontrivial kernel. Let $\phi$ be the 1-cocycle defined by $\phi(g) = g(Q) - Q$ for $g \in G_\mathbb{Q}$. This cocycle has values in $E(F)[3]$, which is easily seen to be just $\Phi = E(\mathbb{Q})[3]$, and factors through $\mathrm{Gal}(F/\mathbb{Q})$. Thus it generates $\mathrm{Hom}(\mathrm{Gal}(F/\mathbb{Q}), \Phi)$ and is certainly contained in $\mathrm{Sel}_E(\mathbb{Q})_3$.

For $n \geq 1$, it turns out that $\ker(\mathrm{Sel}_E(\mathbb{Q}_n)_3 \to \mathrm{Sel}_E(\mathbb{Q}_\infty)_3) = 0$. This can be seen by checking that the local condition at any prime of $\mathbb{Q}_n$ lying above 17 (which will be inert in $\mathbb{Q}_\infty/\mathbb{Q}_n$) fails to be satisfied by a nontrivial element of $\mathrm{Hom}(\mathrm{Gal}(\mathbb{Q}_\infty/\mathbb{Q}_n), \Phi)$. The fact that $E$ has split, multiplicative reduction at 17 helps here. The argument given in [HaMa] shows that $\mathrm{Sel}_E(\mathbb{Q}_\infty)_3$ has no proper $\Lambda$-submodule of finite index. As remarked above, the $\mu$-invariant is 0. A calculation of McCabe for the $p$-adic $L$-function associated to $E$ combined with Kato's theorem implies that the $\lambda$-invariant of $\mathrm{Sel}_E(\mathbb{Q}_\infty)_3$ is 1. It follows that $\mathrm{Sel}_E(\mathbb{Q}_\infty)_3 = E(\mathbb{Q}_\infty) \otimes (\mathbb{Q}_3/\mathbb{Z}_3) \cong \mathbb{Q}_3/\mathbb{Z}_3$, on which $\Gamma$ acts trivially.

**Conductor $= 26$.** Consider 26(B1, B2). These curves are related by isogenies with kernels isomorphic to $\mu_7$ and $\mathbb{Z}/7\mathbb{Z}$. Let $E_1$ be 26(B1). Then $\mathrm{Sel}_{E_1}(\mathbb{Q})_7$ should be zero. From [Cre], we have $c_2 = 7$, $c_{13} = 1$, $a_7 = 1$, and $|E_1(\mathbb{Q})| = 7$. Take $p = 7$. Theorem 4.1 then implies that $f_{E_1}(0) \sim 7$. Thus, $f_{E_1}(T)$ is an irreducible element of $\Lambda$. The only nonzero, proper, $G_\mathbb{Q}$-invariant subgroup of $E_1[7]$ is $E_1(\mathbb{Q}) \cong \mathbb{Z}/7\mathbb{Z}$. Although we haven't verified it, it seems likely that $\mu_{E_1} = 0$. (Conjecture 1.11 would predict this.) If this is so, then $\lambda_{E_1} > 0$. Let $E_2$ be 26(B2). Then $c_2 = c_{13} = 1$, $a_7 = 1$, and $E_2(\mathbb{Q}) = 0$. One can verify that $\mathrm{Sel}_{E_2}(\mathbb{Q})_7 = 0$. Then by Theorem 4.1, we have $f_{E_2}(0) \sim 7^2$. Since $(f_{E_1}(T))$ and $(f_{E_2}(T))$ can differ only by multiplication by a power of 7, it is clear that $f_{E_2}(T) = 7f_{E_1}(T)$, up to a factor in $\Lambda^\times$. Thus, $\mu_{E_2} \geq 1$, which also follows from proposition 5.7 because $E_2[7]$ contains the odd, ramified $G_\mathbb{Q}$-submodule $\mu_7$.

**Conductor = 147**. Consider 147(B1, B2), which we denote by $E_1$ and $E_2$, respectively. They are related by isogenies of degree 13. For $E_1$, one has $c_3 = c_7 = 1$, $a_{13} = 1$, $E_1(\mathbb{Q}) = 0$, and $\mathrm{Sel}_{E_1}(\mathbb{Q}) = 0$. Take $p = 13$. By theorem 4.1, $f_{E_1}(0) \sim 13^2$. For $E_2$, one has $c_3 = 13$, $c_7 = 1$, $a_{13} = 1$, $E_2(\mathbb{Q}) = 0$, and $\mathrm{Sel}_{E_2}(\mathbb{Q}) = 0$. Thus, $f_{E_2}(0) \sim 13^3$. Since an isogeny $E_1 \to E_2$ of degree 13 induces a homomorphism $\mathrm{Sel}_{E_1}(\mathbb{Q}_\infty)_{13} \to \mathrm{Sel}_{E_2}(\mathbb{Q}_\infty)_{13}$ with kernel and cokernel of exponent 13, it is clear that $f_{E_2}(T) = 13 f_{E_1}(T)$, up to a factor in $\Lambda^\times$. Conjecturally, $\mu_{E_1} = 0$ and hence $\mu_{E_2} = 1$. Let $\xi$ be the quadratic character of conductor 7, which is odd. Then $E_1^\xi$ and $E_2^\xi$ are the curves 147(C1, C2). Proposition 5.10 implies that $\mathrm{Sel}_{E_1^\xi}(\mathbb{Q}_\infty)_{13}$ and $\mathrm{Sel}_{E_2^\xi}(\mathbb{Q}_\infty)_{13}$ have $\mu$-invariant equal to zero. In fact, for both $E_1^\xi$ and $E_2^\xi$, we in fact have $c_3 = c_7 = 1$, $a_{13} = -1$, $E_1^\xi(\mathbb{Q}) = E_2^\xi(\mathbb{Q}) = 0$, $\mathrm{Sel}_{E_1^\xi}(\mathbb{Q})_{13} = \mathrm{Sel}_{E_2^\xi}(\mathbb{Q})_{13} = 0$. By proposition 3.8, we have $\mathrm{Sel}_{E_1^\xi}(\mathbb{Q}_\infty)_{13} = \mathrm{Sel}_{E_2^\xi}(\mathbb{Q}_\infty)_{13} = 0$.

**Conductor = 1225**. Consider now $E_1 : y^2 + xy + y = x^3 + x^2 - 8x + 6$ and also $E_2 : y^2 + xy + y = x^3 + x^2 - 208083x - 36621194$. These curves have conductor 1225 and are related by a $\mathbb{Q}$-isogeny of degree 37. They have additive reduction at 5 and 7. Hence the Tamagawa factors are at most 4. The $j$-invariants are in $\mathbb{Z}$ and so these curves have potentially good reduction at 5 and 7. We take $p = 37$. Since $a_{37} = 8$, $E_1$ and $E_2$ have good, ordinary reduction at $p$. Let $\Phi$ be the $G_\mathbb{Q}$-invariant subgroup of $E_2[37]$ and let $\Psi = E_2[37]/\Phi$. Thus, $\Phi$ is the kernel of the isogeny from $E_2$ to $E_1$. The real periods $\Omega_1$, $\Omega_2$ of $E_1$, $E_2$ are given by: $\Omega_1 = 4.1353\ldots$, $\Omega_2 = .11176\ldots$. Since $\Omega_1 = 37\Omega_2$, one finds that $\Phi$ must be odd. Let $\varphi, \psi$ be the $(\mathbb{Z}/37\mathbb{Z})^\times$-valued characters which describe the action of $G_\mathbb{Q}$ on $\Phi$ and $\Psi$. We can regard them as Dirichlet characters. They have conductor dividing $5 \cdot 7 \cdot 37$ and one of them (which we denote by $\theta$) is unramified at 37. By examining the Fourier coefficients of the corresponding modular form, one finds that $\theta$ is characterized by $\theta(2) = 8 + 37\mathbb{Z}$, $\theta(13) = 6 + 37\mathbb{Z}$. The character $\theta$ is even and has order 12 and conductor 35. But since $\varphi$ is odd, we must have $\theta = \psi$. Thus, $\Phi$ is odd and ramified at 37. Therefore, by proposition 5.7, we have $\mu_{E_2} \geq 1$. By using the result given in [Pe2] or [Sch3], one finds that $\mu_{E_2} = \mu_{E_1} + 1$. Conjecturally, $\mu_{E_2} = 1$, $\mu_{E_1} = 0$. In any case, we have $(f_{E_2}(T)) = (37 f_{E_1}(T))$. Now $E_1(\mathbb{Q})$ and $E_2(\mathbb{Q})$ have rank 1. It is interesting to note that the fact that $\mathrm{Sel}_{E_1}(\mathbb{Q})_{37}$, $\mathrm{Sel}_{E_2}(\mathbb{Q})_{37}$ are infinite can be deduced from Theorem 4.1. For if one of these Selmer groups were finite, then so would the other. One would then see that both $f_{E_1}(0)$ and $f_{E_2}(0)$ would have even valuation. This follows from Cassels' theorem that $|\mathrm{Sel}_{E_i}(\mathbb{Q})|$ is a perfect square for $i = 1, 2$ together with the fact that the Tamagawa factors for $E_i$ at 5 and 7 cannot be divisible by 37. But $f_{E_2}(0) \sim 37 f_{3_1}(0)$, which gives a contradiction. Similar remarks apply to even quadratic twists of $E_1$ and $E_2$.

Now we will state and prove the analogues of propositions 5.7 and 5.10 for $p = 2$. It is necessary to define the terms "ramified" and "odd" somewhat more carefully. Assume that $E$ is an elliptic curve/$\mathbb{Q}$ with good, ordinary or multiplicative reduction at 2. Suppose that $\Phi$ is a cyclic $G_{\mathbb{Q}}$-invariant subgroup of $E[2^{\infty}]$. We say that $\Phi$ is "ramified at 2" if $\Phi \subseteq C_2$, where $C_2$ is the subgroup of $E[2^{\infty}]$ which occurs in the description of the image of the local Kummer map for $E$ over $\mathbb{Q}_2$ given in section 2. (It is characterized by $C_2 \cong \mu_{2^{\infty}}$ for the action of $I_{\mathbb{Q}_2}$. Here $I_{\mathbb{Q}_2}$ is the inertia subgroup of $G_{\mathbb{Q}_2}$, identified with a subgroup of $G_{\mathbb{Q}}$ by choosing a prime of $\overline{\mathbb{Q}}$ lying over 2. Then $D_2 = E[2^{\infty}]/C_2$ is an unramified $G_{\mathbb{Q}_2}$-module.) We say that $\Phi$ is "odd" if $\Phi \subseteq C_{\infty}$, where $C_{\infty}$ denotes the maximal divisible subgroup of $E[2^{\infty}]$ on which $\mathrm{Gal}(\mathbb{C}/\mathbb{R})$ acts by $-1$: $C_{\infty} = (E[2^{\infty}]^-)_{\mathrm{div}}$. Then $C_{\infty} \cong \mathbb{Q}_2/\mathbb{Z}_2$ as a group. Here we identify $\mathrm{Gal}(\mathbb{C}/\mathbb{R})$ with a subgroup of $G_{\mathbb{Q}}$ by choosing an infinite prime of $\overline{\mathbb{Q}}$. (We remark that $C_{\infty} \cong \mu_{2^{\infty}}$ as $\mathrm{Gal}(\mathbb{C}/\mathbb{R})$-modules and that $\mathrm{Gal}(\mathbb{C}/\mathbb{R})$ acts trivially on $D_{\infty} = E[2^{\infty}]/C_{\infty}$.) Since $\Phi$ is $G_{\mathbb{Q}}$-invariant, these definitions are easily seen to be independent of the choice of primes of $\overline{\mathbb{Q}}$ lying over 2 and over $\infty$. We now prove the analogue of proposition 5.7.

**Proposition 5.13.** *Suppose that $E$ is an elliptic curve/$\mathbb{Q}$ with good, ordinary or multiplicative reduction at 2. Suppose also that $E[2^{\infty}]$ contains a $G_{\mathbb{Q}}$-invariant subgroup $\Phi$ of order $2^m$ which is ramified at 2 and odd. Then the $\mu$-invariant of $\mathrm{Sel}_E(\mathbb{Q}_{\infty})_2$ is at least $m$.*

*Proof.* The argument is virtually the same as that for proposition 5.7. We consider $\mathrm{Im}(\epsilon)$ where $\epsilon$ is the map

$$\epsilon : H^1(\mathbb{Q}_{\Sigma}/\mathbb{Q}_{\infty}, \Phi) \to H^1(\mathbb{Q}_{\Sigma}/\mathbb{Q}_{\infty}, E[2^{\infty}]).$$

The kernel is finite. Since $\Phi \subseteq C_2$, the elements of $\mathrm{Im}(\epsilon)$ satisfy the local conditions defining $\mathrm{Sel}_E(\mathbb{Q}_{\infty})_2$ at the prime of $\mathbb{Q}_{\infty}$ lying over 2. Also, just as previously, a subgroup of finite index in $\mathrm{Im}(\epsilon)$ satisfies the local conditions for all other nonarchimedean primes of $\mathbb{Q}_{\infty}$. Now we consider the archimedean primes of $\mathbb{Q}_{\infty}$. Note that $H^1(\mathbb{R}, C_{\infty}) = 0$. Since $\Phi \subseteq C_{\infty}$, it is clear that elements in $\mathrm{Im}(\epsilon)$ are locally trivial in $H^1((\mathbb{Q}_{\infty})_{\eta}, E[2^{\infty}])$ for every infinite prime $\eta$ of $\mathbb{Q}_{\infty}$. Therefore, $\mathrm{Im}(\epsilon) \cap \mathrm{Sel}_E(\mathbb{Q}_{\infty})_2$ has finite index in $\mathrm{Im}(\epsilon)$.

It remains to show that the $\Lambda$-module $H^1(\mathbb{Q}_{\Sigma}/\mathbb{Q}_{\infty}, \Phi)$ has $\mu$-invariant equal to $m$. Since the $G_{\mathbb{Q}}$-composition factors for $\Phi$ are isomorphic to $\mathbb{Z}/2\mathbb{Z}$, lemma 5.9 implies that the $\mu$-invariant for $H^1(\mathbb{Q}_{\Sigma}/\mathbb{Q}_{\infty}, \Phi)$ is at most $m$. On the other hand, the Euler characteristic of the $\mathrm{Gal}(\mathbb{Q}_{\Sigma}/\mathbb{Q}_n)$-module $\Phi$ is $\prod_{v|\infty} |\Phi/\Phi^{D_v}|^{-1}$, where $v$ runs over the infinite primes of $\mathbb{Q}_n$ and $D_v = \mathrm{Gal}(\mathbb{C}/\mathbb{R})$ is a corresponding decomposition group. Assume that $m \geq 1$. Then $|\Phi^{D_v}| = 2$ and so this Euler characteristic is $2^{-(m-1)2^n}$ for all $n \geq 0$. Now $H^0(\mathbb{Q}_{\Sigma}/\mathbb{Q}_n, \Phi)$ just has order 2. As for $H^2(\mathbb{Q}_{\Sigma}/\mathbb{Q}_n, \Phi)$, it is known that the map

$$H^2(\mathbb{Q}_{\Sigma}/\mathbb{Q}_n, \Phi) \to \prod_{v|\infty} H^2((\mathbb{Q}_n)_v, \Phi)$$

is surjective. (This is corollary 4.16 in [Mi].) Since $H^2(D_v, \Phi)$ has order 2, it follows that $|H^2(\mathbb{Q}_\Sigma/\mathbb{Q}_n, \Phi)| \geq 2^{2^n}$. Therefore,

$$H^1(\mathbb{Q}_\Sigma/\mathbb{Q}_n, \Phi)| \geq 2^{m2^n+1}.$$

The restriction map $H^1(\mathbb{Q}_\Sigma/\mathbb{Q}_n, \Phi) \to H^1(\mathbb{Q}_\Sigma/\mathbb{Q}_\infty, \Phi)^{\Gamma_n}$ is surjective and has kernel $H^1(\Gamma_n, \mathbb{Z}/2\mathbb{Z})$, which has order 2. Thus,

$$|H^1(\mathbb{Q}_\Sigma/\mathbb{Q}_n, \Phi)^{\Gamma_n}| \geq 2^{m2^n}$$

for all $n$. This implies that $H^1(\mathbb{Q}_\Sigma/\mathbb{Q}_\infty, \Phi)$ has $\mu$-invariant at least $m$. Therefore, the $\mu$-invariant of $H^1(\mathbb{Q}_\Sigma/\mathbb{Q}_\infty, \Phi)$ and hence of $\mathrm{Im}(\epsilon)$ is exactly $m$, proving proposition 5.13.

*Remark.* As we mentioned before (for any $p$), if $E$ admits a $\mathbb{Q}$-isogeny of degree 2 and if $\mathrm{Sel}_E(\mathbb{Q}_\infty)_2$ is $\Lambda$-cotorsion, then $\mathrm{Sel}_E(\mathbb{Q}_\infty)_2$ contains a $\Lambda$-submodule pseudo-isomorphic to $\widehat{\Lambda}[2^{\mu_E}]$. It is known that there are infinitely many elliptic curves/$\mathbb{Q}$ admitting a cyclic $\mathbb{Q}$-isogeny of degree 16, but none with such an isogeny of degree 32. We will give examples below where the assumptions in proposition 5.13 are satisfied and $|\Phi| = 2^m$ with $m = 0, 1, 2, 3$, or 4. For any elliptic curve $E/\mathbb{Q}$, there is a maximal $G_\mathbb{Q}$-invariant subgroup $\Phi$ which is ramified and odd. Define $m_E$ by $|\Phi| = 2^{m_E}$. Conjecturally, $\mu_E = m_E$. Thus the possible values of $\mu_E$ as $E$ varies over elliptic curves/$\mathbb{Q}$ with good, ordinary or multiplicative reduction at 2 should be 0, 1, 2, 3, or 4. Examples where $\mu_E > 0$ are abundant. It suffices to have a point $P \in E(\mathbb{Q})$ of order 2 such that $P \in C_2$ and $P \in C_\infty$, using the notation introduced earlier. If the discriminant of a Weierstrass equation for $E$ is negative, then $E(\mathbb{R})$ has just one component. In this case, $C_\infty[2] = E(\mathbb{R})[2]$ and so if $P \in E(\mathbb{Q})$ has order 2, then $\Phi = \langle P \rangle$ is automatically odd. (Note that in this case $H^1(\mathbb{R}, E[2^\infty]) = 0$ and so the local conditions at the infinite primes of $\mathbb{Q}_\infty$ occurring in the definition of $\mathrm{Sel}_E(\mathbb{Q}_\infty)_2$ are trivially satisfied anyway.) Similarly, if this discriminant is not a square in $\mathbb{Q}_2^\times$, then $\Phi = \langle P \rangle$ is automatically ramified since then $C_2[2] = E(\mathbb{Q}_2)[2]$.

We now prove the analogue of proposition 5.10, which gives a sufficient condition for $\mu_E = 0$ in case $p = 2$.

**Proposition 5.14.** *Suppose that $E$ is an elliptic curve/$\mathbb{Q}$ with good, ordinary or multiplicative reduction at 2. Suppose also that $E(\mathbb{Q})$ contains an element $P$ of order 2 and that $\Phi = \langle P \rangle$ is either ramified at 2 but not odd or odd but not ramified at 2. Then $\mathrm{Sel}_E(\mathbb{Q}_\infty)_2$ is $\Lambda$-cotorsion and $\mu_E = 0$.*

*Proof.* We must show that $\mathrm{Sel}_E(\mathbb{Q}_\infty)_2[2]$ is finite. Consider the map

$$\alpha_E : H^1(\mathbb{Q}_\Sigma/\mathbb{Q}_\infty, E[2^\infty]) \to \mathcal{P}_{E[2^\infty]}^{(\infty)}(\mathbb{Q}_\infty)$$

which occurred in the proof of proposition 5.8. By definition we have

$$\mathrm{Sel}_E(\mathbb{Q}_\infty)_2 \subseteq \ker(\alpha_E).$$

Under the hypothesis that $E$ admits a $\mathbb{Q}$-rational isogeny of degree 2 (i.e., that $E(\mathbb{Q})$ has an element of order 2), we showed earlier that $\ker(\alpha_E)[2]$ has $(\Lambda/2\Lambda)$-corank equal to 1. Consider the map

$$\epsilon : H^2(\mathbb{Q}_\Sigma/\mathbb{Q}_\infty, \Phi) \to H^1(\mathbb{Q}_\Sigma/\mathbb{Q}_\infty, E[2^\infty]).$$

Then $\ker(\epsilon)$ is finite and so, by lemma 5.9, $\mathrm{Im}(\epsilon)$ also has $(\Lambda/2\Lambda)$-corank equal to 1.

Assume first that $\Phi$ is odd but not ramified at 2. Then $\Phi \subseteq C_\infty$. Since we have $H^1(\mathbb{R}, C_\infty) = 0$, it is clear that $\mathrm{Im}(\epsilon) \subseteq \ker(\alpha_E)$. It follows that $\mathrm{Im}(\epsilon)$ has finite index in $\ker(\alpha_E)[2]$. Thus, it suffices to prove that $\mathrm{Im}(\epsilon) \cap \mathrm{Sel}_E(\mathbb{Q}_\infty)_2$ is finite. To do this, consider the composite map $\beta$ defined by the commutative diagram

$$H^1(\mathbb{Q}_\Sigma/\mathbb{Q}_\infty, \Phi) \xrightarrow{\ \epsilon\ } H^1(\mathbb{Q}_\Sigma/\mathbb{Q}_\infty, E[2^\infty]) \longrightarrow H^1(I_\pi, E[2^\infty])$$

$$\searrow{\scriptstyle \beta} \qquad\qquad \downarrow$$

$$H^1(I_\pi, D_2)$$

where $\pi$ is the unique prime of $\mathbb{Q}_\infty$ lying above 2 and $I_\pi$ is the inertia subgroup of $G_{(\mathbb{Q}_\infty)_\pi}$. Let $B = \ker(\beta)$. If $[\sigma] \in H^1(\mathbb{Q}_\Sigma/\mathbb{Q}_\infty, \Phi)$, then the local condition defining $\mathrm{Sel}_E(\mathbb{Q}_\infty)_2$ at the prime $\pi$ is satisfied by $\epsilon([\sigma])$ precisely when $[\sigma] \in B$. Since $\Phi \nsubseteq C_2$, the map $E[2^\infty] \to D_2$ induces an isomorphism of $\Phi$ to $D_2[2]$. Also, since $I_\pi$ acts trivially on $D_2$, the map $H^1(I_\pi, D_2[2]) \to H^1(I_\pi, D_2)$ is injective. Hence

$$B = \ker(H^1(\mathbb{Q}_\Sigma/\mathbb{Q}_\infty, \Phi) \to H^1(I_\pi, \Phi)).$$

If we let $H^1_{\mathrm{unr}}(\mathbb{Q}_\infty, \Phi)$ denote the subgroup of $H^1(\mathbb{Q}_\infty, \Phi)$ consisting of elements which are unramified at all *nonarchimedean* primes of $\mathbb{Q}_\infty$, then $H^1_{\mathrm{unr}}(\mathbb{Q}_\infty, \Phi)$ is a subgroup of $B$ and the index is easily seen to be finite. (Only finitely many nonarchimedean primes $\eta$ of $\mathbb{Q}_\infty$ exist lying over primes in $\Sigma$. $H^1((\mathbb{Q}_\infty)_\eta, \Phi)$ is finite if $\eta \nmid 2$.) Now $G_{\mathbb{Q}_\infty}$ acts trivially on $\Phi$. Let $L^*_\infty$ denote the maximal abelian pro-2 extension of $\mathbb{Q}_\infty$ which is unramified at all nonarchimedean primes of $\mathbb{Q}_\infty$. Then

$$H^1_{\mathrm{unr}}(\mathbb{Q}_\infty, \Phi) = \mathrm{Hom}(\mathrm{Gal}(L^*_\infty/\mathbb{Q}_\infty), \Phi).$$

But it is easy to verify that $L^*_\infty = \mathbb{Q}_\infty$. (For example, one can note that $L^*_\infty \mathbb{Q}_\infty(i)/\mathbb{Q}_\infty(i)$ is *everywhere* unramified. But $\mathbb{Q}_\infty(i) = \mathbb{Q}(\mu_{2^\infty})$. It is known that $\mathbb{Q}(\mu_{2^n})$ has odd class number for all $n \geq 0$. Thus, $\mathbb{Q}_\infty \subseteq L^*_\infty \subseteq \mathbb{Q}_\infty(i)$, from which $L^*_\infty = \mathbb{Q}_\infty$ follows.) Therefore $B$ is finite. Therefore $\mathrm{Im}(\epsilon) \cap \mathrm{Sel}_E(\mathbb{Q}_\infty)_2$ is indeed finite.

Now assume that $\Phi$ is ramified at 2 but not odd. Let $\epsilon$ be as above. Since $\Phi$ is not odd, it follows that $E(\mathbb{R})$ must have two connected components. Hence, by proposition 5.8, $H^1(\mathbb{Q}_\Sigma/\mathbb{Q}_\infty, E[2^\infty])[2]$ has $(\Lambda/2\Lambda)$-corank equal to 2. This implies that the $(\Lambda/2\Lambda)$-corank of $H^1(\mathbb{Q}_\Sigma/\mathbb{Q}_\infty, E[2])$ is 2. On the other hand, $H^1((\mathbb{Q}_\infty)_\pi, E[2])$ also has $(\Lambda/2\Lambda)$-corank equal to 2. Consider the map

$$a : H^1(\mathbb{Q}_\Sigma/\mathbb{Q}_\infty, E[2]) \to H^1((\mathbb{Q}_\infty)_\pi, E[2]).$$

We will show that the kernel is finite. It follows from this that the cokernel is also finite. We have an exact sequence

$$0 \to \Phi \to E[2] \to \Psi \to 0$$

of $G_{\mathbb{Q}_\infty}$-modules, where $\Phi \cong \Psi \cong \mathbb{Z}/2\mathbb{Z}$. The finiteness of the group $B$ introduced earlier in this proof, and the corresponding fact for $\Psi$, implies rather easily that $\ker(a)$ is indeed finite. Consider the map

$$b : H^1((\mathbb{Q}_\infty)_\pi, E[2]) \to H^1((\mathbb{Q}_\infty)_\pi, D_2[2])$$

induced by the map $E[2^\infty] \to D_2$. Since $(\mathbb{Q}_\infty)_\pi$ has 2-cohomological dimension 0, $b$ is surjective. It follows that $\operatorname{coker}(b \circ a)$ is finite. Using the fact that $H^1((\mathbb{Q}_\infty)_\pi, D_2[2])$ has $(\Lambda/2\Lambda)$-corank 1, we see that $\ker(b \circ a)$ also has $(\Lambda/2\Lambda)$-corank 1. Consider the map

$$\gamma_E : H^1(\mathbb{Q}_\Sigma/\mathbb{Q}_\infty, E[2^\infty]) \to H^1((\mathbb{Q}_\infty)_\pi, D_2).$$

The above remarks imply easily that $\ker(\gamma_E)[2]$ has $(\Lambda/2\Lambda)$-corank equal to 1.

The rest of the argument is now rather similar to that for the first case. It is clear that $\operatorname{Im}(\epsilon) \subseteq \ker(\gamma_E)$ since $\Phi \subseteq C_2$. Thus, $\operatorname{Im}(\epsilon)$ has finite index in $\ker(\gamma_E)[2]$. Also, by definition, we have

$$\operatorname{Sel}_E(\mathbb{Q}_\infty)_2 \subseteq \ker(\gamma_E).$$

It then suffices to show that $\operatorname{Im}(\epsilon) \cap \operatorname{Sel}_E(\mathbb{Q}_\infty)_2$ is finite. To do this, we consider the composite map $\delta_\eta$ defined by the following commutative diagram.

$$
\begin{array}{ccc}
H^1(\mathbb{Q}_\Sigma/\mathbb{Q}_\infty, \Phi) \xrightarrow{\ \epsilon\ } & H^1(\mathbb{Q}_\Sigma/\mathbb{Q}_\infty, E[2^\infty]) \longrightarrow & H^1((\mathbb{Q}_\infty), E[2^\infty]) \\
& \searrow{\scriptstyle \delta_\eta} & \downarrow \\
& & H^1((\mathbb{Q}_\infty)_\eta, D_\infty)
\end{array}
$$

where $\eta$ is any infinite prime of $\mathbb{Q}_\infty$. If $[\sigma] \in H^1(\mathbb{Q}_\Sigma/\mathbb{Q}_\infty, \Phi)$, then the local condition defining $\operatorname{Sel}_E(\mathbb{Q}_\infty)_2$ at $\eta$ for the element $\epsilon([\sigma])$ would imply that $\delta_\eta([\sigma]) = 0$. But since $\Phi \not\subseteq C_\infty$, $\Phi$ is identified with $D_\infty[2]$. The map $H^1(\mathbb{R}, D_\infty[2]) \to H^1(\mathbb{R}, D_\infty)$ is injective since $\operatorname{Gal}(\mathbb{C}/\mathbb{R})$ acts trivially on $D_\infty$. Hence

$$\ker(\delta_\eta) = \ker(H^1(\mathbb{Q}_\Sigma/\mathbb{Q}_\infty, \Phi) \to H^1((\mathbb{Q}_\infty)_\eta, \Phi)).$$

By lemma 5.9, we know that $\bigcap_\eta \ker(\delta_\eta)$ is finite (where $\eta$ varies over all the infinite primes of $\mathbb{Q}_\infty$. It follows from this that $\mathrm{Im}(\epsilon) \cap \mathrm{Sel}_E(\mathbb{Q}_\infty)_2$ is also finite. This implies that $\mathrm{Sel}_E(\mathbb{Q}_\infty)_2[2]$ is finite, finishing the proof of proposition 5.14. ∎

We now consider various examples.

**Conductor = 15**. There are eight curves of conductor 15, all related by $\mathbb{Q}$-isogenies whose degrees are powers of 2. We will let $E_i$ denote the curve labeled $A_i$ in [Cre] for $1 \le i \le 8$. The following table summarizes the situation for $p = 2$.

|               | $E_1$ | $E_2$ | $E_3$ | $E_4$ | $E_5$ | $E_6$ | $E_7$ | $E_8$ |
|---------------|-------|-------|-------|-------|-------|-------|-------|-------|
| $|\text{III}|$ = | 1 | 1 | 1 | 1 | 1 | 1 | 1 | 1 |
| $|T|$ = | 8 | 4 | 8 | 8 | 2 | 2 | 4 | 4 |
| $c_3, c_5$ = | 2,4 | 2,2 | 2,2 | 2,8 | 2,1 | 2,1 | 1,1 | 1,1 |
| $f_{E_i}(0) \sim$ | 2 | 4 | 1 | 4 | 8 | 8 | 1 | 1 |
| $\mu_{E_i}$ = | 1 | 2 | 0 | 2 | 3 | 3 | 0 | 0 |

For conductor 15, the Selmer group $\mathrm{Sel}_{E_i}(\mathbb{Q}_\infty)_2$ has $\lambda$-invariant equal to 0 and the $\mu$-invariant varies from 0 to 3. Now $\text{III} = \mathrm{Sel}_E(\mathbb{Q})$. Its order was computed under the assumption of the Birch and Swinnerton-Dyer conjecture by evaluating $L(E_i/\mathbb{Q}, 1)/\Omega_{E_i}$. The real period $\Omega_{E_i}$ was computed using PARI. $|T|$, $c_3$, and $c_5$ are as listed in [Cre]. Using the fact that $a_2 = -1$, and hence $|\tilde{E}_i(\mathbb{F}_2)| = 4$ for each $i$, the fourth row is a consequence of theorem 4.1. In particular, it is clear that $f_{E_3}(T) \in \Lambda^\times$. Hence $\mu_{E_3} = \lambda_{E_3} = 0$. The $\lambda$-invariant of $\mathrm{Sel}_{E_i}(\mathbb{Q}_\infty)_2$ is unchanged by a $\mathbb{Q}$-isogeny. Hence $\lambda_{E_i} = 0$ for all $i$. It is then obvious that $f_{E_i}(0) \sim 2^{\mu_{E_i}}$, which gives the final row.

It is not difficult to reconcile these results with propositions 5.13 and 5.14. For example, consider $E_3 : y^2 + xy + y = x^3 + x^2 - 5x + 2$. We have $E_3[2] \cong (\mathbb{Z}/2\mathbb{Z})^2$. The points of order 2 are $(1, -1)$, $(\frac{3}{4}, -\frac{7}{8})$, and $(-3, 1)$. The second point generates $C_2[2]$; the third point generates $C_\infty[2]$. (Remark: It is not hard to find the generator of $C_\infty[2]$. It is the point in $E(\mathbb{R})[2]$ whose $x$-coordinate is minimal.) Thus, $E_3(\mathbb{Q})[2]$ contains a subgroup which is ramified at 2 but not odd and another subgroup which is odd but not ramified at 2. Proposition 5.14 implies that $\mu_{E_3} = 0$. Similarly, one can verify that $\mu_{E_7} = \mu_{E_8} = 0$. Both $E_7(\mathbb{Q})[2]$ and $E_8(\mathbb{Q})[2]$ have order 2. $E_7(\mathbb{Q})[2]$ is ramified at 2 but not odd. $E_8(\mathbb{Q})[2]$ is odd but not ramified at 2. Proposition 5.14 again applies.

The $\mu$-invariants listed above turn out to be just as predicted by proposition 5.13. One can deduce this from the isogeny data given in [Cre]. One uses the following observation. *Suppose that $\varphi : E \to E'$ is a $\mathbb{Q}$-isogeny such that $\Phi = \ker(\varphi)$ is ramified and odd. Suppose also that $E'[2^\infty]$ contains a $\mathbb{Q}$-rational subgroup $\Phi'$ which is ramified and odd. Then $\varphi^{-1}(\Phi')$ is ramified* *and odd too. Its order is* $|\Phi| \cdot |\Phi'|$. For example, $E_2(\mathbb{Q})$ has three subgroups of order 2, one of which is ramified and odd. There is an $\mathbb{Q}$-isogeny of degree 2 from $E_2$ to $E_1$, $E_5$, and $E_6$. One can verify that $E_1(\mathbb{Q})$, $E_5(\mathbb{Q})$, and $E_6(\mathbb{Q})$ each has a subgroup of order 2 which is ramified and odd. Thus, $E_2[2^\infty]$ must have a subgroup $\Phi'$ of order 4 which is ramified and odd. Now $\Phi = E_5(\mathbb{Q})[2]$ is of order 2, generated by $(-\frac{109}{4}, \frac{105}{8})$. This $\Phi$ is ramified and odd. Since $E_5/\Phi \cong E_2$, it follows that $E_5[2^\infty]$ has a ramified and odd $\mathbb{Q}$-rational subgroup of order 8. Thus, proposition 5.13 implies that $\mu_{E_5} \geq 3$.

**Conductor = 69.** There are two such elliptic curves $E/\mathbb{Q}$. Both should have $|\text{III}_E(\mathbb{Q})| = 1$ and $|E(\mathbb{Q})| = 2$. For one of them, we have $c_3 = 1$, $c_{23} = 2$. For the other, $c_3 = 2$, $c_{23} = 1$. We have $a_2 = 1$ and so $|\widetilde{E}(\mathbb{F}_2)| = 2$. By theorem 4.1, we have $f_E(0) \sim 2$. Hence $f_E(T)$ is an irreducible element of $\Lambda$. Now let $\Phi = E(\mathbb{Q}) \cong \mathbb{Z}/2\mathbb{Z}$. For one of these curves, $\Phi$ is ramified at 2 but not odd. For the other, $\Phi$ is odd but not ramified at 2. Hence proposition 5.14 implies that $\mu_E = 0$. Since $f_E(T) \notin \Lambda^\times$, it follows that $\lambda_E \geq 1$. In fact, it turns out that $\lambda_E = \lambda_E^{M-W} = 1$, $\lambda_E^{\text{III}} = 0$, and $f_E(T) = T + 2$, up to a factor in $\Lambda^\times$. To see this, consider the quadratic twist $E^\xi$, where $\xi$ is the quadratic character corresponding to $\mathbb{Q}(\sqrt{2})$. Now $E^\xi(\mathbb{Q})$ has rank 1. Therefore, $E(\mathbb{Q}(\sqrt{2}))$ has rank 1. But $\mathbb{Q}(\sqrt{2})$ is the first layer in the cyclotomic $\mathbb{Z}_2$-extension $\mathbb{Q}_\infty/\mathbb{Q}$. Therefore, $\text{Sel}_E(\mathbb{Q}_\infty)_2$ contains the image of $E(\mathbb{Q}(\sqrt{2})) \otimes (\mathbb{Q}_2/\mathbb{Z}_2)$ under restriction as a $\Lambda$-submodule. Its characteristic ideal is $(T + 2)$. The assertions made above follow easily.

**Conductor = 195.** We will discuss the isogeny class consisting of A1–A8 in [Cre]. Some of the details below were worked out by Karl Rubin and myself with the help of PARI. We denote these curves by $E_1, \ldots, E_8$, respectively. We will show that $\lambda_{E_i} = \lambda_{E_i}^{M-W} = 3$, $\lambda_{E_i}^{\text{III}} = 0$, and that $\mu_{E_i}$ varies from 0 to 4 for $1 \leq i \leq 8$. Here is a table of the basic data.

|              | $E_1$ | $E_2$ | $E_3$ | $E_4$ | $E_5$ | $E_6$ | $E_7$ | $E_8$ |
|--------------|-------|-------|-------|-------|-------|-------|-------|-------|
| $|\text{III}| =$ | 1 | 1 | 1 | 1 | 1 | 1 | 4 | 1 |
| $|T| =$ | 4 | 8 | 8 | 4 | 4 | 4 | 2 | 2 |
| $c_3, c_5, c_{13} =$ | 4,1,1 | 8,2,1 | 4,4,4 | 16,1,1 | 2,8,2 | 2,2,8 | 1,4,1 | 1,16,1 |
| $f_{E_i}(0) \sim$ | 4 | 8 | 16 | 16 | 32 | 32 | 64 | 64 |
| $\mu_{E_i} =$ | 0 | 1 | 2 | 2 | 3 | 3 | 4 | 4 |

As before, we evaluated $|\text{III}|$ by assuming the Birch and Swinnerton-Dyer conjecture. But one could confirm directly that $|\text{III}_2|$ is as listed, which would be sufficient for us. [Cre] gives $|T|$ and the Tamagawa factors $c_3$, $c_5$, and $c_{13}$. The fourth row is a consequence of theorem 4.1. Since the $\lambda_{E_i}$'s are equal, clearly the $\mu$-invariants must vary. The last row becomes clear if we can show that $\mu_{E_1} = 0$. Unfortunately, this does not follow from proposition 5.14. The problem is that $\Phi = E(\mathbb{Q})[2]$ is of order 2, but is neither ramified at 2 nor odd. In fact, $\Phi$ is generated by $(6, -3)$, which is clearly not in the kernel of reduction modulo 2 and so is not in $C_2$. Also, $E(\mathbb{R})[2]$ has order 4 and $(6, -3)$ is in the connected component of $O_E$. This implies that $(6, -3) \notin C_\infty$.

We will verify that $\mu_{E_1} = 0$ by showing that $f_{E_1}(T)$ is divisible by $g(T) = (T + 2)(T^2 + 2T + 2)$ in $\Lambda$. Since the characteristic ideals of $\mathrm{Sel}_{E_i}(\mathbb{Q}_\infty)_2$ differ only by multiplication by a power of 2, it is equivalent to show that $g(T)$ divides $f_{E_i}(T)$ for any $i$. It then follows that $(f_{E_1}(T)) = (g(T))$ since $f_{E_1}(0)$ and $g(0)$ have the same valuation. Therefore, $\mu_{E_1}$ must indeed be zero. Let $F$ and $K$ denote the first and second layers in the cyclotomic $\mathbb{Z}_2$-extension $\mathbb{Q}_\infty/\mathbb{Q}$. Thus $\mathrm{Gal}(K/\mathbb{Q})$ is cyclic of degree 4 and $F$ is the unique quadratic subfield of $K$. In fact, $F = \mathbb{Q}(\sqrt{2})$, $K = F(\sqrt{2 + \sqrt{2}})$. We will show that $E_2(K) \otimes \mathbb{Q}$, considered as a $\mathbb{Q}$-representation of $\mathrm{Gal}(K/\mathbb{Q})$ contains the two nontrivial, $\mathbb{Q}$-irreducible representations of $\mathrm{Gal}(K/\mathbb{Q})$. One of them has degree 1 and factors through $\mathrm{Gal}(F/\mathbb{Q})$. The other has degree 2 and is faithful. The fact that $g(T)$ divides $f_{E_2}(T)$, and hence $f_{E_1}(T)$, follows easily.

The equation $y^2 + xy = x^3 - 115x + 392$ is the minimal Weierstrass equation defining $E_2$. It is slightly more convenient to calculate with the nonminimal equation $y^2 = (x - 1)(x - 2)(16x + 49)$, obtained by a simple change of variables. We single out the following two points satisfying this equation:

$$P = \left(0, 7\sqrt{2}\right), \qquad Q = \left(10 + 9\sqrt{2}, (123 + 78\sqrt{2})\sqrt{2 + \sqrt{2}}\right).$$

Now $P$ is rational over $F$, $Q$ is rational over $K$. To study $E_2(K)$, it is useful to first determine its torsion subgroup. In fact, we have

$$E_2(\mathbb{Q}_\infty)_{\mathrm{tors}} = E_2(\mathbb{Q})_{\mathrm{tors}} \cong (\mathbb{Z}/2\mathbb{Z}) \times (\mathbb{Z}/4\mathbb{Z}).$$

The structure of $E_2(\mathbb{Q})$ is given in [Cre]. It is easy to see that $E_2(\mathbb{Q}_\infty)_{\mathrm{tors}}$ is a 2-primary group since $E_2$ has good reduction at 2, 2 is totally ramified in $\mathbb{Q}_\infty/\mathbb{Q}$, and $|\widetilde{E}_2(\mathbb{F}_2)| = 4$. Now $\mathbb{Q}_\infty$ is totally real and so $E_2(\mathbb{Q}_\infty)_{\mathrm{tors}} \cong \mathbb{Z}/2\mathbb{Z} \times \mathbb{Z}/2^t\mathbb{Z}$ where $t \geq 2$. Assume $t \geq 3$. Then $E_2(\mathbb{Q}_\infty)_{\mathrm{tors}}$ would have 8 elements of order 8. Since their squares are in $E_2(\mathbb{Q})$, the orbit under $\Gamma = \mathrm{Gal}(\mathbb{Q}_\infty/\mathbb{Q})$ of an element of order 8 has cardinality at most 4. Hence such an element would be rational over $K$. We can rule out this possibility by noting that $E_2$ has good reduction at 31, 31 splits completely in $K/\mathbb{Q}$, and $|\widetilde{E}_2(\mathbb{F}_{31})| = 40$, which is not divisible by 16.

It is now clear that $P$ and $Q$ have infinite order. Also, $\mathrm{Gal}(F/\mathbb{Q})$ acts on $\langle P \rangle$ by $-1$ since $(0, -7\sqrt{2}) = -P$. Thus, $\langle P \rangle \otimes \mathbb{Q}$ is a $\mathrm{Gal}(K/\mathbb{Q})$-invariant subspace of $E_2(K) \otimes \mathbb{Q}$ giving the degree 1, nontrivial representation of $\mathrm{Gal}(K/\mathbb{Q})$. Similarly, $Q$ belongs to $\ker(Tr_{K/F})$, the kernel of the trace map from $E_2(K)$ to $E_2(F)$. Thus, $\ker(Tr_{K/F}) \otimes \mathbb{Q}$ is nonzero and provides at least one copy of the 2-dimensional, irreducible $\mathbb{Q}$-representation of $\mathrm{Gal}(K/\mathbb{Q})$. Therefore, $\mathrm{rank}(E_2(K)) \geq 3$. Considering the action of $\gamma = 1 + T$ on the image of $E_2(K) \otimes (\mathbb{Q}_2/\mathbb{Z}_2)$ in $\mathrm{Sel}_{E_2}(\mathbb{Q}_\infty)_2$ makes it clear that $g(T)$ does indeed divide $f_{E_2}(T)$, as claimed. As noted above, it now follows that $(f_{E_1}(T)) = (g(T))$. This implies that $\lambda_{E_1} = 3$, $\mu_{E_1} = 0$. More precisely, it is clear that $\lambda_{E_1}^{M-W} = 3$, $\lambda_{E_1}^{III} = 0$. For the $\mathbb{Q}$-isogenous curves $E_i$, $1 \le i \le 8$, we also have $\lambda_{E_i}^{M-W} = 3$, $\lambda_{E_i}^{III} = 0$, but $(f_{E_i}(T)) = (2^{\mu_{E_i}}g(T))$.

One can verify that in this example $\mu_{E_i} = m_{E_i}$. It is not hard to prove the existence of a $G_\mathbb{Q}$-invariant subgroup $\Phi_i$ of $E_i[2^\infty]$ with the expected order satisfying the hypotheses of proposition 5.13. Just as for conductor 15, one uses the $\mathbb{Q}$-isogenies between the $E_i$'s. By direct verification, one finds that $E_i[2]$ contains a ramified, odd $G_\mathbb{Q}$-invariant subgroup for $i = 2, \ldots, 8$. The listed isogenies then imply that $E_i[4]$ has a ramified, odd $G_\mathbb{Q}$-invariant subgroup of order 4 for $i = 3, \ldots, 8$. Then one sees that $E_i[8]$ contains such a subgroup of order 8 for $i = 5, \ldots, 8$. Finally, both $E_7$ and $E_8$ admit $\mathbb{Q}$-isogenies of degree 2 to $E_5$. The kernels of these $\mathbb{Q}$-isogenies are ramified and odd. The inverse image of the ramified, odd, $G_\mathbb{Q}$-invariant subgroup $\Phi_5$ of $E_5[8]$ will be the ramified, odd, $G_\mathbb{Q}$-invariant subgroup $\Phi_i$ of $E_i[16]$ of order 16, for $i = 7$ or 8.

Ken Kramer has found a description of the family of elliptic curves/$\mathbb{Q}$ which satisfy the hypotheses of proposition 5.13 for $m = 1, 2, 3$, and 4. Here we will give his description for $m = 1$ and $m = 4$, with the additional condition that $E$ have square-free conductor. For $m = 1$, his family is

$$E : y^2 + xy = x^3 - ax^2 - 4bx + (4a-1)b, \qquad \Phi = \left\langle \left( \frac{4a-1}{4}, \frac{1-4a}{8} \right) \right\rangle$$

where $a, b \in \mathbb{Z}$, $gcd(4a-1, b) = 1$, $(4a-1)^2 > 64b$, and either $a$ or $b$ is negative. The last conditions assure that $\Phi \subseteq C_\infty$. (If $b < 0$, then $E(\mathbb{R})$ has only one connected component. Then $\Phi$ is automatically contained in $C_\infty$. If $b > 0$ and $a < 0$, then the inequality $1 - 4a > 8\sqrt{b}$ implies that the above generator of $\Phi$ is the element of $E(\mathbb{R})[2]$ with minimal $x$-coordinate.) The discriminant of this equation, which is minimal, is $b((4a-1)^2 - 64b)^2$. If $b$ is odd, then $E$ has good, ordinary reduction at 2. Conjecturally, this family should give all elliptic curves/$\mathbb{Q}$ with good, ordinary or multiplicative reduction at 2 and square-free conductor such that $\mathrm{Sel}_E(\mathbb{Q}_\infty)_2$ has positive $\mu$-invariant. Kramer describes the elliptic curves with square-free conductor having a subgroup $\Phi$ of order 16 which is ramified and odd by the following equation:

$$E : y^2 = (x + 2c^4 - d^4)(x^2 + 4(cd)^4 - 4c^8)$$

where $c, d$ are distinct odd, positive integers, $c \equiv d \pmod 4$, and $gcd(c, d) = 1$. This equation is not minimal, but the discriminant of a minimal Weierstrass equation for $E$ is $(c^4 - d^4)c^4d^{16}/16$. Interchanging $c$ and $d$ gives a second elliptic curve, $\mathbb{Q}$-isogenous to $E$, but with discriminant of opposite sign. Thus, there are an even number of such elliptic curves in a $\mathbb{Q}$-isogeny class. A similar statement is true for $m = 2$ or 3, as Kramer shows. We refer to [K] for a more complete discussion.

We will end this article by returning to some of our earlier examples and discussing a few other examples, but now using Kato's theorem in conjunction with some calculations recently carried out by Ted McCabe. Assume that $E$ is a modular elliptic curve/$\mathbb{Q}$ and that $p$ is a prime where $E$ has good, ordinary reduction. Kato's theorem asserts that $f_E(T)$ divides $p^m f_E^{\mathrm{anal}}(T)$ in $\Lambda$ for some $m \geq 0$. Let $\lambda_E^{\mathrm{anal}}$ and $\mu_E^{\mathrm{anal}}$ denote $\lambda(f_E^{\mathrm{anal}})$ and $\mu(f_E^{\mathrm{anal}})$. Kato's theorem implies that $\lambda_E \leq \lambda_E^{\mathrm{anal}}$. McCabe has calculated approximations to the first few coefficients when $f_E^{\mathrm{anal}}(T)$ is written as a power series in $T$, enough to verify that $\mu_E^{\mathrm{anal}} = 0$ and to determine $\lambda_E^{\mathrm{anal}}$ for the examples he considers. These calculations allow us to justify several statements that were made earlier. As previously, we will use the value of $|\mathrm{Sel}_E(\mathbb{Q})_p|$ which is predicted by the Birch and Swinnerton-Dyer conjecture. In [M-SwD], one finds the results of calculations of the $p$-adic $L$-functions for the elliptic curves of conductors 11 and 17 and all primes $< 100$.

Kato's theorem reduces the verification of conjecture 1.13 to showing that $\lambda_E = \lambda_E^{\mathrm{anal}}$ and $\mu_E = \mu_E^{\mathrm{anal}}$. In a number of the following examples, these equalities can be shown. Before discussing the examples, we want to mention two situations which occur rather frequently.

$\boldsymbol{\lambda_E^{\mathrm{anal}} = \mu_E^{\mathrm{anal}} = 0}$. This means that $f_E^{\mathrm{anal}}(T) \in \Lambda^\times$. By Kato's theorem, it follows that $\lambda_E = 0$. Also, $f_E^{\mathrm{anal}}(0) = (1 - \beta_p p^{-1})^2 L(E/\mathbb{Q}, 1)/\Omega_E$ is a $p$-adic unit. Kolyvagin's theorem can then be used to verify the Birch and Swinnerton-Dyer conjecture, i.e., that $\mathrm{Sel}_E(\mathbb{Q})_p$ has the predicted order. Then by theorem 4.1, one would obtain that $f_E(0) \in \mathbb{Z}_p^\times$ too. That is, $f_E(T) \in \Lambda^\times$ and hence $\mu_E = 0$ and conjecture 1.13 is valid for $E$ and $p$.

$\boldsymbol{\lambda_E^{\mathrm{anal}} = 1,\ \mu_E^{\mathrm{anal}} = 0}$. We will also assume that $p$ is odd. Since $\lambda_E^{\mathrm{anal}} = 1$, $f_E^{\mathrm{anal}}(T)$ has exactly one root: $T = a$, where $a \in p\mathbb{Z}_p$. We mentioned in the introduction that $f_E^{\mathrm{anal}}(T^\iota)/f_E^{\mathrm{anal}}(T) \in \Lambda^\times$, where $T^\iota = (1 + T)^{-1} - 1$. Thus $(1 + a)^{-1} - 1$ is also a root of $f_E^{\mathrm{anal}}(T)$. It follows that $(1 + a)^2 = 1$ and, since $p$ is odd, $a = 0$. (For $p = 2$, we would have another possibility: $a = -2$.) Hence $f_E^{\mathrm{anal}}(0) = 0$ and so $T | f_E^{\mathrm{anal}}(T)$. We then must have $f_E^{\mathrm{anal}}(T)/T \in \Lambda^\times$. The $p$-adic $L$-function $L_p(E/\mathbb{Q}, s)$ would have a simple zero at $s = 1$. Assuming that $E$ has good, ordinary reduction at $p$, the complex $L$-function $L(E/\mathbb{Q}, s)$ would have an odd order zero at $s = 1$. (The "signs" in the functional equations for $L_p(E/\mathbb{Q}, s)$ and $L(E/\mathbb{Q}, s)$ are the same. See [M-T-T].) Perrin-Riou's analogue of the Gross-Zagier formula implies that $L'(E/\mathbb{Q}, 1) \neq 0$. Hence rank($E(\mathbb{Q})$) $= 1$. Consequently, $f_E(0) = 0$, $\lambda_E = 1$, and $f_E(T) = p^{\mu_E} T$, up to a factor in $\Lambda^\times$. Furthermore, Perrin-Riou's formula also shows that $h_p(P) \neq 0$, where $P$ is a generator of $E(\mathbb{Q})/E(\mathbb{Q})_{\mathrm{tors}}$, and that

$$h_p(P) \sim p(1 - \beta_p p^{-1})^{-2}(L'(E/\mathbb{Q}, 1)/\Omega_E h_\infty(P))^{-1}.$$

Kolyvagin's theorem should allow one to verify that $L'(E/\mathbb{Q}, 1)/\Omega_E h_\infty(P) \sim \left(\prod_{v \text{ bad}} c_v^{(p)}\right) |\mathrm{III}_E(\mathbb{Q})_p|/|E(\mathbb{Q})_p|^2$. If one then uses Schneider's result (for the case $F = \mathbb{Q}, r = 1$), one would obtain that $f_E(T)/T \in \Lambda^\times$, thus verifying that $\lambda_E = 1$, $\mu_E = 0$, and that conjecture 1.13 is valid for $E$ and $p$.

**Conductor = 67**. We consider $p = 3$. As expected, $\mu_E^{\text{anal}} = 0$. We can't verify that $\mu_E = 0$, as conjecture 1.11 predicts. ($E[3]$ is irreducible as a $G_{\mathbb{Q}}$-module.) McCabe finds that $\lambda_E^{\text{anal}} = 2$. As pointed out earlier, $\text{Sel}_E(\mathbb{Q}_\infty)_3$ is infinite. Hence, assuming that $\mu_E = 0$, we have $\lambda_E > 0$. By proposition 3.10, $\lambda_E$ must be even. Thus, $\lambda_E = 2$. Hence, if $\mu_E = 0$, it is clear from Kato's theorem, that $(f_E(T)) = (f_E^{\text{anal}}(T))$, i.e., conjecture 1.13 holds for $E$ and $p = 3$. In fact, we would have $\lambda_E^{M\text{-}W} = 0$ and $\lambda_E^{\text{III}} = 2$. To see this, suppose that $\lambda_E^{M\text{-}W} > 0$. Now $\Gamma$ acts in the finite-dimensional $\mathbb{Q}$-vector space $E(\mathbb{Q}_\infty) \otimes \mathbb{Q}$. The irreducible $\mathbb{Q}$-representations of $\Gamma$ have degrees $1, 2, 6, \ldots,$ $2 \cdot 3^{n-1}$ for $n \geq 1$. Since $E(\mathbb{Q})$ is finite and $\lambda_E = 2$, we would have $\lambda_E^{M\text{-}W} = 2$ and $\text{rank}(E(\mathbb{Q}_1)) = 2$, where $\mathbb{Q}_1$ is the first layer in $\mathbb{Q}_\infty/\mathbb{Q}$. This would imply that $g(T) = T^2 + 3T + 3$ divides $f_E(T)$. Hence $g(T)$ and $f_E(T)$ would differ by a factor in $\Lambda^\times$, which is impossible since $f_E(0) \sim 3^2$, $g(0) \sim 3$.

**Conductor = 915**. We consider again the elliptic curve $E$ corresponding to 915(A1) in [Cre]. We take $p = 7$ or $p = 43$. In both cases, McCabe finds that $\lambda_E^{\text{anal}} = 2$. Thus, $\lambda_E \leq 2$. It is then clear that $\lambda_E^{M\text{-}W} = 0$. For the only irreducible $\mathbb{Q}$-representation of $\Gamma$ with degree $\leq 2$ is the trivial representation. (The nontrivial irreducible $\mathbb{Q}$-representations have degree divisible by $p - 1$.) But $E(\mathbb{Q})$ is finite in this case. Hence, assuming that $\mu_E = 0$, we have $\lambda_E^{\text{III}} = 2$ for both $p = 7$ and $p = 43$. Also, just as for the preceding example, conjecture 1.13 would hold if $\mu_E = 0$.

**Conductor = 34**. We considered before the elliptic curve $E$ corresponding to 34(A1) in [Cre] and found that $\lambda_E^{M\text{-}W} = 2$, $\lambda_E^{\text{III}} = 0$, and $\mu_E = 0$ for $p = 3$. In this case, McCabe finds that $\lambda_E^{\text{anal}} = 2$ and $\mu_E^{\text{anal}} = 0$. Thus, Kato's theorem again implies conjecture 1.13: $(f_E(T)) = (f_E^{\text{anal}}(T))$ for $p = 3$. There are four elliptic curves of conductor 34, all $\mathbb{Q}$-isogenous. In general, conjecture 1.13 is preserved by $\mathbb{Q}$-isogeny. The power of $p$ dividing $f_E(T)$ changes in a way predicted by the result of [Sch3] or [Pe2]. The power of $p$ dividing $f_E^{\text{anal}}(T)$ changes in a compatible way, determined just by the change in $\Omega_E$. ($\Omega_E$ is the only thing that changes in the definition of $f_E^{\text{anal}}(T)$.) One can verify all of this directly. For $E$, PARI gives $\Omega_E = 4.4956\ldots$. Let $E'$ be 34(A3) in [Cre], which is related to $E$ by a $\mathbb{Q}$-isogeny of degree 3. Using the fact that $\mu_E = 0$, one fines that $\mu_{E'} = 1$. Therefore, $f_{E'}(T) = 3f_E(T)$. But PARI gives $\Omega_{E'} = 1.4985\ldots = \Omega_E/3$. (This must be exact.) Thus, one sees that $f_{E'}^{\text{anal}}(T) = 3f_E^{\text{anal}}(T)$. Conjecture 1.13 is valid for $E'$ too.

**Conductor = 26**. We take $p = 7$. For 26(B1), which we previously denoted by $E_1$, McCabe finds that $\mu_{E_1}^{\text{anal}} = 0$, $\lambda_{E_1}^{\text{anal}} = 4$, and $f_{E_1}^{\text{anal}}(0) \sim 7$. Thus, $f_{E_1}^{\text{anal}}(T)$ is an irreducible element of $\Lambda$. If $\mu_{E_1} = 0$, as conjecturally should be true, then Kato's theorem implies that $f_{E_1}(T) = f_{E_1}^{\text{anal}}(T)$, up to a factor in $\Lambda^\times$. Conjecture 1.13 would then be valid for $E_1$ (and for $E_2$ too). Thus, in this example, if $\mu_{E_1} = 0$, then $\lambda_{E_1} = 4$. Note that proposition 3.10 would tell us only that $\lambda_{E_1}$ is even. Also, just as in the example of conductor 915, we would have $\lambda_{E_1}^{M\text{-}W} = 0$, $\lambda_{E_1}^{\text{III}} = 4$.

**Conductor = 147.** Let $p = 13$. We will denote 147(B1, B2) by $E_1$ and $E_2$ as earlier. McCabe's calculation for 147B1 gives $\mu_{E_1}^{\text{anal}} = 0$, $\lambda_{E_1}^{\text{anal}} = 2$. Proposition 3.10 shows that $\lambda_{E_1}$ is even. If $\mu_{E_1} = 0$, as conjecture 1.11 predicts, then $\lambda_{E_1} > 0$. Hence $\lambda_{E_1} = 2$ and conjecture 1.13 would again follow from Kato's theorem. As in previous examples, we would have $\lambda_{E_1}^{M\text{-}W} = 0$, $\lambda_{E_1}^{\text{III}} = 2$.

**Conductor = 1225.** We consider again the two curves $E_1$ and $E_2$ of conductor 1225 discussed earlier. We take $p = 37$. McCabe finds that $\lambda_{E_1}^{\text{anal}} = 1$, $\mu_{E_1}^{\text{anal}} = 0$. Since $L(E_1/\mathbb{Q}, 1) = 0$, it follows that $f_{E_1}^{\text{anal}}(0) = 0$ and that $f_{E_1}^{\text{anal}}(T)/T \in \Lambda^\times$. As remarked earlier, it then should follow that $\lambda_{E_1} = 1$, $\mu_{E_1} = 0$ and that conjecture 1.13 holds. For $E_2$, we have $\lambda_{E_2} = 1$, $\mu_{E_2} = 1$. Conjecture 1.13 holds for $E_2$ too.

**Conductor = 58.** We consider $E' : y^2 + xy = x^3 - x^2 - x + 1$ and $p = 5$. In this case, $E'(\mathbb{Q}) \cong \mathbb{Z}$ and the predicted order of $\text{III}_{E'}(\mathbb{Q})$ is 1. McCabe finds that $\lambda_{E'}^{\text{anal}} = 1$, $\mu_{E'}^{\text{anal}} = 0$. It then follows that $\lambda_{E'} = 1$, $\mu_{E'} = 0$.

**Conductor = 406.** Consider $E : y^2 + xy = x^3 + x^2 - 2124x - 60592$. This is 406(D1) in [Cre]. We take $p = 5$. We have $c_2 = c_{29} = 2$, $c_7 = 5$, $|E(\mathbb{Q})| = 2$, and $\text{Sel}_E(\mathbb{Q})$ is predicted to have order 1. Thus, by theorem 4.1, $f_E(0) \sim 5$. Now it turns out that $E[5] \cong E'[5]$ as $G_{\mathbb{Q}}$-modules, where $E'$ is the elliptic curve of conductor 58 considered above. One verifies this by comparing the $q$-expansions of the modular forms corresponding to these curves. Since $\mu_{E'} = 0$, it follows that $\mu_E = 0$. Therefore, $\lambda_E \geq 1$. By proposition 3.10 $\lambda_E$ must be even. However, 7 splits completely in $\mathbb{Q}_1/\mathbb{Q}$, where $\mathbb{Q}_1$ denotes the first layer of the cyclotomic $\mathbb{Z}_5$-extension $\mathbb{Q}_\infty$ of $\mathbb{Q}$. (This is because $7^4 \equiv 1 \pmod{5^2}$.) Thus, there are 5 primes of $\mathbb{Q}_1$ lying over 7, each with Tamagawa factor equal to 5. The proof of corollary 5.6 can be used to show that $\lambda_E \geq 5$ and hence, since it is even, we must have $\lambda_E \geq 6$. McCabe finds that $\lambda_E^{\text{anal}} = 6$, $\mu_E^{\text{anal}} = 0$. Therefore, it follows that $\lambda_E = 6$, $\mu_E = 0$, and conjecture 1.13 holds for $E$ and $p = 5$. We also can conclude that $\lambda_E^{M\text{-}W} = 0$. This is so because $E(\mathbb{Q})$ is finite, $E(\mathbb{Q}_\infty) \otimes \mathbb{Q}$ is a finite dimensional $\mathbb{Q}$-representation of $\Gamma$, and the nontrivial irreducible $\mathbb{Q}$-representations of $\Gamma$ have degree divisible by 4. Hence $\text{Sel}_E(\mathbb{Q}_\infty)_5 = \text{III}_E(\mathbb{Q}_\infty)_5$ and $\lambda_E^{\text{III}} = 6$.

# References

[B-D-G-P]   K. Barré-Sirieix, G. Diaz, F. Gramain, G. Philibert, Une preuve de la conjecture de Mahler-Manin, *Invent. Math.* **124** (1996), 1–9.

[Be]        M. Bertolini, Selmer groups and Heegner points in anticyclotomic $\mathbb{Z}_p$-extensions, *Compositio Math.* **99** (1995), 153–182.

[BeDa1]     M. Bertolini, H. Darmon, Heegner points on Mumford-Tate curves, *Invent. Math.* **126** (1996), 413–456.

[BeDa2]     M. Bertolini, H. Darmon, Nontriviality of families of Heegner points and ranks of Selmer groups over anticyclotomic towers, to appear in *Journal of the Ramanujan Math. Society.*

[B-G-S]    D. Bernardi, C. Goldstein, N. Stephens, Notes $p$-adiques sur les courbes elliptiques, *J. reine angew. Math.* **351** (1984), 129–170.

[CoMc]    J. Coates, G. McConnell, Iwasawa theory of modular elliptic curves of analytic rank at most 1, *J. London Math. Soc.* **50** (1994), 243–269.

[CoGr]    J.. Coates, R. Greenberg, Kummer theory for abelian varieties over local fields, *Invent. Math.* **124** (1996), 129–174.

[CoSc]    J. Coates, C.-J. Schmidt, Iwasawa theory for the symmetric square of an elliptic curve, *J. reine angew. Math* **375/376** (1987), 104–156.

[Cre]    J. E. Cremona, *Algorithms for Modular Elliptic Curves*, Cambridge University Press (1992).

[D]    F. Diamond, On deformation rings and Hecke rings, *Ann. of Math.* **144** (1996), 137–166.

[F]    E. Friedman, Ideal class groups in basic $\mathbb{Z}_{p_1} \times \cdots \times \mathbb{Z}_{p_s}$-extensions of abelian number fields, *Invent. Math.* **65** (1982), 425–440.

[FeWa]    B. Ferrero, L. Washington, The Iwasawa invariant $\mu_p$ vanishes for abelian number fields, *Ann. of Math.* **109** (1979), 377–395.

[Gr1]    R. Greenberg, On a certain $l$-adic representation, *Invent. Math.* **21** (1973), 117–124.

[Gr2]    R. Greenberg, Iwasawa theory for $p$-adic representations, *Advanced Studies in Pure Mathematics*, **17** (1989), 97–137.

[Gr3]    R. Greenberg, Iwasawa theory for motives, *LMS Lecture Notes Series* **153** (1991), 211–233.

[Gr4]    R. Greenberg, Trivial zeroes of $p$-adic $L$-functions, *Contemporary Math.* **165** (1994), 149–174.

[Gr5]    R. Greenberg, The structure of Selmer groups, *Proc. Natl. Acad. Sci.* **94** (1997), 11125–11128.

[Gr6]    R. Greenberg, Iwasawa theory for $p$-adic representations II, in preparation.

[GrVa]    R. Greenberg, V. Vatsal, On the Iwasawa invariants of elliptic curves, in preparation.

[Gu1]    L. Guo, On a generalization of Tate dualities with application to Iwasawa theory, *Compositio Math.* **85** (1993), 125–161.

[Gu2]    L. Guo, General Selmer groups and critical values of Hecke $L$-functions, *Math. Ann.* **297** (1993), 221–233.

[HaMa]    Y. Hashimori, K. Matsuno, On the finite $\Lambda$-submodules of the Selmer groups of elliptic curves, preprint.

[Im]    H. Imai, A remark on the rational points of abelian varieties with values in cyclotomic $\mathbb{Z}_p$-extensions, *Proc. Japan Acad.* **51** (1975), 12–16.

[Jo]    J. Jones, Iwasawa $L$-functions of multiplicative Abelian varieties, *Duke Math. J.* **59** (1989), 399–420.

[K]    K. Kramer, Elliptic curves with non-trivial 2-adic Iwasawa $\mu$-invariant, in preparation.

[Man]    J. Manin, Cyclotomic fields and modular curves, *Russian Math. Surveys* **26** (1971), 7–78.

[Maz1]    B. Mazur, Rational points of abelian varieties with values in towers of number fields, *Invent. Math.* **18** (1972), 183–266.

[Maz2]    B. Mazur, Rational isogenies of prime degree, *Invent. Math.* **44** (1978), 129–162.

[Maz3]     B. Mazur, On the arithmetic of special values of *L*-functions, *Invent. Math.* **55** (1979), 207–240.

[Maz4]     B. Mazur, Modular curves and arithmetic, Proceedings of the International Conference of Mathematicians, Warszawa (1983).

[M-SwD]    B. Mazur, P. Swinnerton-Dyer, Arithmetic of Weil curves, *Invent. Math.* **25** (1974), 1–61.

[M-T-T]    B. Mazur, J. Tate, J. Teitelbaum, On *p*-adic analogues of the conjectures of Birch and Swinnerton-Dyer, *Invent. Math.* **84** (1986), 1–48.

[Mi]       J. S. Milne, *Arithmetic Duality Theorems*, Academic Press (1986).

[Mo]       P. Monsky, Generalizing the Birch-Stephens theorem, *Math. Zeit.* **221** (1996), 415–420.

[Pe1]      B. Perrin-Riou, Arithmetique des courbes elliptiques et théorie d'Iwasawa, *Memoire Soc. Math. France* **17** (1984).

[Pe2]      B. Perrin-Riou, Variation de la fonction *L p*-adique par isogénie, *Advanced Studies in Pure Mathematics* **17** (1989), 347–358.

[Pe3]      B. Perrin-Riou, Points de Heegner et dérivées de fonctions *L p*-adiques, *Invent. Math.* **89** (1987), 455–510.

[Pe4]      B. Perrin-Riou, Théorie d'Iwasawa *p*-adique locale et globale, *Invent. Math.* **99** (1990), 247–292.

[Ri]       K. Ribet, Torsion points of abelian varieties in cyclotomic extensions, *Enseign. Math.* **27** (1981), 315–319.

[Ro]       D. Rohrlich, On *L*-functions of elliptic curves and cyclotomic towers, *Invent. Math.* **75** (1984), 409–423.

[Ru1]      K. Rubin, On the main conjecture of Iwasawa theory for imaginary quadratic fields, *Invent. Math.* **93** (1988), 701–713.

[Ru2]      K. Rubin, The "main conjecture" of Iwasawa theory for imaginary quadratic fields, *Invent. Math.* **103** (1991), 25–68.

[R-W]      K. Rubin, A. Wiles, Mordell-Weil groups of elliptic curves over cyclotomic fields, *Progress in Mathematics* **26** (1982), 237–254.

[Sch1]     P. Schneider, Iwasawa *L*-functions of varieties over algebraic number fields, A first approach, *Invent. Math.* **71** (1983), 251–293.

[Sch2]     P. Schneider, *p*-adic height pairings II, *Invent. Math.* **79** (1985), 329–374.

[Sch3]     P. Schneider, The *μ*-invariant of isogenies, *Journal of the Indian Math. Soc.* **52** (1987), 159–170.

[Se1]      J. P. Serre, *Abelian l-adic Representations and elliptic curves*, W. A. Benjamin (1968).

[Se2]      J. P. Serre, Cohomologie Galoisienne, *Lecture Notes in Math.* **5** (1964).

[Si]       J. Silverman, The Arithmetic of Elliptic Curves, *Grad. Texts in Math.* **106** (1986).

[St]       G. Stevens, Stickelberger elements and modular parametrizations of elliptic curves, *Invent. Math.* **98** (1989), 75–106.

[Ta]       J. Tate, The arithmetic of elliptic curves, *Invent. Math.* **23** (1974), 179–206.

[Wa1]      L. Washington, The non-*p*-part of the class number in a cyclotomic $\mathbb{Z}_p$-extension, *Invent. Math.* **49** (1979), 87–97.

[Wa2]      L. Washington, Introduction to Cyclotomic Fields, *Grad. Texts in Math.* **83** (1986).